# The solution cosmological constant problem.Quantum Field Theory in fractal space-time with negative Hausdorff-Colombeau dimensions and dark matter nature.

J. Foukzon
Department of mathematics, Israel Institute of Technology, Haifa, Israel
E-mail: jaykovfoukzon@list.ru
E.Menkova
A.Potapov

**Abstract**. The cosmological constant problem arises because the magnitude of vacuum energy density predicted by quantum field theory is about $120$ orders of magnitude larger than the value implied by cosmological observations of accelerating cosmic expansion. We pointed out that the fractal nature of the quantum space-time with negative Hausdorff- Colombeau dimensions can resolve this tension. The canonical Quantum Field Theory is widely believed to break down at some fundamental high-energy cutoff $\Lambda_*$ and therefore the quantum fluctuations in the vacuum can be treated classically seriously only up to this high-energy cutoff. In this paper we argue that Quantum Field Theory in fractal space-time with negative Hausdorff-Colombeau dimensions gives high-energy cutoff on natural way.We argue that there exists hidden physical mechanism which cancel divergences in canonical $QED_4, QCD_4,$ Higher-Derivative-Quantum-Gravity, etc. In fact we argue that corresponding supermassive Pauli-Villars ghost fields really exists.It means that there exist the ghost- driven acceleration of the univers hidden in cosmological constant.

In order to obtain desired physical result we apply the canonical Pauli-Villars regularization up to $\Lambda_*$.This would fit in the observed value of the dark energy needed to explain the accelerated expansion of the universe if we choose highly symmetric masses distribution between standard matter and ghost matter below that scale $\Lambda_*$, i.e., $f_{s.m}(\mu) \approx -f_{g.m}(\mu),\ \mu = mc, \mu \leq \mu_{\mathbf{eff}}, \mu_{\mathbf{eff}}c < \Lambda_*$ The small value of the cosmological constant explaned by tiny violation of the symmetry between standard matter and ghost matter.Dark matter nature also explaned using a common origin of the dark energy and dark matter phenomena.

**Content**

# 1.Introduction

## 1.1.The cosmological constant problem and Quantum Field Theory in fractal spacetime with negative dimension.

One of the greatest challenges in modern physics is to reconcile general relativity and elementary particles physics into a unified theory. Perhaps the most dramatic clash between the two theories lies in the cosmological constant problem [1-6] and in the problem of the Dark (i.e., non-luminous and non-absorbing) Matter nature is, arguably, the most widely discussed topic in contemporary particle physics.Naive predictions of vacuum energy from canonical quantum field theory predict a magnitude so high that the expansion of the Universe should have accelerated so quickly that no any structure could have formed. The predicted rate of acceleration resulting from vacuum energy is famously $120$ orders of magnitude larger than what is observed. In order to avoid these difficultnes mentioned above we assume that:(i) physics of elementary particles essentially is separated into low/high energy ones, (ii) the standard notion of smooth spacetime is assumed to be altered at a high energy cutoff scale $\Lambda_*$ and a new treatment based on QFT in a fractal spacetime with negative dimension is used above that scale $\Lambda_*$.In this paper we argue that Quantum Field Theory in fractal space-time with negative Hausdorff-Colombeau dimensions [15] gives high-energy cutoff on natural way.No one knows what dark energy is, but we need it to explain the discovered accelerated expansion of the Universe. The most elegant and natural solution is to identify dark energy with the energy of the quantum vacuum predicted by Quantum Field Theory, but the trouble is that QFT predicts the energy density of the vacuum to be orders of magnitude larger than the observed dark energy density:

$$\varepsilon_{\mathbf{de}} \approx 7.5 \times 10^{-27} kg/m^3. \qquad (1.1.1)$$

Recall that it was stressed by Zeldovich [1] that quantum field theory generically demands that cosmological constant or, let us repeat, what is the same, vacuum energy is non-vanishing.Summing the zero-point energies of all normal modes of some quantum field of mass $m$ up to a wave number cut-off $\Lambda_*/c^2 \gg m$, QFT yields [1],[5] a vacuum energy density

$$\varepsilon_{\mathbf{vac}} \sim \int_0^{p_*} d^3p\sqrt{p^2+m^2} \simeq p_*^4. \qquad (1.1.2)$$

# 1.2. Sources of Vacuum Energy

It is not excluded experimentally that the number of fermionic and bosonic species in Nature are the same. Moreover it is practically a necessity, because otherwise vacuum energy density would be infinite. Still the masses of bosons and corresponding fermions are different and, with arbitrary relations between their masses, only the leading term, which diverges as the fourth power of the integration limit, would be canceled out. However in some supersymmetric theories with spontaneous symmetry breaking there may be specific relations between masses of different fields which ensure the compensation not only of the leading term but also quadratically and logarithmically divergent terms. This looks as a very strong argument in favor of such models. However the finite terms are not compensated. Moreover in global supersymmetric theories finite contributions into ρvac must be nonzero and by the order of magnitude they are equal to

$$\varepsilon_{vac}^{susy} = m_{susy}^4 \qquad (1.2.1)$$

where $m^{susy}$ is the scale of supersymmetry breaking. It is known from experiment that $m_{susy} \geq 100$ GeV. Correspondingly $\varepsilon_{vac}^{susy} \geq 10^8$ GeV, i.e. 55 orders of magnitude larger than the permitted upper bound. In more advanced supersymmetric theories which include gravity (the so called supergravity or local supersymmetry) the condition of non-vanishing vacuum energy in the broken symmetry phase is not obligatory. However, if one does not take a special care, the value of vacuum energy in unbroken supergravity models is typically about $m_{Pl}^4 \approx 10^{76}$ GeV. One can choose in principle the parameters in such a way that this contribution into ρvac is compensated down to zero with the accuracy $10^{-123}$ but this demands a fantastic fine-tuning.

One more source of vacuum energy is the energy of the scalar (Higgs) field in the theories with spontaneous symmetry breaking.

## 1.3.New Model of ”Nullification” of Vacuum Energy

Several possible approaches to the problem of vacuum energy have been discussed in the contemporary literature, for the review see ref. [5]. They can be roughly devided into four different groups:

(1) Modification of gravity on large scales.

(2) Anthropic principle.

(3) Symmetry leading to $\rho_{\mathbf{vac}} = 0$.

(4) Adjustment mechanism.

(5) Hidden nonstandard matter sector and corresponding symmetry leading to $\rho_{\mathbf{vac}} \simeq 0$.

1.A modification of gravity at large scales should be done in such a way that the general covariance, which ensures vanishing of the graviton mass, is preserved, energy momentum tensor is covariantly conserved, and simultaneously the vacuum part of this tensor, which is proportional to $g_{\mu\nu}$,does not gravitate. This is definitely not an easy thing to do. Possibly due to this reasons there is no satisfactory model of this kind at the present time.

2.Anthropic principle states that the conditions in the universe must be suitable for life, otherwise there would be no observer that could put a question why the universe is such and not another. With cosmological constant which is as large as predicted by natural estimates in quantum theory, life of our type is definitely impossible. Still this point of view does not look very appealing. The situation is similar to the one that existed in the Friedmann cosmology before inflationary resolution of the fundamental cosmological problems has been proposed. There is one more difficulty in the implimenttion of the anthropic principle. Even if we assume that it is effective, there are no visible building blocks to achieve the necessary compensation of vacuum energy. One can say of course that this compensation is not achieved by a physical field but just by a subtraction constant or in other words by a choice of the position of zero on the energy axis. In other words it is assumed that there is some energy coming from nowhere, which exactly cancels out all the contributions of different physical fields. Though formally this is not excluded, it definitely does not look beautiful.

3.*Probably the most appealing would be a model based on a symmetry principle which forbids a nonzero vacuum energy. Such a symmetry should connect known fields with new unknown ones. Some of those fields should be very light to achieve* *the cancellation on the scale $10^{-3}$eV. Neither such fields are observed, nor such a symmetry is known.*

4. An adjustment mechanism seems the most promising one at the present time. The idea is similar to the mechanism of solving the problem of natural CP-conservation in quantum chromodynamics by the axion field. The axion potential automatically acquires a minimum at the value of the field amplitude that cancels out the CP-odd contribution from the so called theta-term, $\theta G\widetilde{G}$. Similar mechanism can hopefully kill vacuum energy. Let us assume that there is a very light or massless field coupled to gravity in such a way that it is unstable in De Sitter background and develops the condensate whose energy-momentum tensor is equal by magnitude and opposite by sign to the original vacuum energy-momentum tensor. Though it looks rather promising, it is very difficult, if possible at all, to construct a realistic model based on this idea.

5.Hidden nonstandard matter sector and corresponding symmetry leading to $\rho_{\mathbf{vac}} \simeq 0$.

The luminous (light-emitting) components of the universe only comprise about 0.4% of the total energy. The remaining components are dark. Of those, roughly 3.6% are identified: cold gas and dust, neutrinos, and black holes. About 23% is dark matter, and the overwhelming majority is some type of gravitationally self-repulsive dark energy. There is no candidate in the standard model of particle physics.In what way does dark matter extend the standard model?

**Remark 1.3.1.**In order to explain physical nature of dark matter sector we assume that main part of dark matter,i.e., $\simeq 23\% - 4.6\% = 18\%$ (see Fig.2.3.3) formed by supermassive ghost particles vith masess such that $mc^2 > \Lambda_*$.

**Remark 1.3.2.**In order to obtain QFT description of the dark component of matter in natural way we expand now the standard model of particle physics on a sector of ghost particles. QFT in a ghost sector developed in Sect.3.1-3.4 and Sect.4.1-4.8.

**The paper is organized as follows**:

In **Sec.2**, classical Zel'dovich approuch [1] to cosmological constant problem revisited.

In **Sec.2.1**, we summarize the aspects of the cosmological constant problem that are relevant to this work.

In **Sec.2.2**, we summarize the model of cosmological dynamics in the presence of a vacuum energy that was introduced in [1-4], and how it attempts to resolve the problem.

In **Sec.2.3**, dark matter nature is considered. We argue that dark matter sector essentially formed only by super massive ghost particles.The Standard Model of fundamental interactions is extendent on a ghost sector.

In **Sec.3**, Pauli–Villars ghosts as physical dark matter is considered.

In **Sec.3.1**,Pauli– Villars renormalization of the $\lambda\varphi_4^4$ field theory by using Pauli–Villars

host

fields is considered.New physical interpretation of the scalar Pauli–Villars host fields is given.

In **Sec**.**3**.**2**, Pauli–Villars renormalization of the the $QED_3$ by using Pauli–Villars ghost fields is considered. New physical interpretation is given.

In **Sec**.**3**.**3**, High covariant derivatives renormalization as Pauli–Villars renormalization of non-Abelian gauge theories.New physical interpretation is given.

In Sec.**3**.**4**, Pauli–Villars renormalization of $QED_4$ via Colombeau generalized functions is considered successfully.The physical significance of Pauli-Villars renormalization is explained.

In **Sec**.**4**,we construct OFT in a ghost sector via dimensional renormalization supported by Colombeau generalized functions.

In **Sec**.**4**.**1**, Dimensional Regularization via Colombeau generalized functions is given.

In **Sec**.**4**.**2**, the scalar theory $\lambda\varphi_4^4$ in standard sector via Colombeau generalized functions in one-loop approximation is given.

In **Sec**.**4**.**3**, the scalar theory $\lambda\varphi_4^4$ in a ghost sector via Colombeau generalized functions in two-loop approximation is given.

In **Sec**.**4**.**4**, quantum electrodynamics in a ghost sector via Colombeau generalized functions is given.

In **Sec**.**4**.**5**, quantum chromodynamics in a ghost sector via Colombeau generalized functions is given.

In **Sec**.**4**.**6**, the general structure of the $\mathfrak{R}$-operation in a ghost fields sector via Colombeau generalized functions is given.

In **Sec**.**4**.**7**, the renormalization Group in a ghost sector is considered.

In **Sec**.**4**.**8**, dimensional regularization and the $\overline{MS}$ scheme in a ghost sector is given.

In **Sec**.**5**, the higher-derivative-quantum-gravity is considered as physical quantum-gravity theory below high energy cutoff $\Lambda_*$The renormalizable models of quantum-gravity which we have considered in this section, many years mistakenly regarded only as constructs for a study of the ultraviolet problem of quantum gravity. The difficulties with unitarity appear to preclude their direct acceptability as canonical physical theories in locally Minkowski space-time. In canonical case they do have only some promise as phenomenological models.However, for their unphysical behavior may be restricted to arbitrarily large energy

scales mentioned above by an appropriate limitation on the renormalized masses $m_2$ and

$m_0$.Actually, it is only the massive spin-two excitations of the field which give the trouble

with unitarity and thus require a very large mass. The limit on the mass $m_0$ is determined only by the observational constraints on the static field.

In Sec.**6**, Hausdorff-Colombeau measure and associated negative Hausdorff-Colombeau dimensions is considered successfully.

In Sec.**6**.**1**, we provide fractional integration in negative dimensions on natural way via Colombeau generalized functions.

In Sec.**6**.**2**, Using Hausdorff measure with associated positive Hausdorff dimension the rigorous definition of the Colombeau-Feynman path integral in $D = 4$ from dimensional regularization is given.

In Sec.**6**.**3**, we provide Hausdorff-Colombeau measure and associated negative Hausdorff-Colombeau dimensions.

In Sec.**6**.**4**, we provide the main properties of the Hausdorff- Colombeau metric measures with associated negative Hausdorff-Colombeau dimensions.

In Sec.**7**, we provide scalar quantum field theory in spacetime with Hausdorff-Colombeau negative dimensions.

In Sec.**7**.**1**,the equation of motion and Hamiltonian in spacetime with Hausdorff-Colombeau negative dimensions is considerd.

In Sec.**7**.**2**,propagator of a free scalar quantum field in configuration space with Hausdorff-Colombeau negative dimensions is considerd.

In Sec.**7**.**3**,Green’s functions corresponding to a self-interecting scalar quantum field in spacetime with Hausdorff-Colombeau negative dimensions is considerd.

In Sec.**7**.**4**,saddle-point evaluation of the Colombeau-Feynman path integral corresponding to a self-interecting scalar quantum field in negative dimensions is considerd successfully.

In Sec.**7**.**5**,an criteria of the power-counting renormalizability of $P(\varphi)_{D^-}$ scalar quantum field theory in negative dimensions $D^- < 0$ is considerd successfully.

In Sec.**7**.**6**,we have proved power-counting renormalizability of Einstein gravity in negative dimensions.

In Sec.**7**.**7**,an criteria of thepower-counting renormalizability of Hŏrava gravity in negative dimensions.

In Sec.**8**,the solution cosmological constant problem is considerd successfully.

In Sec.**8**.**1**,Zeropoint energy density corresponding to Einstein-Gliner-Zel’dovich vacuum with tiny Lorentz invariance violation is considerd.

In Sec.**8**.**2**,Zeropoint energy density corresponding to a non-singular Gliner cosmology is considerd.

In Sec.**8**.**3**, Zeropoint energy density in models with supermassive physical ghost fields is considerd.

In Sec.**9**, we compare the classical and non classical assumptions that are made in the different formulation of the cosmological constant problem.

In Sec.**9**.**1**,we briefly review the canonical assumptions that are made in the usual formulation of the cosmological constant problem.

In Sec.**9**.**2**,we list the modified assumptions that are made in this paper.

In Sec.**9**.**3**,

In Sec.**9**.**4**,

In Sec.**9**.**5**,

In Sec.**9**.**6**,semiclassical Möller-Rosenfeld gravity via aprouch proposed in this paper is considerd. We conclude that Moller-Rosenfeld equation holds again in a good approximation.

In Sec.**9**.**7**,we briefly discussed higher-derivative quantum gravity at energy scale $\Lambda \leq \Lambda_*$ and corresponding controlable tiny violetion of the unitarity condition.

We conclude with the physical significance of the new results in Sec.**9**-**10**.

# 2. Zel'dovich approach to cosmological constant problem by using Pauli-Villars regularization revisited.Ghost particles as physical dark matter.

## 2.1.The formulation of the cosmological constant problem.

The cosmological constant problem arises at the intersection between general relativity and quantum field theory, and is regarded as a fundamental unsolved problem in modern physics. Remind that a peculiar and truly quantum mechanical feature of the quantum fields is that they exhibit zero-point fluctuations everywhere in space, even in regions which are otherwise 'empty' (i.e.devoid of matter and radiation).This vacuum energy density is believed to act as a contribution to the cosmological constant $\lambda$ appearing in Einstein's field equations from 1917,

$$R_{\mu\nu} - \frac{1}{2}g_{\mu\nu}R = \frac{8\pi G}{c^4}T'_{\mu\nu} \qquad (2.1.1)$$

where $R_{\mu\nu}$ and $R$ refer to the curvature of space-time, $g_{\mu\nu}$ is the metric, $T'_{\mu\nu}$ the energy-momentum tensor,

$$T'_{\mu\nu} = T_{\mu\nu} + \frac{c^4\lambda}{8\pi G}\begin{pmatrix} 1 & 0 & 0 & 0 \\ 0 & -1 & 0 & 0 \\ 0 & 0 & -1 & 0 \\ 0 & 0 & 0 & 1 \end{pmatrix} \qquad (2.1.2)$$

where $T_{\mu\nu}$ is the energy-momentum tensor of matter. Thus $T'_{00} = T_{00} + \varepsilon_\lambda$, $T'_{\alpha\beta} = T_{\alpha\beta} + \delta_{\alpha\beta}P_\lambda$, where

$$\varepsilon_\lambda = -P_\lambda = c^4\lambda/8\pi G. \qquad (2.1.3)$$

Remind that under Lorentz transformations $(\varepsilon_\Lambda, P_\Lambda) \to \varepsilon'_\Lambda, (\varepsilon_\Lambda, P_\Lambda) \to P'_\Lambda$ the quantities $\varepsilon_\Lambda$ and $P_\lambda$ are changes by law

$$\varepsilon'_\lambda = \frac{\varepsilon_\lambda + \beta^2 P_\lambda}{1-\beta^2}, P'_\lambda = \frac{P_\lambda + \beta^2\varepsilon_\lambda}{1-\beta^2}. \qquad (2.1.4)$$

Thus for the quantities $\varepsilon_\lambda$ and $P_\lambda$ Lorentz invariance holds by Eq.(2.1.3) [1].

In modern cosmology it is assumed that the observable universe was initially vacuumlike, i.e., the cosmological medium was non-singular and Lorentz invariant. In the earlier, non-singular Friedmann cosmology the Friedmann universe comes into being during the phase transition of an initial vacuumlike state to the state of 'ordinary' matter [2],[3].

The Friedmann equations start with the simplifying assumption that the universe is spatially homogeneous and isotropic, i.e. the cosmological principle; empirically, this is justified on scales larger than ~100 Mpc. The cosmological principle implies that the metric of the universe must be of the form Robertson-Walker metric [2].
Robertson-Walker metric reads

$$ds^2 = dt^2 - a^2(t)\left[\frac{dr^2}{1-kr^2} + r^2(d\theta^2 + \sin^2\theta d\varphi^2)\right]. \quad (2.1.5)$$

For such a metric, the Ricci curvature scalar is $R = -6k$ and it is said that space has the curvature $k$.The scaling factor $a(t)$ rescales this curvature for a given time $t$, producing a curvature $k(t) = k/a(t)$.The scaling factor $a(t)$ is given by two independent Friedmann equations for modeling a homogeneous, isotropic universe reads

$$\dot{a}^2 = \frac{G}{3}\varepsilon a^2 - k, \ddot{a} = -\frac{G}{6}(\varepsilon + 3p) \quad (2.1.6)$$

and the equation of state

$$p = p(\varepsilon), \quad (2.1.7)$$

where $p$ is pressure and $\varepsilon$ is a density of the cosmological medium. For the case of the vacuumlike cosmological medium equation of state reads [1],[2],[3],[4]:

$$p = -\varepsilon. \quad (2.1.8)$$

By virtue of Friedman's equations (2.1.6) in the Universe filled with a vacuum-like medium, the density of the medium is preserved, i.e. $\varepsilon = const$, but the scale factor $a(t)$ grows exponentially. By virtue of continuity, it can be assumed that the admixture of a substance does not change the nature of the growth of the latter, and the density of the medium hardly changes. This growth, interpreted by analogy with the Friedmann models as an expansion of the universe, but almost without changing the density of the medium! - was named inflation. The idea of inflation is the basis of inflation scenarios [2].

Non-singular cosmology [2],[3],[4] suggests that the initial state of the observable universe was vacuum-like, but unstable with respect to the phase transition to the ordinary non-Lorentz-invariant medium. This, for example, takes place if, by virtue of the equations of state of the medium, a fluctuation decrease in its density $d$ violates the condition of vacuum-like degeneration, $p = -\varepsilon$ or, which is the same, $3p + \varepsilon = -2\varepsilon < 0$, replacing it with

$$-2\varepsilon < 3p + \varepsilon < 0. \quad (2.1.9)$$

According to Friedman's equations, it corresponds to an accelerated expansion of the cosmological medium, accompanied by a drop in its density, which makes the process irreversible [2]. The impulse for expansion in this scenario, the vacuum-like environment, is not reported to itself (bloating), but to the emerging Friedmann environment.

In review [5], Weinberg indicates that the first published discussion of the contribution of quantum fluctuations to the cosmological constant was a 1967 paper by Zel'dovich [6].In his article [1] Zel'dovich emphasizes that zeropoint energies of particle physics theories cannot be ignored when gravitation is taken into account, and since he explicitly discusses the discrepancy between estimates of vacuum energy and observations, he is clearly pointing to a cosmological constant problem. As well known zeropoint energy density of scalar quantum field,etc.is divergent

$$\varepsilon_{\mathbf{vac}}(m) = \frac{2\pi c}{(2\pi\hbar)^3}\int_0^\infty \sqrt{p^2+m^2c^2}\,p^2dp = \infty. \tag{2.1.10}$$

In order avoid difficultnes mentioned above, in article [1] Zel'dovich has applied canonical Pauli-Villars regularization [7],[8] and formally has obtained an finite result (his formulas [1], Eqs. (VIII.12)-(VIII.13) p.228)

$$\varepsilon_{\mathbf{vac}} = -p_{\mathbf{vac}} = \frac{1}{8}\int_0^{\mu_{\mathbf{eff}}} f(\mu)\mu^4(\ln\mu)d\mu = \frac{c^4\lambda}{8\pi G}, \tag{2.1.11}$$

where

$$\int_0^{\mu_{\mathbf{eff}}} f(\mu)d\mu = \int_0^{\mu_{\mathbf{eff}}} f(\mu)\mu^2d\mu = \int_0^{\mu_{\mathbf{eff}}} f(\mu)\mu^4d\mu = 0. \tag{2.1.12}$$

**Remark 2.1.1.**Unfortunately the Eq(2.1.11)-Eq(2.1.12) gives nothing in order to obtain desired small numerical values of the zero-point energy density $\varepsilon_{\mathbf{vac}}$.It is clear that aditional physical assumptions is needed.

In his paper [1], Zel'dovich arrives at a zero-point energy (his formula [1],Eq.(IX.1))

$$\varepsilon_{\mathbf{vac}} = m\left(\frac{mc}{\hbar}\right)^3 \sim 10^{17}g/cm^3, \lambda \sim 10^{-10}cm^{-2}, \tag{2.1.13}$$

where $m$ (the ultra-violet cut-of ) is taken equal to the proton mass. Zel'dovich notes that since this estimate exceeds observational bounds by 46 orders of magnitude it is clear that "...such an estimate has nothing in common with reality".

In his paper [1], Zel'dovich wroted:" Recently A.D. Sakharov proposed a theory of gravitation, or, more precisely, a justification GR equations based on consideration of vacuum fluctuations.In this theory, the essential assumption is that there is some elementary length $L$ or the corresponding limiting momentum $p_0 = \hbar/L$. Shorter lengths or for large impulses theory is not applicable. Sakharov gets the expression of gravitational constant $G$ through $L$ or $p_0$ (his formula [1],Eq.(IX.6))

$$G = \frac{c^3L^2}{\hbar} = \frac{\hbar c^3}{p_0^2}. \tag{2.1.14}$$

This expression has been known since the days of Planck, but it was read "from right to left": gravity determines the length $L$ and the momentum $p_0$. According to Sakharov, $L$ and $p_0$ are primary. Substitute Eq.(IX. 6) in the expression Eq.(IX.4) (see [1]), we get

$$\rho_\Lambda = \frac{m^6c^5}{p_0^2\hbar^3}, \varepsilon_{vac} = \frac{m^6c^7}{p_0^2\hbar^3}. \tag{2.1.15}$$

That is expressions that the first members (in the formulas [1],Eqs.(VIII.10)-(VIII. 11)) which are vanishes (with $p_0 \to \infty$).Thus, we can suggest the following interpretation of the cosmological constant: there is a theory of elementary particles, which would give (according to the mechanism that has not been revealed at the present time) identically zero vacuum energy, if this theory were applicable infinitely, up to arbitrarily large momentum; there is a momentum $p_0$, beyond which the theory is nont aplicable; along with other implications, modifying the theory gives different from zero vacuum energy; general considerations make it likely that the effect is portional $p_0^{-2}$.Clarification of the question of the existence and magnitude of the cosmological constant will also be of fundamental importance for the theory of elementary particles".

In contrast with Zel'dovich paper [1] we assume that Poincaré group is deformed at some fundamental high-energy cutoff $\Lambda_*$ [9],[10],[11] in accordance on the basis of the following deformed Poisson brackets

$$\begin{gathered}\{x^\mu, x^\nu\} = \varkappa^{-1}(x^\mu \eta^{0\nu} - x^\nu \eta^{\mu 0}), \{p^\mu, p^\nu\} = 0, \\ \{x^\mu, p^\nu\} = -\eta^{\mu\nu} + \varkappa^{-1}\eta^{\mu 0}p^\nu\end{gathered} \tag{2.1.16}$$

where $\mu, \nu, = 0, 1, 2, 3,\ \eta^{\mu\nu} = (+1, -1, -1, -1)$ and $\varkappa$ is a parameter identified as the ratio between the high-energy cutoff $\Lambda_*$ and the light speed. The corresponding to (2.1.16) momentum transformation reads [11]

$$\begin{aligned} p_0' &= \frac{\gamma(p_0 - up_x)}{1 + (c\varkappa)^{-1}[(\gamma-1)p_0 - \gamma u p_x]}, p_x' = \frac{\gamma(p_x - up_0/c^2)}{1 + (c\varkappa)^{-1}[(\gamma-1)p_0 - \gamma u p_x]}, \\ p_y' &= \frac{p_y}{1 + (c\varkappa)^{-1}[(\gamma-1)p_0 - \gamma u p_x]}, p_z' = \frac{p_z}{1 + (c\varkappa)^{-1}[(\gamma-1)p_0 - \gamma u p_x]}, \end{aligned} \tag{2.1.17}$$

and coordinate transformation reads [11]

$$\begin{aligned} t' &= \frac{\gamma(t - ux/c^2)}{1 + (c\varkappa)^{-1}[(\gamma-1)p_0 - \gamma u p_x]}, x' = \frac{\gamma(x - ut)}{1 + (c\varkappa)^{-1}[(\gamma-1)p_0 - \gamma u p_x]}, \\ y' &= \frac{y}{1 + (c\varkappa)^{-1}[(\gamma-1)p_0 - \gamma u p_x]}, z' = \frac{z}{1 + (c\varkappa)^{-1}[(\gamma-1)p_0 - \gamma u p_x]}, \end{aligned} \tag{2.1.18}$$

where $\gamma = \sqrt{1 - u^2/c^2}$. It is easy to check that the energy $E = c\varkappa$ , identified as the high-energy cutoff $\Lambda_*$, is an invariant as it is also the case for the fundamental length $l_{\Lambda_*} = \hbar c/E = \hbar/\varkappa$.

**Remark 2.1.2**. Note that the transformation (2.1.17) defined in $p$-space and the transformation (2.1.18) defined in $x$-space becomes Lorentz for small energies and momenta and defines a large invariant energy $l_{\Lambda_*}^{-1}$. The high-energy cutoff $\Lambda_*$ is preserved by the modified action of the Lorentz group [9],[10].

This meant that the canonical concept of metric as quadratic invariant collapses at high energies, being replaced by the non-quadratic invariant [9]:

$$\|p\|^2 = \frac{\eta^{ab}p_a p_b}{(1 + l_{\Lambda_*}p_0)}, \tag{2.1.19}$$

or by the non-quadratic invariant

$$\|p\|^2 = \frac{\eta^{ab}p_a p_b}{(1 - l_{\Lambda_*}p_0)}, \tag{2.1.20}$$

where $l_{\Lambda_*} = \Lambda_*^{-1}, a, b = 0, 1, 2, 3.$

**Remark 2.1.3**.Note that:

(i) the invariant (2.1.16) is infinite for the new negative invariant energy scale of the theory $\Lambda_* = -l_{\Lambda_*}^{-1}$, and it's not quadratic for energies close or above and

(ii) the invariant (2.1.17) is infinite for the new positive invariant energy scale of the theory $\Lambda_* = l_{\Lambda_*}^{-1}$.

**Remark 2.1.4**.It is also clear from Eq.(2.1.16) and Eq.(2.1.17) that the symmetry of positive and negative values of the energy is broken.The two theories with the two signs of $l_\Lambda$ obviously are physically distinct; and we know of no theoretical argument

which fixes

the sign of $l_\Lambda$

The massive particles have a positive invariant $\|p\|^2 > 0$ which can be identified with the

square of the mass $\|p\|^2 = m^2, (c = 1)$.Thus in the case of the invariant (2.1.19) we

obtain

$$\frac{p_0^2 - p^2}{(1 + l_{\Lambda *}p_0)^2} = m^2, p_0 \in (-l_{\Lambda *}^{-1}, \infty) \qquad (2.1.21)$$

From Eq.(2.1.21) we obtain

$$p_0 = \frac{m^2 l_{\Lambda *}}{1 - m^2 l_{\Lambda *}^2} + \frac{1}{\sqrt{1 - m^2 l_{\Lambda *}^2}} \sqrt{\frac{m^4 l_{\Lambda *}^2}{1 - m^2 l_{\Lambda *}^2} + (p^2 + m^2)}\,. \qquad (2.1.22)$$

In the case of the invariant (2.1.20) we obtain

$$\frac{p_0^2 - p^2}{(1 - l_{\Lambda *}p_0)^2} = m^2, p_0 \in (-\infty, l_{\Lambda *}^{-1}). \qquad (2.1.23)$$

From Eq.(2.1.23) we obtain

$$p_0 = -\frac{m^2 l_{\Lambda *}}{1 - m^2 l_{\Lambda *}^2} - \frac{1}{\sqrt{1 - m^2 l_{\Lambda *}^2}} \sqrt{\frac{m^4 l_{\Lambda *}^2}{1 - m^2 l_{\Lambda *}^2} + (p^2 + m^2)} \qquad (2.1.24)$$

The action for a scalar field $\varphi$ must be invariant under the deformed Lorentz transformations.The invariant action reads [10]

$$S = \frac{1}{2}\int d^4x \left( \frac{\eta^{ab}(\partial_a\varphi)(\partial_b\varphi)}{[1 + l_{\Lambda *}\partial_0\varphi]} + \frac{m^2}{2}\varphi^2 + V(\varphi) \right). \qquad (2.1.25)$$

Thus there is no linear field equation even if $V(\varphi) = 0$.

**Remark 2.1.5**.Throughout this paper, we use the perturbative expansion

$$S = \frac{1}{2}\int d^4x \left( \eta^{ab}(\partial_a\varphi)(\partial_b\varphi) + \frac{m^2}{2}\varphi^2 \right) + O(l_{\Lambda *}). \qquad (2.1.26)$$

and dealing in Lorentz invariant approximation

$$S \simeq \frac{1}{2}\int d^4x \left( \eta^{ab}(\partial_a\varphi)(\partial_b\varphi) + \frac{m^2}{2}\varphi^2 \right). \qquad (2.1.27)$$

since for $l_{\Lambda *} \ll 1$ the expansion (2.1.26) holds.

## 2.2. Zel'dovich approuch by using Pauli-Villars renormalization revisited.What is wrong with Pauli-Villars renormalization. Ghost particles as physical dark matter.

### 2.2.1. Zel'dovich approuch by using Pauli-Villars renormalization revisited.

Remind that vacuum energy density for free scalar quantum field is [1]:

$$\varepsilon(\mu) = \frac{1}{2}\frac{1}{(2\pi\hbar)^3}\int_0^\infty 4\pi c\sqrt{p^2 + \mu^2}\,p^2 dp = K\int_0^\infty \sqrt{p^2 + \mu^2}\,p^2 dp = KI(\mu), \qquad (2.2.1)$$

where $\mu = m_0c$. From the basic definitions [1]

$$p = T_{xx}, p(\mu) = \frac{1}{2}\frac{1}{(2\pi\hbar)^3}\int_0^\infty u_x p_x 4\pi p^2 dp, \mathbf{u} = \frac{c\mathbf{p}}{\sqrt{p^2+\mu^2}}, \overline{u_x p_x} = \frac{1}{3}\langle \mathbf{u}, \mathbf{p}\rangle$$

one obtains

$$p(\mu) = \frac{K}{3}\int_0^\infty \frac{p^4 dp}{\sqrt{p^2+\mu^2}} = KF(\mu). \tag{2.2.2}$$

**Remark 2.2.1.** Note that the integral in RHS of the Eq.(2.2.1) and the Eq.(2.2.2) divergent and ultraviolet cutoff is needed.

Thus in accordance with [1] we set

$$\varepsilon(\mu, p_0) = KI(\mu, p_0), p(\mu, p_0) = KF(\mu, p_0), \tag{2.2.3}$$

where

$$I(\mu, p_0) = \int_0^{p_0}\sqrt{p^2+\mu^2}\, p^2 dp, F(\mu, p_0) = \int_0^{p_0}\frac{p^4 dp}{\sqrt{p^2+\mu^2}}, \tag{2.2.4}$$

where $p_0 c < \Lambda_*$. For fermionic quantum field similarly one obtains [1]

$$\varepsilon(\mu, p_0) = -4KI(\mu, p_0), p(\mu) = -4KF(\mu, p_0). \tag{2.2.5}$$

Thus from Eqs.(2.2.3)-(2.2.5) by using formally Pauli-Villars regularization [7],[8] and regularization by high-energy cutoff the expression for free vacuum energy density $\varepsilon$ reads

$$\varepsilon = \sum_{i=0}^{2M} f_i I(\mu_i, p_0) \tag{2.2.6}$$

and the expression for pressure $p$ reads

$$p = \sum_{i=0}^{2M} f_i F(\mu_i, p_0). \tag{2.2.7}$$

Here $\mu_i$ is a finite positive sequence $\mu_i \in \mathbb{R}_+, i = 1,2,\ldots,2M$ and $f_i$ is a finite sequence $f_i \in \mathbb{R}, f_i \leq 1, i = 1,2,\ldots,2M$.

**Definition 2.2.1.** We define now discrete distribution $f_{PV} : \mathbb{R}_+ \to \mathbb{R}$ by formula

$$f_{PV}(\mu_i) = f_i, \tag{2.2.8}$$

and we will call it as a full discrete Pauli-Villars masses distribution.

**Remark 2.2.2.** We assum now that in Eqs.(2.2.6)-(2.2.7): (i) the quantities $\mu_i^{s.m} = \mu_i, i = 1,2,\ldots,M$ is a masses of a physical particles corresponding to standard matter and (ii) the quantities $\mu_i^{g.m} = \mu_i, i = M+1,2,\ldots,2M$ is a masses of ghost particles with a wrong kinetic term and wrong statistics corresponding to a physical dark matter.

**Remark 2.2.3.** We recall that the Euler-Maclaurin summation formula reads

$$\sum_{i=1}^{2M} g(\mu_1 + (i-1)h) = \int_{\mu_1}^{\mu_{2M}} f(\mu)d\mu + A_1[g(\mu_{2M}) - g(\mu_1)] + A_2 h[g'(\mu_{2M}) - g'(\mu_1)] + O(h^2), \tag{2.2.9}$$
$$f(\mu) = \frac{1}{h}g(\mu)$$

Let $g(\mu)$ be an appropriate continuous function such that: (i) $g(\mu_i) = f_i, i = 1,2,\ldots,2M$,
(ii) $g'(\mu_{2M}) = 0, g'(\mu_1) = 0$.
Thus from Eqs.(2.2.6)-(2.2.7) and Eqs.(2.2.9) we obtain

$$\varepsilon = \sum_{i=0}^{2M} f_i I(\mu_i, p_0) = \int_{\mu_1}^{\mu_{2M}} f(\mu) I(\mu, p_0) d\mu + A_1 h[f(\mu_{2M}) I(\mu_{2M}, p_0) - f(\mu_1) I(\mu_1, p_0)] + O(h^2) \qquad (2.2.10)$$

and

$$p = \sum_{i=0}^{2M} f_i F(\mu_i, p_0) = \int_{\mu_1}^{\mu_{2M}} f(\mu) F(\mu, p_0) d\mu + A_1 h[f(\mu_{2M}) F(\mu_{2M}, p_0) - f(\mu_1) F(\mu_1, p_0)] + O(h^2). \qquad (2.2.11)$$

**Definition 2.2.2.** We will call the function $f(\mu)$ as a full continuous Pauli-Villars masses distribution.

**Definition 2.2.3.** We define now: (i) discrete distribution $f_{PV}^{s.m} : \mathbb{R}_+ \to \mathbb{R}$ by formula

$$f_{PV}^{s.m}(\mu_i^{s.m}) = f_i, i = 1,2,\ldots,M \qquad (2.2.12)$$

and we will call it as discrete Pauli-Villars masses distribution of the standard matter and
(ii) discrete distribution $f_{PV}^{g.m} : \mathbb{R}_+ \to \mathbb{R}$ by formula

$$f_{PV}^{g.m}(\mu_i) = f_i, i = M+1,2,\ldots,2M \qquad (2.2.13)$$

and we will call it as discrete Pauli-Villars masses distribution of the ghost matter.

**Remark 2.2.4.** We rewrite now the Eqs.(2.2.6)-(2.2.7) in the following equivalent form

$$\varepsilon = \sum_{i=1}^{M} f_{PV}^{s.m}(\mu_i^{s.m}) I(\mu_i^{s.m}, p_0) + \sum_{j(i)=M+1}^{2M} f_{PV}^{g.m}\left(\mu_{j(i)}^{g.m}\right) I\left(\mu_{j(i)}^{g.m}, p_0\right) \qquad (2.2.14)$$

and

$$p = \sum_{i=1}^{M} f_{PV}^{s.m}(\mu_i^{s.m}) F(\mu_i^{s.m}, p_0) + \sum_{j(i)=M+1}^{2M} f_{PV}^{g.m}\left(\mu_{j(i)}^{g.m}\right) F\left(\mu_{j(i)}^{g.m}, p_0\right), \qquad (2.2.15)$$

where $j(i) = i + M, i = 1+1,2,\ldots,M$.

**Remark 2.2.5.** We assume now that: (i) $\mu_i^{s.m} \approx \mu_{j(i)}^{g.m}$, (ii) $\left|f_{PV}^{s.m}(\mu_i^{s.m}) + f_{PV}^{g.m}\left(\mu_{j(i)}^{g.m}\right)\right| \ll 1$, i.e.,

$$f_{PV}^{s.m}(\mu_i^{s.m}) \approx -f_{PV}^{g.m}\left(\mu_{j(i)}^{g.m}\right). \qquad (2.2.16)$$

Note that Eq.(2.2.16) meant highly symmetric discrete Pauli-Villars masses distribution
between standard matter and ghost matter below that scale $\Lambda_*$
Thus from Eqs.(2.2.14)-(2.2.15) and Eqs.(2.2.16) we obtain

$$\varepsilon = \sum_{i=1}^{M} f_{PV}^{s.m}(\mu_i^{s.m}) I(\mu_i^{s.m}, p_0) + \sum_{j(i)=M+1}^{2M} f_{PV}^{g.m}\left(\mu_{j(i)}^{g.m}\right) I\left(\mu_{j(i)}^{g.m}, p_0\right) = \sum_{i=1}^{M}\left[f_{PV}^{s.m}(\mu_i^{s.m}) + f_{PV}^{g.m}\left(\mu_{j(i)}^{g.m}\right)\right] I(\mu_i, p_0) \tag{2.2.17}$$

and

$$p = \sum_{i=1}^{M} f_{PV}^{s.m}(\mu_i^{s.m}) F(\mu_i^{s.m}, p_0) + \sum_{j(i)=M+1}^{2M} f_{PV}^{g.m}\left(\mu_{j(i)}^{g.m}\right) F\left(\mu_{j(i)}^{g.m}, p_0\right) = \sum_{i=1}^{M}\left[f_{PV}^{s.m}(\mu_i^{s.m}) + f_{PV}^{g.m}\left(\mu_{j(i)}^{g.m}\right)\right] F(\mu_i, p_0). \tag{2.2.18}$$

From Eqs.(2.2.17)-(2.2.18) and Eqs.(2.2.9) finally we obtain

$$\varepsilon = \sum_{i=1}^{M}\left[f_{PV}^{s.m}(\mu_i^{s.m}) + f_{PV}^{g.m}\left(\mu_{j(i)}^{g.m}\right)\right] I(\mu_i, p_0) = \int_{\mu_1}^{\mu_{\mathbf{eff}}}\left[f_{PV}^{s.m}(\mu) + f_{PV}^{g.m}(\mu)\right] I(\mu, p_0) d\mu \tag{2.2.19}$$

and

$$p = \sum_{i=1}^{M}\left[f_{PV}^{s.m}(\mu_i^{s.m}) + f_{PV}^{g.m}\left(\mu_{j(i)}^{g.m}\right)\right] F(\mu_i^{s.m}) = \int_{\mu_1}^{\mu_{\mathbf{eff}}}\left[f_{PV}^{s.m}(\mu) + f_{PV}^{g.m}(\mu)\right] F(\mu, p_0) d\mu, \tag{2.2.20}$$

where obviously

$$f_{PV}^{s.m}(\mu) + f_{PV}^{g.m}(\mu) = f_{PV}(\mu) \approx 0. \tag{2.2.21}$$

**Definition 2.2.5**.We will call $f_{PV}^{s.m}(\mu)$ and $f_{PV}^{g.m}(\mu)$ as continuous Pauli-Villars masses distribution of the standard matter and ghost matter correspondingly.
Thus finally we obtain

$$\varepsilon(\mu_{\mathbf{eff}}, p_0) = \int_{\mu_1}^{\mu_{\mathbf{eff}}} f_{PV}(\mu) I(\mu, p_0) d\mu, \tag{2.2.22}$$

and

$$p(\mu_{\mathbf{eff}}, p_0) = \int_{\mu_1}^{\mu_{\mathbf{eff}}} f_{PV}(\mu) F(\mu, p_0) d\mu, \tag{2.2.23}$$

In order to calculate $\varepsilon(\mu_{\mathbf{eff}}, p_0)$ and $p(\mu_{\mathbf{eff}}, p_0)$ let us evaluate now the following quantities defined above by Eqs.(2.2.4)

$$I(\mu, p_0) = \int_0^{p_0} p^2\sqrt{p^2+\mu^2}\,dp = \int_0^{p_\mu} p^2\sqrt{p^2+\mu^2}\,dp + \int_{p_\mu}^{p_0} p^2\sqrt{p^2+\mu^2}\,dp = \int_0^{p_\mu} p^3\sqrt{1+\frac{\mu^2}{p^2}}\,dp + \int_{p_\mu}^{p_0} p^3\sqrt{1+\frac{\mu^2}{p^2}}\,dp \tag{2.2.24}$$

and

$$F(\mu,p_0) = \frac{1}{3}\int_0^{p_0} \frac{p^4dp}{\sqrt{p^2+\mu^2}} = \frac{1}{3}\int_0^{p_\mu} \frac{p^4dp}{\sqrt{p^2+\mu^2}} + \frac{1}{3}\int_{p_\mu}^{p_0} \frac{p^4dp}{\sqrt{p^2+\mu^2}} = \frac{1}{3}\int_0^{p_\mu} \frac{p^3dp}{\sqrt{1+\frac{\mu^2}{p^2}}} + \frac{1}{3}\int_{p_\mu}^{p_0} \frac{p^3dp}{\sqrt{1+\frac{\mu^2}{p^2}}}, \quad (2.2.25)$$

where $p_\mu = r\mu, r > 1, \mu/p < 1/r < 1$. Note that:

$$\sqrt{1+\frac{\mu^2}{p^2}} = 1+\frac{1}{2}\frac{\mu^2}{p^2}-\frac{1}{8}\frac{\mu^4}{p^4}+\frac{1}{16}\frac{\mu^6}{p^6}+\ldots.$$
$$p^2\sqrt{p^2+\mu^2} = p^3\sqrt{1+\frac{\mu^2}{p^2}} = p^3+\frac{1}{2}\mu^2p-\frac{1}{8}\frac{\mu^4}{p}+\frac{1}{16}\frac{\mu^6}{p^3}+\ldots. \quad (2.2.26)$$

By inserting Eq.(2.2.26) into Eq.(2.2.24) one obtains

$$I(\mu,p_0) = C_1\mu^4+\frac{1}{4}p_0^4+\frac{1}{4}\mu^2p_0^2-\frac{1}{8}\mu^4\ln\left(\frac{p_0}{\mu}\right)-\frac{1}{32}\frac{\mu^6}{p_0^2}+p_0^{-5}O(\mu^8), \quad (2.2.27)$$

where $C_1\mu^4 = \int_0^{p_\mu} p^2\sqrt{p^2+\mu^2}\,dp$. Note that:

$$\left(\sqrt{1+\frac{\mu^2}{p^2}}\right)^{-1} = 1-\frac{1}{2}\frac{\mu^2}{p^2}+\frac{3}{8}\frac{\mu^4}{p^4}-\frac{5}{16}\frac{\mu^6}{p^6}+\ldots. \quad (2.2.28)$$

By inserting Eq.(2.2.28) into Eq.(2.2.25) one obtains

$$F(\mu,p_0) = C_2\mu^4+\frac{1}{12}p_0^4-\frac{1}{12}\mu^2p_0^2+\frac{1}{8}\mu^4\ln\left(\frac{p_0}{\mu}\right)+\frac{5}{32}\frac{\mu^6}{p_0^2}+p_0^{-5}O(\mu^8). \quad (2.2.29)$$

By inserting Eq.(2.2.27) into Eq.(2.2.22) one obtains

$$\varepsilon(\mu_{\mathbf{eff}},p_0) = \frac{1}{4}p_0^4\int_0^{\mu_{\mathbf{eff}}} f_{PV}(\mu)d\mu+\frac{1}{4}p_0^2\int_0^{\mu_{\mathbf{eff}}} f_{PV}(\mu)\mu^2d\mu+\left(C_1-\frac{1}{8}\ln p_0\right)\int_0^{\mu_{\mathbf{eff}}} f_{PV}(\mu)\mu^4d\mu+ +\frac{1}{8}\int_0^{\mu_{\mathbf{eff}}} f_{PV}(\mu)\mu^4(\ln\mu)d\mu-\left(\frac{1}{p_0^2}\right)\frac{1}{32}\int_0^{\mu_{\mathbf{eff}}} f_{PV}(\mu)\mu^6d\mu+O\left(\int_0^{\mu_{\mathbf{eff}}} f_{PV}(\mu)\mu^8\right)p_0^{-5}. \quad (2.2.30)$$

By inserting Eq.(2.2.29) into Eq.(2.2.23) one obtains

$$p(\mu_{\mathbf{eff}},p_0) = \frac{1}{12}p_0^4\int_0^{\mu_{\mathbf{eff}}} f_{PV}(\mu)d\mu-\frac{1}{12}p_0^2\int_0^{\mu_{\mathbf{eff}}} f_{PV}(\mu)\mu^2d\mu+\left(C_2+\frac{1}{8}\ln p_0\right)\int_0^{\mu_{\mathbf{eff}}} f_{PV}(\mu)\mu^4d\mu- -\frac{1}{8}\int_0^{\mu_{\mathbf{eff}}} f_{PV}(\mu)\mu^4(\ln\mu)d\mu+\left(\frac{5}{p_0^2}\right)\frac{1}{32}\int_0^{\mu_{\mathbf{eff}}} f_{PV}(\mu)\mu^6d\mu+O\left(\int_0^{\mu_{\mathbf{eff}}} f_{PV}(\mu)\mu^8\right)p_0^{-5}. \quad (2.2.31)$$

We formally choose now [1]

$$\int_0^{\mu_{\mathbf{eff}}} f_{PV}(\mu)d\mu = \int_0^{\mu_{\mathbf{eff}}} f_{PV}(\mu)\mu^2 d\mu = \int_0^{\mu_{\mathbf{eff}}} f_{PV}(\mu)\mu^4 d\mu = 0. \quad (2.2.32)$$

By inserting Eq.(2.2.32) into Eqs.(2.2.30)-(2.2.31) one obtains

$$\begin{aligned} \varepsilon(\mu_{\mathbf{eff}}, p_0) &= \frac{1}{8}\int_0^{\mu_{\mathbf{eff}}} f_{PV}(\mu)\mu^4(\ln\mu)d\mu + O(p_0^{-2}), \\ p(\mu_{\mathbf{eff}}, p_0) &= -\frac{1}{8}\int_0^{\mu_{\mathbf{eff}}} f_{PV}(\mu)\mu^4(\ln\mu)d\mu + O(p_0^{-2}). \end{aligned} \quad (2.2.33)$$

Taking the limit $p_0 \to \infty$ in Eq.(2.2.33) gives

$$\begin{aligned} \varepsilon(\mu_{\mathbf{eff}}) &= \frac{1}{8}\int_0^{\mu_{\mathbf{eff}}} f_{PV}(\mu)\mu^4(\ln\mu)d\mu, \\ p(\mu_{\mathbf{eff}}) &= -\frac{1}{8}\int_0^{\mu_{\mathbf{eff}}} f_{PV}(\mu)\mu^4(\ln\mu)d\mu. \end{aligned} \quad (2.2.34)$$

Thus finally we obtain well known Zel'dovich formulas [1]:

$$\varepsilon(\mu_{\mathbf{eff}}) = -p(\mu_{\mathbf{eff}}) = \frac{1}{8}\int_0^{\mu_{\mathbf{eff}}} f_{PV}(\mu)\mu^4(\ln\mu)d\mu = \frac{c^4\lambda}{8\pi G}. \quad (2.2.35)$$

**Remark 2.2.6**. If we assume that

$$f_{PV}^{s.m}(\mu) + f_{PV}^{g.m}(\mu) = f_{PV}(\mu) \equiv 0 \quad (2.2.36)$$

instead Eq.(2.2.21) we obtain zero value of the cosmological constant $\lambda$. In this paper a small value of the cosmological constant explained by tiny violation of the simmetry required by Eq.(2.2.16).

## 2.2.2 Pauli-Villars renormalization.What is wrong with canonical Pauli-Villars renormalization.Ghost particles as physical dark matter.

**Remark 2.2.7**. Remind that canonical Pauli-Villars regularization consists of introducing a fictitious mass term. For example, we would replace a propagator $1/\left(k^2 - m_0^2 + i\epsilon\right)$, by the regulated propagator

$$\Delta(k^2) = \sum_{i=0}^N \frac{a_i}{k^2 - m_i^2 + i\epsilon} = \frac{1}{k^2 - m_0^2 + i\epsilon} - \sum_{i=1}^N \frac{a_i}{k^2 - m_i^2 + i\epsilon}, \quad (2.2.37)$$

where $a_0 = 1$ and $m_i, i = 1, 2, \ldots N$ can be thought of as the mass of a fictitious heavy particle, whose contribution is subtracted from that of an ordinary particle. Assume that $m_i^2/k^2 < 1$, if we expand each term of this sum (2.2.37) as a power series in $k^2 + i\epsilon$, i.e.,

$$\frac{a_i}{k^2 - m_i^2 + i\epsilon} = \frac{a_i}{k^2 + i\epsilon}\frac{1}{1 - \frac{m_i^2}{k^2+i\epsilon}} = \frac{a_i}{k^2 + i\epsilon} \times$$
$$\left[1 + \frac{m_i^2}{k^2 + i\epsilon} + \frac{m_i^4}{\left(k^2 + i\epsilon\right)^2} + \ldots\right] = \frac{a_i}{k^2 + i\epsilon} + \frac{a_i m_i^2}{\left(k^2 + i\epsilon\right)^2} + \frac{a_i m_i^4}{\left(k^2 + i\epsilon\right)^3} + \ldots, \qquad (2.2.38)$$

where $m_i^2/k^2 < 1$, we get

$$\Delta(k^2) = \sum_{i=0}^{N} \frac{a_i}{k^2 + i\epsilon} + \sum_{i=0}^{N} \frac{a_i m_i^2}{\left(k^2 + i\epsilon\right)^2} + \sum_{i=0}^{N} O\left(\frac{1}{\left(k^2 + i\epsilon\right)^3}\right). \qquad (2.2.39)$$

For a renormalizable theory the maximum supercriticial power of divergence of any Feinman integral is quadratic, so that the $O(1/k^6)$ terms are ultraviolet finite. The finiteness of the regulated Feinman integral is then guaranteed by requiring that

$$\sum_{i=0}^{N} a_i = 0, \sum_{i=0}^{N} a_i m_i^2 = 0. \qquad (2.2.40)$$

**Remark 2.2.8**. Note that in order to obtain renormalized physical quantities canonical procedure required take the limits $\lim_{m_i \to \infty} m_i, i = 1,\ldots,N$ of the regulated Feinman integral.

Unfortunately under these limits the expansion (2.2.38) obviously breaks down since the inequalities $m_i^2/k^2 < 1, i = 1,\ldots,N$ in the limits $\lim_{m_i \to \infty} m_i, i = 1,\ldots,N$ obviously no longer holds!Thus canonical Pauli-Villars procedure does not make any rigorous mathematical sense.In fact that is formal deletion of the divergences by hands.

**Remark 2.2.9**. In order to avoid these difficultness mentioned above we have choose Pauli-Villars masses $m_{PV}$ of order $m_{PV} \asymp \Lambda_*/c^2$ or $m \gg \Lambda_*/c^2$ and therefore there is no problems arises with unitarity condition,see section V.

**Remark 2.2.10**.We claim that such sufficiently larges Pauli-Villars masses $m_{PV}$ corresponds to a physical ghost particles which formed Dark Matter sector of the Universe.

**Remark 2.2.11**.Note that in order to aply modified Pauli-Villars renormalization mentioned above (see Remarks 2.2.9-2.2.10) to QFT with Lagrangian $\mathcal{L}(\varphi,\psi,\partial_\mu\varphi,\partial_\mu\psi)$ we would replace the Lagrangian $\mathcal{L}(\varphi,\psi,\partial_\mu\varphi,\partial_\mu\psi)$ by Lagrangian $\underline{\mathcal{L}}\left(\underline{\varphi},\underline{\psi},\partial_\mu\underline{\varphi},\partial_\mu\underline{\psi}\right)$,where [7]:

$$\underline{\varphi}(x) = \varphi(x) + \sum_n b_n\breve{\varphi}_n(x,\mu_{PV,n}^2), \underline{\psi}(x) = \psi(x) + \sum_n c_n\breve{\psi}_n(x,\varkappa_{PV,n}^2), \qquad (2.2.41)$$

and where commutator for $\breve{\varphi}_n$ and anticommutator for $\breve{\psi}_n$ reads

$$\begin{gathered}[\breve{\varphi}_m(x,\mu_{PV,m}^2),\breve{\varphi}_n(x',\mu_{PV,n}^2)] = -i\rho_n\Delta(x - x',\mu_{PV,n}^2)\delta_{mn},\\ \{\breve{\psi}_m(x,\varkappa_{PV,m}^2),\breve{\psi}_n(x',\varkappa_{PV,n}^2)\} = -i\varepsilon_n S(x - x',\varkappa_{PV,n}^2)\delta_{mn}.\end{gathered} \qquad (2.2.42)$$

From the Eqs.(2.2.41)-Eqs.(2.2.42) one obtains

$$\begin{gathered}\left[\underline{\varphi}(x),\underline{\varphi}(x')\right] = i\sum_{n=0}^{N} \rho_n b_n^2\Delta(x - x',\mu_{PV,n}^2),\\ \left[\underline{\psi}(x),\underline{\psi}(x')\right] = -i\sum_{n=0}^{N} \varepsilon_n\bar{c}_n c_n S(x - x',\varkappa_{PV,n}^2).\end{gathered} \qquad (2.2.43)$$

Assume now that

$$\sum_{n=0}^{N}\rho_n b_n^2 = 0, \sum_{n=0}^{N}\rho_n b_n^2 \mu_{PV,n}^2 = 0, \sum_{n=0}^{N}\varepsilon_n \bar{c}_n c_n = 0, \sum_{n=0}^{N}\varepsilon_n \bar{c}_n c_n \varkappa_{PV,n}^2 = 0. \quad (2.2.44)$$

From Eqs.(2.2.44) it follows directly that QFT with Lagrangian $\underline{\mathcal{L}}\left(\underline{\varphi}, \underline{\psi}, \partial_\mu \underline{\varphi}, \partial_\mu \underline{\psi}\right)$ is finite QFT with indefinite metric [7]. It meant that a "bad ghosts" with Pauli-Villars masses $\mu_{PV,n}$, $\varkappa_{PV,n}$ appears in Lagrangian as real quantum fields corresponding to real ghost particles.

**Remark 2.2.12**.Thus we argue that UV divergence in fact arises from physically wrong Lagrangian $\mathcal{L}(\varphi, \psi, \partial_\mu \varphi, \partial_\mu \psi)$ in which Pauli-Villars ghosts $\breve{\varphi}$ are missing.

**Remark 2.2.13**.Note that "bad ghosts" represent general meaning of the word "ghost" in theoretical physics: states of negative norm [7] or fields with the wrong sign of the kinetic term, such as Pauli–Villars ghosts $\breve{\varphi}$, whose existence allows the probabilities to be negative thus violating unitarity.The quadratic lagrangian $\mathcal{L}_{\breve{\varphi}}^2$ for $\breve{\varphi}$ begins with a wrong sign kinetic term [in $(+---)$ signature]

$$\mathcal{L}_{\breve{\varphi}}^2 = -\frac{1}{2}\partial^\mu \breve{\varphi} \partial_\mu \breve{\varphi} + \ldots \quad (2.2.45)$$

**Remark 2.2.14**.Note that in order to obtain Eqs.(2.2.34), the standard quantum fields do not need to couple alwais directly to the ghost sector. In this paper the ghost sector is considered as real physical mechanism which acts on a function $f_{PV}(\mu)$ in Eqs.(2.2.34).It means that there exist the ghost-driven acceleration of the univers hidden in cosmological constant $\lambda$.

**Remark 2.2.11**.As pointed out in paper [12] even if the standard model fields have no direct couplings to the ghost sector, they will indirectly interact with it through gravity, and the propagation of gravity through the ghost condensate gives rise to a fascinating modification of gravity in the IR. However,no modifications of gravity can be seen directly, and no cosmological experiment can distinguish the ghost-driven acceleration from a cosmological constant.

**Remark 2.2.12**.In order to obtain desired physical result from Eqs.(2.2.35),i.e.,

$$\varepsilon_{\mathbf{vac}} = 0.7 \times 10^{-29} gcm^{-3} = 2.8 \times 10^{-47} Gev^4/\hbar^3 c^5 \quad (2.2.45)$$

we assume that

$$f_{PV}(\mu) = f_{PV}^{s.m}(\mu) + f_{PV}^{g.m.}(\mu) \approx 0, \quad (2.2.46)$$

where $f_{s.m.}(\mu)$ corresponds to standard matter and where $f_{g.m.}(\mu)$ corresponds to a physical ghost matter. Obviously Eq.(2.2.46) required tiny violation of the symmetry between standard matter and ghost matter.

**Remark 2.2.13**.In additional we assume now that

$$|f_{PV}(\mu)| = \begin{cases} O\left(\mu_{\mathbf{eff}}^{-n}\right), n > 1 & m_0 \leq \mu \leq \mu_{\mathbf{eff}} \\ 0 & \mu > \mu_{\mathbf{eff}} \end{cases} \tag{2.2.47}$$

where $\mu_{\mathbf{eff}} \ll \Lambda_*/c = p_0$.

For definiteness we have chosen now that

$$f_{PV}(\mu) = \begin{cases} O\left(\mu_{\mathbf{eff}}^{-n}\right)\cos\mu, n > 1 & m_0 \leq \mu \leq \mu_{\mathbf{eff}} \\ 0 & \mu > \mu_{\mathbf{eff}} \end{cases} \tag{2.2.48}$$

where $n \in \mathbb{N}$ and

$$p_0^4 \sin(\gamma\mu_{\mathbf{eff}}) \ll p_0^2, \cos\gamma\mu_{\mathbf{eff}} \simeq 1, \tag{2.2.49}$$

where $\gamma \in \mathbb{R}$. Note that

$$\frac{1}{4}p_0^4 \int_0^{\mu_{\mathbf{eff}}} f_{PV}(\mu)d\mu = \frac{\mu_{\mathbf{eff}}^{-n}}{4\gamma} p_0^4 \sin(\gamma\mu_{\mathbf{eff}}) \ll \frac{\mu_{\mathbf{eff}}^{-n}}{4\gamma} p_0^2,$$

$$\frac{1}{4}p_0^2 \int_0^{\mu_{\mathbf{eff}}} f_{PV}(\mu)\mu^2 d\mu =$$

$$\frac{1}{4\gamma^3} p_0^2[(\gamma^2\mu^2\sin(\gamma\mu) - 2\sin(\gamma\mu) + 2\gamma\mu\cos(\gamma\mu))|_0^{\mu_{\mathbf{eff}}}] \simeq \frac{\mu_{\mathbf{eff}}^{-n+1}}{2\gamma^2} p_0^2, \tag{2.2.50}$$

$$\ln p_0 \int_0^{\mu_{\mathbf{eff}}} f_{PV}(\mu)\mu^4 d\mu = \frac{\ln p_0}{\gamma^5} \times$$

$$[(24\sin(\gamma\mu) + 4\gamma^3\mu^3\cos(\gamma\mu) - 12\gamma^2\mu^2\sin(\gamma\mu) + \gamma^4\mu^4\sin(\gamma\mu) - 24\gamma\mu\cos(\gamma\mu))|_0^{\mu_{\mathbf{eff}}}]$$

$$\simeq \frac{4\mu_{\mathbf{eff}}^{-n+3}}{\gamma^2} \ln p_0.$$

From Eqs.(2.2.30)-(2.2.31) and Eqs.(2.2.50) one obtains

$$\varepsilon(\mu_{\mathbf{eff}}, p_0) \simeq -\frac{\mu_{\mathbf{eff}}^{-n+1}}{2\gamma^2} p_0^2 + \frac{4\mu_{\mathbf{eff}}^{-n+3}}{\gamma^2} \ln p_0 \tag{2.2.50}$$

and

$$p(\mu_{\mathbf{eff}}, p_0) \simeq \frac{\mu_{\mathbf{eff}}^{-n+1}}{6\gamma^2} p_0^2 - \frac{4\mu_{\mathbf{eff}}^{-n+3}}{\gamma^2} \ln p_0. \tag{2.2.51}$$

**Remark 2.2.14.** Note that from Eqs.(2.2.50)-(2.2.51) it follows

$$3p(\mu_{\mathbf{eff}}, p_0) + \varepsilon(\mu_{\mathbf{eff}}, p_0) \simeq -\frac{8\mu_{\mathbf{eff}}^{-n+3}}{\gamma^2} \ln p_0, \tag{2.2.52}$$

and according to Friedman's equations (2.1.6), it corresponds to an accelerated expansion
of the cosmological medium.

**Remark 2.2.15.** Note that from Eq.(2.2.47) and Eqs.(2.2.33) it follows directly that

$$|p(\mu_{\mathbf{eff}})| = |\varepsilon(\mu_{\mathbf{eff}})| \simeq \frac{1}{8}\left|\int\limits_0^{\mu_{\mathbf{eff}}} f(\mu)\mu^4(\ln\mu)d\mu\right| \leq O\left(\mu_{\mathbf{eff}}^{-n+5}\ln\mu_{\mathbf{eff}}\right). \quad (2.2.53)$$

**Remark 2.2.16**.However the Eqs.(2.2.32) required in Zel'dovich paper [1] that is mathematical abstraction, which arises from Pauli-Villars renormalization procedure and nothing more. The Eqs.(2.2.32) in fact demand the unphysical fine tuning of the mass distribution $f_{PV}(\mu)$,i.e.,

$$p_0^4 \int\limits_0^{\mu_{\mathbf{eff}}} f_{PV}(\mu)d\mu \equiv 0, p_0^2 \int\limits_0^{\mu_{\mathbf{eff}}} f_{PV}(\mu)\mu^2 d\mu \equiv 0, \int\limits_0^{\mu_{\mathbf{eff}}} f_{PV}(\mu)\mu^4 d\mu \equiv 0. \quad (2.2.54)$$

for any $p_0$. Obviously the Eqs.(2.2.54) unstable relative to arbitrarily small perturbations of the mass distribution $f_{PV}(\mu)$.

**Remark 2.2.17**.Note that:

(i) any Lorentz invariant theory of elementary particles breaks down under physical conditions of the stability of the Eqs.(2.2.54) relative to arbitrarily small perturbations of the mass distribution $f_{PV}(\mu)$ because for Lorentz invariant theory the limit $p_*^2 \to \infty$ is required (ii) in order to obtain physical model of the cosmological constant without fine tuning fundamental high-energy cutoff $\Lambda_*$is required.

**Remark 2.2.18**.In order to avoid these difficultnes mentioned above we assume that

(i) physics of elementary particles is separated into low/high energy ones,

(ii) the standard notion of smooth spacetime is assumed to be altered at a high energy cutoff scale and a new treatment based on QFT in a fractal spacetime with negative dimension is used above that scale

(iii) instead Eqs.(2.2.54) we assume now that

$$p_0^4 \int\limits_0^{\mu_{\mathbf{eff}}(p_0)} f_{PV}(\mu)d\mu \approx 0, p_*^2 \int\limits_0^{\mu_{\mathbf{eff}}(p_0)} f_{PV}(\mu)\mu^2 d\mu \approx 0, \int\limits_0^{\mu_{\mathbf{eff}}(p_0)} f_{PV}(\mu)\mu^4 d\mu \approx 0. \quad (2.2.54)$$

required by Eq.(2.2.16).

## 2.2.3. Problems from non-renormalizability of canonical quantum gravity with Einstein-Hilbert action.

**Remark 2.2.19**.However serious problems arises from non-renormalizability of canonical quantum gravity with Einstein-Hilbert action

$$S_{EH} = \frac{1}{16\pi G}\int d^4x\sqrt{-g}R. \quad (2.2.55)$$

For example taking $\Lambda_*^3$ particles of energy a per unit volume gives the gravitational self-energy density of order $\Lambda_*^6$,i.e.,the density $\varepsilon_{\Lambda_*}$ diverges as $\Lambda_*^6$

$$\varepsilon_{\Lambda_*} \simeq G\Lambda_*^6, \quad (2.2.56)$$

where $\Lambda_*$ is a high-energy cutoff [5].

In order to avoid these difficulties we apply instead Einstein-Hilbert action (2.2.55) the gravitational action which include terms quadratic in the curvature tensor

$$\Im = -\int d^4x\sqrt{-g}\,(\alpha R_{\mu\nu}R^{\mu\nu} - \beta R^2 + 2\kappa^{-2}R), \tag{2.2.57}$$

**Remark 2.2.20**.Gravitational actions (2.2.57) which include terms quadratic in the curvature tensor are renormalizable [13]. The requirement that the graviton propagator behave like $p^{-4}$ for large momenta makes it necessary to choose the indefinite-metric vector space over the negative-energy states.These negative-norm states cannot be excluded from the physical sector of the vector space without destroying the unitarity of the $\mathbf{S}$ matrix, however, for their unphysical behavior may be restricted to arbitrarily large energy scales $\Lambda_*$ by an appropriate limitation on the renormalized masses $m_2$ and $m_0$.

**Remark 2.2.21**.We assum that $m_0c \gg \mu_{\mathbf{eff}}, m_2c \gg \mu_{\mathbf{eff}}$.

**Remark 2.2.22**.The canonical Quantum Field Theory is widely believed to break down at some fundamental high-energy cutoff $\Lambda_*$ and therefore the quantum fluctuations in the vacuum can be treated classically seriously only up to this high-energy cutoff, see for example [14]. In this paper we argue that Quantum Field Theory in fractal space-time with negative Hausdorff-Colombeau dimensions [15] gives high-energy cutoff on natural way.

## 2.3.Dark matter nature. A common origin of the dark energy and dark matter phenomena.

Dark matter is a hypothetical form of matter that is thought to account for approximately 85% of the matter in the universe, and about a quarter of its total energy density. The majority of dark matter is thought to be non-baryonic in nature, possibly being composed of some as-yet undiscovered subatomic particles.Its presence is implied in a variety of astrophysical observations, including gravitational effects that cannot be explained unless more matter is present than can be seen. For this reason, most experts think dark matter to be ubiquitous in the universe and to have had a strong influence on its structure and evolution. Dark matter is called dark because it does not appear to interact with observable electromagnetic radiation, such as light, and is thus invisible to the entire electromagnetic spectrum, making it extremely difficult to detect using usual astronomical equipment

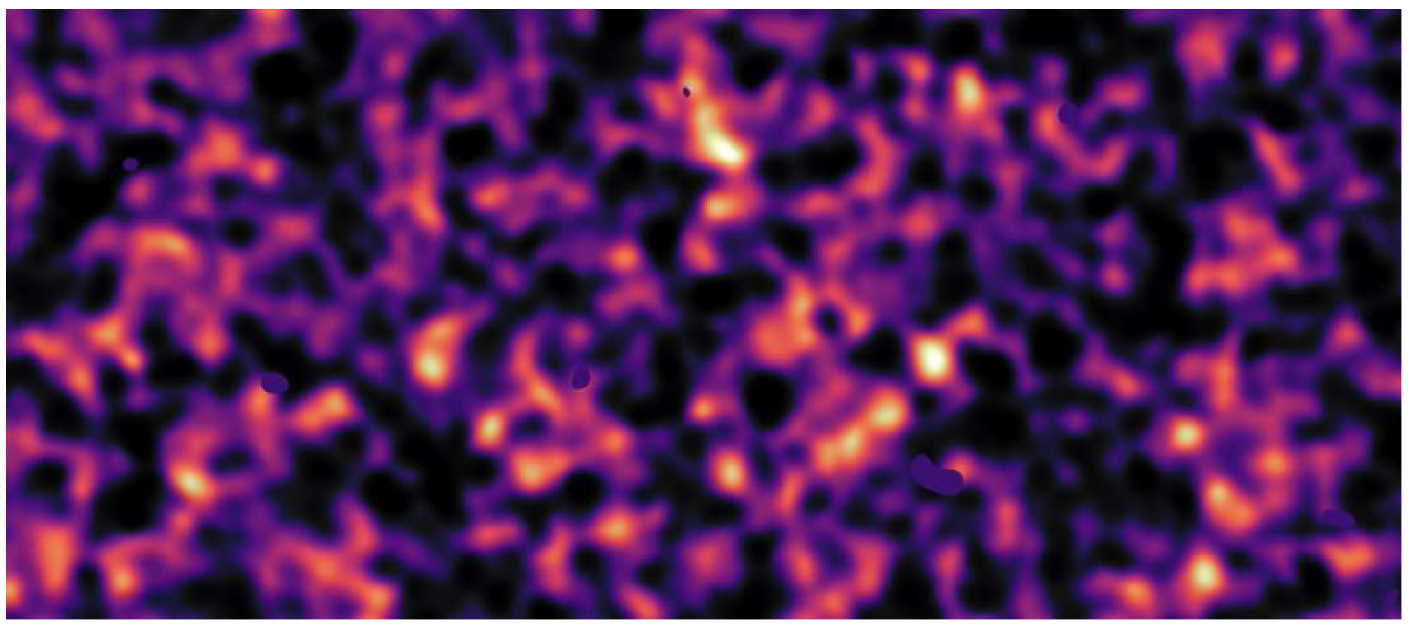

Fig.2.3.1.Dark matter map for a patch of sky
based on gravitational lensing analysis [18].

Fig.2.3.1.Analysis of a giant new galaxy survey, made with ESO’s VLT Survey Telescope in Chile, suggests that dark matter may be less dense and more smoothly distributed throughout space than previously thought. An international team used data from the Kilo Degree Survey (KiDS) to study how the light from about 15 million distant galaxies was affected by the gravitational influence of matter on the largest scales in the Universe. The results appear to be in disagreement with earlier results from the Planck satellite.

This map of dark matter in the Universe was obtained from data from the KiDS survey, using the VLT Survey Telescope at ESO’s Paranal Observatory in Chile. It reveals an expansive web of dense (light) and empty (dark) regions. This image is one out of five patches of the sky observed by KiDS. Here the invisible dark matter is seen rendered in pink, covering an area of sky around 420 times the size of the full moon. This image reconstruction was made by analysing the light collected from over three million distant galaxies more than 6 billion light-years away. The observed galaxy images were warped by the gravitational pull of dark matter as the light travelled through the Universe. Some small dark regions, with sharp boundaries, appear in this image. They are the locations of bright stars and other nearby objects that get in the way of the observations of more distant galaxies and are hence masked out in these maps as no weak-lensing signal can be measured in these areas [16-17].

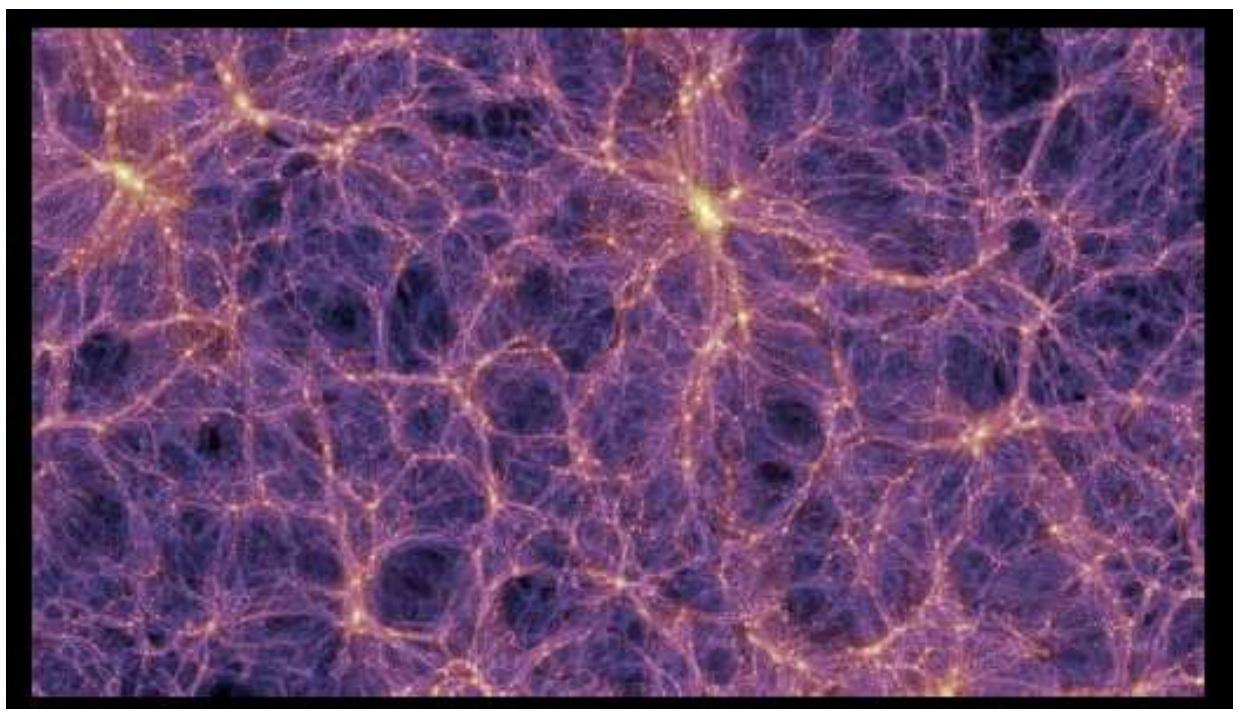

Fig.2.3.2.A simulation of the dark matter distribution
in the universe13.6 billion years ago.

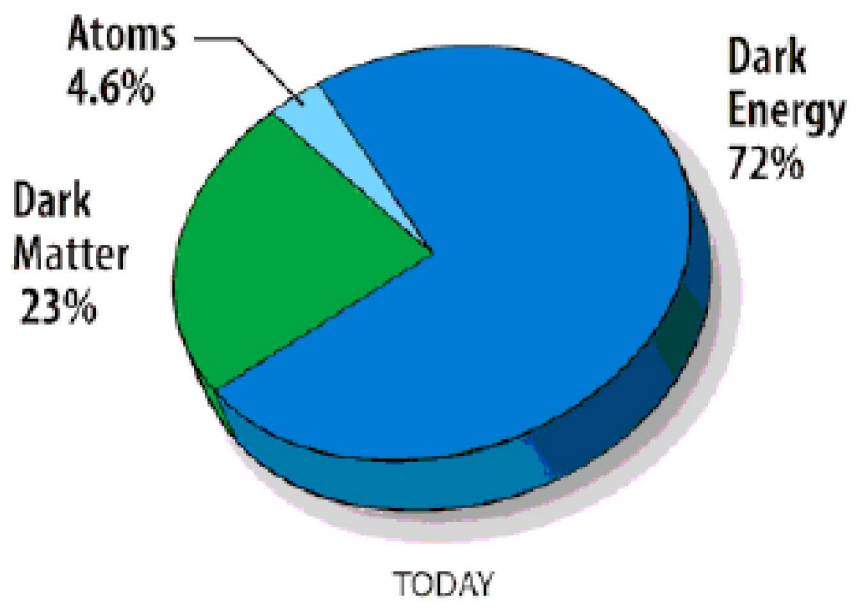

Fig.2.3.3.Matter and energy distribution in
the universe today.The luminous (light-emitting)
components of the universe only comprise
about 0.4% of the total energy.
The remaining components are dark.

The luminous (light-emitting) components of the universe only comprise about 0.4% of the total energy. The remaining components are dark. Of those, roughly 3.6% are identified: cold gas and dust, neutrinos, and black holes. About 23% is dark matter, and the overwhelming majority is some type of gravitationally self-repulsive dark energy.

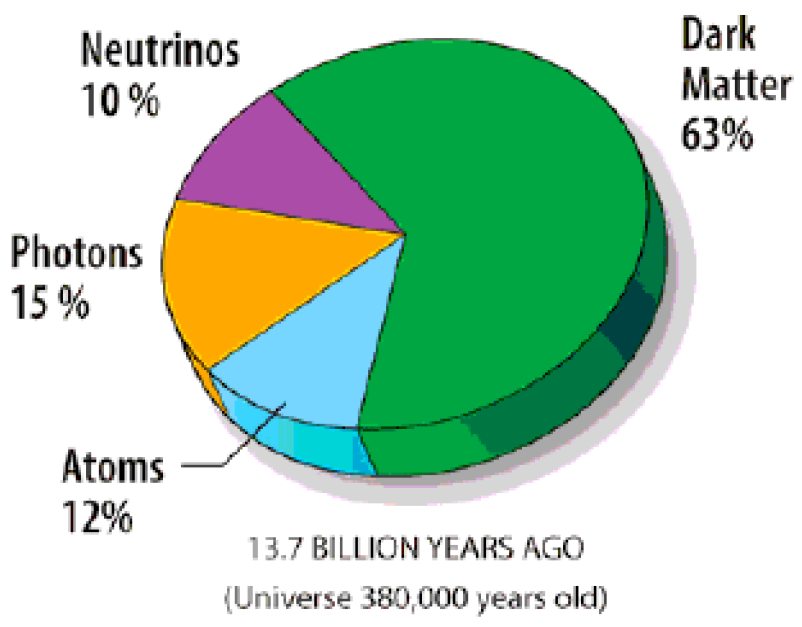

Fig.2.3.4.Matter and energy distribution in
the universe 13.7 bullion years ago.

There is no candidate in the standard model of particle physics.In what way does dark matter extend the standard model?

**Remark 2.3.1.**In order to explain physical nature of dark matter sector we assume that main part of dark matter,i.e., $\simeq 23\% - 4.6\% = 18\%$ (see Fig.2.3.3) formed by supermassive ghost particles vith masess such that $mc^2 > \Lambda_*$.

**Remark 2.3.2.**In order to obtain QFT description of the dark component of matter in natural way we expand now the standard model of particle physics on a sector of ghost
particles. QFT in a ghost sector developed in Sect.3.1-3.4 and Sect.4.1-4.8.

## 2.3.1.The Standard Model of fundamental interactions in standard matter sector below fundamental high-energy cutoff $\Lambda_*$.

Let's remind that in the Standard Model (SM) of fundamental interactions besides the gauge interactions and the quartic interaction of the Higgs fields there are also Yukawa type interactions of the fermion fields with the Higgs field. These interactions are also renormalizable and is characterized by the Yukawa coupling constants, one for each fermion field. The peculiarity of the SM is that the masses of the fields appear as a result of spontaneous symmetry breaking when the Higgs field develops a vacuum expectation value. As a result the masses are not independent but are expressed via the coupling constant multiplied by the vacuum expectation value. Another property of the Standard Model is that it has the gauge group $SU_c(3) \times SU_L(2) \times U_Y(1)$ which is spontaneously broken to $SU_c(3) \times U_{EM}(1)$. In the theories with spontaneously broken symmetry, according to the Goldstone theorem there are massless particles, the goldstone bosons. The Lagrangian of the Standard Model in standard matter sector consists of the following three parts:

$$\mathcal{L}^{s.m.} = \mathcal{L}^{s.m.}_{gauge} + \mathcal{L}^{s.m.}_{Yukawa} + \mathcal{L}^{s.m.}_{Higgs}, \tag{2.3.1}$$

The gauge part is totally fixed by the requirement of the gauge invariance leaving only the values of the couplings as free parameters

$$\begin{aligned} \mathcal{L}^{s.m.}_{gauge} &= -\frac{1}{4}G^a_{\mu\nu}G^a_{\mu\nu} - \frac{1}{4}W^i_{\mu\nu}W^i_{\mu\nu} - \frac{1}{4}B_{\mu\nu}B_{\mu\nu} \\ &+ i\overline{L}_\alpha\gamma^\mu D_\mu L_\alpha + i\overline{Q}_\alpha\gamma^\mu D_\mu Q_\alpha + i\overline{E}_\alpha\gamma^\mu D_\mu E_\alpha \\ &+ i\overline{U}_\alpha\gamma^\mu D_\mu U_\alpha + i\overline{D}_\alpha\gamma^\mu D_\mu D_\alpha + (D_\mu H)^\dagger(D_\mu H), \end{aligned} \tag{2.3.2}$$

where the following notation for the covariant derivatives is used

$$\begin{aligned} G^a_{\mu\nu} &= \partial_\mu G^a_\nu - \partial_\nu G^a_\mu + g_s f^{abc}G^b_\mu G^c_\nu, W^i_{\mu\nu} = \partial_\mu W^i_\nu - \partial_\nu W^i_\mu + g\epsilon^{ijk}W^j_\mu W^k_\nu, \\ B_{\mu\nu} &= \partial_\mu B_\nu - \partial_\nu B_\mu, D_\mu L_\alpha = \left(\partial_\mu - i\frac{g}{2}\tau^i W^i_\mu + i\frac{g'}{2}B_\mu\right)L_\alpha, \\ D_\mu E_\alpha &= (\partial_\mu + ig'B_\mu)E_\alpha, D_\mu Q_\alpha = \left(\partial_\mu - i\frac{g}{2}\tau^i W^i_\mu - i\frac{g'}{6}B_\mu - i\frac{g_s}{2}\lambda^a G^a_\mu\right)Q_\alpha, \\ D_\mu U_\alpha &= \left(\partial_\mu - i\frac{2}{3}g'B_\mu - i\frac{g_s}{2}\lambda^a G^a_\mu\right)U_\alpha, D_\mu D_\alpha = \left(\partial_\mu + i\frac{1}{3}g'B_\mu - i\frac{g_s}{2}\lambda^a G^a_\mu\right)D_\alpha. \end{aligned} \tag{2.3.3}$$

The Yukawa part of the Lagrangian which is needed for the generation of the quark and lepton masses is also chosen in the gauge invariant form and contains arbitrary Yukawa couplings (we ignore the neutrino masses, for simplicity)

$$\mathcal{L}^{s.m.}_{Yukawa} = y^L_{\alpha\beta}\overline{L}_\alpha E_\beta H + y^D_{\alpha\beta}\overline{Q}_\alpha D_\beta H + y^U_{\alpha\beta}\overline{Q}_\alpha U_\beta \tilde{H} + h.c., \tag{2.3.4}$$

where $\tilde{H} = i\tau_2 H^\dagger$.At last the Higgs part of the Lagrangian contains the Higgs potential which is chosen in such a way that the Higgs field acquires the vacuum expectation value and the potential itself is stable

$$\mathcal{L}^{s.m.}_{Higgs} = -V = m^2H^\dagger H - \frac{\lambda}{2}(H^\dagger H)^2. \tag{2.3.5}$$

Here there are two arbitrary parameters: $m^2$ and $\lambda$. The ghost fields and the gauge fixing terms are omitted.The Lagrangian of the SM contains the following set of free parameters:

(1) 3 gauge couplings $g_s, g, g'$, (2) 3 Yukawa matrices $y^L_{\alpha\beta}, y^D_{\alpha\beta}, y^U_{\alpha\beta}$,

(3) Higgs coupling constant $\lambda$, (4) Higgs mass parameter $m^2$,

(4) the number of the matter fields (generations).

All particles obtain their masses due to spontaneous breaking of the $SU_{left}(2)$

symmetry group via a nonzero vacuum expectation value (v.e.v.) of the Higgs field

$$< H >= \begin{pmatrix} v \\ 0 \end{pmatrix}, \quad v = \frac{m}{\sqrt{\lambda}}. \tag{2.3.6}$$

As a result, the gauge group of the SM is spontaneously broken down to

$$SU_c(3) \otimes SU_L(2) \otimes U_Y(1) \Rightarrow SU_c(3) \otimes U_{EM}(1). \tag{2.3.7}$$

The physical weak intermediate bosons are linear combinations of the gauge ones

$$W^{\pm}_{\mu} = \frac{W^1_{\mu} \mp iW^2_{\mu}}{\sqrt{2}}, \quad Z_{\mu} = -\sin\theta_W B_{\mu} + \cos\theta_W W^3_{\mu} \tag{2.3.8}$$

with masses

$$m_W = \frac{1}{\sqrt{2}} g v, \quad m_Z = \frac{m_W}{\cos\theta_W}, \quad \tan\theta_W = \frac{g'}{g}, \tag{2.3.9}$$

while the photon field

$$\gamma_{\mu} = \cos\theta_W B_{\mu} + \sin\theta_W W^3_{\mu} \tag{2.3.10}$$

remains massless. The matter fields acquire masses proportional to the corresponding Yukawa couplings:

$$M^u_{\alpha\beta} = y^u_{\alpha\beta} v, \; M^d_{\alpha\beta} = y^d_{\alpha\beta} v, \; M^l_{\alpha\beta} = y^l_{\alpha\beta} v, \; m_H = \sqrt{2}\, m. \tag{2.3.11}$$

The mass matrices have to be diagonalized to get the quark and lepton masses. The explicit mass terms in the Lagrangian are forbidden because they are not $SU_{\mathbf{left}}(2)$ symmetric. They would destroy the gauge invariance and, hence, the renormalizability of the Standard Model. To preserve the gauge invariance we use the mechanism of spontaneous symmetry breaking which, as was explained above, allows one to get the renormalizable theory with massive fields.

The Feynman rules in the SM include the ones for QED and QCD with additional new vertices corresponding to the $SU(2)$ group and the Yukawa interaction, as well as the vertices with goldstone particles if one works in the renormalizable gauge. We will not write them down due to their complexity, though the general form is obvious.

Let us consider now the one-loop divergent diagrams in the SM. Besides the familiar diagrams in $QED$ and $QCD$ discussed below in section IV.5 one has the diagrams presented in Fig.2.3.5. The diagrams containing the goldstone bosons are omitted. The calculation of these diagrams is similar to what we have done below in section IV.5. Therefore, we show only the results for the renormalization constants of the fields and the coupling constants. They have the form (for the gauge fields we use the Feynman gauge)

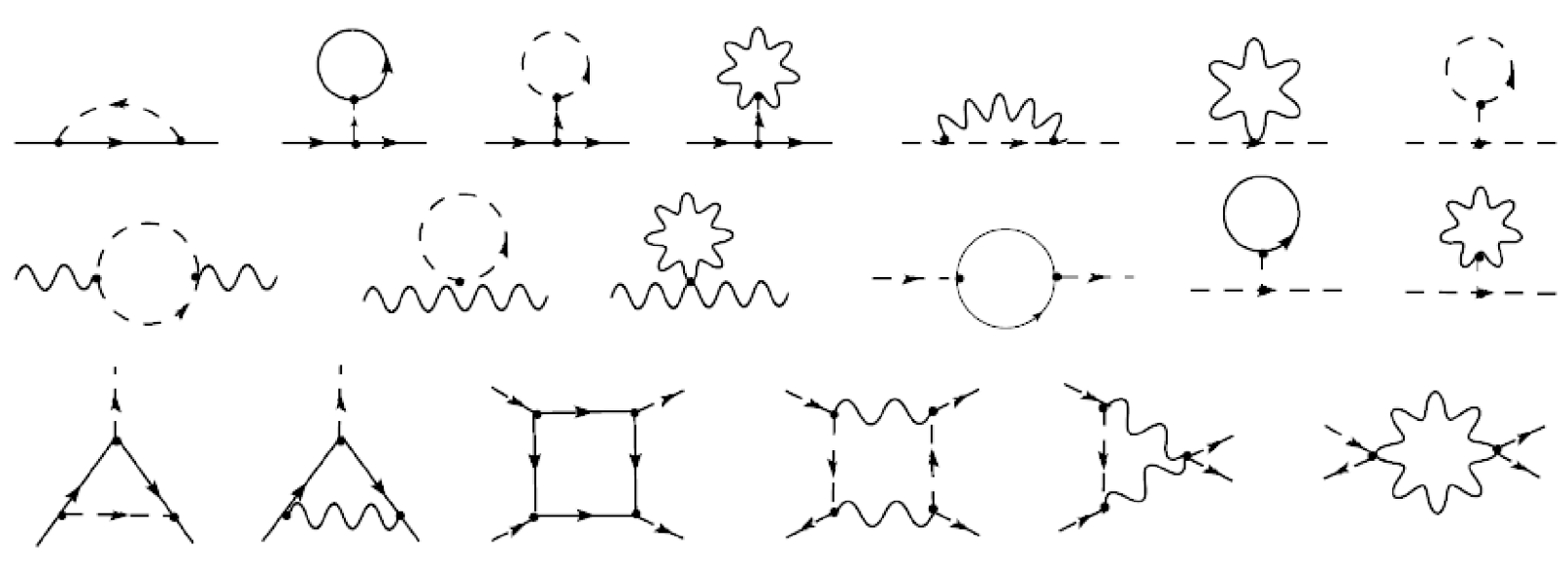

Fig.2.3.5.

Fig.2.3.5.Some divergent one-loop diagrams in the SM. The dotted line denotes the Higgs
field, the solid line - the quark and lepton fields, and the wavy line - the gauge fields

**Remark 2.3.3**.In standard sector the renormalization constants of the fields and the coupling constants reads (see Sect.IV.1-IV.8):

$$
\begin{aligned}
&Z_{2Q_L}(\varepsilon, g'(\varepsilon), g(\varepsilon), g_s(\varepsilon)) = \\
&1 - \frac{1}{\varepsilon}\frac{1}{16\pi^2}\Big[\frac{1}{36}g'^2(\varepsilon) + \frac{3}{4}g^2(\varepsilon) + \frac{4}{3}g_s^2(\varepsilon) + \frac{1}{2}y_U^2 + \frac{1}{2}y_D^2\Big],\\
&Z_{2u_R}(\varepsilon, g'(\varepsilon), g_s(\varepsilon)) = 1 - \frac{1}{\varepsilon}\frac{1}{16\pi^2}\Big[\frac{4}{9}g'^2(\varepsilon) + \frac{4}{3}g_s^2(\varepsilon) + y_U^2\Big],\\
&Z_{2d_R}(\varepsilon, g'(\varepsilon), g_s(\varepsilon)) = 1 - \frac{1}{\varepsilon}\frac{1}{16\pi^2}\Big[\frac{1}{9}g'^2(\varepsilon) + \frac{4}{3}g_s^2(\varepsilon) + y_D^2\Big],\\
&Z_{2L_L}(\varepsilon,, g'(\varepsilon), g(\varepsilon)) = 1 - \frac{1}{\varepsilon}\frac{1}{16\pi^2}[\frac{1}{4}g'^2(\varepsilon) + \frac{3}{4}g^2(\varepsilon) + \frac{1}{2}y_L^2],\\
&Z_{2e_R}(\varepsilon) = 1 - \frac{1}{\varepsilon}\frac{1}{16\pi^2}[g'^2(\varepsilon) + y_L^2],\\
&Z_{2H}(\varepsilon, g'(\varepsilon), g(\varepsilon)) = 1 + \frac{1}{\varepsilon}\frac{1}{16\pi^2}\Big[\frac{1}{2}g'^2(\varepsilon) + \frac{3}{2}g^2(\varepsilon) - 3y_U^2 - 3y_D^2 - y_L^2\Big],\\
&Z_{3B}(\varepsilon, g'(\varepsilon)) = 1 - \frac{1}{\varepsilon}\frac{1}{16\pi^2}\Big[\frac{20}{9}N_F + \frac{1}{6}N_H\Big]g'^2(\varepsilon) \qquad \text{U(1)}_\text{Y}\ \text{boson}\\
&Z_{3A}(\varepsilon, e(\varepsilon)) = 1 + \frac{1}{\varepsilon}\frac{1}{16\pi^2}\Big[3 - \frac{32}{9}N_F\Big]e^2(\varepsilon) \qquad \text{photon}\\
&Z_{3W}(\varepsilon, g(\varepsilon)) = 1 + \frac{1}{\varepsilon}\frac{1}{16\pi^2}\Big[\frac{10}{3} - \frac{1}{3}(N_F + 3N_F) - \frac{1}{6}N_H\Big]g^2(\varepsilon),\\
&Z_{3G}(\varepsilon, g_s(\varepsilon)) = 1 + \frac{1}{\varepsilon}\frac{1}{16\pi^2}\Big[5 - \frac{4}{3}N_F\Big]g_s^2(\varepsilon),\\
&Z_{g_3^2}(\varepsilon, g_s(\varepsilon)) = 1 + \frac{1}{\varepsilon}\frac{1}{16\pi^2}[-11 + \frac{4}{3}N_F]g_s^2(\varepsilon), \qquad (2.3.12)\\
&Z_{g_2^2}(\varepsilon, g(\varepsilon)) = 1 + \frac{1}{\varepsilon}\frac{1}{16\pi^2}\Big[-\frac{22}{3} + \frac{4}{3}N_F + \frac{1}{6}N_H\Big]g^2(\varepsilon),\\
&Z_{g'^2}(\varepsilon, g'(\varepsilon)) = 1 + \frac{1}{\varepsilon}\frac{1}{16\pi^2}\Big[\frac{20}{9}N_F + \frac{1}{6}N_H\Big]g'^2(\varepsilon),\\
&Z_{y_U^2}(\varepsilon, g'(\varepsilon), g(\varepsilon), g_s(\varepsilon)) = \\
&1 + \frac{1}{\varepsilon}\frac{1}{16\pi^2}\Big[-\frac{17}{12}g'^2(\varepsilon) - \frac{9}{4}g^2(\varepsilon) - 8g_s^2(\varepsilon) + \frac{9}{2}y_U^2 + \frac{3}{2}y_D^2 + y_L^2\Big],\\
&Z_{y_D^2}(\varepsilon, g'(\varepsilon), g(\varepsilon), g_s(\varepsilon)) = \\
&1 + \frac{1}{\varepsilon}\frac{1}{16\pi^2}\Big[-\frac{5}{12}g'^2(\varepsilon) - \frac{9}{4}g^2(\varepsilon) - 8g_s^2(\varepsilon) + \frac{3}{2}y_U^2 + \frac{9}{2}y_D^2 + y_L^2\Big],\\
&Z_{y_L^2}(\varepsilon, g'(\varepsilon), g(\varepsilon)) = \\
&1 + \frac{1}{\varepsilon}\frac{1}{16\pi^2}\Big[-\frac{15}{4}g'^2(\varepsilon) - \frac{9}{4}g^2(\varepsilon) + \frac{9}{4}y_L^2 + 3y_U^2 + 3y_D^2\Big],\\
&Z_\lambda(\varepsilon, g'(\varepsilon), g(\varepsilon), \lambda(\varepsilon)) = \\
&1 + \frac{1}{\varepsilon}\frac{1}{16\pi^2}\Big[-\frac{3}{2}g'^2(\varepsilon) - \frac{9}{2}g^2(\varepsilon) + 2(3y_U^2 + 3y_D^2 + y_L^2) + 6\lambda(\varepsilon)\\
&-2(3y_U^4 + 3y_D^4 + y_L^4) \times \lambda^{-1}(\varepsilon) + \\
&\Big(\frac{3}{8}g'^4(\varepsilon) + \frac{9}{8}g^4(\varepsilon) + \frac{3}{4}g^2(\varepsilon)g'^2(\varepsilon)\Big) \times \lambda^{-1}(\varepsilon)\Big],
\end{aligned}
$$

where, for simplicity, we ignored the mixing between the generations and assumed the Since the masses of all the particles are equal to the product of the gauge or Yukawa couplings and the vacuum expectation value of the Higgs field, in the minimal subtraction scheme the mass ratios are renormalized the same way as the ratio of couplings. To find the renormalization of the mass itself, one should know how the v.e.v. is renormalized or find explicitly the mass counter-term from Feynman diagrams. In this case, one has also to take into account the tad-pole diagrams shown in Fig.2.3.5, including the diagrams with goldstone bosons. For illustration we present the renormalization constant of the $b$-quark mass in the SM

$$Z_{m_b}(\varepsilon,g'(\varepsilon),g(\varepsilon),\lambda(\varepsilon)) =$$

$$1+\frac{1}{\varepsilon}\frac{1}{16\pi^2}\left[\sum\nolimits_l \frac{y_l^4}{\lambda(\varepsilon)}+3\sum\nolimits_q \frac{y_q^4}{\lambda(\varepsilon)}-\frac{3}{2}\lambda(\varepsilon)+\frac{3}{4}(y_b^2-y_t^2)\right. \tag{2.3.13}$$

$$\left.-\frac{3}{16}\frac{(g^2(\varepsilon)+g'^2(\varepsilon))^2}{\lambda}-\frac{3}{8}\frac{g^4(\varepsilon)}{\lambda(\varepsilon)}-3Q_b(Q_b-T_b^3)g'^2(\varepsilon)-4g_s^2(\varepsilon)\right].$$

The result for the t-quark can be obtained by replacing b by t. For the light quarks the Yukawa constants are very small and can be ignored in Eq.(2.3.15).

**Remark 2.3.4.** In the standard sector the bare Lagrangian reads

$$\mathcal{L}^{s.m.Bare} = \mathcal{L}^{s.m.Bare}_{gauge}+\mathcal{L}^{s.m.Bare}_{Yukawa}+\mathcal{L}^{s.m.Bare}_{Higgs}, \tag{2.3.14}$$

with renormalization constants of the fields,coupling constants (2.3.12) and masess (2.3.13) satisfies the following conditions: (see Sect.IV.1-IV.8)

$$\begin{gathered}
0<Z^{s.m.}_{2Q_L}(\bar\varepsilon,g'(\bar\varepsilon),g(\bar\varepsilon),g_s(\bar\varepsilon)),0<Z^{s.m.}_{2u_R}(\bar\varepsilon,g'(\bar\varepsilon),g_s(\bar\varepsilon)),0<Z^{s.m.}_{2d_R}(\varepsilon,g'(\varepsilon),g_s(\varepsilon)),\\
0<Z^{s.m.}_{2L_L}(\bar\varepsilon,,g'(\bar\varepsilon),g(\bar\varepsilon)),0<Z^{s.m.}_{2H}(\bar\varepsilon,g'(\bar\varepsilon)),0<Z^{s.m.}_{3B}(\bar\varepsilon,g'(\bar\varepsilon)),0<Z^{s.m.}_{3A}(\bar\varepsilon,e(\bar\varepsilon)),\\
0<Z^{s.m.}_{3W}(\bar\varepsilon,g(\bar\varepsilon)),0<Z^{s.m.}_{3G}(\bar\varepsilon,g_s(\bar\varepsilon)),0<Z^{s.m.}_{g_3^2}(\bar\varepsilon,g_s(\bar\varepsilon)),0<Z^{s.m.}_{g_2^2}(\bar\varepsilon,g(\bar\varepsilon)),\\
0<Z^{s.m.}_{g'^2}(\bar\varepsilon,g'(\bar\varepsilon)),0<Z^{s.m.}_{y_U^2}(\bar\varepsilon,g'(\bar\varepsilon),g'(\bar\varepsilon),g_s(\bar\varepsilon)),0<Z^{s.m.}_{y_D^2}(\bar\varepsilon,g'(\bar\varepsilon),g(\bar\varepsilon),g_s(\bar\varepsilon)),\\
0<Z^{s.m.}_{y_L^2}(\bar\varepsilon,g'(\bar\varepsilon),g(\bar\varepsilon)),0<Z^{s.m.}_{\lambda}(\bar\varepsilon,g'(\bar\varepsilon),g(\bar\varepsilon),\lambda(\bar\varepsilon)),0<Z^{s.m.}_{m_b}(\bar\varepsilon,g'(\bar\varepsilon),g(\bar\varepsilon),\lambda(\bar\varepsilon)),\\
\text{etc.},
\end{gathered} \tag{2.3.15}$$

for some $\bar\varepsilon \in (0,1]$, see Sec.4.1-4.8.

## 2.3.2.The Standard Model of fundamental interactions in a ghost matter sector below fundamental high-energy cutoff $\Lambda_*$.

In the extendended Standard Model of fundamental interactions besides the gauge interactions and the quartic interaction of the Higgs ghost fields there are also Yukawa type interactions of the fermion ghost fields with the Higgs ghost field. These interactions are also renormalizable and is characterized by the Yukawa coupling constants, one for each fermion field. Another property of the Standard Model is that it has the gauge group $SU_c(3)\times SU_L(2)\times U_Y(1)$ which is spontaneously broken to $SU_c(3)\times U_{EM}(1)$. In the theories with spontaneously broken symmetry, according to the Goldstone theorem there are massless particles, the goldstone ghost bosons.

The Lagrangian of the extendended Standard Model in a ghost sector consists of the following three parts:

$$\mathcal{L}^{g.s.} = \mathcal{L}^{g.s.}_{gauge}+\mathcal{L}^{g.s.}_{Yukawa}+\mathcal{L}^{g.s.}_{Higgs}, \tag{2.3.16}$$

The Lagrangian (2.3.16) is obtained from the Lagrangian (2.3.1) of the Standard Model as corresponding bare Lagrangian with the renormalization constants of the fields and the coupling constants given by Eq.(2.3.12).

$$\mathcal{L}^{g.s.} = \mathcal{L}^{Bare}_{gauge}+\mathcal{L}^{Bare}_{Yukawa}+\mathcal{L}^{Bare}_{Higgs}, \tag{2.3.17}$$

The gauge part is totally fixed by the requirement of the gauge invariance leaving only the values of the couplings as free parameters

$$\mathcal{L}^{g.s.} = \mathcal{L}^{Bare}_{gauge} = -\frac{1}{4}\breve{G}^a_{\mu\nu}\breve{G}^a_{\mu\nu} - \frac{1}{4}\breve{W}^i_{\mu\nu}\breve{W}^i_{\mu\nu} - \frac{1}{4}\breve{B}_{\mu\nu}\breve{B}_{\mu\nu}$$
$$+i\overline{\breve{L}}_\alpha\gamma^\mu\breve{D}_\mu\breve{L}_\alpha + i\overline{\breve{Q}}_\alpha\gamma^\mu\breve{D}_\mu\breve{Q}_\alpha + i\overline{\breve{E}}_\alpha\gamma^\mu\breve{D}_\mu\breve{E}_\alpha \tag{2.3.18}$$
$$+i\overline{\breve{U}}_\alpha\gamma^\mu\breve{D}_\mu\breve{U}_\alpha + i\overline{\breve{D}}_\alpha\gamma^\mu\breve{D}_\mu\breve{D}_\alpha + (\breve{D}_\mu\breve{H})^\dagger(\breve{D}_\mu\breve{H}),$$

where we abraviate $\breve{G}^a_{\mu\nu}, \breve{W}^i_{\mu\nu}, \breve{B}_{\mu\nu}, \breve{L}_\alpha, \breve{Q}_\alpha, \breve{E}_\alpha, \breve{U}_\alpha, \breve{H}$ for a short instead $G^{a,Bare}_{\mu\nu}, W^{i,Bare}_{\mu\nu}, B^{Bare}_{\mu\nu}, \ldots,\ H^{Bare}$ and where the following notations for the covariant derivatives are used

$$\breve{G}^a_{\mu\nu} = \partial_\mu\breve{G}^a_\nu - \partial_\nu\breve{G}^a_\mu + g_s f^{abc}\breve{G}^b_\mu\breve{G}^c_\nu, \breve{W}^i_{\mu\nu} = \partial_\mu\breve{W}^i_\nu - \partial_\nu\breve{W}^i_\mu + g\epsilon^{ijk}\breve{W}^j_\mu\breve{W}^k_\nu,$$
$$\breve{B}_{\mu\nu} = \partial_\mu\breve{B}_\nu - \partial_\nu\breve{B}_\mu, \breve{D}_\mu\breve{L}_\alpha = \left(\partial_\mu - i\frac{g}{2}\tau^i\breve{W}^i_\mu + i\frac{g^{'}}{2}\breve{B}_\mu\right)\breve{L}_\alpha,$$
$$\breve{D}_\mu\breve{E}_\alpha = (\partial_\mu + ig'\breve{B}_\mu)\breve{E}_\alpha, \breve{D}_\mu\breve{Q}_\alpha =$$
$$\left(\partial_\mu - i\frac{g}{2}\tau^i\breve{W}^i_\mu - i\frac{g^{'}}{6}\breve{B}_\mu - i\frac{g_s}{2}\lambda^a\breve{G}^a_\mu\right)\breve{Q}_\alpha, \tag{2.3.19}$$
$$\breve{D}_\mu\breve{U}_\alpha = \left(\partial_\mu - i\frac{2}{3}g'\breve{B}_\mu - i\frac{g_s}{2}\lambda^a\breve{G}^a_\mu\right)\breve{U}_\alpha, \breve{D}_\mu\breve{D}_\alpha =$$
$$\left(\partial_\mu + i\frac{1}{3}g'\breve{B}_\mu - i\frac{g_s}{2}\lambda^a\breve{G}^a_\mu\right)\breve{D}_\alpha.$$

The Yukawa part of the Lagrangian in a ghost sector, which is needed for the generation of the ghost quark and ghost lepton masses is also chosen in the gauge invariant form and contains arbitrary Yukawa couplings (we ignore the neutrino masses, for simplicity)

$$\mathcal{L}^{Bare}_{\mathbf{Yukawa}} = y^L_{\alpha\beta}\overline{\breve{L}}_\alpha\breve{E}_\beta\breve{H} + y^D_{\alpha\beta}\overline{\breve{Q}}_\alpha\breve{D}_\beta\breve{H} + y^U_{\alpha\beta}\overline{\breve{Q}}_\alpha\breve{U}_\beta\widetilde{\breve{H}} + h.c., \tag{2.3.20}$$

where $\tilde{H} = i\tau_2\breve{H}^\dagger$.At last the ghost Higgs part of the Lagrangian contains the Higgs potential which is chosen in such a way that the ghost Higgs field acquires the vacuum expectation value and the potential itself is stable

$$\mathcal{L}^{Bare}_{\mathbf{Higgs}} = -V^{Bare} = \breve{m}^2\breve{H}^\dagger\breve{H} - \frac{\breve{\lambda}}{2}(\breve{H}^\dagger\breve{H})^2. \tag{2.3.21}$$

Here there are two bare parameters $\breve{m}^2_\varepsilon, \breve{\lambda}_\varepsilon$.

The Lagrangian of the SM in a ghost sector contains the following set of the bare parameters:

(1) 3 gauge couplings $\breve{g}_s, \breve{g}, \breve{g}'$, (2) 3 Yukawa matrices $y^L_{\alpha\beta}, y^D_{\alpha\beta}, y^U_{\alpha\beta}$,

(3) Higgs coupling constant $\breve{\lambda}$, (4) ghost Higgs mass parameter $\breve{m}^2$,

(4) the number of the ghost matter fields (generations).

All particles obtain their masses due to spontaneous breaking of the $SU_{left}(2)$ symmetry group via a nonzero vacuum expectation value (v.e.v.) of the Higgs field

$$< \breve{H} >= \begin{pmatrix} \breve{v} \\ 0 \end{pmatrix}, \quad \breve{v} = \frac{\breve{m}}{\sqrt{\breve{\lambda}}}. \tag{2.3.22}$$

As a result, the gauge group is spontaneously broken down to

$$SU_c(3) \otimes SU_L(2) \otimes U_Y(1) \Rightarrow SU_c(3) \otimes U_{EM}(1). \tag{2.3.23}$$

The physical weak intermediate bosons are linear combinations of the gauge ones

$$\breve{W}_{\mu}^{\pm Bare} = \frac{\breve{W}_{\mu}^{1Bare} \mp i\breve{W}_{\mu}^{2Bare}}{\sqrt{2}}, \quad \breve{Z}_{\mu}^{Bare} = -\sin\theta_{\breve{W}}\breve{B}_{\mu}^{Bare} + \cos\theta_{\breve{W}}\breve{W}_{\mu}^{3Bare} \qquad (2.3.24)$$

with masses

$$\breve{m}_{\breve{W}} = \frac{1}{\sqrt{2}}\breve{g}v, \quad \breve{m}_{\breve{Z}} = \frac{\breve{m}_{\breve{W}}}{\cos\theta_{\breve{W}}}, \quad \tan\theta_{\breve{W}} = \frac{\breve{g}'}{\breve{g}}, \qquad (2.3.25)$$

while the ghost photon field

$$\breve{\gamma}_{\mu} = \cos\theta_{\breve{W}}\breve{B}_{\mu} + \sin\theta_{\breve{W}}\breve{W}_{\mu}^{3} \qquad (2.3.26)$$

remains massless. The matter fields acquire masses proportional to the corresponding Yukawa couplings:

$$\breve{M}_{\alpha\beta}^{u} = y_{\alpha\beta}^{u}\breve{v}, \ \breve{M}_{\alpha\beta}^{d} = y_{\alpha\beta}^{d}\breve{v}, \ \breve{M}_{\alpha\beta}^{l} = y_{\alpha\beta}^{l}\breve{v}, \ \breve{m}_{\breve{H}} = \sqrt{2}\,\breve{m}. \qquad (2.3.27)$$

**Remark 2.3.5.** In a ghost standard sector the bare Lagrangian reads

$$\mathcal{L}^{g.m.Bare} = \mathcal{L}_{gauge}^{g.m.Bare} + \mathcal{L}_{Yukawa}^{g.m.Bare} + \mathcal{L}_{Higgs}^{g.m.Bare}, \qquad (2.3.28)$$

with renormalization constants of the fields,coupling constants (2.3.12) and masess (2.3.13) satisfies the following conditions: (see Sect.IV.1-IV.8)

$$\begin{gathered}
Z_{2Q_L}^{g.m.}(\bar{\varepsilon}, g'(\bar{\varepsilon}), g(\bar{\varepsilon}), g_s(\bar{\varepsilon})) < 0, Z_{2u_R}^{g.m.}(\bar{\varepsilon}, g'(\bar{\varepsilon}), g_s(\bar{\varepsilon})) < 0, \\
Z_{2d_R}^{g.m.}(\varepsilon, g'(\varepsilon), g_s(\varepsilon)) < 0, \\
Z_{2L_L}^{g.m.}(\bar{\varepsilon},, g'(\bar{\varepsilon}), g(\bar{\varepsilon})) < 0, Z_{2H}^{g.m.}(\bar{\varepsilon}, g'(\bar{\varepsilon})) < 0, Z_{3B}^{g.m.}(\bar{\varepsilon}, g'(\bar{\varepsilon})) < 0, \\
Z_{3A}^{g.m.}(\bar{\varepsilon}, e(\bar{\varepsilon})) < 0, \\
Z_{3W}^{g.m.}(\bar{\varepsilon}, g(\bar{\varepsilon})) < 0, Z_{3G}^{g.m.}(\bar{\varepsilon}, g_s(\bar{\varepsilon})) < 0, Z_{g_3^2}^{g.m.}(\bar{\varepsilon}, g_s(\bar{\varepsilon})) < 0, \\
Z_{g_2^2}^{g.m.}(\bar{\varepsilon}, g(\bar{\varepsilon})) < 0, \qquad (2.3.29) \\
Z_{g'^2}^{g.m.}(\bar{\varepsilon}, g'(\bar{\varepsilon})) < 0, Z_{y_U^2}^{g.m.}(\bar{\varepsilon}, g'(\bar{\varepsilon}), g'(\bar{\varepsilon}), g_s(\bar{\varepsilon})) < 0, \\
Z_{y_D^2}^{g.m.}(\bar{\varepsilon}, g'(\bar{\varepsilon}), g(\bar{\varepsilon}), g_s(\bar{\varepsilon})) < 0, \\
Z_{y_L^2}^{g.m.}(\bar{\varepsilon}, g'(\bar{\varepsilon}), g(\bar{\varepsilon})) < 0, Z_{\lambda}^{g.m.}(\bar{\varepsilon}, g'(\bar{\varepsilon}), g(\bar{\varepsilon}), \lambda(\bar{\varepsilon})) < 0, \\
Z_{m_b}^{g.m.}(\bar{\varepsilon}, g'(\bar{\varepsilon}), g(\bar{\varepsilon}), \lambda(\bar{\varepsilon})) < 0, \\
\text{etc.,}
\end{gathered}$$

for some $\bar{\varepsilon} \in (0,1]$.

Remind that vacuum energy density for free scalar quantum field with a wrong statistic is:

$$\varepsilon(\mu) = -\frac{1}{2}\frac{1}{(2\pi\hbar)^3}\int_0^{\infty} 4\pi c\sqrt{p^2+\mu^2}\,p^2 dp = K'\int_0^{\infty}\sqrt{p^2+\mu^2}\,p^2 dp = K'I(\mu), \qquad (2.3.30)$$

where $\mu = m_0 c$. From the basic definitions [1]:

$$p = T_{xx}, p(\mu) = -\frac{1}{2}\frac{1}{(2\pi\hbar)^3}\int_0^{\infty} u_x p_x 4\pi p^2 dp, \mathbf{u} = \frac{c\mathbf{p}}{\sqrt{p^2+\mu^2}}, \overline{u_x p_x} = \frac{1}{3}\langle\mathbf{u},\mathbf{p}\rangle$$

one obtains

$$p(\mu) = \frac{K'}{3}\int_0^{\infty}\frac{p^4 dp}{\sqrt{p^2+\mu^2}} = K'F(\mu). \qquad (2.3.32)$$

**Remark 2.3.6.** Note that the integral in RHS of the Eq.(2.3.30) and in the Eq.(2.3.32)

divergent and ultraviolet cutoff is needed.
Thus in accordance with [1] we set

$$\varepsilon(\mu,p_0) = K'I(\mu,p_0), p(\mu,p_0) = K'F(\mu,p_0), \qquad (2.3.33)$$

where

$$I(\mu,p_0) = \int_0^{p_0}\sqrt{p^2+\mu^2}\,p^2dp, F(\mu,p_0) = \int_0^{p_0}\frac{p^4dp}{\sqrt{p^2+\mu^2}}, \qquad (2.3.34)$$

where $p_0 < \Lambda_*/c$. For fermionic quantum field with a wrong statistic, similarly one obtains

$$\varepsilon(\mu,p_0) = -4K'I(\mu,p_0), p(\mu) = -4K'F(\mu,p_0). \qquad (2.3.35)$$

Thus from Eqs.(2.3.33)-(2.3.35) by using formally Pauli-Villars regularization [7],[8] and regularization by high-energy cutoff the expression for free vacuum energy density $\varepsilon$ reads

$$\varepsilon = \sum_{i=0}^{2M} f_i I(\mu_i,p_0) \qquad (2.3.36)$$

and the expression for pressure $p$ reads

$$p = \sum_{i=0}^{2M} f_i F(\mu_i,p_0). \qquad (2.3.37)$$

Here $\mu_i$ is a finite positive sequence $\mu_i \in \mathbb{R}_+, i = 1,2,\ldots,2M$ and $f_i$ is a finite sequence $f_i \in \mathbb{R}, f_i \leq 1, i = 1,2,\ldots,2M$.

**Definition 2.3.1.** We define now discrete distribution $f_{PV} : \mathbb{R}_+ \to \mathbb{R}$ by formula

$$f_{PV}(\mu_i) = f_i, \qquad (2.3.38)$$

and we will call it as a full discrete Pauli-Villars masses distribution.

**Remark 2.3.7.** We assum now that in Eqs.(2.3.36)-(2.3.37): (i) the quantities $\mu_i^{s.m} = \mu_i, i = 1,2,\ldots,M$ is a masses of a physical particles corresponding to standard matter and (ii) the quantities $\mu_i^{g.m} = \mu_i, i = M+1,2,\ldots,2M$ is a masses of ghost particles with a wrong kinetic term and wrong statistics corresponding to a physical dark matter.

**Remark 2.3.8.** We recall that the Euler-Maclaurin summation formula reads

$$\begin{gathered}\sum_{i=1}^{2M} g(\mu_1+(i-1)h) = \int_{\mu_1}^{\mu_{2M}} f(\mu)d\mu + A_1[g(\mu_{2M})-g(\mu_1)] + \\ A_2h[g'(\mu_{2M})-g'(\mu_1)] + O(h^2), \\ f(\mu) = \frac{1}{h}g(\mu)\end{gathered} \qquad (2.3.39)$$

Let $g(\mu)$ be an appropriate continuous function such that: (i) $g(\mu_i) = f_i, i = 1,2,\ldots,2M$, (ii) $g'(\mu_{2M}) = 0, g'(\mu_1) = 0$.

Thus from Eqs.(2.3.36)-(2.3.37) and Eqs.(2.3.39) we obtain

$$\begin{gathered}\varepsilon = \sum_{i=0}^{2M} f_i I(\mu_i,p_0) = \\ \int_{\mu_1}^{\mu_{2M}} f(\mu)I(\mu,p_0)d\mu + A_1h[f(\mu_{2M})I(\mu_{2M},p_0) - f(\mu_1)I(\mu_1,p_0)] + O(h^2)\end{gathered} \qquad (2.3.40)$$

and

$$p = \sum_{i=0}^{2M} f_i F(\mu_i, p_0) = \int_{\mu_1}^{\mu_{2M}} f(\mu) F(\mu, p_0) d\mu + A_1 h [f(\mu_{2M}) F(\mu_{2M}, p_0) - f(\mu_1) F(\mu_1, p_0)] + O(h^2). \tag{2.3.41}$$

**Definition 2.3.2.** We will call the function $f(\mu)$ as a full continuous Pauli-Villars masses distribution.

**Definition 2.3.3.** We define now: (i) discrete distribution $f_{PV}^{b.g.m} : \mathbb{R}_+ \to \mathbb{R}$ by formula

$$f_{PV}^{b.g.m}(\mu_i^{s.m}) = f_i, i = 1, 2, \ldots, M \tag{2.3.42}$$

and we will call it as discrete Pauli-Villars masses distribution of the bosonic ghost matter

and

(ii) discrete distribution $f_{PV}^{f.g.m} : \mathbb{R}_+ \to \mathbb{R}$ by formula

$$f_{PV}^{f.g.m}(\mu_i) = f_i, i = M+1, 2, \ldots, 2M \tag{2.3.43}$$

and we will call it as discrete Pauli-Villars masses distribution of the fermionic ghost matter.

**Remark 2.3.9.** We rewrite now the Eqs.(2.3.36)-(2.3.37) in the following equivalent form

$$\varepsilon = \sum_{i=1}^{M} f_{PV}^{b.g.m}(\mu_i^{s.m}) I\left(\mu_i^{b.g.m}, p_0\right) + \sum_{j(i)=M+1}^{2M} f_{PV}^{f.g.m}\left(\mu_{j(i)}^{f.g.m}\right) I\left(\mu_{j(i)}^{f.g.m}, p_0\right) \tag{2.3.44}$$

and

$$p = \sum_{i=1}^{M} f_{PV}^{b.g.m}\left(\mu_i^{b.g.m}\right) F\left(\mu_i^{b.g.m}, p_0\right) + \sum_{j(i)=M+1}^{2M} f_{PV}^{f.g.m}\left(\mu_{j(i)}^{f.g.m}\right) F\left(\mu_{j(i)}^{f.g.m}, p_0\right), \tag{2.3.45}$$

where $j(i) = i + M, i = 1+1, 2, \ldots, M.$

**Remark 2.3.10.** We assume now that:(i) $\mu_i^{b.g.m} \approx \mu_{j(i)}^{f.g.m}$,

(ii) $\left| f_{PV}^{b.g.m}\left(\mu_i^{b.g.m}\right) + f_{PV}^{f.g.m}\left(\mu_{j(i)}^{f.g.m}\right) \right| \ll 1$, i.e.,

$$f_{PV}^{b.g.m}\left(\mu_i^{b.g.m}\right) \approx -f_{PV}^{f.g.m}\left(\mu_{j(i)}^{f.g.m}\right). \tag{2.3.46}$$

Note that Eq.(2.3.46) meant highly symmetric discrete Pauli-Villars masses distribution

between bosonic ghost matter and fermionic ghost matter above that scale $\Lambda_*$

Thus from Eqs.(2.3.44)-(2.3.45) and Eqs.(2.3.46) we obtain

$$\varepsilon = \sum_{i=1}^{M} f_{PV}^{b.g.m}\left(\mu_i^{b.g.m}\right) I\left(\mu_i^{b.g.m}, p_0\right) + \sum_{j(i)=M+1}^{2M} f_{PV}^{f.g.m}\left(\mu_{j(i)}^{f.g.m}\right) I\left(\mu_{j(i)}^{f.g.m}, p_0\right) = \sum_{i=1}^{M} \left[ f_{PV}^{b.g.m}\left(\mu_i^{b.g.m}\right) + f_{PV}^{f.g.m}\left(\mu_{j(i)}^{f.g.m}\right) \right] I(\mu_i, p_0) \tag{2.3.47}$$

and

$$p = \sum_{i=1}^{M} f_{PV}^{b.g.m}\left(\mu_i^{b.g.m}\right) F\left(\mu_i^{b.g.m}, p_0\right) + \sum_{j(i)=M+1}^{2M} f_{PV}^{f.g.m}\left(\mu_{j(i)}^{f.g.m}\right) F\left(\mu_{j(i)}^{f.g.m}, p_0\right) = \sum_{i=1}^{M}\left[f_{PV}^{b.g.m}\left(\mu_i^{b.g.m}\right) + f_{PV}^{f.g.m}\left(\mu_{j(i)}^{f.g.m}\right)\right] F(\mu_i, p_0). \quad (2.3.48)$$

From Eqs.(2.3.47)-(2.3.48) and Eqs.(2.3.39) finally we obtain

$$\varepsilon = \sum_{i=1}^{M}\left[f_{PV}^{b.g.m}\left(\mu_i^{b.g.m}\right) + f_{PV}^{f.g.m}\left(\mu_{j(i)}^{f.g.m}\right)\right] I(\mu_i, p_0) = \int_{\mu_1}^{\mu_{\mathbf{eff}}}\left[f_{PV}^{b.g.m}(\mu) + f_{PV}^{f.g.m}(\mu)\right] I(\mu, p_0) d\mu \quad (2.2.49)$$

and

$$p = \sum_{i=1}^{M}\left[f_{PV}^{b.g.m}(\mu_i^{s.m}) + f_{PV}^{f.g.m}\left(\mu_{j(i)}^{f.g.m}\right)\right] F(\mu_i, p_0) = \int_{\mu_1}^{\mu_{\mathbf{eff}}}\left[f_{PV}^{b.g.m}(\mu) + f_{PV}^{f.g.m}(\mu)\right] F(\mu, p_0) d\mu, \quad (2.3.50)$$

where obviously

$$f_{PV}^{b.g.m}(\mu) + f_{PV}^{f.g.m}(\mu) = f_{PV}^{g.m.}(\mu) \approx 0. \quad (2.3.51)$$

Thus finally we obtain

$$\varepsilon^{g.m.}\left(\mu_{\mathbf{eff}}^{(1)}, \mu_{\mathbf{eff}}^{(2)}, p_0\right) = \int_{\mu_{\mathbf{eff}}^{(1)}}^{\mu_{\mathbf{eff}}^{(2)}} f_{PV}^{g.m.}(\mu) I(\mu, p_0) d\mu, \quad (2.3.52)$$

and

$$p^{g.m.}\left(\mu_{\mathbf{eff}}^{(1)}, \mu_{\mathbf{eff}}^{(2)}, p_0\right) = \int_{\mu_{\mathbf{eff}}^{(1)}}^{\mu_{\mathbf{eff}}^{(2)}} f_{PV}^{g.m.}(\mu) F(\mu, p_0) d\mu, \quad (2.3.53)$$

where $\mu_{\mathbf{eff}}^{(1)}, \mu_{\mathbf{eff}}^{(2)} \gg p_0$. In order to calculate $\varepsilon^{g.m.}\left(\mu_{\mathbf{eff}}^{(1)}, \mu_{\mathbf{eff}}^{(2)}, p_0\right)$ and $p^{g.m.}\left(\mu_{\mathbf{eff}}^{(1)}, \mu_{\mathbf{eff}}^{(2)}, p_0\right)$ let us evaluate now the following quantities defined above by Eqs.(2.3.34)

$$I(\mu, p_0) = \int_0^{p_0} p^2 \sqrt{p^2 + \mu^2}\, dp = \int_0^{p_0} \mu p^2 \sqrt{1 + \frac{p^2}{\mu^2}}\, dp \quad (2.3.54)$$

and

$$F(\mu, p_0) = \frac{1}{3}\int_0^{p_0} \frac{p^4 dp}{\sqrt{p^2 + \mu^2}} = \frac{1}{3}\int_0^{p_0} \frac{p^4 \mu^{-1} dp}{\sqrt{1 + \frac{p^2}{\mu^2}}}, \quad (2.3.55)$$

where $p_0/\mu \ll 1$. Note that

$$\sqrt{1+\frac{p^2}{\mu^2}} = 1+\frac{1}{2}\frac{p^2}{\mu^2}-\frac{1}{8}\frac{p^4}{\mu^4}+\frac{1}{16}\frac{p^6}{\mu^6}+\ldots.$$

$$p^2\sqrt{p^2+\mu^2} = p^2\mu\sqrt{1+\frac{p^2}{\mu^2}} = p^2\mu\left(1+\frac{1}{2}\frac{p^2}{\mu^2}-\frac{1}{8}\frac{p^4}{\mu^4}+\frac{1}{16}\frac{p^6}{\mu^6}+\ldots.\right) = \tag{2.3.56}$$

$$p^2\mu+\frac{1}{2}\frac{p^4}{\mu}-\frac{1}{8}\frac{p^6}{\mu^3}+\frac{1}{16}\frac{p^8}{\mu^5}+\ldots.$$

By inserting Eq.(2.3.56) into Eq.(2.3.54) one obtains

$$I(\mu,p_0) = \int_0^{p_0}\left(p^2\mu+\frac{1}{2}\frac{p^4}{\mu}-\frac{1}{8}\frac{p^6}{\mu^3}+\frac{1}{16}\frac{p^8}{\mu^5}+\ldots.\right)dp = \tag{2.3.56}$$

$$\frac{1}{3}p_0^3\mu+\frac{1}{10}\frac{p_0^5}{\mu}-\frac{1}{7\bullet 8}\frac{p_0^7}{\mu^3}+\frac{1}{9\bullet 16}\frac{p_0^9}{\mu^5}+\ldots.$$

Note that

$$\left(1+\frac{p^2}{\mu^2}\right)^{-1/2} = 1-\frac{1}{2}\frac{p^2}{\mu^2}+\frac{3}{8}\frac{p^4}{\mu^4}+\ldots$$

$$p^4\mu^{-1}\left(1+\frac{p^2}{\mu^2}\right)^{-1/2} = \frac{p^4}{\mu}-\frac{1}{2}\frac{p^6}{\mu^3}+\frac{3}{8}\frac{p^8}{\mu^5}+\ldots \tag{2.3.57}$$

By inserting Eq.(2.3.57) into Eq.(2.3.55) one obtains

$$F(\mu,p_0) = \frac{1}{3}\int_0^{p_0}\frac{p^4\mu^{-1}dp}{\sqrt{1+\frac{p^2}{\mu^2}}} = \frac{1}{3}\int_0^{p_0}\left(\frac{p^4}{\mu}-\frac{1}{2}\frac{p^6}{\mu^3}+\frac{3}{8}\frac{p^8}{\mu^5}+\ldots\right)dp = \tag{2.3.58}$$

$$\frac{p_0^5}{3\bullet 5\mu}-\frac{1}{2\bullet 3\bullet 7}\frac{p_0^7}{\mu^3}+\frac{1}{8\bullet 9}\frac{p_0^9}{\mu^5}+\ldots$$

By inserting Eq.(2.3.56) into Eq.(2.3.52) one obtains

$$\varepsilon^{g.m.}\left(\mu_{\mathbf{eff}}^{(1)},\mu_{\mathbf{eff}}^{(2)},p_0\right) = \int_{\mu_{\mathbf{eff}}^{(1)}}^{\mu_{\mathbf{eff}}^{(2)}} f_{PV}^{g.m.}(\mu)I(\mu,p_0)d\mu =$$

$$\int_{\mu_{\mathbf{eff}}^{(1)}}^{\mu_{\mathbf{eff}}^{(2)}} f_{PV}^{g.m.}(\mu)\left(\frac{1}{3}p_0^3\mu+\frac{1}{10}\frac{p_0^5}{\mu}-\frac{1}{7\bullet 8}\frac{p_0^7}{\mu^3}+\frac{1}{9\bullet 16}\frac{p_0^9}{\mu^5}+\ldots.\right)d\mu = \tag{2.3.59}$$

$$\frac{p_0^3}{3}\int_{\mu_{\mathbf{eff}}^{(1)}}^{\mu_{\mathbf{eff}}^{(2)}} f_{PV}^{g.m.}(\mu)\mu d\mu+\frac{p_0^5}{10}\int_{\mu_{\mathbf{eff}}^{(1)}}^{\mu_{\mathbf{eff}}^{(2)}}\frac{f_{PV}^{g.m.}(\mu)d\mu}{\mu}-\frac{p_0^7}{7\bullet 8}\int_{\mu_{\mathbf{eff}}^{(1)}}^{\mu_{\mathbf{eff}}^{(2)}}\frac{f_{PV}^{g.m.}(\mu)d\mu}{\mu^3}+\ldots.$$

By inserting Eq.(2.3.58) into Eq.(2.3.52) one obtains

$$p^{g.m.}\left(\mu_{\mathbf{eff}}^{(1)},\mu_{\mathbf{eff}}^{(2)},p_0\right) = \int_{\mu_{\mathbf{eff}}^{(1)}}^{\mu_{\mathbf{eff}}^{(2)}} f_{PV}^{g.m.}(\mu)F(\mu,p_0)d\mu =$$

$$\int_{\mu_{\mathbf{eff}}^{(1)}}^{\mu_{\mathbf{eff}}^{(2)}} f_{PV}^{g.m.}(\mu)\left(\frac{p_0^5}{3\cdot 5\mu} - \frac{1}{2\cdot 3\cdot 7}\frac{p_0^7}{\mu^3} + \frac{1}{8\cdot 9}\frac{p_0^9}{\mu^5} + \ldots\right)d\mu = \tag{2.3.53}$$

$$\frac{p_0^5}{3\cdot 5}\int_{\mu_{\mathbf{eff}}^{(1)}}^{\mu_{\mathbf{eff}}^{(2)}} \frac{f_{PV}^{g.m.}(\mu)d\mu}{\mu} - \frac{p_0^7}{2\cdot 3\cdot 7}\int_{\mu_{\mathbf{eff}}^{(1)}}^{\mu_{\mathbf{eff}}^{(2)}} \frac{f_{PV}^{g.m.}(\mu)d\mu}{\mu^3} + \frac{p_0^9}{8\cdot 9}\int_{\mu_{\mathbf{eff}}^{(1)}}^{\mu_{\mathbf{eff}}^{(2)}} \frac{f_{PV}^{g.m.}(\mu)d\mu}{\mu^5} + \ldots$$

**Remark 2.3.11.**We assume now that

$$\left|f_{PV}^{g.m.}(\mu)\right| = \begin{cases} O\left(\left(\mu_{\mathbf{eff}}^{(1)}\right)^{-n}\right), n > 7 & \mu_{\mathbf{eff}}^{(1)} \leq \mu \leq \mu_{\mathbf{eff}}^{(2)} \\ 0 & \mu > \mu_{\mathbf{eff}}^{(2)} \end{cases} \tag{2.3.60}$$

Note that under assumption (2.3.60) the quantities $\varepsilon^{g.m.}\left(\mu_{\mathbf{eff}}^{(1)},\mu_{\mathbf{eff}}^{(2)},p_0\right)$ and $p^{g.m.}\left(\mu_{\mathbf{eff}}^{(1)},\mu_{\mathbf{eff}}^{(2)},p_0\right)$can not contribute in the value of the cosmological constant.

# 3. Pauli-Villars ghosts as physical dark matter.

## 3.1.Pauli–Villars renormalization of the $\lambda\varphi^4$ field theory.What is wrong with Pauli-Villars renormalization of the $\lambda\varphi^4$.New physical interpretation Pauli–Villars ghost fields.

Before explaining the role of PV ghosts,etc.as physical dark matter remind the idea of PV renormalization as a conventional $UV$ renormalization.We consider, as an example, the scalar field theory with the interaction $\lambda\varphi^4$. Lagrangian density of this theory reads

$$\mathcal{L} = \frac{1}{2}\partial^\mu\varphi\partial_\mu\varphi - \frac{m_0^2}{2}\varphi^2 + \lambda\varphi^4. \tag{3.1.1}$$

This theory requires UV renormalization (e.g. in (2+1) and (3+1) dimensions). Let us show that it is sufficient to introduce $N$ extra fields with large mass playing the role of the regularization parameter. Lagrangian density can be rewritten as follows

$$\begin{gathered}\mathcal{L} = \sum_{i=0}^{N}(-1)^i\left(\frac{1}{2}\partial^\mu\varphi\partial_\mu\varphi - \frac{m_i^2}{2}\varphi_i^2\right) + \lambda : \varphi^4 :, \\ \varphi = \varphi_0 + \breve{\varphi} = \sum_{i=0}^{N}\varphi_i, \breve{\varphi} = \sum_{i=1}^{N}a_i\varphi_i.\end{gathered} \tag{3.1.2}$$

Here the symbol ”::” means that in perturbation theory we drop Feynman diagrams with loops containing only one vertex. The $\varphi_0$ is usual field with mass $m_0$ and the $\varphi_i, i = 1,\ldots,N$ is the extra field with mass $m_i, i = 1,\ldots,N.$ It can be shown that in (3+1)-dimensional theory the introduction of one PV field is sufficient for the ultraviolet regularization of perturbation theory in $\lambda$. One can show that momentum space Feynman diagrams in the original theory with Lagrangian density (2.1.1) diverge no more than quadratically [19]-[20] (beside of vacuum diagrams) shown in Fig.3.1.1.

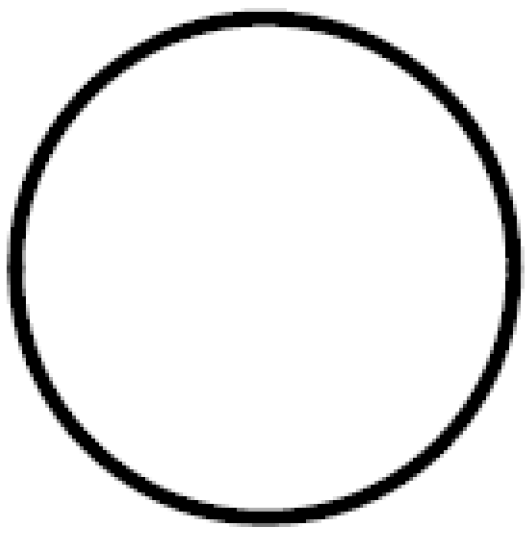

Fig.3.1.1.One-loop massive vacuum diagram.

If we consider now Feynman diagrams in the theory with Lagrangian density (3.1.2) we see that propagators of fields $\varphi_0$ and $\breve{\varphi}$ sum up in corresponding diagrams so that we obtain the following expression which plays the role of regularized propagator

$$\Delta(k^2) = \sum_{j=0}^{N} \frac{a_j}{k^2 - m_j^2 + i0} = \frac{1}{k^2 - m_0^2 + i0} - \sum_{j=1}^{N} \frac{a_j}{k^2 - m_j^2 + i0}, \tag{3.1.3}$$

where $k^2 = k_0^2 - k_1^2 + k_2^2 + k_3^2$. Integral corresponding to vacuum diagram is

$$\mathfrak{I} = \int \frac{d^4k}{(2\pi)^4} \Delta(k^2) = \int \frac{d^4k}{(2\pi)^4} \sum_{j=0}^{N} \frac{a_j}{k^2 - m_j^2 + i0}. \tag{3.1.4}$$

To do this integral, since it is convergent, we can Wick rotate.

**Remark 3.1.1**. All the integrals in quantum field theory are written in Minkowski space, however, the ultraviolet divergence appears for large values of modulus of momentum and it is useful to regularize it in Euclidean space [20].Transition to Euclidean space can be achieved by replacing the zeroth component of momentum $k_0 \to ik_4$, so that the squares of all momenta and the scalar products change the sign $k^2 = k_0^2 - \vec{k}^2 \to -k_4^2 - \vec{k}^2 = -k_E^2$ and the measure of integration becomes equal to $d^4k \to id^4k_E$, where the integration over the fourth component of momenta goes along the imaginary axis. To go to the integration along the real axis, one has to perform the (Wick) rotation of the integration contour by $90^\circ$ (see.Fig.3.1.2). This is possible since the integral over the big circle vanishes and during the transformation of the contour it does not cross the poles.

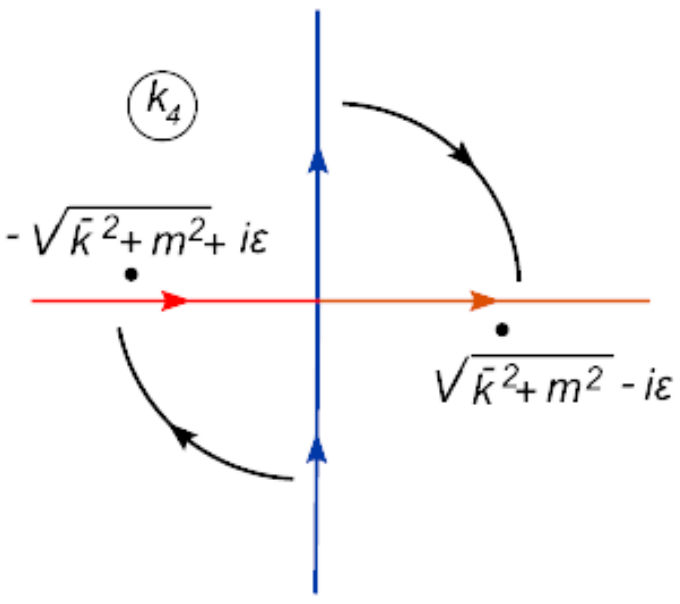

Fig.3.1.2.The Wick rotation of the integration contour.

Then we get

$$\Im_E = \frac{i}{8\pi^2}\int_0^\infty dk_E \sum\nolimits_{j=0}^N \frac{a_j k_E^3}{k_E^2+m_j^2}. \tag{3.1.5}$$

To do this integral, since it is convergent, we can dealing with regularized integral

$$\Im(\varepsilon,\Lambda) = \frac{i}{8\pi^2}\int_\varepsilon^\Lambda dk_E \sum\nolimits_{j=0}^N \frac{a_j k_E^3}{k_E^2+m_j^2}, \tag{3.1.6}$$

where $\varepsilon \asymp 0, \Lambda \asymp \infty$,i.e. $\Im(\varepsilon,\Lambda) \approx \Im_E$.We assume now that Pauli-Villars conditions given by Eqs.(2.2.39) holds.Let us consider now the quantity

$$\Im_\eta \triangleq \Im_\eta(\varepsilon,\Lambda) = \frac{i}{8\pi^2}\int_\varepsilon^\Lambda dk_E \sum\nolimits_{j=0}^N \frac{a_j k_E^3}{k_E^2+\eta m_j^2}, \tag{3.1.7}$$

where $\eta \in (0,1]$,and therefore from Eq.(3.1.7) we obtain

$$\Im_\eta|_{\eta=0} = \frac{i}{8\pi^2}\int_\varepsilon^\Lambda dk_E \sum\nolimits_{j=0}^N a_j k_E = \frac{i}{8\pi^2}\sum\nolimits_{j=0}^N a_j \int_\varepsilon^\Lambda k_E dk_E \equiv 0, \tag{3.1.8}$$

since Eqs.(2.2.39) holds.From Eq.(3.1.7) by differentiation we obtain

$$\frac{d}{d\eta}\Im_\eta = \frac{i}{8\pi^2}\int_\varepsilon^\Lambda dk_E \sum\nolimits_{j=0}^N \frac{a_j m_j^2 k_E^3}{(k_E^2+\eta m_j^2)^2}, \tag{3.1.9}$$

and therefore from Eq.(3.1.9) by using Eq.(2.2.39) we obtain

$$\begin{gathered}\left.\frac{d}{d\eta}\Im_\eta\right|_{\eta=0} = \frac{i}{8\pi^2}\int_\varepsilon^\Lambda dk_E \sum\nolimits_{j=0}^N \left.\frac{a_j m_j^2 k_E^3}{(k_E^2+\eta m_j^2)^2}\right|_{\eta=0} = \\ = \frac{i}{8\pi^2}\sum\nolimits_{j=0}^N a_j m_j^2 \int_\varepsilon^\Lambda k_E^{-1} dk_E \equiv 0,\end{gathered} \tag{3.1.10}$$

since Eqs.(2.2.39) holds.From Eq.(3.1.9) by differentiation we obtain

$$\begin{gathered}\frac{d^2}{d\eta^2}\Im_\eta = \sum\nolimits_{j=0}^N \Re_j(\eta) = \frac{i}{4\pi^2}\int_\varepsilon^\Lambda dk_E \sum\nolimits_{j=0}^N \frac{a_j m_j^4 k_E^3}{(k_E^2+\eta m_j^2)^3}, \\ \Re_j(\eta) = \frac{i a_j m_j^4}{4\pi^2}\int_\varepsilon^\Lambda dk_E \frac{k_E^3}{(k_E^2+\eta m_j^2)^3}.\end{gathered} \tag{3.1.11}$$

Note that

$$\Re_j(\eta) \simeq \frac{i a_j m_j^4}{4\pi^2}\int_0^\infty dk_E \frac{k_E^3}{(k_E^2+\eta m_j^2)^3} = \frac{i a_j m_j^4}{4\pi^2}\frac{-i}{4\eta m_j^2} = \frac{a_j m_j^2}{16\pi^2 \eta}. \tag{3.1.12}$$

Thus

$$\frac{d}{d\eta}\Im_\eta = \sum\nolimits_{j=0}^N \int_0^1 \Re_j(\eta)d\eta = \sum\nolimits_{j=0}^N \frac{a_j m_j^2}{16\pi^2}\ln\eta \tag{3.1.13}$$

and

$$\Im_\eta = \sum\nolimits_{j=0}^N \frac{a_j m_j^2}{16\pi^2}(\eta\ln\eta - \eta), \tag{3.1.14}$$

Therefore

$$\Im(\varepsilon,\Lambda) = \Im_\eta|_{\eta=1} = -\sum\nolimits_{j=0}^N \frac{a_j m_j^2}{16\pi^2} \equiv 0, \tag{3.1.15}$$

since Eqs.(2.2.39) holds.Thus integral (3.1.4) corresponding to vacuum diagram by using Pauli-Villars renormalization identically equal zero,i.e.

$$\mathbf{Ren}_{PV}(\Im) = \int \frac{d^4k}{(2\pi)^4}\Delta(k^2) = \int \frac{d^4k}{(2\pi)^4} \sum_{j=0}^{N} \frac{a_j}{k^2 - m_j^2 + i0} \equiv 0. \quad (3.1.16)$$

Let us consider now how this method works in the case of the simplest scalar diagram shown in Fig.3.1.3. The corresponding Feinman integral has the form

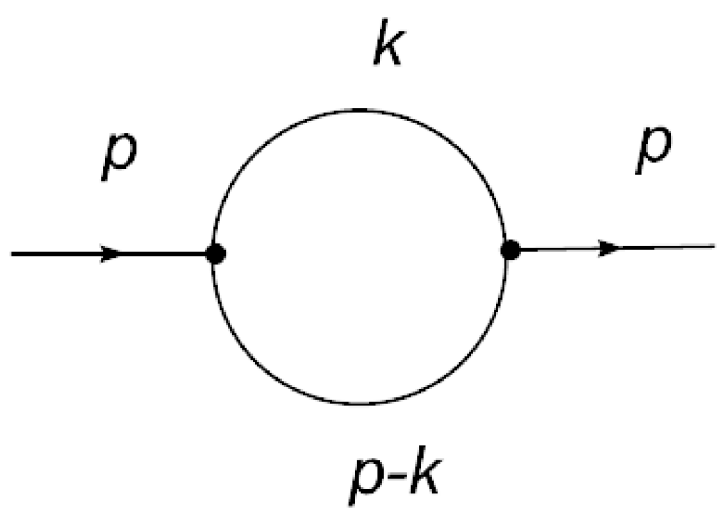

Fig.3.1.3

$$\Im(p^2) = \frac{1}{(2\pi)^4} \int \frac{d^4k}{(k^2 - m_0^2 + i0)[(p^2 - k^2) - m_0^2 + i0]}. \quad (3.1.17)$$

Regularized Feinman integral (3.1.17) reads

$$\Im_{reg}(p^2) = \frac{1}{(2\pi)^4} \int \sum_{j=0}^{N} \frac{a_j d^4k}{(k^2 - m_j^2 + i0)[(p^2 - k^2) - m_j^2 + i0]}, \quad (3.1.18)$$

where $N = 1$.To do this integral, since it is convergent, we can Wick rotate [20]. Then we get

$$\Im_{reg}(p^2) = \frac{i}{(2\pi)^4} \int \sum_{j=0}^{N} \frac{a_j d^4k}{(k^2 + m_j^2)[(p^2 - k^2) + m_j^2]}. \quad (3.1.19)$$

The integral (3.1.19) can be written as

$$\Im_{reg}(p^2) = \frac{i}{(2\pi)^4} \int_0^1 dx \int \sum_{j=0}^{N} \frac{a_j d^4k}{[k^2 + p^2x(1-x) + m_j^2]^2} =$$
$$\frac{i}{8\pi^2} \int_0^1 dx \int \sum_{j=0}^{N} \frac{a_j k_E^3 dk_E}{[k_E^2 + p^2x(1-x) + m_j^2]^2}. \quad (3.1.20)$$

To do this integral, since it is convergent, we can dealing with regularized integral

$$\Im_{reg}(p^2, \varepsilon, \Lambda) = \frac{i}{8\pi^2} \int_0^1 dx \int_{\varepsilon}^{\Lambda} \sum_{j=0}^{N} \frac{a_j k_E^3 dk_E}{[k_E^2 + p^2x(1-x) + m_j^2]^2}. \quad (3.1.21)$$

Let us consider now the quantity

$$\Im_{\eta}(p^2, \varepsilon, \Lambda) = \frac{i}{8\pi^2} \int_0^1 dx \int_{\varepsilon}^{\Lambda} \sum_{j=0}^{N} \frac{a_j k_E^3 dk_E}{[k_E^2 + p^2x(1-x) + \eta m_j^2]^2}. \quad (3.1.22)$$

where $\eta \in (0,1]$,and therefore from Eq.(3.1.22) we obtain $\Im_0(p^2, \varepsilon, \Lambda) \equiv 0$,since Eqs.(2.2.39) holds.From Eq.(3.1.22) by differentiation we obtain

$$\frac{d}{d\eta}\Im_\eta(p^2,\varepsilon,\Lambda) =$$

$$-\frac{i}{4\pi^2}\int_0^1 dx\int_\varepsilon^\Lambda \sum\nolimits_{j=0}^N \frac{a_j m_j^2 k_E^3 dk_E}{[k_E^2+p^2x(1-x)+\eta m_j^2]^3} \simeq$$

$$-\frac{i}{4\pi^2}\sum\nolimits_{j=0}^N a_j m_j^2 \Re_j(p^2,\eta,\Lambda,\varepsilon), \tag{3.1.23}$$

$$\Re_j(p^2,\eta,\Lambda,\varepsilon) \simeq \int_0^1 dx \int \frac{k_E^3 dk_E}{[k_E^2+p^2x(1-x)+\eta m_j^2]^3} = \frac{1}{4}\int_0^1 \frac{dx}{p^2x(1-x)+\eta m_j^2}.$$

From Eq.(3.1.23) we obtain

$$\frac{d}{d\eta}\Im_\eta(p^2,\varepsilon,\Lambda) \simeq -\frac{i}{4\pi^2}\sum\nolimits_{j=0}^N a_j m_j^2 \Re_j(p^2,\eta,\varepsilon,\Lambda) =$$

$$-\frac{i}{16\pi^2}\sum\nolimits_{j=0}^N a_j \int_0^1 \frac{dx}{m_j^{-2}p^2x(1-x)+\eta}. \tag{3.1.24}$$

From Eq.(2.1.24) we obtain

$$\Im_{reg}(p^2) = -\frac{i}{16\pi^2}\sum\nolimits_{j=0}^N a_j \int_0^1 dx \int_0^1 \frac{d\eta}{m_j^{-2}p^2x(1-x)+\eta}. \tag{3.1.25}$$

Note that

$$\int_0^1 \frac{d\eta}{m_j^{-2}p^2x(1-x)+\eta} =$$

$$[m_j^{-2}p^2x(1-x)+\eta]\ln[m_j^{-2}p^2x(1-x)+\eta]\Big|_0^1 - 1 = \tag{3.1.26}$$

$$[m_j^{-2}p^2x(1-x)+1]\ln[m_j^{-2}p^2x(1-x)+1] -$$

$$-[m_j^{-2}p^2x(1-x)]\ln[m_j^{-2}p^2x(1-x)] - 1.$$

Thus

$$\Im_{reg}(p^2) = -\frac{i}{16\pi^2}\sum_{j=0}^{N=1} a_j \int_0^1 dx \int_0^1 \frac{d\eta}{m_j^{-2}p^2x(1-x)+\eta} =$$

$$-\frac{i}{16\pi^2}\sum_{j=0}^{N=1} a_j \int_0^1 dx\{[m_j^{-2}p^2x(1-x)+1]\ln[m_j^{-2}p^2x(1-x)+1]-$$
$$-[m_j^{-2}p^2x(1-x)]\ln[m_j^{-2}p^2x(1-x)]\} + \frac{i}{16\pi^2}\sum_{j=0}^{N=1} a_j =$$
$$-\frac{i}{16\pi^2}\sum_{j=0}^{N=1} a_j \int_0^1 dx\{[m_j^{-2}p^2x(1-x)+1]\ln[m_j^{-2}p^2x(1-x)+1]-$$
$$-[m_j^{-2}p^2x(1-x)]\ln[m_j^{-2}p^2x(1-x)]\} = \tag{3.1.27}$$
$$-\frac{i}{16\pi^2}\int_0^1 dx\{[m_0^{-2}p^2x(1-x)+1]\ln[m_0^{-2}p^2x(1-x)+1]-$$
$$-[m_0^{-2}p^2x(1-x)]\ln[m_0^{-2}p^2x(1-x)]\} +$$
$$+\frac{i}{16\pi^2}\int_0^1 dx\{[m_1^{-2}p^2x(1-x)+1]\ln[m_1^{-2}p^2x(1-x)+1]-$$
$$-[m_1^{-2}p^2x(1-x)]\ln[m_1^{-2}p^2x(1-x)]\}.$$

From Eq.(3.1.27) we obtain

$$\Im_{reg}(p^2, m_0, m_1) =$$
$$-\frac{i}{16\pi^2}\int_0^1 dx\{[m_0^{-2}p^2x(1-x)+1]\ln[m_0^{-2}p^2x(1-x)+1]-$$
$$-[m_0^{-2}p^2x(1-x)]\ln[m_0^{-2}p^2x(1-x)]\} + \tag{3.1.28}$$
$$+\frac{i}{16\pi^2}\int_0^1 dx\{[m_1^{-2}p^2x(1-x)+1]\ln[m_1^{-2}p^2x(1-x)+1]-$$
$$-[m_1^{-2}p^2x(1-x)]\ln[m_1^{-2}p^2x(1-x)]\}.$$

We assume now that $m_1^{-2}p^2 \ll 1$ and from Eq.(3.1.28) finally we obtain

$$\Im_{reg}(p^2, m_0, m_1) =$$
$$-\frac{i}{16\pi^2}\int_0^1 dx\{[m_0^{-2}p^2x(1-x)+1]\ln[m_0^{-2}p^2x(1-x)+1]- \tag{3.1.29}$$
$$-[m_0^{-2}p^2x(1-x)]\ln[m_0^{-2}p^2x(1-x)]\} + O(m_1^{-2}p^2).$$

**Remark 3.1.2.** Note that by taking the limit: $\lim_{m_1\to\infty}\Im_{reg}(p^2, m_0, m_1)$ from Eq.(3.1.28) one

obtains instead Eq.(3.1.29) the following purely formal result:

$$\Im_{reg}(p^2,m_0) \triangleq \lim_{m_1\to\infty}\Im_{reg}(p^2,m_0,m_1) =$$

$$-\frac{i}{16\pi^2}\int_0^1 dx\{[m_0^{-2}p^2x(1-x)+1]\ln[m_0^{-2}p^2x(1-x)+1]- \tag{3.1.30}$$

$$-[m_0^{-2}p^2x(1-x)]\ln[m_0^{-2}p^2x(1-x)]\}.$$

However this result completely wrong and mathematically is not sound!!!

**Remark 3.1.3**.Note that in the limit: $m_i \to \infty, i = 1,..,N$ the expansions

$$\frac{a_i}{k^2-m_i^2+i\epsilon} = \frac{a_i}{k^2+i\epsilon}\frac{1}{1-\dfrac{m_i^2}{k^2+i\epsilon}} = \frac{a_i}{k^2+i\epsilon}\times$$

$$\left[1+\frac{m_i^2}{k^2+i\epsilon}+\frac{m_i^4}{\left(k^2+i\epsilon\right)^2}+\ldots\right] = \tag{3.1.31}$$

$$\frac{a_i}{k^2+i\epsilon}+\frac{a_im_i^2}{\left(k^2+i\epsilon\right)^2}+\frac{a_im_i^4}{\left(k^2+i\epsilon\right)^3}+\ldots,$$

obviously breaks down (see Remark 2.2.8) and therefore the Eq.(3.1.32) is not holds

$$\Delta(k^2,,m_0,m_1,\ldots,m_N) =$$

$$\sum\nolimits_{i=0}^{N}\frac{a_i}{k^2+i\epsilon}+\sum\nolimits_{i=0}^{N}\frac{a_im_i^2}{\left(k^2+i\epsilon\right)^2}+\sum\nolimits_{i=0}^{N}O\left(\frac{1}{\left(k^2+i\epsilon\right)^3}\right). \tag{3.1.32}$$

But therefore the Eq.(3.1.29) also is not holds and Pauli–Villars procedure completely breaks down.

**Remark 3.1.4**.It follows from consideration above that we cannot deleted from final result

Pauli–Villars masses $m_i, i = 1,\ldots,N$ which appeas in Lagrangian density (3.1.1).

**Remark 3.1.5**. Note that one can dealing instead regularized integral (3.1.21) with the following Colombeau integral

$$(\Im_{reg}(p^2,\varepsilon,\delta))_\varepsilon = \frac{i}{8\pi^2}\left(\int_0^1 dx\int_\varepsilon^{1/\delta}\sum\nolimits_{j=0}^{N}\frac{a_jk_E^3dk_E}{[k_E^2+p^2x(1-x)+m_j^2]^2}\right)_\delta, \tag{3.1.32}$$

$\delta \in (0,1]$, see [21]-[26].

## 3.2.Pauli–Villars renormalization of $QED_3$.What is wrong with Pauli–Villars renormalization of $QED_3$.New physical interpretation Pauli–Villars ghost fields.

### 3.2.1.Pauli–Villars renormalization of $QED_3$.What is wrong with Pauli–Villars renormalization of $QED_3$.

Let us consider now the conditions that must be required on the masses and coupling constants of the regulator fields such that a regularized closed fermion loop in $2+1$ dimensions is rendered finite in the calculations.Thus corresponding Feynman integral is proportional to

$$\int \mathcal{L}(p,m)d^3p = $$
$$\int d^3p \frac{\mathbf{Tr}[\gamma_{\mu_1}(m+\not p)\gamma_{\mu_2}(m+\not p+\not k_1)\ldots\gamma_{\mu_n}(m+\not p+\ldots+\not k_{n-1})]}{(m^2-p^2+i\epsilon)\left[m^2-(p+k_1)^2+i\epsilon\right]\ldots\left[m^2-(p+k_{n-1})^2+i\epsilon\right]}, \quad (3.2.1)$$

so,for large $|p|$,its integrand behaves like $|p|^{-n}$, whereas for $n < 4$ the integral (3.2.1) diverges as $\int_0^\infty \frac{p^2dp}{p^n} \sim \int_0^\infty \frac{dp}{p^{n-2}}$. We apply now the momentum $\Lambda$-cutoff regularization $|p| \leq \Lambda$ and have to replaced ill defined formal expression (3.2.1) by the well defined Integral

$$\Im_\Lambda(k_1,\ldots,k_{n-1},m) = \int_\Lambda \mathcal{L}(p,m)d^3p \triangleq$$
$$\int_{|p|\leq\Lambda} d^3p \frac{\mathbf{Tr}[\gamma_{\mu_1}(m+\not p)\gamma_{\mu_2}(m+\not p+\not k_1)\ldots\gamma_{\mu_n}(m+\not p+\ldots+\not k_{n-1})]}{(m^2-p^2+i\epsilon)\left[m^2-(p+k_1)^2+i\epsilon\right]\ldots\left[m^2-(p+k_{n-1})^2+i\epsilon\right]}. \quad (3.2.2)$$

**Remark 3.2.1.**Note that $\Lambda$-cutoff regularization meant lattice $QED_3$ on a lattice at length
scales $a \sim \Lambda^{-1}$. However quantity $\Im_\Lambda$ behaves as $\sim \int_\Lambda p^{-n+2}dp$ and therefore gives unphysical result which strictly depends on parametr $\Lambda$.

**Remark 3.2.2.**In order obtain physical result we apply now Pauli–Villars renormalization.

The integrand $\mathcal{L}(p,m)$ in (3.2.2) can be written as

$$\mathcal{L}(p,m) = \sum_k m^k a_{n+k} \sim p^{-(n+k)}. \quad (3.2.3)$$

Therefore, in making the canonical substitution $\mathcal{L}(p,m) \to \sum_{i=0}^{n_s} c_i\mathcal{L}(p,M_i)$ where

$$M_i^2/|p|^2 \ll 1, \quad (3.2.4)$$

and where $n_s$ is the number of auxiliary spinor fields, we must impose in the vacuum polarization case ($n = 2$) the following conditions

$$(a)\sum_{i=0}^{n_s} c_i \equiv 0, (b)\sum_{i=0}^{n_s} c_iM_i \equiv 0. \quad (3.2.5)$$

in order to eliminate the linear and logarithmic divergences, respectively. Let us calculate of the vacuum polarization tensor in spinor $QED_3$. In the standard notation, the regularized expression for the vacuum polarization tensor reads

$$\prod\nolimits_{\mu\nu}^M(k,\Lambda) = \frac{ie^2}{(2\pi)^3}\sum_{i=0}^{n_s} c_i \int_\Lambda d^3p \frac{P(M_i)}{(M_i^2-p_1^2)(M_i^2-p_2^2)}, \quad (3.2.6)$$

where $c_0 = 1, M_0 = m, M_i = m\lambda_i, i = 1,\ldots,n_f, p_{1,2} = p \mp \frac{1}{2}k$ and

$$P(M_i) = \mathbf{Tr}[\gamma_\mu(\not p_1 + M_i)\gamma_\nu(\not p_2 + M_i)] =$$
$$2[M_i^2 g_{\mu\nu} + p_{1\mu}p_{2\nu} + p_{1\nu}p_{2\mu} - g_{\mu\nu}(p_1 \cdot p_2) - iM_i\epsilon_{\mu\nu\alpha}k^\alpha]. \quad (3.2.7)$$

We choose now both the electron mass and that of the auxiliary field $M_1$ to be positive. Using now the canonical Feynman parametrization $[(M_i^2-p_1^2)(M_i^2-p_2^2)]^{-1} = \int_0^1 d\xi[M_i^2 - p_1^2 - (p_2^2-p_1^2)\xi]^{-1}$and performing the momentum shift $p_\mu \to p_\mu(1/2-\xi)k_\mu$ we obtain

$$\Pi^M_{\mu\nu}(k,\Lambda) = \left(g_{\mu\nu} - \frac{k_\mu k_\nu}{k^2}\right)\Pi^M_1(k^2,\Lambda) + im\epsilon_{\mu\nu\alpha}k^\alpha\Pi^M_2(k^2,\Lambda) + \Pi^M_3(k^2,\Lambda), \tag{3.2.8}$$

where

$$\begin{gathered}
\Pi^M_1(k^2,\Lambda) = 4ie^2k^2\sum_{i=0}^{n_s} c_i\int_0^1 d\xi\xi(1-\xi)\int_\Lambda \frac{d^3p}{(2\pi)^3}\frac{1}{\left(Q_i^2(\xi)-p^2\right)^2},\\
\Pi^M_2(k^2,\Lambda) = -\frac{2ie^2}{m}\sum_{i=0}^{n_s} c_i\int_0^1 d\xi\int_\Lambda \frac{d^3p}{(2\pi)^3}\frac{1}{\left(Q_i^2(\xi)-p^2\right)^2},\\
\Pi^M_3(k^2,\Lambda) = \frac{2}{3}ie^2g_{\mu\nu}\sum_{i=0}^{n_s} c_i\left(I^{(1)}_{i,\Lambda} + I^{(2)}_{i,\Lambda}\right),\\
I^{(1)}_{i,\Lambda} = 3\int_0^1 d\xi\int_\Lambda \frac{d^3p}{(2\pi)^3}\frac{1}{\left(Q_i^2(\xi)-p^2\right)}, I^{(2)}_i = 2\int_0^1 d\xi\int_\Lambda \frac{d^3p}{(2\pi)^3}\frac{1}{\left(Q_i^2(\xi)-p^2\right)^2},\\
Q_i^2(\xi) = M_i^2 - \xi(1-\xi)k^2.
\end{gathered} \tag{3.2.9}$$

If we carry out the momentum integrations in Eqs.(3.2.9), it is straightforward to arrive at $\Pi^M_3(k^2) = 0$, as expected by gauge invariance.Let us take now $c_1 = \alpha - 1, c_2 = -\alpha, c_{j=0}, j > 2,$ where the parameter $\alpha$ can assume any real value except zero and unity, so that condition given by Eq.(3.2.5.a) is satisfied.From Eqs.(3.2.9) for instance we find for a sufficiently large value of the parametr $\Lambda$ :

$$\begin{gathered}
\Pi^M_2(k^2,\Lambda) \simeq\\
\frac{e^2}{4\pi m}\int_0^1 d\xi\frac{m}{\left[m^2-\xi(1-\xi)k^2\right]^{1/2}} + \frac{e^2}{4\pi m}\int_0^1 d\xi\frac{(\alpha-1)M_1}{\left[M_1^2-\xi(1-\xi)k^2\right]^{1/2}} +\\
+\frac{e^2}{4\pi m}\int_0^1 d\xi\frac{\alpha M_2}{\left[M_2^2-\xi(1-\xi)k^2\right]^{1/2}} = \frac{e^2}{4\pi m}\int_0^1 d\xi\frac{m}{\left[m^2-\xi(1-\xi)k^2\right]^{1/2}} +\\
\mathcal{F}_1(M_1,\Lambda) + \mathcal{F}_2(M_2,\Lambda),
\end{gathered} \tag{3.2.10}$$

where

$$\begin{gathered}
\mathcal{F}_1(M_1,\Lambda) \simeq \frac{e^2}{4\pi m}\int_0^1 d\xi\frac{(\alpha-1)M_1}{\left[M_1^2-\xi(1-\xi)k^2\right]^{1/2}},\\
\mathcal{F}_2(M_2,\Lambda) \simeq -\frac{e^2}{4\pi m}\int_0^1 d\xi\frac{\alpha M_2}{\left[M_2^2-\xi(1-\xi)k^2\right]^{1/2}}.
\end{gathered} \tag{3.2.11}$$

Note that for $M_1, M_2 \gg m$ from Eq.(3.2.5.b) follows that $M_2 \simeq M_1(1-\alpha^{-1})$ and therefore from Eq.(3.2.11) for $k^2/M_1 \ll 1$ we obtain

$$\begin{gathered}
\mathcal{F}_1(M_1,\Lambda) \simeq \frac{e^2}{4\pi m}\int_0^1 d\xi\frac{(\alpha-1)}{\left[1-\xi(1-\xi)\left(\frac{k}{M_1}\right)^2\right]^{1/2}} =\\
= \frac{e^2(\alpha-1)}{4\pi m} + \frac{e^2}{4\pi m}\frac{(\alpha-1)}{2}\left(\frac{k}{M_1}\right)^2\int_0^1 d\xi\xi(1-\xi) + O\left(\frac{k^4}{M_1^4}\right)
\end{gathered} \tag{3.2.12}$$

and

$$\mathcal{F}_2(M_2,\Lambda) \simeq -\frac{e^2}{4\pi m}\int_0^1 d\xi \frac{\alpha(sign(1-\alpha^{-1}))}{\left[1-\xi(1-\xi)\left(\frac{k}{M_2}\right)^2\right]^{1/2}} =$$

$$= -\frac{e^2\alpha(sign(1-\alpha^{-1}))}{4\pi m} -$$

$$-\frac{e^2\alpha(sign(1-\alpha^{-1}))}{8\pi m}\left(\frac{k}{M_1(1-\alpha^{-1})}\right)^2\int_0^1 d\xi\xi(1-\xi) \tag{3.2.13}$$

$$+O\left(\frac{k^4}{M_1^4(1-\alpha^{-1})^4}\right).$$

From Eq.(3.2.8) and Eq.(3.2.10)-(3.2.13) we obtain

$$\Pi_2^M(k^2,\Lambda) \simeq$$

$$\frac{e^2}{4\pi m}\int_0^1 d\xi \frac{m}{[m^2-\xi(1-\xi)k^2]^{1/2}} + \frac{e^2(\alpha-1)}{4\pi m} - \frac{e^2\alpha(sign(1-\alpha^{-1}))}{4\pi m} +$$

$$+\frac{e^2}{4\pi m}\frac{(\alpha-1)}{2}\left(\frac{k}{M_1}\right)^2\int_0^1 d\xi\xi(1-\xi) - \tag{3.2.14}$$

$$-\frac{e^2\alpha(sign(1-\alpha^{-1}))}{8\pi m}\left(\frac{k}{M_1(1-\alpha^{-1})}\right)^2\int_0^1 d\xi\xi(1-\xi) + O\left(\frac{k^4}{M_1^4}\right).$$

For $\alpha = 1/2$ we obtain

$$\Pi_2^M(k^2,\Lambda) \simeq$$

$$\frac{e^2}{4\pi m}\int_0^1 d\xi \frac{m}{[m^2-\xi(1-\xi)k^2]^{1/2}} - \frac{1}{2}\frac{e^2}{4\pi m} + \frac{1}{2}\frac{e^2}{4\pi m} +$$

$$-\frac{e^2}{16\pi m}\left(\frac{k}{M_1}\right)^2\int_0^1 d\xi\xi(1-\xi) - \tag{3.2.15}$$

$$+\frac{e^2}{16\pi m}\left(\frac{k}{M_1}\right)^2\int_0^1 d\xi\xi(1-\xi) + O\left(\frac{k^4}{M_1^4}\right) =$$

For instance for $k = 0$ and $\alpha = 1/2$ one obtains

$$\Pi_2^M(0) = \frac{e^2}{4\pi m}. \tag{3.2.16}$$

**Remark 3.2.3.**Note that when using canonical QFT in order to recover the original theory
without physical ghost matterthe absolute values of the coefficients $\lambda_i$ ultimately alwais go
to infinity $\lambda_i \to \infty$. However under limit: $\lambda_i \to \infty$ the inequality (3.2.4) is not holds and therefore Pauli–Villars procedure completely breaks down!

## 3.2.2.New physical interpretation Pauli–Villars ghost fields.

**Remark 3.2.4.**These PV-renormalizable models with finite mases $m_i, i = 1,\ldots,N$, which
we have considered in this section many years regarded only as constructs for a study of
the ultraviolet problem of QFT. The difficulties with unitarity appear to preclude their

direct
acceptability as canonical physical theories in locally Minkowski space-time. However, for
their unphysical behavior may be restricted to *arbitrarily large energy scales* $\Lambda_*$
mentioned above by an appropriate limitation on the finite masses $m_i$ such that

$$m \ll m_i < \Lambda_*, i = 1,\ldots,N. \qquad (3.2.17)$$

**Remark 3.2.5.** We have aplied now the Colombeau approach [21]-[26] and replaced ill defined formal expression (3.2.1) by well defined Colombeau generalized function [21]-[26]:

$$(\Im_\varepsilon(k_1,\ldots,k_{n-1},m))_\varepsilon = \left(\int_{-1/\varepsilon}^{1/\varepsilon} \mathcal{L}(p,k_1,\ldots,k_{n-1},m)d^3p\right)_\varepsilon =$$
$$\int_0^{1/\varepsilon} d^3p \frac{\mathbf{Tr}[\gamma_{\mu_1}(m+\not{p})\gamma_{\mu_2}(m+\not{p}+\not{k}_1)\ldots\gamma_{\mu_n}(m+\not{p}+\ldots+\not{k}_{n-1})]}{(m^2-p^2+i\epsilon)\left[m^2-(p+k_1)^2+i\epsilon\right]\ldots\left[m^2-(p+k_{n-1})^2+i\epsilon\right]}, \qquad (3.2.18)$$

where $\varepsilon \in (0,1]$.

**Remark 3.2.6.** Note that $[(\mathcal{L}_\varepsilon(k_1,\ldots,k_{n-1},m))_\varepsilon] \in \mathcal{G}(\mathbb{R}^{4(n-1)})$ and for any (n-1)-tuple $(\bar{k}_1,\ldots,\bar{k}_{n-1})$ obviously $(\Im_\varepsilon(\bar{k}_1,\ldots,\bar{k}_{n-1},m))_\varepsilon \in \widetilde{\mathbb{R}}$ (see [23]) and any $(\Im_\varepsilon(\bar{k}_1,\ldots,\bar{k}_{n-1},m))_\varepsilon$ are infinite large Colombeau generalized numbers [23],[26].

**Remark 3.2.7.** The integrand $\mathcal{L}(p,k_1,\ldots,k_{n-1},m)$ in (3.2.18) for any $m = (m_{i,\varepsilon})_\varepsilon \in \widetilde{\mathbb{R}}$ can be
written as

$$\mathcal{L}(p,k_1,\ldots,k_{n-1},m_{i,\varepsilon}) =$$
$$\frac{P_n(p)+m_{i,\varepsilon}P_{n-1}(p)+m_{i,\varepsilon}^2P_{n-2}(p)+\ldots+m_{i,\varepsilon}^n}{P_{2n}(p)+m_{i,\varepsilon}^2P_{n-2}(p)+\ldots+m_{i,\varepsilon}^{2n}}, \qquad (3.2.19)$$

where $P_i(p)$ stands for a polynomial of degree $i$ in the components of $p$. We can write the
denominator of $\mathcal{L}(p,k_1,\ldots,k_{n-1},m_{i,\varepsilon})$ in the following form

$$P_{n-2}(p)\left(1+m_{i,\varepsilon}^2\frac{P_{2n-2}(p)}{P_{2n}(p)}+\ldots+m_{i,\varepsilon}^{2n}\frac{1}{P_{2n}(p)}\right) \qquad (3.2.20)$$

and, for sufficiently large $p = p(\varepsilon) \gg m_{i,\varepsilon}, i = 0,1,\ldots,n_s$ and for any $\varepsilon \in (0,1]$ perform the
expansions

$$\left(1+m_{i,\varepsilon}^2\frac{P_{2n-2}(p)}{P_{2n}(p)}+\ldots+m_{i,\varepsilon}^{2n}\frac{1}{P_{2n}(p)}\right)^{-1} =$$
$$1-\left(m_{i,\varepsilon}^2\frac{P_{2n-2}(p)}{P_{2n}(p)}+\ldots+m_{i,\varepsilon}^{2n}\frac{1}{P_{2n}(p)}\right) \qquad (3.2.21)$$
$$-\left(m_i^2\frac{P_{2n-2}(p)}{P_{2n}(p)}+\ldots+m_i^{2n}\frac{1}{P_{2n}(p)}\right)^2-\ldots$$

so that the integrand $\mathcal{L}(p,k_1,\ldots,k_{n-1},m_{i,\varepsilon})$ behaves like

$$\mathcal{L}(p,k_1,\ldots,k_{n-1},m_{i,\varepsilon}) \asymp \frac{P_n(p)}{P_{2n}(p)} + m_{i,\varepsilon}\frac{P_{n-1}(p)}{P_{2n}(p)} +$$
$$+m_{i,}^2\frac{P_n(p)}{P_{2n}(p)}\left[\frac{P_{n-2}(p)}{P_n(p)} - \frac{P_{2n-2}(p)}{P_{2n}(p)}\right] + m_{i,\varepsilon}^3\frac{P_{n-1}(p)}{P_{2n}(p)}\left[\frac{P_{n-3}(p)}{P_{n-1}(p)} - \frac{P_{2n-2}(p)}{P_{2n}(p)}\right] \asymp \tag{3.2.22}$$
$$\asymp \sum\nolimits_k m_{i,\varepsilon}^k p^{-(n+k)}.$$

Therefore, in making the substitution in Colombeau integral (3.2.18)

$$\mathcal{L}(p,k_1,\ldots,k_{n-1},m) \to \sum\nolimits_{i=0}^{n_s} \mathcal{L}(p,k_1,\ldots,k_{n-1},m_{i,\varepsilon}), \tag{3.2.23}$$

where $(m_{0,\varepsilon})_\varepsilon \equiv m \in \mathbb{R}$ for any $\varepsilon \in (0,1]$ and where $n_s$ is the number of auxiliary spinor fields, we obtain

$$(\Im_\varepsilon(k_1,\ldots,k_{n-1},m))_\varepsilon \to (\Im_\varepsilon(k_1,\ldots,k_{n-1},\{m_{i,\varepsilon}\}_{i=0}^{n_s}))_\varepsilon, \tag{3.2.24}$$

where

$$(\Im_\varepsilon(k_1,\ldots,k_{n-1},\{m_i\}_{i=0}^{n_s}))_\varepsilon =$$
$$= \sum\nolimits_{i=0}^{n_s} c_i(\Im_\varepsilon(k_1,\ldots,k_{n-1},m_{i,\varepsilon}))_\varepsilon = \left(\sum\nolimits_{i=0}^{n_s}\int_{-1/\varepsilon}^{1/\varepsilon}\mathcal{L}(p,k_1,\ldots,k_{n-1},m_{i,\varepsilon})d^3p\right)_\varepsilon = \tag{3.2.25}$$
$$\left(\int_0^{1/\varepsilon}\mathcal{L}(p,k_1,\ldots,k_{n-1},m_0)d^3p\right)_\varepsilon + \left(\int_{-1/\varepsilon}^{1/\varepsilon}\sum\nolimits_{i=1}^{n_s}\mathcal{L}(p,k_1,\ldots,k_{n-1},m_{i,\varepsilon})d^3p\right)_\varepsilon$$

**Remark 3.2.8.**Note that Colombeau integral $(\Im_\varepsilon(\bar{k}_1,\ldots,\bar{k}_{n-1},\{m_i\}_{i=0}^{n_s}))_\varepsilon$ for any (n-1)-tuple
$(\bar{k}_1,\ldots,\bar{k}_{n-1})$ is infinite large Colombeau generalized number, however we can impose the
following conditions:

$$(a)\sum\nolimits_{i=0}^{n_s} c_i \equiv 0, (b)\sum\nolimits_{i=0}^{n_s} c_i(m_{i,\varepsilon})_\varepsilon \equiv 0,\ldots,\sum\nolimits_{i=0}^{n_s} c_i(m_{i,\varepsilon}^q)_\varepsilon \equiv 0,$$
$$c_0 = 1, (m_{0,\varepsilon})_\varepsilon \equiv m \in \mathbb{R}, \tag{3.2.26}$$

in order to eliminate the infinite large Colombeau generalized quantities ~ $(\ln\varepsilon^{-1})_\varepsilon$, $(\varepsilon^{-1})_\varepsilon,\ldots,(\varepsilon^{-q})$, respectively from Colombeau integral (3.2.25),i.e., make it finite in canonical sense.

Using now the canonical Feynman parametrization in the following Colombeau form

$$\left([(M_{i,\varepsilon}^2 - p_1^2)(M_{i,\varepsilon}^2 - p_2^2)]^{-1}\right)_\varepsilon = \left(\int_0^1 d\xi[M_{i,\varepsilon}^2 - p_1^2 - (p_2^2 - p_1^2)\xi]^{-1}\right)_\varepsilon \tag{3.2.27}$$

and performing the momentum shift $p_\mu \to p_\mu(1/2 - \xi)k_\mu$ we obtain instead Eq.(3.2.8)

$$(\Pi_{\mu\nu}^M(k,\varepsilon))_\varepsilon = \left(g_{\mu\nu} - \frac{k_\mu k_\nu}{k^2}\right)(\Pi_1^M(k^2,\varepsilon))_\varepsilon + im\epsilon_{\mu\nu\alpha}k^\alpha(\Pi_2^M(k^2,\varepsilon))_\varepsilon + (\Pi_3^M(k^2,\varepsilon))_\varepsilon, \tag{3.2.28}$$

where

$$
\begin{aligned}
&(\Pi_1^M(k^2,\varepsilon))_\varepsilon = 4ie^2k^2\sum_{i=0}^{n_s} c_i\left(\int_0^1 d\xi\xi(1-\xi)\int_{-1/\varepsilon}^{1/\varepsilon}\frac{d^3p}{(2\pi)^3}\frac{1}{(Q_{i,\varepsilon}^2(\xi)-p^2)^2}\right)_\varepsilon,\\
&(\Pi_2^M(k^2,\varepsilon))_\varepsilon = -\frac{2ie^2}{m}\sum_{i=0}^{n_s} c_i\left(\int_0^1 d\xi\int_{-1/\varepsilon}^{1/\varepsilon}\frac{d^3p}{(2\pi)^3}\frac{1}{(Q_{i,\varepsilon}^2(\xi)-p^2)^2}\right)_\varepsilon,\\
&(\Pi_3^M(k^2,\varepsilon))_\varepsilon = \frac{2}{3}ie^2g_{\mu\nu}\sum_{i=0}^{n_s} c_i\left(\left(I_{i,\varepsilon}^{(1)}\right)_\varepsilon+\left(I_{i,\varepsilon}^{(2)}\right)_\varepsilon\right),\\
&\left(I_{i,\varepsilon}^{(1)}\right)_\varepsilon = 3\left(\int_0^1 d\xi\int_{-1/\varepsilon}^{1/\varepsilon}\frac{d^3p}{(2\pi)^3}\frac{1}{(Q_{i,\varepsilon}^2(\xi)-p^2)}\right)_\varepsilon,\\
&\left(I_{i,\varepsilon}^{(2)}\right)_\varepsilon = 2\left(\int_0^1 d\xi\int_{-1/\varepsilon}^{1/\varepsilon}\frac{d^3p}{(2\pi)^3}\frac{1}{(Q_{i,\varepsilon}^2(\xi)-p^2)^2}\right)_\varepsilon,\\
&(Q_{i,\varepsilon}^2(\xi))_\varepsilon = (M_{i,\varepsilon}^2)_\varepsilon-\xi(1-\xi)k^2.
\end{aligned}
\tag{3.2.29}
$$

Carry out the momentum integrations in Eqs.(3.2.29), it is straightforward to arrive at $(\Pi_3^M(k^2,\varepsilon))_\varepsilon = 0,$ as expected by gauge invariance.Let us take now $c_1 = \alpha-1, c_2 = -\alpha,$ $c_{j=0}, j > 2,$ where the parameter $\alpha$ can assume any real value except zero and unity, so that condition given by Eq.(3.2.26.a) is satisfied.From Eqs.(3.2.29) for instance we find

$$
\begin{aligned}
&(\Pi_2^M(k^2,\varepsilon))_\varepsilon \approx_{\widetilde{\mathbb{R}}}\\
&\frac{e^2}{4\pi m}\int_0^1 d\xi\frac{m}{[m^2-\xi(1-\xi)k^2]^{1/2}}+\frac{e^2}{4\pi m}\left(\int_0^1 d\xi\frac{(\alpha-1)M_{1,\varepsilon}}{[M_{1,\varepsilon}^2-\xi(1-\xi)k^2]^{1/2}}\right)_\varepsilon+\\
&+\frac{e^2}{4\pi m}\left(\int_0^1 d\xi\frac{\alpha M_{2,\varepsilon}}{[M_{2,\varepsilon}^2-\xi(1-\xi)k^2]^{1/2}}\right)_\varepsilon = \frac{e^2}{4\pi m}\int_0^1 d\xi\frac{m}{[m^2-\xi(1-\xi)k^2]^{1/2}}+\\
&(\mathcal{F}_1(M_{1,\varepsilon},\varepsilon))_\varepsilon+(\mathcal{F}_2(M_{2,\varepsilon},\varepsilon))_\varepsilon,
\end{aligned}
\tag{3.2.30}
$$

where

$$
\begin{aligned}
&(\mathcal{F}_1(M_{1,\varepsilon},\varepsilon))_\varepsilon \approx_{\widetilde{\mathbb{R}}} \frac{e^2}{4\pi m}\left(\int_0^1 d\xi\frac{(\alpha-1)M_{1,\varepsilon}}{[M_{1,\varepsilon}^2-\xi(1-\xi)k^2]^{1/2}}\right)_\varepsilon,\\
&(\mathcal{F}_2(M_{2,\varepsilon},\varepsilon))_\varepsilon \approx_{\widetilde{\mathbb{R}}} -\frac{e^2}{4\pi m}\left(\int_0^1 d\xi\frac{\alpha M_{2,\varepsilon}}{[M_{2,\varepsilon}^2-\xi(1-\xi)k^2]^{1/2}}\right)_\varepsilon.
\end{aligned}
\tag{3.2.31}
$$

If we chose now infinite large Pauli-Villars masses $(M_{1,\varepsilon})_\varepsilon \in \widetilde{\mathbb{R}}, (M_{2,\varepsilon})_\varepsilon \in \widetilde{\mathbb{R}}\backslash\mathbb{R}$ we obtain

$$
\begin{aligned}
&(\mathcal{F}_1(M_{1,\varepsilon},\varepsilon))_\varepsilon \approx_{\widetilde{\mathbb{R}}}\\
&\frac{e^2}{4\pi m}\left(\int_0^1 d\xi(\alpha-1)\left[1-\xi(1-\xi)\left(\frac{k}{M_{1,\varepsilon}}\right)^2\right]^{-1/2}\right)_\varepsilon \approx_{\widetilde{\mathbb{R}}} \frac{e^2(\alpha-1)}{4\pi m}
\end{aligned}
\tag{3.2.32}
$$

and

$$(\mathcal{F}_2(M_2,\varepsilon))_\varepsilon \approx_{\widetilde{\mathbb{R}}}$$

$$-\frac{e^2}{4\pi m}\left(\int_0^1 d\xi\alpha(sign(1-\alpha^{-1}))\left[1-\xi(1-\xi)\left(\frac{k}{M_{2,\varepsilon}}\right)^2\right]^{-1/2}\right)_\varepsilon \approx_{\widetilde{\mathbb{R}}} \tag{3.2.33}$$

$$\approx_{\widetilde{\mathbb{R}}} -\frac{e^2\alpha(sign(1-\alpha^{-1}))}{4\pi m}.$$

From Eq.(3.2.30) and Eq.(3.2.32)-(3.2.33) we obtain canonical result

$$ren(\Pi_2^M(k^2)) = (\Pi_2^M(k^2,\varepsilon))_\varepsilon \approx_{\widetilde{\mathbb{R}}}$$

$$\frac{e^2}{4\pi}\left(\int_0^1 \frac{d\xi}{[m^2-\xi(1-\xi)k^2]^{1/2}}\right) + \frac{e^2(\alpha-1)}{4\pi m} - \frac{e^2\alpha(sign(1-\alpha^{-1}))}{4\pi m} \in \mathbb{R}. \tag{3.2.34}$$

**Remark 3.2.9.** Note that in order to obtain finite physical result (3.2.33) using infinite large Pauli-Villars masses $(M_{1,\varepsilon})_\varepsilon, (M_{2,\varepsilon})_\varepsilon \in \widetilde{\mathbb{R}}\backslash\mathbb{R}$ mentioned above, one needs write down a Pauli-Villars Lagrangian density for $QED_3$, which works for instance at the $1$-loop level, as

$$(\mathcal{L}_{ren}^{PV}(\varepsilon))_\varepsilon = -\frac{1}{4}(F_{\mu\nu,\varepsilon}^2)_\varepsilon +$$

$$\left(\overline{\psi}_\varepsilon\left(i\partial\!\!\!/ - eA\!\!\!/_\varepsilon - eA_{\mu,\varepsilon}^{gh} - m\right)\psi_\varepsilon\right)_\varepsilon + \frac{1}{4}\left(\left(F_{\mu\nu,\varepsilon}^{gh}\right)_\varepsilon\right)^2 - \frac{1}{2}\left(\left(M_{1,\varepsilon}A_{\mu,\varepsilon}^{gh}\right)_\varepsilon\right)^2 + \tag{3.2.35}$$

$$\left(\overline{\psi}_\varepsilon^{gh}\left(i\partial\!\!\!/ - eA\!\!\!/_\varepsilon - eA_{\mu,\varepsilon}^{gh} - M_{2,\varepsilon}\right)\psi_\varepsilon^{gh}\right)_\varepsilon$$

with $\left(A_{\mu,\varepsilon}^{gh}\right)_\varepsilon$ the ghost photon and $\left(\psi_\varepsilon^{gh}\right)_\varepsilon$ the ghost electron and $\left(F_{\mu\nu,\varepsilon}^{gh}\right)_\varepsilon = \left(\partial_\mu A_{\nu,\varepsilon}^{gh}\right)_\varepsilon - \left(\partial_\nu A_{\mu,\varepsilon}^{gh}\right)_\varepsilon$. We assume in (3.2.35) that both the ghost photon and ghost electron have bosonic statistics and the ghost photon has a wrong-sign kinetic term,

**Remark 3.2.10.** In contast with canonical interpetation of the renormalization constant as formal infinite quantities, in this paper we argue that only finite physical quantity can appears in physical Pauli-Villars Lagrangian density and therefore Pauli-Villars masses $(M_{1,\varepsilon})_\varepsilon, (M_{2,\varepsilon})_\varepsilon$ can be choosen arbitraly large but finite. Thus true physical result is given by Eq.(3.2.15).

## 3.3. High covariant derivatives renormalization as Pauli–Villars renormalization of non-Abelian gauge theories. New physical interpretation.

The standard approach to regularization of non-Abelian gauge theories is dimensional regularization but this of course is inherently perturbative [27]. However, the ordinary PV renormalization of non-Abelian theories fails. Gauge invariance is violated, is blocking any hope of BRST invariance, which confounds proofs of renormalizability.

**Remark 3.3.1.** Note that the existence of interesting non-perturbative phenomena in gauge theories requires the introduction of a non-perturbative regularization. Discretization of space-time leads in a natural way to lattice regularizations which preserve gauge invariance and have a non-perturbative meaning. The construction of a non-perturbative gauge invariant regularization of gauge theories in a continuum space-time has been a challenging problem in gauge theories. A natural candidate has always been a gauge invariant generalization of Pauli-Villars methods involving high derivatives in the action.

Recall that the euclidean action of Yang-Mills theory is given by

$$S(A) = \frac{1}{g^2}\int d^4x F^a_{\mu\nu}F^{\mu\nu}_a, \tag{3.3.1}$$

where $F^a_{\mu\nu} = \partial_\mu A^a_\nu - \partial_\nu A^a_\mu + f^{abc}A^b_\mu A^c_\nu$ is a strength of the gauge field $A^a_\mu$. Recall that the high covariant derivatives method proposed in papers [28]-[30] proceeds by two steps. The Yang-Mills action is replaced by its regularized version

$$S(A,\Lambda_*) = \frac{1}{g^2}\int d^4x F^a_{\mu\nu}[(I+\Delta_\lambda\Lambda_*^{-2})^n]^{a'\mu}_{a\mu'}F^{\mu\nu}_a, \tag{3.3.2}$$

where $\Delta_\lambda \triangleq (\Delta_\lambda)^{a\mu}_{a'\mu'} = -D^{2a}_{a'}\delta^\mu_{\mu'} + 2\lambda f^{\ a}_{a'c}F^{c\mu}_{\ \mu'}$ is the covariant differential operator given in terms of the covariant derivative $D^a_{\mu b} = \partial_\mu\delta^a_b + f^a_{bc}A^c_\mu$ and $\lambda \in \mathbb{R}$ is an arbitrary real constant. Then the partition function for the regularized action in $\alpha_0$–gauge reads

$$Z(A,\Lambda_*) = \int\prod_x \mathbf{D}[A(x)]\det(\partial^\mu D_\mu)\exp\left\{-S(A,\Lambda_*) - \frac{1}{2\alpha_0}\int d^4x\partial^\mu A^a_\mu(I-\partial^2\Lambda_*^{-2})^n\partial^\nu A^a_\nu\right\}. \tag{3.3.3}$$

In this way, provided $n \geq 2$, all 1PI diagrams with more than one loop acquire a negative degree of divergence by power counting. However, the degree of divergence of one-loop 1PI diagrams is unchanged by the addition of covariant derivatives. Therefore the theory is not completely regularized by the simple fact of adding higher covariant derivatives to the action as for the case of scalar field theories, however, that due to the regular behaviour of the gluonic propagator the contributions in the effective action to the ghost two point function and gluon-ghost vertex are finite at one loop. This implies that one loop divergences exclusively arise in diagrams with only external gluon lines, and are given by the following product of determinants

$$Z_{\text{div}} = \det(-\partial^\mu D_\mu)\det{}^{-1/2}(\Im), \tag{3.3.4}$$

where

$$\det{}^{-1/2}(\Im) = \int\prod_x \mathbf{D}[q(x)]\exp\left\{-\frac{1}{2}d^4xd^4yq^a_\mu(x)\frac{\delta^2 S(A,\Lambda_*)}{\delta A^a_\mu(x)\delta A^b_\nu(y)}q^b_\nu(y) - \right. \\ \left. -\frac{1}{2\alpha_0}\int d^4x\partial^\mu q^a_\mu(I-\partial^2\Lambda_*^{-2})^n\partial^\nu q^a_\nu \right. \tag{3.3.5}$$

Since Faddeev-Popov ghost fields only get finite renormalizations at one loop, the divergences of $Z$ can be written in a gauge invariant way.Recall that one loop divergences of Yang-Mills theory $Z_{\text{div}}$ are formally equal to those of [28]

$$Z_{\text{div}} = \det(-\partial^{\mu} D_{\mu}) \det{}^{-1/2}(\Im_0^L), \tag{3.3.6}$$

with

$$\det{}^{-1/2}(\Im_0^L) = \int \prod_x \mathbf{D}[q(x)] \delta(D^{\mu} q(x)) \exp\left\{ -\frac{1}{2} d^4x d^4y q_{\mu}^a(x) \frac{\delta^2 S(A, \Lambda_*)}{\delta A_{\mu}^a(x) \delta A_{\nu}^b(y)} q_{\nu}^b(y) \right\} \tag{3.3.7}$$

**Remark 3.3.2**.Note that the that all the determinants in (3.3.6) are explicitly gauge invariant. This fact can be understood as a consequence of the absence of divergent radiative corrections to the interaction of Faddeev-Popov ghost fields, which also implies that the BRST symmetry is only renormalized by finite radiative corrections.

**Remark 3.3.3**.Note that gauge invariance is not lost when we add mass terms in (3.3.6).

Then, it seems natural to use these determinants as the Pauli-Villars counterterms that

subtract divergences at one loop in a gauge invariant way. This is the Slavnov approach

introduced in Ref.[28] where the Slavnov introduced the following Pauli-Villars regulator

$$\det{}^{-1/2}(\Im_m) = \det(\Lambda_*^2 m^2 - D^2) \det{}^{-1/2}(\Im_m^L) \tag{3.3.8}$$

with

$$\det{}^{-1/2}(\Im_m) = \int \prod_x \mathbf{D}[q(x)] \delta(D^{\mu} q(x)) \exp\left\{ -\frac{1}{2} d^4x d^4y q_{\mu}^a(x) \frac{\delta^2 S(A, \Lambda_*)}{\delta A_{\mu}^a(x) \delta A_{\nu}^b(y)} q_{\nu}^b(y) - -\frac{1}{2} m^2 \Lambda_*^2 \int d^4x q^2(x). \right. \tag{3.3.9}$$

The regularized partition function reads [28]

$$Z(A, \Lambda_*) = \int \prod_x \mathbf{D}[A(x)] \exp\left\{ -S(A, \Lambda_*) - \frac{1}{2\alpha_0} \int d^4x \partial^{\mu} A_{\mu}^a (I - \partial^2 \Lambda_*^{-2})^n \partial^{\nu} A_{\nu}^a \right\} \times \det(-\partial^{\mu} D_{\mu}) \prod_j \det{}^{-s_j/2}(\Im_{m_j}). \tag{3.3.10}$$

is, then, free of divergences at one loop provided the $s_j$ parameters are chosen so that

$$1 + \sum_j s_j = 0. \tag{3.3.11}$$

**Remark 2.2.4**.Note that Pauli-Villars conditions do not involve the masses as it is usually the case. This is due to gauge invariance and the high derivative terms in the action that make finite the terms depending on $m$.

The problem is that Pauli-Villars determinants $\det^{-1/2}(\Im_m)$ do not converge formally to a

constant, as they should, when the cutoff is removed. In fact, we have that [31]

$$\lim_{\Lambda_*\to\infty} \det^{-1/2}(\Im_m) = \int \prod_x \mathbf{D}[q(x)]\delta(D^\mu q(x)) \exp\left\{-\frac{1}{2}\int d^4xq^2(x)\right\}. \quad (3.3.12)$$

that depends on $A$ through the delta functional $\delta(D^\mu q(x))$.These difficulties mentioned above has been resolved in paper [31].

# 3.4.Pauli–Villars renormalization of $QED_4$via Colombeau generalized functions.What is the physical significance of Pauli-Villars renormalization?

The regularization method of Pauli–Villars (PV) subtraction is of long standing in quantum field theory.In the more common dimensional regularization the properties of for instance the Dirac algebra are dimension dependent (and in particular the treatment of $\gamma_5$ is not unambiguous), and hence problems may arise in the study of chiral phenomena.

Recall that Pauli-Villars regularization requires that for each particle of mass $m$ a new ghost particle of mass $M_{PV}$ be added with either the wrong statistics or the wrong-sign kinetic term. These new particles are designed to cancel exactly loop amplitudes with physical particles at asymptotically large loop momentum. For example, one can write down a Pauli-Villars Lagrangian for $QED_4$, which works at the $1$-loop level, as

$$\begin{gathered}\mathcal{L}_{ren}^{PV} = -\frac{1}{4}F_{\mu\nu}^2 + \overline{\psi}\left(i\partial\!\!\!/ - eA\!\!\!/ - eA_\mu^{gh} - m\right)\psi + \frac{1}{4}\left(F_{\mu\nu}^{gh}\right)^2 - \frac{1}{2}M_{1,PV}^2\left(A_\mu^{gh}\right)^2 + \\ \overline{\psi}\left(i\partial\!\!\!/ - eA\!\!\!/ - eA_\mu^{gh} - M_{2,PV}^2\right)\psi\end{gathered} \quad (3.4.1)$$

with $A_\mu^{gh}$ the ghost photon and $\psi^{gh}$ the ghost electron and $F_{\mu\nu}^{gh} = \partial_\mu A_\nu^{gh} - \partial_\nu A_\mu^{gh}$. We assume that both the ghost photon and ghost electron have bosonic statistics; the ghost photon has a wrong-sign kinetic term.For example, $\mathcal{L}_{ren}^{PV}$ leads to a Feynman-gauge ghost-photon propagator of the form

$$\left\langle 0\left|T\left\{A_\mu^{gh}(x), A_\nu^{gh}(y)\right\}\right|0\right\rangle = \int \frac{d^4p}{(2\pi)^4}\exp[ip(x-y)]\frac{ig^{\mu\nu}}{p^2 - M_{PV}^2 + i\epsilon}. \quad (3.4.2)$$

**Remark 3.4.1.**(i) Since this has the opposite sign from the photon propagator, it will cancel the photon's contribution, for example, to the electron self-energy loop for loop momenta $k^\mu \gg M$. The ghost electron propagator is the same as the regular electron propagator; however, ghost electron loops do not get a factor of $-1$ (since they are bosonic) and therefore cancel regular electron loops when $k^\mu \gg M_{PV}$.
(ii) At the end of the calculation the limit $M_{PV} \to \infty$ is implied.
For example from Eq.(3.2.13) by taking the limit $M_{PV} \to \infty$ we get the canonical result independent on Pauli-Villars mas $M_{PV}$ :

$$\begin{gathered}\Pi_2(k^2) = \lim_{M_{PV}\to\infty}\Pi_2^{M_{PV}}(k^2) = \\ \frac{e^2}{4\pi m}\int_0^1 d\xi\frac{m}{[m^2 - \xi(1-\xi)k^2]^{1/2}} + \frac{e^2(\alpha-1)}{4\pi m} - \frac{e^2\alpha(sign(1-\alpha^{-1}))}{4\pi m}.\end{gathered} \quad (3.4.3)$$

**Remark 3.4.2.**The sign of the residue of the propagator is normally dictated by

unitarity -

a particle whose propagator has the sign in Eq.(2.3.2) has negative norm, and would generate probabilities greater than $1$. So, $A_\mu^{gh}$ cannot create or destroy physical on-shell particles. Thus, fields such as $A_\mu^{gh}$ are said to be associated with Pauli-Villars ghosts.

**Remark 3.4.3.** Indeed, the introduction of Pauli-Villars ghosts is much more clearly a deformation in the $UV$, relevant at energy scales $\Lambda$ of order the Pauli-Villars mass $M_{PV}$ or larger, than analytically continuing to $4-\varepsilon$ dimensions.

**Remark 3.4.4.** In order to avoid difficulties with unitarity mentioned above, we assume that:(i) physics of elementary particles is separated into low/high energy ones,

(ii) the standard notion of smooth spacetime is assumed to be altered at a high energy cutoff scale and a new treatment based on QFT in a fractal spacetime with negative dimension is used above that energy scale $\Lambda_* \gg M_{PV} \gg m$

(iii) at the end of the calculation the limit $M_{PV} \to \infty$ is not implied. For example instead Eq.(3.4.3) we set

$$\begin{gathered}\Pi_2(k^2) \triangleq \Pi_2^{M_{PV}}(k^2) = \\ \frac{e^2}{4\pi m}\int_0^1 d\xi \frac{m}{[m^2-\xi(1-\xi)k^2]^{1/2}} + \frac{e^2(\alpha-1)}{4\pi m} - \frac{e^2\alpha(sign(1-\alpha^{-1}))}{4\pi m} + \\ + \frac{e^2}{4\pi m}\frac{(\alpha-1)}{2}\left(\frac{k}{M_{PV}}\right)^2\int_0^1 d\xi\xi(1-\xi) - \\ -\frac{e^2\alpha(sign(1-\alpha^{-1}))}{8\pi m}\left(\frac{k}{M_{PV}(1-\alpha^{-1})}\right)^2\int_0^1 d\xi\xi(1-\xi) + O\left(\frac{k^4}{M_{PV}^4}\right).\end{gathered} \tag{3.4.4}$$

In N. N. Bogoliubov handbook [8] Pauli–Villars renormalization of $QED_4$ is considered.

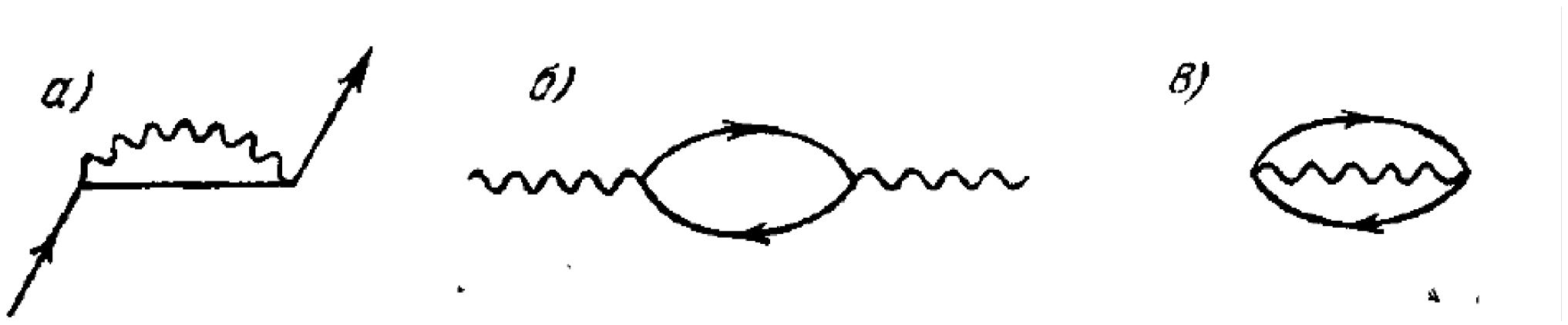

Fig.3.4.1.

The term of the scattering matrix corresponding to the diagram of Fig.3.4.1a) reads [8]:

$$\begin{gathered}-e^2 : \overline{\psi}(x)\gamma_n S^c(x-y)\gamma^n D_c^0(x-y)\psi(y) := \\ -i : \overline{\psi}(x)\Sigma(x-y)\psi(y) :\end{gathered} \tag{3.4.5}$$

where

$$\Sigma(x-y) = -ie^2\gamma^n S^c(x-y)\gamma_n D_c^0(x-y). \tag{3.4.6}$$

The formal Fourier transform $\widehat{\Sigma}(p)$ of the operator (3.4.6) reads

$$\hat{\Sigma}(p) = \frac{e^2}{i(2\pi)^4}\int dk\hat{D}_c^0(k)\gamma^n\hat{S}^c(p-k)\gamma_n =$$

$$\frac{e^2}{i(2\pi)^4}\int\frac{d^4k}{k^2+i\epsilon}\gamma^n\frac{\hat{p}-\hat{k}+m}{(p-k)^2-m^2+i\epsilon}\gamma_n, \quad (3.4.7)$$

$$\Sigma(x-y) = \frac{1}{(2\pi)^4}\int\exp[-ip(x-y)]\hat{\Sigma}(p)d^4p,$$

where

$$D_c^0(x-y) = -\frac{1}{(2\pi)^4}\int\frac{d^4k\exp[-ip(x-y)]}{k^2+i\epsilon}$$

$$S_{\alpha\beta}^c(x-y) = \frac{1}{(2\pi)^4}\int\frac{d^4k\left[(m+\hat{p})_{\alpha\beta}\right]\exp[-ip(x-y)]}{m^2-p^2-i\epsilon} \quad (3.4.8)$$

So that the integrand in (3.4.7) behaves like $|k|^{-3}$ the integral (3.4.7) diverges.

**Remark 3.4.5.**In fact, the causal Green function (3.4.8) of the $QED_4$ is the classical Schwartz distribution which is defined on a test smooth functions. It has the $\delta$-function like singularities and needs an additional definition for the product of several such functions at a single point.The discussed above diagram (see Fig.3.4.1) is precisely this product.

**Remark 3.4.6.**In his handbook [8] N. N. Bogoliubov argue that a problem of the ultraviolet divergences arises exactly from the Schwartz Impossibility Theorem [32].In fact N. N.Bogoliubov argue that these problems has only purely mathematical nature.However this Bogoliubov statement completely wrong but holds from Bogoliubov time until nowodays.

**Remark 3.4.7.**In particular N. N. Bogoliubov wroted [8]: "We thus see that the purely formal rules for dealing with products of causal functions, which we adopted earlier, lead to a meaningless result in this case.This is *essentially a manifestation* of the fact that we did not define the product of singular functions as an integrable singular function. In order to solve the problem of determining the coefficients of the chronological product $T(\mathcal{L}(x_1))\mathcal{L}((x_2))$ as integrable improper functions, we use the method of transition to the limit similar to that used in §18 (see [8],§18). In order to do this, we first consider an auxiliary fictitious case in which the field operator functions satisfy commutation relations in which the causal $\Delta^c$-functions are replaced by $reg(\Delta^c)$."

**Remark 3.4.8.**Recall that classical Schwartz distribution is defined as linear functionals on a test smooth functions [32]. Schwartz distributions may be multiplied by real numbers and added together, so they form a real vector space. Schwartz distributions may also be
multiplied by infinitely differentiable functions, but it is not possible to define a product of
general distributions that extends the usual pointwise product of functions and has the same algebraic properties. This result was shown by Schwartz (1954), and is usually referred to as the Schwartz Impossibility Theorem [32].

**Remark 3.4.9**.Note that:(i) by using linear homomorphism (a) $D_c^0(x^2) \to D_c^0(x^2;\varepsilon) = D_c^0(x^2 + i\varepsilon)$ and (b) $S_{\alpha\beta}^c(x^2) \to S_{\alpha\beta}^c(x^2;\varepsilon) = S_{\alpha\beta}^c(x^2 + i\varepsilon)$
adapted to the Lorentz invariance of the Schwartz distributions $D_c^0(x^2)$ and $S_{\alpha\beta}^c(x^2)$
we can embed these distributions into Colombeau algebra $\mathcal{G}(\mathbb{C}_x^4)$ :

$$D_c^0(x^2) \hookrightarrow (D_c^0(x^2;\varepsilon))_{\varepsilon\in(0,1]} \in \mathcal{G}(\mathbb{C}_x^4),$$
$$S_{\alpha\beta}^c(x^2) \hookrightarrow (S_{\alpha\beta}^c(x^2;\varepsilon))_{\varepsilon\in(0,1]} \in \mathcal{G}(\mathbb{C}_x^4). \quad (3.4.9)$$

**Remark 3.4.10**. Note that in contrast with Schwartz distributions $D_c^0(x^2)$ and $S_{\alpha\beta}^c(x^2)$
Colombeau distributions $(D_c^0(x^2;\varepsilon))_\varepsilon$ and $(S_{\alpha\beta}^c(x^2;\varepsilon))_\varepsilon$ on the light cone are well defined
in the sense of Colombeau generalized numbers,i.e., $(D_c^0(0;\varepsilon))_\varepsilon, (S_{\alpha\beta}^c(0;\varepsilon))_\varepsilon \in \widetilde{\mathbb{R}}$.
Thus there is no any mathematical problems for dealing with products and convolution of
Colombeau causal generalized functions $(D_c^0(x^2;\varepsilon))_\varepsilon$ and $(S_{\alpha\beta}^c(x^2;\varepsilon))_\varepsilon$.
We rewrite now the Eqs.(3.4.8) in the following equivalent form by using Colombeau integration

$$(D_c^0(x-y;\varepsilon))_\varepsilon =$$
$$-\frac{1}{(2\pi)^4}\left(\int_{-1/\varepsilon}^{1/\varepsilon} dk^{(0)} \int_{-1/\varepsilon}^{1/\varepsilon} dk^{(1)} \int_{-1/\varepsilon}^{1/\varepsilon} dk^{(2)} \int_{-1/\varepsilon}^{1/\varepsilon} dk^{(3)} \frac{\exp[-ip(x-y)]}{k^2 + i\epsilon}\right)_\varepsilon \triangleq \quad (3.4.10)$$
$$-\frac{1}{(2\pi)^4}\left(\int_{-1/\varepsilon}^{1/\varepsilon} \frac{d^4k \exp[-ip(x-y)]}{k^2 + i\epsilon}\right)_\varepsilon,$$

and

$$(S_{\alpha\beta}^c(x-y;\varepsilon))_\varepsilon = \frac{1}{(2\pi)^4} \times$$
$$\left(\int_{-1/\varepsilon}^{1/\varepsilon} dk^{(0)} \int_{-1/\varepsilon}^{1/\varepsilon} dk^{(1)} \int_{-1/\varepsilon}^{1/\varepsilon} dk^{(2)} \int_{-1/\varepsilon}^{1/\varepsilon} dk^{(3)} \frac{\left[(m+\widehat{p})_{\alpha\beta}\right] \exp[-ip(x-y)]}{m^2 - p^2 - i\epsilon}\right)_\varepsilon \quad (3.4.11)$$
$$\triangleq \frac{1}{(2\pi)^4}\left(\int_{-1/\varepsilon}^{1/\varepsilon} \frac{d^4k\left[(m+\widehat{p})_{\alpha\beta}\right] \exp[-ip(x-y)]}{m^2 - p^2 - i\epsilon}\right)_\varepsilon.$$

Therefore term of the scattering matrix corresponding to the diagram of Fig.3.4.1a) reads:

$$-e^2 : \overline{\psi}(x)\gamma_n[(S^c(x-y;\varepsilon))_\varepsilon]\gamma^n[(D_c^0(x-y;\varepsilon))_\varepsilon]\psi(y) :=$$
$$-i : \overline{\psi}(x)[(\Sigma(x-y;\varepsilon))_\varepsilon]\psi(y) : \quad (3.4.11)$$

where

$$(\Sigma(x-y;\varepsilon)_\varepsilon)_\varepsilon = -ie^2\gamma^n[(S^c(x-y;\varepsilon))_\varepsilon]\gamma_n[(D_c^0(x-y;\varepsilon))_\varepsilon]. \tag{3.4.12}$$

From Eqs.(3.4.10)-(3.4.11) and Eq.(3.4.12) we find that

$$(\Sigma(x-y;\varepsilon))_\varepsilon = \frac{1}{(2\pi)^4}\left(\int_{-1/\varepsilon}^{1/\varepsilon} \exp[-ip(x-y)]\widehat{\Sigma}(p;\varepsilon)d^4p\right)_\varepsilon, \tag{3.4.13}$$

where

$$\left(\widehat{\Sigma}(p;\varepsilon)\right)_\varepsilon = \frac{e^2}{i(2\pi)^4}\left(\int_{-1/\varepsilon}^{1/\varepsilon} dk\hat{D}_c^0(k;\varepsilon)\gamma^n\left[\widehat{S}^c(p-k;\varepsilon)\right]\gamma_n\right)_\varepsilon = \frac{e^2}{i(2\pi)^4}\left(\int_{-1/\varepsilon}^{1/\varepsilon}\frac{d^4k}{k^2+i\epsilon}\gamma^n\frac{\widehat{p}-\widehat{k}+m}{(p-k)^2-m^2+i\epsilon}\gamma_n\right)_\varepsilon, \tag{3.4.14}$$

Note that in contrast with ill defined formal expressions (3.4.7) the expressions (3.4.14) gives a well defined Colombeau generalized functions: $(\Sigma(x-y;\varepsilon))_\varepsilon \in \mathcal{G}(\mathbb{R}^4)$ and $\left(\widehat{\Sigma}(p;\varepsilon)\right)_\varepsilon \in \mathcal{G}(\mathbb{R}_p^4)$. Note that for any $\bar{p} \in \mathbb{R}_p^4$ obviously $\left(\widehat{\Sigma}(\bar{p};\varepsilon)\right)_\varepsilon \in \widetilde{\mathbb{R}}$ and any $\left(\widehat{\Sigma}(\bar{p};\varepsilon)\right)_\varepsilon$ is infinite large Colombeau generalized number [26].In order to eliminate the infinite large Colombeau generalized quantities $\sim (\ln\varepsilon^{-1})_\varepsilon, (\varepsilon^{-1})_\varepsilon$, respectively from Colombeau integral (3.4.14),i.e., make it finite, we apply Pauli–Villars renormalization via Colombeau generalized functions,see sect.3.2. Finally we get

$$\mathbf{ren}\left\{\left(\widehat{\Sigma}(p;\varepsilon)\right)_\varepsilon\right\} = \frac{e^2}{8\pi^2}\int_0^1 d\xi(2m-\widehat{p}\xi)\ln\left[\frac{m^2}{m^2-\xi p^2}\,\frac{\xi M_{PV}^2-(1-\xi)m^2-\xi(1-\xi)p^2}{\xi[M_{PV}^2-\xi^2(1-\xi)p^2]}\right] \simeq \frac{e^2}{8\pi^2}\int_0^1 d\xi(2m-\widehat{p}\xi)\ln\left[\frac{m^2}{m^2-\xi p^2}\right]+O\left(\frac{p^2}{M_{PV}^2}\right). \tag{3.4.15}$$

where $M_{PV} \in \mathbb{R}_+$ is arbitraly large but *finite* Pauli–Villars mass.

**Remark 3.4.11.**Note that in contrast with canonical formal approuch [8],[33-36] we cannot taking the limit $M_{PV} \to \infty$ in (3.4.15) ,see sect.3.2.

**Remark 3.4.12**. It is clear from consideration above that problem with ultraviolet divergences arises not from the Schwartz Impossibility Theorem [32], as Bogoliubov mistakenly has claimem many years ago [8],but exactly from physically wrong canonical Lagrangian of $QED_4$ in which physical ghost fields was missing.This is essentially a manifestation of the fact that Pauli–Villars renormalization via Colombeau generalized functions (see sect.3.2) signals about real physical nature of the ghost fields.

# 4.OFT in a ghost sector via dimensional renormalization.

## 4.1.Dimensional Regularization via Colombeau generalized functions.

The most popular in gauge theories is the so-called dimensional regularization. In this case, one modifies the integration measure: $d^4q \to (\mu^2)^{\varepsilon} d^{4-2\varepsilon}q$ where $\mu$ is a parameter of dimensional regularization with dimension of a mass.In this case, all the ultraviolet and infrared singularities manifest themselves as pole terms in $\varepsilon$.Consider the earlier discussed example see Fig.3.1.3 and using the Euclidean representation rewrite it formally in $D$-dimensional space

$$I(p^2,m^2;D) = \int_0^1 dx \int \frac{d^Dk}{(k^2+\Delta)^2} =$$
$$\frac{\Omega_D}{2}\int_0^1 dx \int_0^\infty \frac{(k^2)^{\frac{D}{2}-1}dk}{(k^2+\Delta)^2} = (\Delta^2)^{\frac{D}{2}-2}\frac{\Gamma\left(\frac{D}{2}\right)\Gamma\left(\frac{D}{2}-2\right)}{\Gamma(2)}, \tag{4.1.1}$$
$$\Delta = \Delta(p^2,m^2) = p^2x(1-x)+m^2,$$

where we assume that the dimension $D$ is such that the integral exists. In this case this is 2 and 3. The main formula (4.1.1) allows one to perform the analytical continuation over $D = 4$ into the region $D = 4-2\varepsilon, \varepsilon \in (0,1]$. For $\varepsilon = 0$, i.e., in 4 dimensions, the integral does not exist since the $\Gamma$-function has a pole at zero argument. However, in the vicinity of zero we get a regularized expression.From Eq.(4.1.1) we get

$$I(p^2,m^2;D) = \frac{i}{(2\pi)^D}\frac{\Omega_D}{2}\int_0^1 dx \frac{\Gamma\left(\frac{D}{2}\right)\Gamma\left(\frac{D}{2}-2\right)}{[p^2x(1-x)+m^2]^{\frac{D}{2}-2}}. \tag{4.1.2}$$

Substituting now $D = 4-2\varepsilon$ in RHS of the Eq.(4.1.2) and transforming back into the pseudo-Euclidean space we get

$$I_\varepsilon(p^2,m^2) \triangleq I_\varepsilon(p^2,m^2;4-2\varepsilon) = \frac{i(-\pi)^{2-\varepsilon}}{(2\pi)^{4-2\varepsilon}}\Gamma(\varepsilon)\int_0^1 dx \frac{(\mu^2)^\varepsilon}{[p^2x(1-x)+m^2]^\varepsilon}. \tag{4.1.3}$$

The formula (4.1.3) allows one to define the integral $I(p^2,m^2;D=4)$ as Colombeau generalized function (see [21]-[22]) $I(p^2,m^2;D=4) \triangleq (I_\varepsilon(p^2,m^2))_{\varepsilon\in(0,1]} \in \mathcal{G}(\mathbb{R}_x^D)$ :

$$I(p^2,m^2;D=4) \triangleq (I_\varepsilon(p^2,m^2))_\varepsilon \triangleq I_\varepsilon(p^2,m^2;4-2\varepsilon) =$$
$$\frac{i\left((-\pi)^{2-\varepsilon}\right)_\varepsilon}{\left((2\pi)^{4-2\varepsilon}\right)_\varepsilon}(\Gamma(\varepsilon))_\varepsilon\left(\int_0^1 dx \frac{(\mu^2)^\varepsilon}{[p^2x(1-x)+m^2]^\varepsilon}\right)_\varepsilon. \tag{4.1.4}$$

Expanding the denominator of the integrand into the series over $\varepsilon$, finally we get

$$I(p^2,m^2;D=4) \triangleq (I_\epsilon(p^2,m^2))_\epsilon$$
$$\frac{i}{16\pi^2}(\Gamma(1+\epsilon))_\epsilon\left((\epsilon^{-1})_\epsilon - \int_0^1 dx \ln\left[\frac{p^2x(1-x)+m^2}{-\mu^2}\right] + \ln(4\pi)\right) \approx_{\widetilde{\mathbb{R}}} \tag{4.1.5}$$
$$\approx_{\widetilde{\mathbb{R}}} \frac{i}{16\pi^2}(\epsilon^{-1})_\epsilon + \frac{i}{16\pi^2}\left(\ln(4\pi) - \int_0^1 dx \ln\left[\frac{p^2x(1-x)+m^2}{-\mu^2}\right]\right).$$

We see that the classical ultraviolet divergence now takes the rigorous mathematical form of the infinite large Colombeau generalized number $(\varepsilon^{-1})_\varepsilon \in \widetilde{\mathbb{R}}$,see [23].

We present below the main Colombeau integrals needed for the one-loop calculations. They can be obtained via the analytical continuation from the integer values of $D$. We will write them down directly in the pseudo-Euclidean space. First note that

$$\int \frac{d^D p}{[p^2 - 2kp + m^2]^\alpha} = i\frac{\Gamma(\alpha - D/2)}{\Gamma(\alpha)} \frac{(-\pi)^{D/2}}{[m^2 - k^2]^{\alpha - D/2}} \tag{4.1.6}$$

and

$$\widetilde{\int}\frac{d^4 p}{[p^2 - 2kp + m^2]^2} \triangleq \left(\int \frac{d^{4-2\varepsilon} p}{[p^2 - 2kp + m^2]^2}\right)_\varepsilon \triangleq i\frac{(\Gamma(\varepsilon))_\varepsilon}{\Gamma(2)} \frac{\left((-\pi)^{2-\varepsilon}\right)_\varepsilon}{([m^2 - k^2]^\varepsilon)_\varepsilon},$$

$$\widetilde{\int}\frac{p_\mu d^4 p}{[p^2 - 2kp + m^2]^2} \triangleq \left(\int \frac{p_\mu d^{4-2\varepsilon} p}{[p^2 - 2kp + m^2]^2}\right)_\varepsilon \triangleq i\frac{(\Gamma(\varepsilon))_\varepsilon}{\Gamma(2)} \frac{\left((-\pi)^{2-\varepsilon}\right)_\varepsilon k_\mu}{([m^2 - k^2]^\varepsilon)_\varepsilon},$$

$$\widetilde{\int}\frac{p_\mu p_\nu d^4 p}{[p^2 - 2kp + m^2]^2} \triangleq \left(\int \frac{p_\mu p_\nu d^{4-2\varepsilon} p}{[p^2 - 2kp + m^2]^2}\right)_\varepsilon \triangleq \tag{4.1.7}$$

$$i\left((-\pi)^{2-\varepsilon}\right)_\varepsilon \left[\frac{(\Gamma(\varepsilon))_\varepsilon}{\Gamma(2)} \frac{k_\mu k_\nu}{([m^2 - k^2]^\varepsilon)_\varepsilon} + \frac{g^{\mu\nu}}{2} \frac{(\Gamma(\varepsilon - 1))_\varepsilon}{\Gamma(2)} \frac{1}{\left([m^2 - k^2]^{\varepsilon - 1}\right)_\varepsilon}\right],$$

$$(\Gamma(\varepsilon))_\varepsilon = (\varepsilon^{-1})_\varepsilon, \varepsilon \in (0,1].$$

The key formula is (4.1.6). These integrals remain Colombeau generalized functions from
$\mathcal{G}(\mathbb{R}^4)$ and $\mathcal{G}(\mathbb{R}^4 \times \mathbb{R}^4)$.

## 4.2.The scalar theory $\lambda\varphi^4_{D=4}$in a ghost sector via Colombeau generalized functions.The one-loop approximation.

Let us consider the theory described by the Lagrangian

$$\mathcal{L} = \frac{1}{2}(\partial_\mu \varphi)^2 - \frac{m^2}{2}\varphi^2 - \frac{\lambda}{4!}\varphi^4. \tag{4.2.1}$$

**The propagator**: In the first order there is only one diagram of the tad-pole type shown
in Fig.4.2.1.

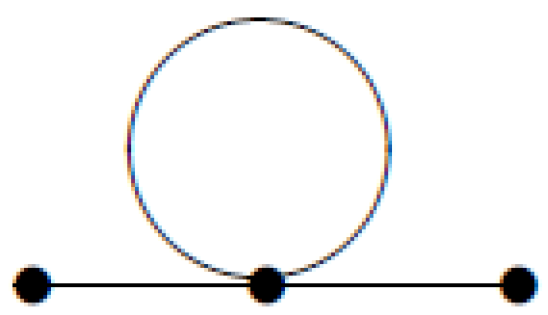

Fig.4.2.1.The one-loop propagator diagram.

The corresponding Colombeau integral is

$$(J_1(p^2, \varepsilon))_\varepsilon = \frac{-i\lambda}{(2\pi)^{4-2\varepsilon}} \frac{i}{2}\left(\int \frac{d^{4-2\varepsilon} k (\mu^2)^\varepsilon}{k^2 - m^2}\right)_\varepsilon, \tag{4.2.2}$$

where $1/2$ is the combinatoric factor and $\varepsilon \in (0,1]$. Calculating the Colombeau integral (4.2.2), according to (4.1.7), we get

$$(J_1(p^2,\varepsilon))_\varepsilon = \frac{-i\lambda}{((4\pi)^{2-\varepsilon})_\varepsilon}\frac{(\Gamma(-1+\varepsilon))_\varepsilon}{2\Gamma(1)}m^2\left(\left(\frac{\mu^2}{m^2}\right)^\varepsilon\right)_\varepsilon = \frac{i\lambda}{32\pi^2}m^2\left[\left(\frac{1}{\varepsilon}\right)_\varepsilon + 1 - \gamma_E + \log(4\pi) - \log\frac{m^2}{\mu^2}\right]. \tag{4.2.2}$$

The fact that the Colombeau integral (4.2.2) well defined but contains infinite Colombeau generalized number $(\varepsilon^{-1})_\varepsilon \in \widetilde{\mathbb{R}}$.

**The vertex**: Here one also has only one diagram but the external momenta can be adjusted in several ways (see Fig.4.2.2). As a result the total contribution to the vertex function consists of three parts $I(s,t,u) = I_1(s) + I_1(t) + I_1(u)$, where we introduced the commonly accepted notation for the Mandelstam variables (we assume here that the momenta $p_1$ and $p_2$ are incoming and the momenta $p_3$ and $p_4$ are outgoing) $s = (p_1+p_2)^2 = (p_3+p_4)^2$, $t = (p_1-p_3)^2 = (p_2-p_4)^2$, $u = (p_1-p_4)^2 = (p_2-p_3)^2$, and the integral equals

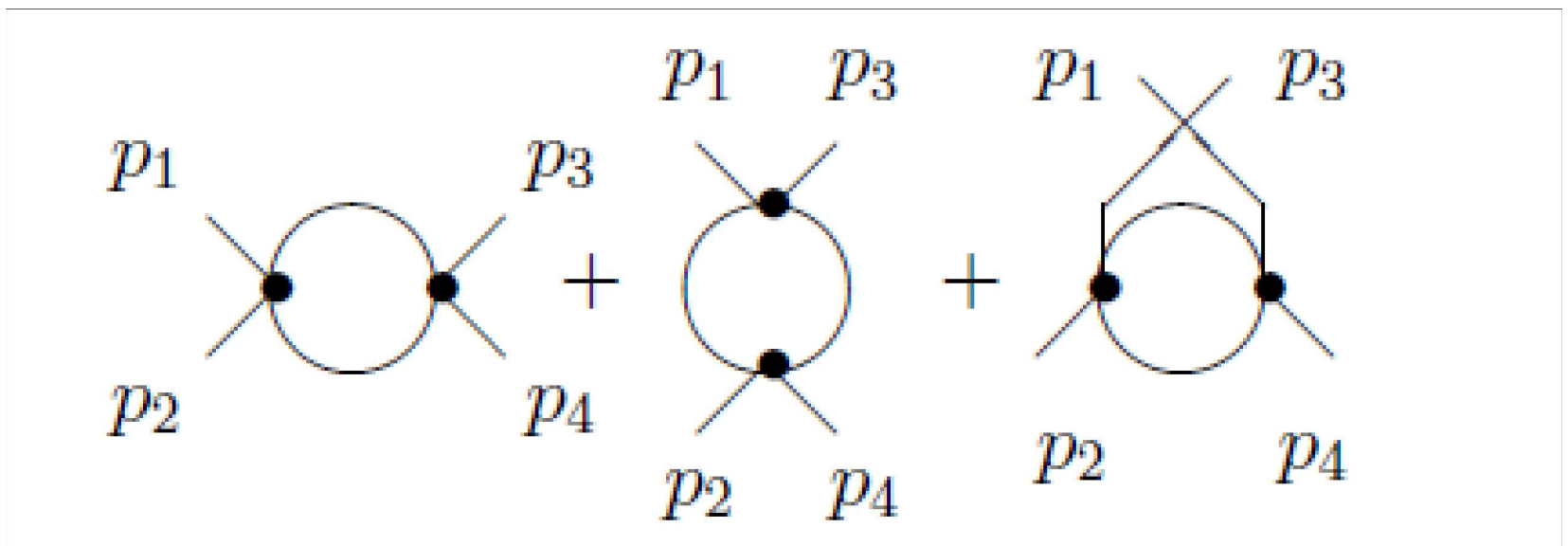

Fig.4.2.2.The one-loop vertex diagram

and the Colombeau integral reads

$$I_1(s) = \frac{(-i\lambda)^2}{48}\frac{((\mu^2)^\varepsilon)_\varepsilon}{((2\pi)^{4-2\varepsilon})_\varepsilon}i^2\left(\int\frac{d^{4-2\varepsilon}k}{[k^2-m^2][(p-k)^2-m^2]}\right)_\varepsilon, \tag{4.2.3}$$

where $1/48$ is the combinatoric coefficient.

Recall that the key formula is the Fourier-transformation (in the sense of generalized functions) of the propagator of a massless particle for $D = 4$ reads [32]

$$\int\frac{d^4p\, e^{ipx}}{p^2} = \frac{i\pi^2}{x^2}, \tag{4.2.4}$$

which holds in arbitrary noncritical dimension $D$ and any power of the propagator as follows [32]:

$$\int\frac{d^Dp\, e^{ipx}}{(p^2)^\alpha} = i(-\pi)^{D/2}\frac{\Gamma(D/2-\alpha)}{\Gamma(\alpha)}\frac{1}{(x^2)^{D/2-\alpha}}. \tag{4.2.5}$$

In the case of the integral (4.2.3) for $m = 0$ one first has to mentally transform both the propagators into coordinate space which, according to (4.2.6), gives the factor $\left(\frac{\Gamma(1-\varepsilon)}{\Gamma(1)}\right)^2$, then multiply the obtained propagators (this gives $1/(x^2)^{2-2\varepsilon}$)) and transform the obtained result back into momentum space that gives the factor $\frac{\Gamma(\varepsilon)}{\Gamma(2-2\varepsilon)}$ and the power of momenta $1/(p^2)^\varepsilon$ (the same as in the argument of the last $\Gamma$-function). Besides this, each loop contains the factor $i(-\pi)^{2-\varepsilon}$. Collecting all together one obtains

$$(I_1(s,\varepsilon))_\varepsilon = \frac{(-i\lambda)^2}{48}\frac{((\mu^2)^\varepsilon)_\varepsilon i^2}{((2\pi)^{4-2\varepsilon})_\varepsilon}\left(\int\frac{d^{4-2\varepsilon}k}{k^2(p-k)^2}\right)_\varepsilon =$$
$$\frac{\lambda^2}{48}\frac{i(\pi^{2-\varepsilon})_\varepsilon}{((2\pi)^{4-2\varepsilon})_\varepsilon}\left(\left(\frac{\mu^2}{-s}\right)^\varepsilon\right)_\varepsilon\frac{(\Gamma(1-\varepsilon)\Gamma(1-\varepsilon)\Gamma(\varepsilon))_\varepsilon}{\Gamma(1)\Gamma(1)(\Gamma(2-2\varepsilon))_\varepsilon}$$
$$= \frac{i}{48}\frac{\lambda^2}{((4\pi)^{2-\varepsilon})_\varepsilon}\left(\left[\frac{\mu^2}{-s}\right]^\varepsilon\right)_\varepsilon\frac{1}{(\varepsilon(1-2\varepsilon))_\varepsilon}\left(\frac{\Gamma^2(1-\varepsilon)\Gamma(1+\varepsilon)}{\Gamma(1-2\varepsilon)}\right)_\varepsilon = \tag{4.2.6}$$
$$\frac{i}{48}\frac{\lambda^2}{16\pi^2}\left[\left(\frac{1}{\varepsilon}\right)_\varepsilon + 2 - \gamma_E + \log 4\pi + \ln\frac{\mu^2}{-s}\right].$$

The four-point vertex in the one-loop approximation reads

$$(\Gamma_{4,\varepsilon})_\varepsilon = -i\lambda\left\{1 - \frac{\lambda}{16\pi^2}\left(\left(\frac{3}{2\varepsilon}\right)_\varepsilon + 3 - \frac{3}{2}\gamma_E + \right.\right.$$
$$\left.\left.\frac{3}{2}\log 4\pi + \frac{1}{2}\ln\frac{\mu^2}{-s} + \frac{1}{2}\ln\frac{\mu^2}{-t} + \frac{1}{2}\ln\frac{\mu^2}{-u}\right)\right\}. \tag{4.2.7}$$

We are interested now in the Colombeau singular parts, i.e.,infinite part in rigorous Colombeau sense.They are given by Eqs.4.2.7 and 4.2.8

$$\mathbf{Sing}\{(J_1(p^2,\varepsilon))_\varepsilon\} = -im^2\left(\frac{\lambda}{16\pi^2}\right)\left(-\frac{1}{2\varepsilon}\right)_\varepsilon,$$
$$\mathbf{Sing}\{((\Gamma_4(s,t,u;\varepsilon)))_\varepsilon\} = -i\lambda\left(\frac{\lambda}{16\pi^2}\right)\left(-\frac{3}{2\varepsilon}\right)_\varepsilon. \tag{4.2.8}$$

Note that the singular parts do not depend on momenta, i.e. their Fourier-transform has the form of the $\delta$-function in coordinate space.In order to remove the obtained Colombeau singularities we add to the Lagrangian (4.2.1) extra terms, the counter-terms equal to the Colombeau singular parts with the opposite sign (the factor $i$ belongs to the S-matrix and does not enter into the Lagrangian),namely,

$$(\Delta\mathcal{L}_\varepsilon)_\varepsilon = \left(\frac{1}{2\varepsilon}\right)_\varepsilon\frac{\lambda}{16\pi^2}\left(-\frac{m^2}{2}\phi^2\right) + \frac{\lambda}{16\pi^2}\left(\frac{3}{2\varepsilon}\right)_\varepsilon\left(-\frac{\lambda}{4!}\phi^4\right). \tag{4.2.9}$$

These counter-terms correspond to additional vertices shown in Fig.4.2.3.

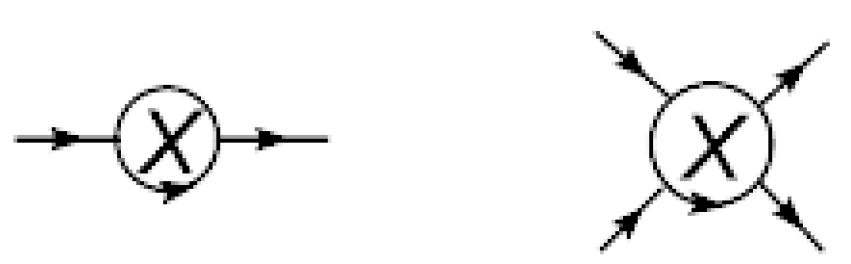

Fig.4.2.3.The one-loop counter-terms
in the scalar theory $\lambda\varphi_4^4$.

With account taken of the these new diagrams the expressions for the propagator (4.2.3) and the vertex (4.2.7) become

$$(J_1(p^2,\varepsilon))_\varepsilon = \frac{i\lambda}{32\pi^2}m^2(1-\gamma_E+\log(4\pi)-\log(m^2/\mu^2)) \tag{4.2.10}$$

and

$$(\Gamma_{4,\varepsilon})_\varepsilon = i\lambda\left\{\frac{\lambda}{16\pi^2}\left(3-\frac{3}{2}\gamma_E+\frac{3}{2}\log(4\pi)+\frac{1}{2}\ln\frac{\mu^2}{-s}+\frac{1}{2}\ln\frac{\mu^2}{-t}+\frac{1}{2}\ln\frac{\mu^2}{-u}\right)\right\} \tag{4.2.11}$$

correspondingly.Notice that the obtained expressions have no Colombeau infinities but contain the dependence on the regularization parameter $\mu^2$ which was absent in the initial theory. The appearance of this dependence on a dimensional parameter is inherent in any regularization and is called the dimensional transmutation,i.e., an appearance of a new scale in a theory. We write the counter-term in the following way

$$\left(\Delta\mathcal{L}_\varepsilon^{(1)}\right)_\varepsilon = -((Z_\varepsilon)_\varepsilon - 1)\frac{m^2}{2}\phi^2 - ((Z_{4,\varepsilon})_\varepsilon - 1)\frac{\lambda}{4!}\phi^4, \qquad (4.2.12)$$

where for different subtraction schemes one has

$$\begin{gathered}(Z_\varepsilon^{MS})_\varepsilon = 1 + \left(\frac{1}{2\varepsilon}\right)_\varepsilon \frac{\lambda}{16\pi^2}, \left(Z_\varepsilon^{\overline{MS}}\right)_\varepsilon = 1 + \left[\left(\frac{1}{2\varepsilon}\right)_\varepsilon + 1 - \gamma_E + \log(4\pi)\right]\frac{\lambda}{16\pi^2},\\ (Z_{4,\varepsilon}^{MS})_\varepsilon = 1 + \left(\frac{3}{2\varepsilon}\right)_\varepsilon \frac{\lambda}{16\pi^2}, \left(Z_{4,\varepsilon}^{\overline{MS}}\right)_\varepsilon = 1 + \left[\left(\frac{3}{2\varepsilon}\right)_\varepsilon - 3\gamma_E + 3\log(4\pi)\right]\frac{\lambda}{16\pi^2}, \\ (Z_{4,\varepsilon}^{MOM})_\varepsilon = 1 + \left[\left(\frac{3}{2\varepsilon}\right)_\varepsilon + 3 - 3\gamma_E + 3\log(4\pi) + \frac{3}{2}\ln\frac{\mu^2}{t^2}\right]\frac{\lambda}{16\pi^2}.\end{gathered} \qquad (4.2.13)$$

The Lagrangian (4.2.1) together with the Colombeau counter-terms (4.2.12) can be written as

$$\mathcal{L} + \left(\Delta_\varepsilon^{(1)}\right)_\varepsilon = (Z_{2,\varepsilon})_\varepsilon\frac{1}{2}(\partial_\mu\varphi)^2 - (Z_\varepsilon)_\varepsilon\frac{m^2}{2}\varphi^2 - (Z_{4,\varepsilon})_\varepsilon\frac{\lambda}{4!}\varphi^4 = (\mathcal{L}_\varepsilon^{Bare})_\varepsilon, \qquad (4.2.14)$$

where the renormalization Colombeau constants $(Z_\varepsilon)_\varepsilon \in \widetilde{\mathbb{R}}$ and $(Z_{4,\varepsilon})_\varepsilon \in \widetilde{\mathbb{R}}$ are given by Eqs.(4.2.13) and the renormalization Colombeau constant $(Z_{2,\varepsilon})_\varepsilon$ in the one-loop approximation equals $1$,i.e.,

$$\begin{gathered}(\mathcal{L}_\varepsilon^{\mathbf{ren}})_\varepsilon = \mathcal{L} + \left(\Delta\mathcal{L}_\varepsilon^{(1)}\right)_\varepsilon \\ = \frac{1}{2}(Z_{2,\varepsilon})_\varepsilon(\partial_\mu\varphi)^2 - (Z_\varepsilon)_\varepsilon\frac{m^2}{2}\varphi^2 - (Z_{4,\varepsilon})_\varepsilon\frac{\lambda}{4!}\varphi^4 = (\mathcal{L}_\varepsilon^{Bare})_\varepsilon,\end{gathered} \qquad (4.2.15)$$

Writing the "bare" Lagrangian in the same form as the initial one but in terms of the "bare" fields and couplings

$$(\mathcal{L}_\varepsilon^{Bare})_\varepsilon = \frac{1}{2}(\partial_\mu\varphi_B)^2 - \frac{(m_{B,\varepsilon}^2)_\varepsilon}{2}\varphi_B^2 - \frac{(\lambda_{B,\varepsilon})_\varepsilon}{4!}\varphi_B^4, \qquad (4.2.16)$$

where $(m_{B,\varepsilon})_\varepsilon, (\lambda_{B,\varepsilon})_\varepsilon \in \widetilde{\mathbb{R}}$ are infinite Colombeau constants, and comparing it with (4.2.15), we get the connection between the "bare" and renormalized Colombeau quantities

$$\varphi_B = \varphi,\ \ (m_{B,\varepsilon}^2)_\varepsilon = m^2(Z_\varepsilon)_\varepsilon,\ \ (\lambda_{B,\varepsilon})_\varepsilon = \lambda(Z_{4,\varepsilon})_\varepsilon, \qquad (4.2.17)$$

where $(Z_\varepsilon)_\varepsilon = 1 + O((\varepsilon^{-1})_\varepsilon), (Z_{4,\varepsilon})_\varepsilon = 1 + O((\varepsilon^{-1})_\varepsilon)$ positive infinite Colombeau constants. Equations (4.2.16) and (4.2.17) imply that the one-loop radiative corrections calculated from the Lagrangian (4.2.16) with parameters chosen according to (4.2.17) and (4.2.13) are finite.

**Remark 4.2.1**. Note that: (i) in one-loop approximation the all renormalization Colombeau
constants are strictly positive:

$$(Z_\varepsilon)_\varepsilon > 0, (Z_{2,\varepsilon})_\varepsilon > 0, (Z_{4,\varepsilon})_\varepsilon > 0, \qquad (4.2.18)$$

(ii) it follows from (4.2.18) in one-loop approximation a ghost sector is absent completely.

## 4.3. The scalar theory $\lambda\varphi^4_{D=4}$ in a ghost sector via Colombeau generalized functions. The two-loop approximation

Consider now the two-loop diagrams. The propagator: In this order of PT there is only one diagram shown in Fig.4.3.1.

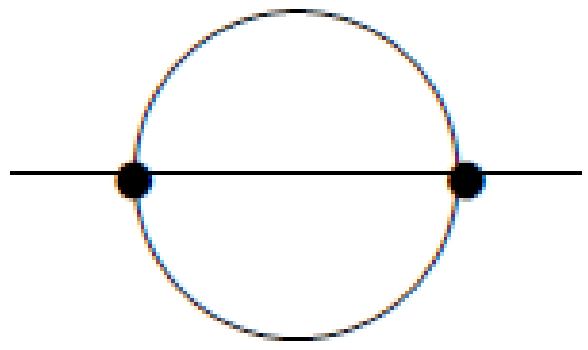

Fig.4.3.1. The two-loop propagator type diagram.

The corresponding Colombeau integral reads

$$(J_2(p^2;\varepsilon))_\varepsilon = \frac{(-i\lambda)^2}{3!}\frac{i^3((\mu^2)^{2\varepsilon})_\varepsilon}{((2\pi)^{8-4\varepsilon})_\varepsilon}\left(\int\int\frac{d^{4-2\varepsilon}k\,d^{4-2\varepsilon}q}{q^2(k-q)^2(p-k)^2}\right)_\varepsilon, \quad (4.3.1)$$

One has to transform each of the propagators into coordinate space, multiply them and transform back to momentum space. This reduces to writing down the corresponding transformation factors. Thus

$$\begin{gathered}(J_2(p^2;\varepsilon))_\varepsilon = \frac{i\lambda^2}{6}\frac{((i\pi^2)^{2-\varepsilon})_\varepsilon}{((2\pi)^{8-4\varepsilon})_\varepsilon}p^2\left(\left(\frac{\mu^2}{-p^2}\right)^{2\varepsilon}\right)_\varepsilon \times \\ \frac{(\Gamma(1-\varepsilon)\Gamma(1-\varepsilon)\Gamma(1-\varepsilon)\Gamma(-1+2\varepsilon))_\varepsilon}{\Gamma(1)\Gamma(1)\Gamma(1)(\Gamma(3-3\varepsilon))_\varepsilon} = \\ \frac{i}{6}\frac{\lambda^2}{(16\pi^2)^2}\left(\left[\frac{\mu^2}{-p^2}\right]^{2\varepsilon}\right)_\varepsilon\frac{p^2}{((2-3\varepsilon)(1-3\varepsilon)(1-2\varepsilon)2\varepsilon)_\varepsilon} = \\ \frac{i}{24}\frac{\lambda^2}{(16\pi^2)^2}p^2\left[\left(\frac{1}{\varepsilon}\right)_\varepsilon + \frac{13}{2} + 2\ln\frac{\mu^2}{-p^2}\right],\end{gathered} \quad (4.3.2)$$

where the Euler constant and $\ln 4\pi$ are omitted. The appeared ultraviolet divergence, the pole in $\varepsilon$, can be removed via the introduction of the (quasi) local Colombeau counter-term

$$\left(\Delta\mathcal{L}^{(2)}_\varepsilon\right)_\varepsilon = \frac{1}{2}((Z_{2,\varepsilon})_\varepsilon - 1)(\partial\phi)^2, \quad (4.3.3)$$

where the wave function renormalization constant $(Z_{2,\varepsilon})_\varepsilon$ in the $\overline{MS}$ scheme is obtained by taking the infinite large part of the Colombeau integral with the opposite sign

$$(Z_{2,\varepsilon})_\varepsilon = 1 - \left(\frac{1}{24\varepsilon}\right)_\varepsilon\left(\frac{\lambda}{16\pi^2}\right)^2. \quad (4.3.4)$$

After that the propagator in the massless case reads

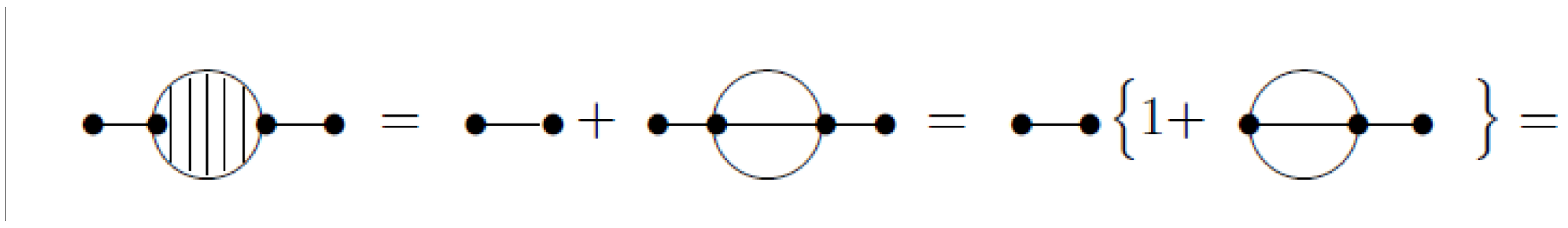

$$= \frac{i}{p^2}\left\{1 - \frac{1}{24}\frac{\lambda^2}{(16\pi^2)^2}\left(\frac{13}{2} + 2\ln\frac{\mu^2}{-p^2}\right)\right\}. \quad (4.3.5)$$

The vertex: In the given order there are two diagrams (remind that in the massless case the tad-poles equal to zero) shown in Fig.4.3.2.

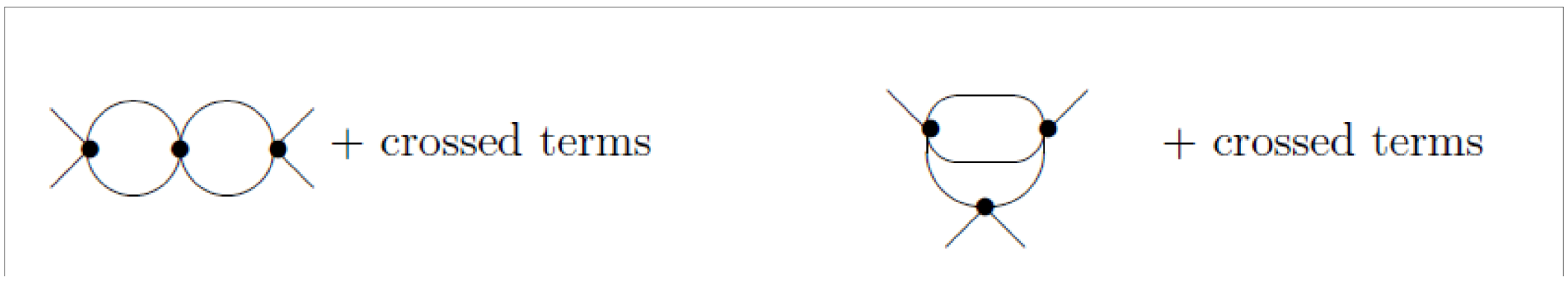

Fig.4.3.2.The two-loop vertex diagrams

The first diagram by analogy with the one-loop case equals the sum of $s,t$ and $u$ channels $(I_{21}(s,t,u;\varepsilon))_\varepsilon = (I_{21}(s;\varepsilon))_\varepsilon + (I_{21}(t;\varepsilon))_\varepsilon + (I_{21}(u;\varepsilon))_\varepsilon$, where each integral is nothing else but the square of the one-loop integral

$$(I_{21}(s;\varepsilon))_\varepsilon = \frac{(-i\lambda)^3}{96}\left(\frac{((\mu^2)^\varepsilon)_\varepsilon}{((2\pi)^{4-2\varepsilon})_\varepsilon} i^2 \left(\int \frac{d^{4-2\varepsilon}k}{k^2(p-k)^2}\right)_\varepsilon\right)^2 = \\ -\frac{i}{96}\frac{\lambda^3}{(16\pi^2)^2}\left(\left(\frac{1}{\varepsilon}\right)_\varepsilon + 2 + \ln\frac{\mu^2}{-s}\right)^2. \tag{4.3.6}$$

In the same order of $\lambda^3$ one gets additional diagrams presented in Fig.4.3.3.

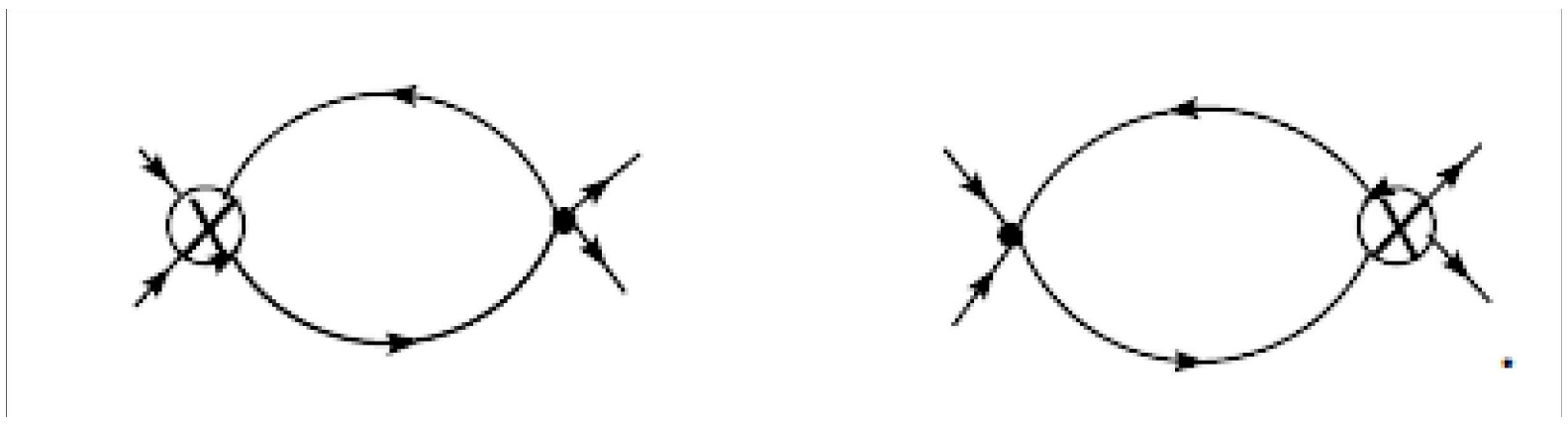

Fig.4.3.3.The diagrams with the counter-terms in
the two-loop approximation

These diagrams lead to the subtraction of divergences in the subgraphs (left and right) in the first diagram of Fig.4.3.2. The subtraction of divergent subgraphs (the $\mathcal{R}$-operation without the last subtraction called the $\mathcal{R}'$-operation) looks like where the subgraph surrounded with the dashed line means its singular part, and the rest of the graph is obtained by shrinking down the singular subgraph to a point. The result has the following form

$$\mathcal{R}'[(I_{21}(s;\varepsilon))_\varepsilon] =$$
$$-\frac{i}{4}\frac{\lambda^3}{(16\pi^2)^2}\left\{\left(\left(\frac{1}{\varepsilon}\right)_\varepsilon + 2 + \ln\frac{\mu^2}{-s}\right)^2 - \left(\frac{2}{\varepsilon}\right)_\varepsilon\left(\left(\frac{1}{\varepsilon}\right)_\varepsilon + 2 + \ln\frac{\mu^2}{-s}\right)\right\} \tag{4.3.7}$$
$$= -\frac{i}{4}\frac{\lambda^3}{(16\pi^2)^2}\left(-\left(\frac{1}{\varepsilon^2}\right)_\varepsilon + 4 + \ln^2\frac{\mu^2}{-s} + 4\ln\frac{\mu^2}{-s}\right).$$

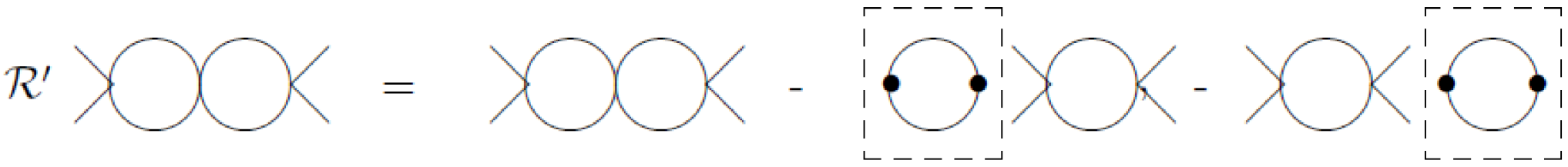

Fig.4.3.4.

Notice that after the subtractions of subgraphs the Colombeau singular part is local, i.e. in momentum space does not contain $\ln p^2$. The terms with the single pole $(1/\varepsilon)_\varepsilon$ are absent since the diagram can be factorized into two diagrams of the lower order. The contribution of a given diagram to the vertex function equals

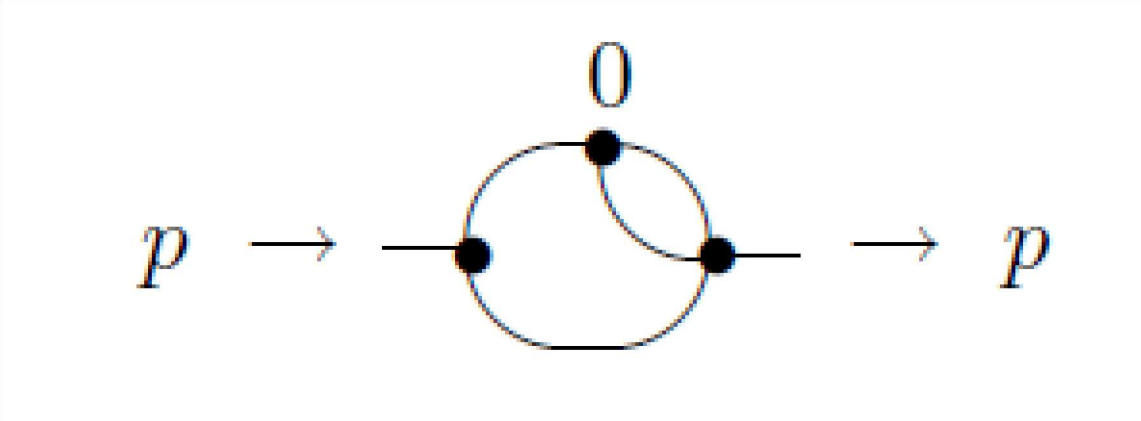

Fig.4.3.5.

$$(\Delta\Gamma_{4,\varepsilon})_\varepsilon = -i\lambda\left\{\frac{1}{4}\frac{\lambda^2}{(16\pi^2)^2}\left(-\left(\frac{3}{\varepsilon^2}\right)_\varepsilon + 12\right.\right.$$
$$\left.\left. + \ln^2\frac{\mu^2}{-s} + 4\ln\frac{\mu^2}{-s} + \ln^2\frac{\mu^2}{-t} + 4\ln\frac{\mu^2}{-t} + \ln^2\frac{\mu^2}{-u} + 4\ln\frac{\mu^2}{-u}\right)\right\} \tag{4.3.8}$$

The contribution to the renormalization constant of the four-point vertex in the $\overline{MS}$ scheme is equal to the singular part with the opposite sign

$$(\Delta Z_{4,\varepsilon})_\varepsilon = \left(\frac{3}{4\varepsilon^2}\right)_\varepsilon\left(\frac{\lambda}{16\pi^2}\right)^2. \tag{4.3.9}$$

The second diagram with the crossed terms contains 6 different cases. Consider one of them. Since we are interested here in the singular parts contributing to the renormalization constants, we perform some simplification of the original integral. We use a very important property of the minimal subtraction scheme that the renormalization constants depend only on dimensionless coupling constants and do not depend on the masses and the choice of external momenta. Therefore, we put all the masses equal to zero, and to avoid artificial infrared divergences, we also put equal to zero one of the external momenta. Then the diagram becomes the propagator type one:

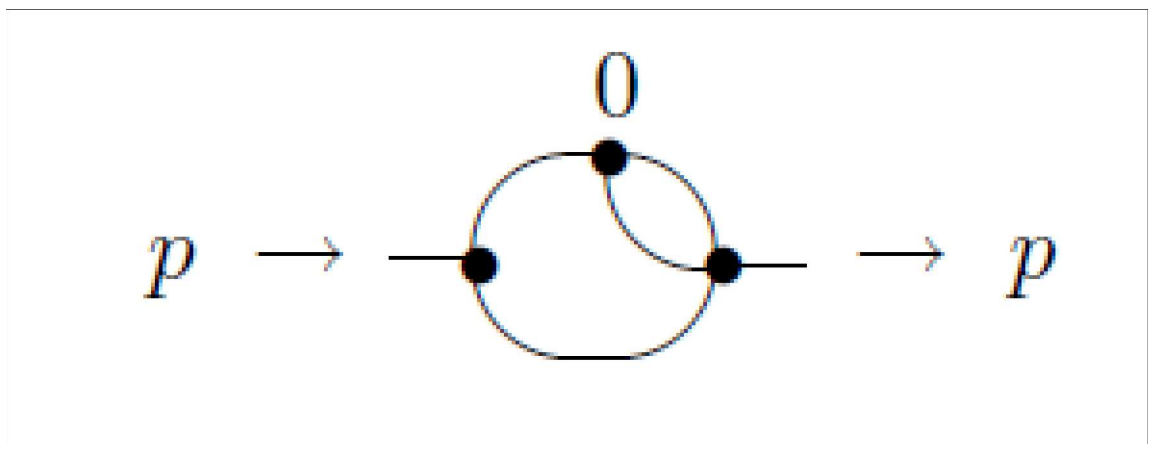

Fig.4.3.6.

The corresponding Colombeau integral is:

$$(I_{22}(p^2;\varepsilon))_\varepsilon = \frac{(-i\lambda)^3}{48}\frac{((\mu^2)^{2\varepsilon})_\varepsilon}{((2\pi)^{8-4\varepsilon})_\varepsilon} i^4 \left(\int \frac{d^{4-2\varepsilon}q d^{4-2\varepsilon}k}{q^2(k-q)^2k^2(p-k)^2}\right)_\varepsilon, \tag{4.3.10}$$

(1/48 is the combinatorial coefficient). Since putting one of the momenta equal to zero we reduced the diagram to the propagator type, we can again use the advocated method to calculate the massless integral. Therefore

$$\begin{aligned}(I_{22}(p^2;\varepsilon))_\varepsilon &= \frac{i\lambda^3}{48}\frac{((\mu^2)^{2\varepsilon})_\varepsilon}{((2\pi)^{8-4\varepsilon})_\varepsilon} i\pi^2 \times \\ &\frac{(\Gamma(1-\varepsilon)\Gamma(1-\varepsilon)\Gamma(\varepsilon))_\varepsilon}{\Gamma(1)\Gamma(1)(\Gamma(2-2\varepsilon))_\varepsilon}\left(\int \frac{d^{4-2\varepsilon}k}{(k^2)^{1+\varepsilon}(p-k)^2}\right)_\varepsilon \\ &= -\frac{i}{48}\frac{\lambda^3}{(16\pi^2)^2}\left(\left(\frac{\mu^2}{-p^2}\right)^{2\varepsilon}\right)_\varepsilon \times \\ &\frac{(\Gamma(1-\varepsilon)\Gamma(1-\varepsilon)\Gamma(\varepsilon)\Gamma(1-2\varepsilon)\Gamma(1-\varepsilon)\Gamma(2\varepsilon))_\varepsilon}{\Gamma(1)\Gamma(1)(\Gamma(2-2\varepsilon)\Gamma(1+\varepsilon)\Gamma(1)\Gamma(2-3\varepsilon))_\varepsilon} \\ &= -\frac{i}{48}\frac{\lambda^3}{(16\pi^2)^2}\left(\left(\frac{\mu^2}{-p^2}\right)^{2\varepsilon}\right)_\varepsilon \frac{1}{2(\varepsilon^2(1-2\varepsilon)(1-3\varepsilon))_\varepsilon} \\ &= -\frac{i}{48}\frac{\lambda^3}{(16\pi^2)^2} \times \\ &\left\{\left(\frac{1}{2\varepsilon^2}\right)_\varepsilon + \left(\frac{5}{2\varepsilon}\right)_\varepsilon + 2 + \left(\frac{\ln(-\mu^2/p^2)}{\varepsilon}\right)_\varepsilon + \ln^2\frac{\mu^2}{-p^2} + 5\ln\frac{\mu^2}{-p^2}\right\}\end{aligned} \tag{4.3.11}$$

As one can see, in this case we again have the second order pole in $\varepsilon$ and, accordingly, the single pole with the logarithm of momentum. The reason of their appearance is the presence of the divergent subgraph. Here we again have to look at the counter-terms of the previous order which eliminate the divergence from the one-loop subgraph. The subtraction of divergent subgraphs (the $\mathcal{R}$-operation without the last subtraction) looks like

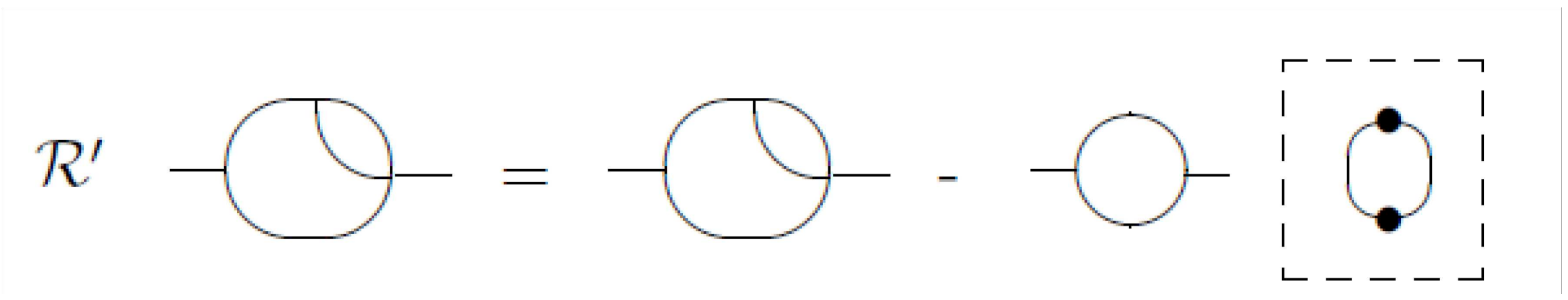

Fig.4.3.7.

$$(\mathcal{R}'I_2(s;\varepsilon))_\varepsilon =$$

$$-\frac{i}{2}\frac{\lambda^3}{(16\pi^2)^2}\left\{\left(\frac{\mu^2}{-p^2}\right)^{2\varepsilon}\frac{1}{2(\varepsilon^2(1-2\varepsilon)(1-3\varepsilon))_\varepsilon}-\left(\frac{\mu^2}{-p^2}\right)^{\varepsilon}\frac{1}{(\varepsilon^2(1-2\varepsilon))_\varepsilon}\right\}$$

$$= -\frac{i}{2}\frac{\lambda^3}{(16\pi^2)^2}\times$$

$$\left\{\left(\left(\frac{1}{2\varepsilon^2}\right)_\varepsilon+\left(\frac{5}{2\varepsilon}\right)_\varepsilon+2+\left(\frac{\ln(-\mu^2/p^2)}{\varepsilon}\right)_\varepsilon+\ln^2\frac{\mu^2}{-p^2}+5\ln\frac{\mu^2}{-p^2}\right)\right. \quad (4.3.12)$$

$$\left.-\left(\left(\frac{1}{\varepsilon^2}\right)_\varepsilon+\left(\frac{2}{\varepsilon}\right)_\varepsilon+4+\left(\frac{\ln(-\mu^2/p^2)}{\varepsilon}\right)_\varepsilon+\frac{1}{2}\ln^2\frac{\mu^2}{-p^2}+2\ln\frac{\mu^2}{-p^2}\right)\right\}=$$

$$= -\frac{i}{2}\frac{\lambda^3}{(16\pi^2)^2}\left\{-\left(\frac{1}{2\varepsilon^2}\right)_\varepsilon+\left(\frac{1}{2\varepsilon}\right)_\varepsilon-2+\frac{1}{2}\ln^2\frac{\mu^2}{-p^2}+3\ln\frac{\mu^2}{-p^2}\right\}.$$

Once again, after the subtraction of the divergent subgraph the singular part is local, i.e. in momentum space does not depend on $\ln p^2$.The contribution to the vertex function from this diagram is:

$$(\Delta\Gamma_{4,\varepsilon})_\varepsilon =$$

$$-i\lambda\left\{\frac{1}{2}\frac{\lambda^2}{(16\pi^2)^2}\left(-\left(\frac{3}{\varepsilon^2}\right)_\varepsilon+\left(\frac{3}{\varepsilon}\right)_\varepsilon-12+\frac{1}{2}\ln^2\frac{\mu^2}{-p^2}+3\ln\frac{\mu^2}{-p^2}+\ldots\right)\right\} \quad (4.3.13)$$

and, accordingly,

$$(\Delta Z_{4,\varepsilon})_\varepsilon = \left(\left(\frac{3}{2\varepsilon^2}\right)_\varepsilon-\left(\frac{3}{2\varepsilon}\right)_\varepsilon\right)\left(\frac{\lambda}{16\pi^2}\right)^2. \quad (4.3.14)$$

Thus, due to (4.2.13) and (4.3.14) in the two-loop approximation the quartic vertex renormalization constant in the $\overline{MS}$ scheme reads:

$$(Z_{4,\varepsilon})_\varepsilon = 1+\left(\frac{3}{2\varepsilon}\right)_\varepsilon\frac{\lambda}{16\pi^2}+\left(\frac{\lambda}{16\pi^2}\right)^2\left[\left(\frac{9}{4\varepsilon^2}\right)_\varepsilon-\left(\frac{3}{2\varepsilon}\right)_\varepsilon\right]. \quad (4.3.15)$$

With taking account of the two-loop renormalization of the propagator (4.3.4) one has:

$$(Z_{\lambda,\varepsilon})_\varepsilon = \left(Z_{4,\varepsilon}Z_{2,\varepsilon}^{-2}\right)_\varepsilon = 1+\left(\frac{3}{2\varepsilon}\right)_\varepsilon\frac{\lambda}{16\pi^2}+\left(\frac{\lambda}{16\pi^2}\right)^2\left[\left(\frac{9}{4\varepsilon^2}\right)_\varepsilon-\left(\frac{17}{12\varepsilon}\right)_\varepsilon\right]. \quad (4.3.16)$$

The Lagrangian (4.2.1) together with the counter-terms (4.3.13)-(4.3.13) can be written as

$$(\mathcal{L}_\varepsilon^{\mathbf{ren}})_\varepsilon = (Z_{2,\varepsilon})_\varepsilon\frac{1}{2}(\partial_\mu\varphi)^2-(Z_{4,\varepsilon})_\varepsilon\frac{\lambda}{4!}\varphi^4. \quad (4.3.17)$$

**Remark 4.3.1**. Note that: (i) if $(Z_{2,\varepsilon})_\varepsilon > 0$ the "bare" Lagrangian reads

$$(\mathcal{L}_\varepsilon^{\mathbf{Bare}})_\varepsilon = \frac{1}{2}(\partial_\mu(\varphi_{B,\varepsilon})_\varepsilon)^2-\frac{1}{4!}(\lambda_{B,\varepsilon})_\varepsilon(\varphi_{B,\varepsilon}^4)_\varepsilon \quad (4.3.18)$$

and by using "bare" Lagrangian (4.3.18) we obtain the scalar theory $\lambda\varphi_4^4$ in standard sector such that that the two-loop radiative corrections calculated from the Lagrangian (4.3.18) with Colombeau parameters chosen according to (4.3.20),
(ii) if $(Z_{2,\varepsilon})_\varepsilon < 0$ the "bare" Lagrangian reads

$$(\mathcal{L}_\varepsilon^{\mathbf{Bare}})_\varepsilon = -\frac{1}{2}(\partial_\mu(\varphi_{B,\varepsilon})_\varepsilon)^2-\frac{1}{4!}(\lambda_{B,\varepsilon})_\varepsilon(\varphi_{B,\varepsilon}^4)_\varepsilon. \quad (4.3.19)$$

Lagrangian (4.3.19) contain wrong kinetic term,i.e.,kinetic term with a wrong sign corresponding to a "bad" ghosts and therefore by using "bare" Lagrangian (4.3.19) we obtain the scalar theory $\lambda\varphi_4^4$ in a ghost sector such that the two-loop radiative corrections calculated from the Lagrangian (4.3.19) with Colombeau parameters chosen according to (4.3.20), where

$$(\varphi_{B,\varepsilon})_\varepsilon = \begin{cases} \left[\left(\sqrt{Z_{2,\varepsilon}}\right)_\varepsilon\right]\varphi & \text{if } (Z_{2,\varepsilon})_\varepsilon > 0 \\ \left[\left(\sqrt{-Z_{2,\varepsilon}}\right)_\varepsilon\right]\varphi & \text{if } (Z_{2,\varepsilon})_\varepsilon < 0 \end{cases}$$

$$(\lambda_{B,\varepsilon})_\varepsilon = \left(Z_{4,\varepsilon}Z_{2,\varepsilon}^{-2}\lambda(\varepsilon)\right)_\varepsilon,$$

$$\left(Z_\varepsilon^{\overline{MS}}\right)_\varepsilon = 1 + \left[\left(\frac{1}{2\varepsilon}\right)_\varepsilon + 1 - \gamma_E + \log(4\pi)\right]\frac{(\lambda(\varepsilon))_\varepsilon}{16\pi^2}, \tag{4.3.20}$$

$$(Z_{2,\varepsilon})_\varepsilon = 1 - \left(\frac{1}{24\varepsilon}\right)_\varepsilon\left(\frac{(\lambda(\varepsilon))_\varepsilon}{16\pi^2}\right)^2,$$

$$(Z_{4,\varepsilon})_\varepsilon = 1 + \left(\frac{3}{2\varepsilon}\right)_\varepsilon\frac{(\lambda(\varepsilon))_\varepsilon}{16\pi^2} + \left(\frac{(\lambda(\varepsilon))_\varepsilon}{16\pi^2}\right)^2\left[\left(\frac{9}{4\varepsilon^2}\right)_\varepsilon - \left(\frac{3}{2\varepsilon}\right)_\varepsilon\right].$$

**Remark 4.3.2**. Note that: (i) if $(Z_{2,\varepsilon})_\varepsilon < 0$, i.e.,

$$\left(\frac{1}{24\varepsilon}\right)_\varepsilon\left(\frac{(\lambda(\varepsilon))_\varepsilon}{16\pi^2}\right)^2 > 1 \tag{4.3.21}$$

and therefore we get

$$(\lambda_{B,\varepsilon})_\varepsilon = -\left\{(\lambda(\varepsilon))_\varepsilon + \left(\frac{3}{2\varepsilon}\right)_\varepsilon\frac{(\lambda^2(\varepsilon))_\varepsilon}{16\pi^2} + \frac{(\lambda^3(\varepsilon))_\varepsilon}{16^2\pi^4}\left[\left(\frac{9}{4\varepsilon^2}\right)_\varepsilon - \left(\frac{3}{2\varepsilon}\right)_\varepsilon\right]\right\}\times \left[\left|1 - \left(\frac{1}{24\varepsilon}\right)_\varepsilon\left(\frac{(\lambda(\varepsilon))_\varepsilon}{16\pi^2}\right)^2\right|\right]^{-2} \tag{4.3.22}$$

with $(\lambda(\varepsilon))_\varepsilon > 8\pi\sqrt{6}\,(\varepsilon^{1/2})_\varepsilon, \varepsilon \in (0,1]$,

(ii) ander condition (4.3.21) Lagrangian (4.3.17) obviously contain wrong kinetic term, i.e.,kinetic term with a wrong sign corresponding to a "bad" ghosts and therefore in two-loop approximation we obtain the scalar theory $\lambda\varphi_4^4$ in a ghost sector.

**Remark 4.3.2**. Note that: (i) if $(Z_{2,\varepsilon})_\varepsilon > 0$, i.e.,

$$\left(\frac{1}{24\varepsilon}\right)_\varepsilon\left(\frac{(\lambda(\varepsilon))_\varepsilon}{16\pi^2}\right)^2 < 1 \tag{4.3.23}$$

and therefore we get

$$(\lambda_{B,\varepsilon})_\varepsilon = \left\{(\lambda(\varepsilon))_\varepsilon + \left(\frac{3}{2\varepsilon}\right)_\varepsilon\frac{(\lambda^2(\varepsilon))_\varepsilon}{16\pi^2} + \frac{(\lambda^3(\varepsilon))_\varepsilon}{16^2\pi^4}\left[\left(\frac{9}{4\varepsilon^2}\right)_\varepsilon - \left(\frac{3}{2\varepsilon}\right)_\varepsilon\right]\right\}\times \left[\left|1 - \left(\frac{1}{24\varepsilon}\right)_\varepsilon\left(\frac{(\lambda(\varepsilon))_\varepsilon}{16\pi^2}\right)^2\right|\right]^{-2} \tag{4.3.24}$$

with $(\lambda(\varepsilon))_\varepsilon < 8\pi\sqrt{6}\,(\varepsilon^{1/2})_\varepsilon, \varepsilon \in (0,1]$,

(ii) ander condition (4.3.23) we obtain the scalar theory $\lambda\varphi_4^4$ in standard sector.

The statement is that the counter-terms introduced this way eliminate all the ultraviolet infinite large Colombeau objects up to two-loop order and make the Green functions

and

hence the radiative corrections finite. In the case of nonzero mass, one should also add the mass counter-term.

The Lagrangian (4.2.1) together with the counter-terms in this case can be written as

$$(\mathcal{L}_{\varepsilon}^{ren})_{\varepsilon} = (Z_{2,\varepsilon})_{\varepsilon}\frac{1}{2}(\partial_{\mu}\varphi)^2 - (Z_{\varepsilon})_{\varepsilon}\frac{m^2}{2}\varphi^2 - (Z_{4,\varepsilon})_{\varepsilon}\frac{\lambda}{4!}\varphi^4 \quad (4.3.25)$$

**Remark 4.3.3**. Note that: (i) if we assume that $(Z_{2,\varepsilon})_{\varepsilon} > 0$ we obtain the Lagrangian $(\mathcal{L}_{\varepsilon}^{ren,s.m.s.})_{\varepsilon}$ corresponding to standard matter sector and the "bare" Lagrangian corresponding to standard matter sector reads

$$(\mathcal{L}_{\varepsilon}^{Bare,s.m.s})_{\varepsilon} = \frac{1}{2}(\partial_{\mu}(\varphi_{B,\varepsilon}^{s.m.s})_{\varepsilon})^2 - \frac{\left((m_{B,\varepsilon}^{s.m.s.})^2\right)_{\varepsilon}}{2}\left((\varphi_{B,\varepsilon}^{s.m.s.})^2\right)_{\varepsilon} - \frac{1}{4!}(\lambda_{B,\varepsilon}^{s.m.s})_{\varepsilon}\left((\varphi_{B,\varepsilon}^{s.m.s})^4\right)_{\varepsilon}, \quad (4.3.26.a)$$

and by using "bare" Lagrangian (4.3.26) we obtain the scalar theory $\lambda\varphi_4^4$ in standard sector, where

$$\begin{aligned}(\varphi_{B,\varepsilon}^{s.m.s.})_{\varepsilon} &= \left[\left(\sqrt{Z_{2,\varepsilon}}\right)_{\varepsilon}\right]\varphi \quad \text{if } (Z_{2,\varepsilon})_{\varepsilon} > 0\\ \left((m_{B,\varepsilon}^{s.m.s.})^2\right)_{\varepsilon} &= \left[(Z_{\varepsilon})_{\varepsilon}(Z_{2,\varepsilon}^{-1})_{\varepsilon}\right]m^2 \quad \text{if } (Z_{2,\varepsilon})_{\varepsilon} > 0 \quad (4.3.26.b)\\ (\lambda_{B,\varepsilon})_{\varepsilon} &= \left(Z_{4,\varepsilon}Z_{2,\varepsilon}^{-2}\lambda(\varepsilon)\right)_{\varepsilon}.\end{aligned}$$

(ii) if if we assume that $(Z_{2,\varepsilon})_{\varepsilon} < 0$ we obtain the Lagrangian $(\mathcal{L}_{\varepsilon}^{ren,g.m.s.})_{\varepsilon}$ corresponding to

a ghost matter sector and the "bare" Lagrangian corresponding to a ghost matter sector

reads

$$\left(\mathcal{L}_{\varepsilon}^{Bare,g.m.s.}\right)_{\varepsilon} = -\frac{1}{2}(\partial_{\mu}(\varphi_{B,\varepsilon}^{g.m.s.})_{\varepsilon})^2 - \frac{\left((m_{B,\varepsilon}^{g.m.s.})^2\right)_{\varepsilon}}{2}\left((\varphi_{B,\varepsilon}^{g.m.s.})^2\right)_{\varepsilon} - \frac{1}{4!}(\lambda_{B,\varepsilon}^{g.m.s.})_{\varepsilon}\left((\varphi_{B,\varepsilon}^{g.m.s.})^4\right)_{\varepsilon} \quad (4.3.27)$$

and by using "bare" Lagrangian (4.3.27) we obviously obtain the scalar theory $\lambda\varphi_4^4$ in a ghost sector,where

$$\begin{aligned}(\varphi_{B,\varepsilon}^{g.m.s.})_{\varepsilon} &= \left[\left(\sqrt{-Z_{2,\varepsilon}}\right)_{\varepsilon}\right]\varphi \quad \text{if } (Z_{2,\varepsilon})_{\varepsilon} < 0\\ \left((m_{B,\varepsilon}^{g.m.s.})^2\right)_{\varepsilon} &= \left[(Z_{\varepsilon})_{\varepsilon}(-Z_{2,\varepsilon}^{-1})_{\varepsilon}\right]m^2 \quad \text{if } (Z_{2,\varepsilon})_{\varepsilon} < 0 \quad (4.3.28)\\ (\lambda_{B,\varepsilon}^{g.m.s.})_{\varepsilon} &= \left(Z_{4,\varepsilon}Z_{2,\varepsilon}^{-2}\lambda(\varepsilon)\right)_{\varepsilon}.\end{aligned}$$

**Remark 4.3.4**. Recall that classical Schwartz distribution is defined as linear functionals on a test smooth functions [32]. Schwartz distributions may be multiplied by real numbers and added together, so they form a real vector space. Schwartz distributions may also be multiplied by infinitely differentiable functions, but it is not possible to define a product of
general distributions that extends the usual pointwise product of functions and has the same algebraic properties. This result was shown by Schwartz (1954), and is usually referred to as the Schwartz Impossibility Theorem [32].

**Remark 4.3.5.** In coordinate space the large values of momenta correspond to the small
distances. Hence, the ultraviolet divergences allow for the singularities at small distances.

Indeed, the simplest divergent loop diagram (Fig.3.1.3) in coordinate space is the product
of two propagators. Each Euclidean propagator $\Delta(\vec{x}) = \Delta(r) \in D'(\mathbb{R}^4), \vec{x} \in \mathbb{R}^4, r = \|x\|$ is
uniquely defined in momentum as well as in coordinate space, but the square of the propagator has already an ill-defined Fourier transform, it is ultraviolet divergent. The reason is that the square of the propagator is singular as $r^2 \to 0$ and behaves like

$$[\Delta(r)]^2 \sim 1/r^4. \tag{4.3.29}$$

In fact, the causal Green function of the QFT is the classical Schwartz distribution which
is defined on a test smooth functions. It has the $\delta$-function like singularities and needs an
additional definition for the product of several such functions at a single point.The discussed above diagram (see Fig.3.1.3) is precisely this product.

**Remark 4.3.6**. In handbook [8] N. N. Bogoliubov argue that a problem of the ultraviolet divergences arises exactly from the Schwartz Impossibility Theorem [32].

**Remark 4.3.7**. Note that:(i) by gomomorfism $\Delta(r^2) \to \Delta(r^2;\varepsilon) = \Delta(r^2+\varepsilon)$ we can embed
the distribution $\Delta(r^2)$ into the Colombeau algebra $\mathcal{G}(\mathbb{R}^4_x)$ :
$\Delta(r^2) \hookrightarrow (\Delta(r^2;\varepsilon))_{\varepsilon\in(0,1]} \in \mathcal{G}(\mathbb{R}^4_x)$,

(ii) in Colombeau algebra $\mathcal{G}(\mathbb{R}^4_x)$ the square of the propagator is $(\Delta^2(r^2;\varepsilon))_\varepsilon \in \mathcal{G}(\mathbb{R}^4_x)$

$$(\Delta^2(r^2;\varepsilon))_\varepsilon \sim 1/(r^2+\varepsilon)^2 \in \mathcal{G}(\mathbb{R}^4_x), \tag{4.3.30}$$

(iii) Colombeau Fourier transform $C\mathcal{F}[\Delta^2(r^2;\varepsilon)_\varepsilon] \triangleq (\mathcal{F}[\Delta^2(r^2;\varepsilon)])_\varepsilon$ (see [21-22]) of the square of the propagator well defined and $C\mathcal{F}[\Delta^2(r^2;\varepsilon)_\varepsilon] \in \mathcal{G}(\mathbb{R}^4_p)$.

**Remark 4.3.8**. Note that in sect.IV Colombeau Fourier transform of the square of the propagator has been defined directly in momentum space by using dimensional regularization see Eq.(4.1.7).

**Remark 4.3.9.** Note that:(i) by using Colombeau algebra $\mathcal{G}(\mathbb{R}^4_x)$ and Colombeau Fourier transform $C\mathcal{F} : \mathcal{G}(\mathbb{R}^4_x) \to \mathcal{G}(\mathbb{R}^4_p)$ there is no any problem arises from the Schwartz Impossibility Theorem, (ii) classical ill-defined ultraviolet divergences replased by well defined infinite large Colombeau generalized numbers.

**Remark 4.3.10**. Thus by using Colombeau algebra $\mathcal{G}(\mathbb{R}^4_x)$ (see [21-22]) of the Colombeau
generalized functions instead classical Schwartz distribution and Colombeau generalized numbers $\widetilde{\mathbb{R}}$ (see [21]-[22],[23]-[26]) there is no any problem arises from the Schwartz Impossibility Theorem.

## 4.4.Quantum electrodynamics in a ghost sector via Colombeau generalized functions.

Let us consider now the calculation of the diagrams in the gauge theories. We start with quantum electrodynamics. The $QED_4$ Lagrangian has the form

$$\mathcal{L}_{QED} = -\frac{1}{4}F_{\mu\nu}^2 + \bar{\psi}(i\gamma^\mu\partial_\mu - m)\psi + e\bar{\psi}\gamma^\mu A_\mu\psi - \frac{1}{2\xi}(\partial_\mu A_\mu)^2, \qquad (4.4.1)$$

where the electromagnetic stress tensor is $F_{\mu\nu} = \partial_\mu A_\nu - \partial_\nu A_\mu$, and the last term in (4.4.1) fixes the gauge. In what follows we choose the Feynman or the diagonal gauge $(\xi = 1)$.The Feynman rules corresponding to the Lagrangian (4.4.1) are shown in Fig.4.4.1.

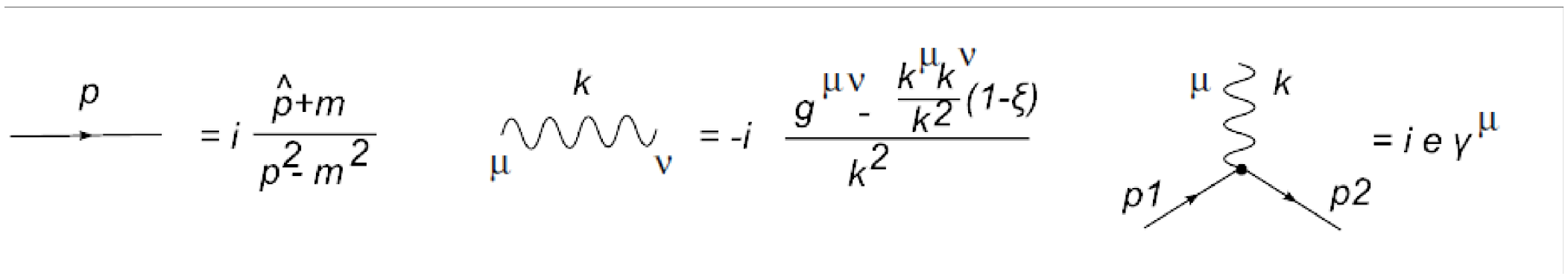

Fig.4.4.1.

In quantum electrodynamics the divergences appear only in the photon propagator, the electron propagator, and the triple vertex. The one-loop divergent diagrams are shown in Fig.4.4.2.

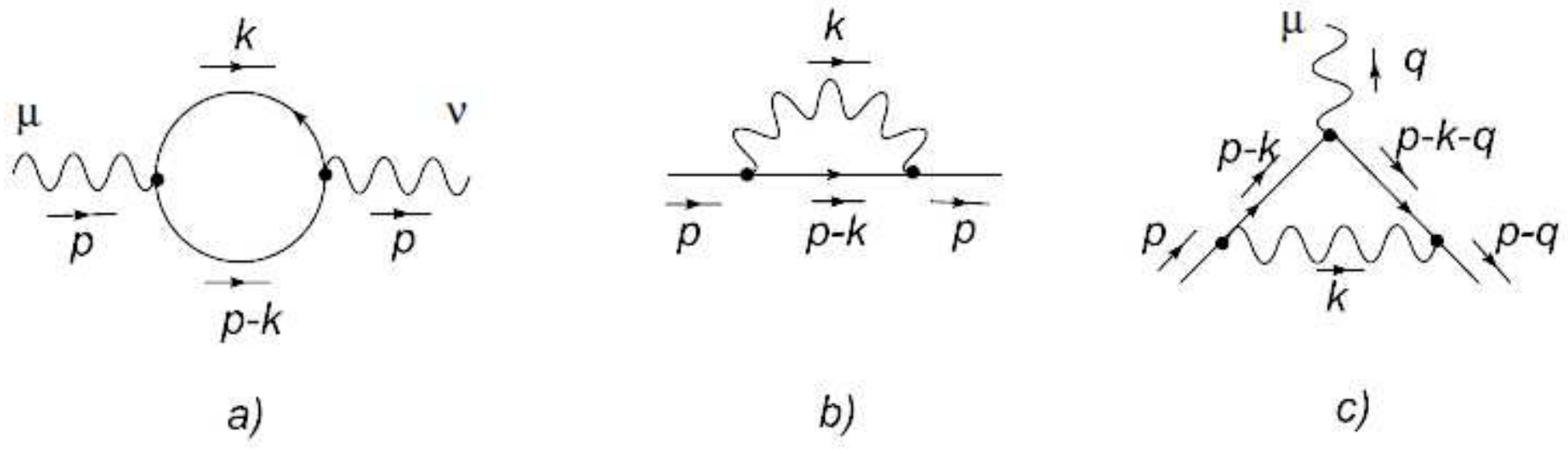

Fig.4.4.2.Theone-loop divergent diagrams in $QED$

We begin with the vacuum polarization graph. It is given by the diagram shown in Fig.4.4.2a). The corresponding formal expression is

$$\Pi_{\mu\nu}(p) = (-)\frac{e^2}{(2\pi)^4}\int d^4k\frac{Tr[\gamma^\mu(m+\hat{k})\gamma^\nu(m+\hat{k}-\hat{p})]}{[m^2-k^2][m^2-(k-p)^2]}, \qquad (4.4.2)$$

where the "-" sign comes from the fermion loop and $\hat{q} \equiv \gamma^\mu q_\mu$. We set now $D = 4-2\varepsilon$, $\varepsilon \in (0,1]$.Then the integral (4.4.2) becomes to Colombeau

$$(\Pi_{\mu\nu}^{Dim}(p;\varepsilon))_\varepsilon = (-)\frac{e^2((\mu^2)^\varepsilon)_\varepsilon}{((2\pi)^{4-2\varepsilon})_\varepsilon}\left(\int d^{4-2\varepsilon}k\frac{Tr[\gamma^\mu(m+\hat{k})\gamma^\nu(m+\hat{k}-\hat{p})]}{[m^2-k^2][m^2-(k-p)^2]}\right)_\varepsilon, \quad (4.4.3)$$

where $(\Pi_{\mu\nu}^{Dim}(p;\varepsilon))_\varepsilon \in \mathcal{G}(\mathbb{R}_p^4)$.We put now $m = 0$ for simplicity. From (4.4.3) one obtains

$$(\Pi_{\mu\nu}^{Dim}(p;\varepsilon))_\varepsilon =$$

$$i\frac{4e^2}{16\pi^2}((4\pi)^\varepsilon)_\varepsilon\left(\left(-\frac{\mu^2}{p^2}\right)^\varepsilon\right)_\varepsilon\left(\frac{\Gamma^2(2-\varepsilon)\Gamma(\varepsilon)}{\Gamma(4-2\varepsilon)}\right)_\varepsilon[2p^\mu p^\nu - g^{\mu\nu}p^2 - g^{\mu\nu}p^2] \quad (4.4.4)$$

$$= -i\frac{8e^2}{16\pi^2}((4\pi)^\varepsilon)_\varepsilon\left(\left(-\frac{\mu^2}{p^2}\right)^\varepsilon\right)_\varepsilon(g^{\mu\nu}p^2 - p^\mu p^\nu)\left(\frac{\Gamma^2(2-\varepsilon)\Gamma(\varepsilon)}{\Gamma(4-2\varepsilon)}\right)_\varepsilon.$$

Expanding now over $\varepsilon$ with the help of Eqs.(4.4.5)

$$(\Gamma(\varepsilon))_\varepsilon = \left(\frac{1}{\varepsilon}\right)_\varepsilon(\Gamma(1+\varepsilon))_\varepsilon,\ (\Gamma(2-\varepsilon))_\varepsilon = ((1-\varepsilon)\Gamma(1-\varepsilon))_\varepsilon,$$

$$(\Gamma(4-2\varepsilon))_\varepsilon = ((3-2\varepsilon)(2-2\varepsilon)(1-2\varepsilon))_\varepsilon(\Gamma(1-2\varepsilon))_\varepsilon, \quad (4.4.5)$$

$$\varepsilon \in (0,1],$$

we obtain

$$(\Pi_{\mu\nu}^{Dim}(p;\varepsilon))_\varepsilon =$$

$$-i\frac{e^2}{16\pi^2}((4\pi)^\varepsilon)_\varepsilon\left(\left(-\frac{\mu^2}{p^2}\right)^\varepsilon\right)_\varepsilon(g^{\mu\nu}p^2 - p^\mu p^\nu)\left(\frac{4(1+5/3\varepsilon)}{3\varepsilon}\right)_\varepsilon e^{-\gamma\varepsilon}$$

$$= -ie^2\frac{g^{\mu\nu}p^2 - p^\mu p^\nu}{16\pi^2}\frac{4}{3}\left[\left(\frac{1}{\varepsilon}\right)_\varepsilon - \gamma_E + \log 4\pi + \log\frac{-\mu^2}{p^2} + \frac{5}{3}\right] \quad (4.4.6)$$

$$= i(g^{\mu\nu}p^2 - p^\mu p^\nu)(\Pi^{Dim}(p^2;\varepsilon))_\varepsilon,$$

where

$$(\Pi^{Dim}(p^2;\varepsilon))_\varepsilon = -\frac{e^2}{16\pi^2}\frac{4}{3}\left[\left(\frac{1}{\varepsilon}\right)_\varepsilon - \gamma_E + \log 4\pi + \log\frac{-\mu^2}{p^2} + \frac{5}{3}\right]. \quad (4.4.7)$$

Given the expression for the vacuum polarization one can construct the photon propagator as shown in Fig.4.4.3.

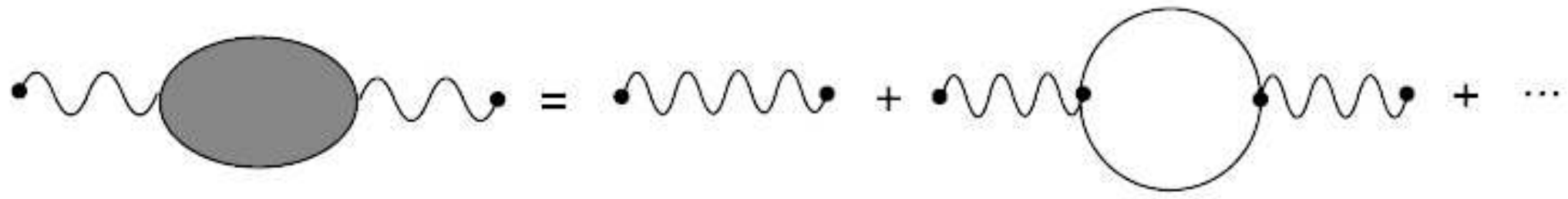

Fig.4.4.3.The photon propagator in $QED$

One has

$$(G_{\mu\nu}(p;\varepsilon))_\varepsilon = \frac{-i}{p^2}g^{\mu\nu} + \frac{-i}{p^2}g^{\mu\rho}(\Pi_{\rho\sigma,\varepsilon})_\varepsilon\frac{-i}{p^2}g^{\sigma\nu} + \cdots$$

$$= \frac{-i}{p^2}g^{\mu\nu} - \frac{(\Pi_\varepsilon^{\mu\nu})_\varepsilon}{p^4} + \cdots = \frac{-i}{p^2}g^{\mu\nu} - \frac{i(g^{\mu\nu} - p^\mu p^\nu/p^2)}{p^2}(\Pi(p^2;\varepsilon))_\varepsilon + \cdots \quad (4.4.8)$$

$$= \frac{-i}{p^2}(g^{\mu\nu} - \frac{p^\mu p^\nu}{p^2})(1 + (\Pi(p^2;\varepsilon))_\varepsilon + \cdots) - \frac{i}{p^2}\frac{p^\mu p^\nu}{p^2},$$

where $(\Pi(p^2;\varepsilon))_\varepsilon$ is given by eq.(4.4.7). Notice that the radiative corrections are always proportional to the transverse tensor $P_{\mu\nu} = g_{\mu\nu} - p_\mu p_\nu/p^2$. This is a consequence of the gauge invariance and follows from the Ward identities

Let us consider now the electron self-energy graph Fig.4.4.2.b). The corresponding formal expression is

$$\Sigma(\hat{p}) = -\frac{e^2}{(2\pi)^4}\int d^4k \frac{\gamma^\mu(\hat{p}-\hat{k}+m)\gamma^\mu}{k^2[(p-k)^2-m^2]}. \tag{4.4.9}$$

Acting in a usual way we go to fractal dimension $D = 4-2\varepsilon, \varepsilon \in (0,1]$ convert the indices of the $\gamma$-matrices and introduce the Feynman parametrization. The result is

$$(\Sigma^{Dim}(\hat{p};\varepsilon))_\varepsilon = -\frac{e^2((\mu^2)^\varepsilon)_\varepsilon}{((2\pi)^{4-2\varepsilon})_\varepsilon}\int_0^1 dx\left(\int \frac{d^{4-2\varepsilon}k[-2(1-\varepsilon)(\hat{p}-\hat{k})+(4-2\varepsilon)m]}{[k^2-2kpx+p^2x-m^2x]^2}\right)_\varepsilon \tag{4.4.10}$$

The integral over $k$ can now be evaluated according to the standard formulas

$$(\Sigma^{Dim}(\hat{p};\varepsilon))_\varepsilon = -i\frac{e^2}{16\pi^2}\frac{((-\mu^2)^\varepsilon)_\varepsilon}{((4\pi)^{-\varepsilon})_\varepsilon}(\Gamma(\varepsilon))_\varepsilon\left(\int_0^1 dx \frac{-2(1-\varepsilon)\hat{p}(1-x)+(4-2\varepsilon)m}{[p^2x(1-x)-m^2x]^\varepsilon}\right)_\varepsilon. \tag{4.4.11}$$

This expression can be expanded in series in $\varepsilon, \varepsilon \in (0,1]$

$$\begin{aligned}(\Sigma^{Dim}(\hat{p};\varepsilon))_\varepsilon = -i\frac{e^2}{16\pi^2}\Bigg[&-\left(\frac{\hat{p}-4m}{\varepsilon}\right)_\varepsilon + \hat{p}-2m-(\hat{p}-4m)(-\gamma_E+\log(4\pi)) \\ &+\int_0^1 dx[2\hat{p}(1-x)-4m]\log\frac{p^2x(1-x)-m^2x}{-\mu^2}\Bigg].\end{aligned} \tag{4.4.12}$$

At last, consider the vertex function Fig.4.4.2c). The corresponding formal expressionis is

$$\Gamma_1(p,q) = \frac{e^3}{(2\pi)^4}\int d^4k \frac{\gamma^\nu(\hat{p}-\hat{k}-\hat{q}+m)\gamma^\mu(\hat{p}-\hat{k}+m)\gamma^\nu}{[(p-k-q)^2-m^2][(p-k)^2-m^2]k^2}. \tag{4.4.13}$$

Transfer (4.4.13) to dimension $4-2\varepsilon, \varepsilon \in (0,1]$ and introduce the Feynman parametrization. This gives corresponding Colombeau integral

$$\begin{aligned}(\Gamma_1(p,q;\varepsilon))_\varepsilon &= \frac{(e^3(\mu^2)^\varepsilon)_\varepsilon}{((2\pi)^{4-2\varepsilon})_\varepsilon}\Gamma(3)\int_0^1 dx\int_0^x dy \\ &\times\left(\int \frac{d^{4-2\varepsilon}k[\gamma^\nu(\hat{p}-\hat{k}-\hat{q}+m)\gamma^\mu(\hat{p}-\hat{k}+m)\gamma^\nu]}{[((p-k-q)^2-m^2)y+((p-k)^2-m^2)(x-y)+k^2(1-x)]^3}\right)_\varepsilon.\end{aligned} \tag{4.4.14}$$

The integral over $k$ is straightforward and reads

$$\begin{aligned}(\Gamma_1^{Dim}(p,q;\varepsilon))_\varepsilon &= ie\frac{e^2}{16\pi^2}\frac{((-\mu^2)^\varepsilon)_\varepsilon}{((4\pi)^{-\varepsilon})_\varepsilon}\int_0^1 dx\int_0^x dy \\ &\Bigg\{(\Gamma(1+\varepsilon))_\varepsilon\frac{[\gamma^\nu(\hat{p}(1-x)-\hat{q}(1-y)+m)\gamma^\mu(\hat{p}(1-x)+\hat{q}y+m)\gamma^\nu]}{([(p-q)^2y(1-x)+p^2(1-x)(x-y)+q^2y(x-y)-m^2x]^{1+\varepsilon})_\varepsilon} \\ &+\frac{(\Gamma(\varepsilon))_\varepsilon}{2}\frac{\gamma^\nu\gamma^\rho\gamma^\mu\gamma^\rho\gamma^\nu}{([(p-q)^2y(1-x)+p^2(1-x)(x-y)+q^2y(x-y)-m^2x]^\varepsilon)_\varepsilon}\Bigg\}.\end{aligned} \tag{4.4.15}$$

As one can see, the first integral is finite Colombeau quantity and the second one is logarithmically divergent. Expanding in series in $\varepsilon, \varepsilon \in (0,1]$ we obtain

$$(\Gamma_1(p,q;\varepsilon))_\varepsilon = ie\frac{e^2}{16\pi^2}\left\{\left(\frac{\gamma^\mu}{\varepsilon}\right)_\varepsilon - 2\gamma^\mu - \gamma^\mu(\gamma_E - \ln(4\pi))\right.$$
$$-2\gamma^\mu\int_0^1 dx\int_0^x dy\ \ln\left[\frac{(p-q)^2y(1-x)+p^2(1-x)(x-y)+q^2y(x-y)-m^2x}{-\mu^2}\right] \quad (4.4.16)$$
$$\left.+\int_0^1 dx\int_0^x dy\frac{\gamma^\nu(\hat{p}(1-x)-\hat{q}(1-y)+m)\gamma^\mu(\hat{p}(1-x)+\hat{q}y+m)\gamma^\nu}{(p-q)^2y(1-x)+p^2(1-x)(x-y)+q^2y(x-y)-m^2x}\right\}$$

Quantum electrodynamics (4.4.1) is a renormalizable theory; hence, the Colombeau counter-terms repeat the structure of the Lagrangian. They can be written as

$$\left(\Delta\mathcal{L}_\varepsilon^{QED}\right)_\varepsilon =$$
$$-\frac{(Z_{3,\varepsilon})_\varepsilon - 1}{4}F_{\mu\nu}^2 + ((Z_{2,\varepsilon})_\varepsilon - 1)i\bar{\psi}\hat{\partial}\psi - m((Z_\varepsilon)_\varepsilon - 1)\bar{\psi}\psi + e((Z_{1,\varepsilon})_\varepsilon - 1)\bar{\psi}\hat{A}\psi. \quad (4.4.17)$$

The term that fixes the gauge is not renormalized. In the leading order of perturbation theory we calculated the corresponding diagrams with the help of dimensional regularization mentioned above. Their infinite large Colombeau parts with the opposite sign give the proper Colombeau renormalization constants. They are, respectively,

$$(Z_{1,\varepsilon})_\varepsilon = 1 - \frac{e^2}{16\pi^2}\left(\frac{1}{\varepsilon}\right)_\varepsilon, (Z_{2,\varepsilon})_\varepsilon = 1 - \frac{e^2}{16\pi^2}\left(\frac{1}{\varepsilon}\right)_\varepsilon,$$
$$(Z_{3,\varepsilon})_\varepsilon = 1 - \frac{e^2}{16\pi^2}\left(\frac{4}{3\varepsilon}\right)_\varepsilon, (Z_\varepsilon)_\varepsilon = 1 - \frac{e^2}{16\pi^2}\left(\frac{4}{\varepsilon}\right)_\varepsilon. \quad (4.4.18)$$

**Remark 4.4.1**. We assume now that

$$(Z_{1,\varepsilon})_\varepsilon > 0, (Z_{2,\varepsilon})_\varepsilon > 0, (Z_{3,\varepsilon})_\varepsilon > 0, (Z_\varepsilon)_\varepsilon > 0. \quad (4.4.19)$$

Adding (4.4.1) with (4.4.17) from (4.4.19) we obtain the $QED_4$ bare Lagrangian in standard matter sector

$$\left(\mathcal{L}_{QED,\varepsilon}^{Bare,s.m.s.}\right)_\varepsilon = \mathcal{L}_{QED} + \left(\Delta\mathcal{L}_\varepsilon^{QED}\right)_\varepsilon = -\frac{(Z_{3,\varepsilon})_\varepsilon}{4}F_{\mu\nu}^2 + (Z_{2,\varepsilon})_\varepsilon i\bar{\psi}\hat{\partial}\psi - m(Z_\varepsilon)_\varepsilon\bar{\psi}\psi +$$
$$e(Z_{1,\varepsilon})_\varepsilon\bar{\psi}\hat{A}\psi - \frac{1}{2\xi}(\partial_\mu A_\mu)^2 =$$
$$= -\frac{1}{4}\left((F_{\mu\nu B,\varepsilon}^{s.m.s.})^2\right)_\varepsilon + i(\bar{\psi}_{B,\varepsilon}^{s.m.s.})_\varepsilon\hat{\partial}(\psi_{B,\varepsilon}^{s.m.s.})_\varepsilon - m(Z_\varepsilon Z_{2,\varepsilon}^{-1})_\varepsilon(\bar{\psi}_{B,\varepsilon}^{s.m.s.})_\varepsilon(\psi_{B,\varepsilon}^{s.m.s.})_\varepsilon +$$
$$e\left(Z_{1,\varepsilon}Z_{2,\varepsilon}^{-1}Z_{3,\varepsilon}^{-1/2}\right)_\varepsilon(\bar{\psi}_{B,\varepsilon}^{s.m.s.})_\varepsilon\left(\hat{A}_{B,\varepsilon}^{s.m.s.}\right)_\varepsilon(\psi_{B,\varepsilon}^{s.m.s.})_\varepsilon \ - \frac{(Z_{3,\varepsilon}^{-1})_\varepsilon}{2\xi}(\partial_\mu(A_{\mu B,\varepsilon}^{s.m.s.})_\varepsilon)^2 = \quad (4.4.20)$$
$$-\frac{1}{4}\left((F_{\mu\nu B,\varepsilon}^{s.m.s})^2\right)_\varepsilon + i(\bar{\psi}_{B,\varepsilon}^{s.m.s})_\varepsilon\hat{\partial}(\psi_{B,\varepsilon}^{s.m.s})_\varepsilon - (m_{B,\varepsilon}^{s.m.s})_\varepsilon(\bar{\psi}_{B,\varepsilon}^{s.m.s.})_\varepsilon(\psi_{B,\varepsilon}^{s.m.s.})_\varepsilon +$$
$$(e_{B,\varepsilon}^{s.m.s.})_\varepsilon(\bar{\psi}_{B,\varepsilon}^{s.m.s.})_\varepsilon\left(\hat{A}_{B,\varepsilon}^{s.m.s.}\right)_\varepsilon(\psi_{B,\varepsilon})_\varepsilon \ - \frac{1}{2\left(\xi_{B,\varepsilon}^{s.m.s.}\right)_\varepsilon}(\partial_\mu(A_{\mu B,\varepsilon}^{s.m.s.})_\varepsilon)^2,$$

where

$$(\psi_{B,\varepsilon}^{s.m.s})_\varepsilon = \left(Z_{2,\varepsilon}^{1/2}\right)_\varepsilon\psi,\ \ (A_{B,\varepsilon}^{s.m.s})_\varepsilon = \left(Z_{3,\varepsilon}^{1/2}\right)_\varepsilon A,$$
$$(m_{B,\varepsilon}^{s.m.s})_\varepsilon = (Z_\varepsilon Z_{2,\varepsilon}^{-1})_\varepsilon m,\ \ (e_{B,\varepsilon}^{s.m.s})_\varepsilon = \left(Z_{1,\varepsilon}Z_{2,\varepsilon}^{-1}Z_{3,\varepsilon}^{-1/2}\right)_\varepsilon e, \quad (4.4.21)$$
$$\left(\xi_{B,\varepsilon}^{s.m.s}\right)_\varepsilon = (Z_{3,\varepsilon})_\varepsilon\xi.$$

**Remark 4.4.2**. We assume now that

$$(Z_{1,\varepsilon})_\varepsilon < 0, (Z_{2,\varepsilon})_\varepsilon < 0, (Z_{3,\varepsilon})_\varepsilon < 0, (Z_\varepsilon)_\varepsilon < 0. \qquad (4.4.22)$$

Adding (4.4.1) with (4.4.17) from (4.4.19) we obtain the $QED_4$ bare Lagrangian in a ghost

matter sector

$$\begin{gathered}
\left(\mathcal{L}^{Bare,g.m.s.}_{QED,\varepsilon}\right)_\varepsilon = \mathcal{L}_{QED} + \left(\Delta\mathcal{L}^{QED}_\varepsilon\right)_\varepsilon = -\frac{(Z_{3,\varepsilon})_\varepsilon}{4} F^2_{\mu\nu} + (Z_{2,\varepsilon})_\varepsilon i\bar{\psi}\hat{\partial}\psi - m(Z_\varepsilon)_\varepsilon \bar{\psi}\psi + \\
e(Z_{1,\varepsilon})_\varepsilon \bar{\psi}\hat{A}\psi - \frac{1}{2\xi}(\partial_\mu A_\mu)^2 = \\
= \frac{1}{4}\left((F^{g.m.s.}_{\mu\nu B,\varepsilon})^2\right)_\varepsilon - i(\bar{\psi}^{g.m.s.}_{B,\varepsilon})_\varepsilon \hat{\partial}(\psi^{g.m.s.}_{B,\varepsilon})_\varepsilon - m(Z_\varepsilon Z^{-1}_{2,\varepsilon})_\varepsilon (\bar{\psi}^{g.m.s.}_{B,\varepsilon})_\varepsilon (\psi^{g.m.s.}_{B,\varepsilon})_\varepsilon + \\
e\left(Z_{1,\varepsilon} Z^{-1}_{2,\varepsilon} Z^{-1/2}_{3,\varepsilon}\right)_\varepsilon (\bar{\psi}^{g.m.s.}_{B,\varepsilon})_\varepsilon \left(\hat{A}^{g.m.s.}_{B,\varepsilon}\right)_\varepsilon (\psi^{g.m.s.}_{B,\varepsilon})_\varepsilon \; - \frac{(Z^{-1}_{3,\varepsilon})_\varepsilon}{2\xi} (\partial_\mu (A^{g.m.s.}_{\mu B,\varepsilon})_\varepsilon)^2 = \\
-\frac{1}{4}\left((F^{g.m.s}_{\mu\nu B,\varepsilon})^2\right)_\varepsilon + i(\bar{\psi}^{g.m.s}_{B,\varepsilon})_\varepsilon \hat{\partial}(\psi^{g.m.s}_{B,\varepsilon})_\varepsilon - (m^{s.m.s}_{B,\varepsilon})_\varepsilon (\bar{\psi}^{g.m.s.}_{B,\varepsilon})_\varepsilon (\psi^{g.m.s.}_{B,\varepsilon})_\varepsilon + \\
(e^{g.m.s.}_{B,\varepsilon})_\varepsilon (\bar{\psi}^{g.m.s.}_{B,\varepsilon})_\varepsilon \left(\hat{A}^{g.m.s.}_{B,\varepsilon}\right)_\varepsilon (\psi^{g.m.s.}_{B,\varepsilon})_\varepsilon \; - \frac{1}{2\left(\xi^{g.m.s.}_{B,\varepsilon}\right)_\varepsilon} (\partial_\mu (A^{g.m.s.}_{\mu B,\varepsilon})_\varepsilon)^2
\end{gathered} \qquad (4.4.23)$$

From (4.4.23) we finally obtain the $QED_4$ bare Lagrangian in a ghost matter sector

$$\begin{gathered}
\left(\mathcal{L}^{Bare,g.m.s.}_{QED,\varepsilon}\right)_\varepsilon = \\
\frac{1}{4}\left((F^{g.m.s.}_{\mu\nu B,\varepsilon})^2\right)_\varepsilon - i(\bar{\psi}^{g.m.s.}_{B,\varepsilon})_\varepsilon \hat{\partial}(\psi^{g.m.s.}_{B,\varepsilon})_\varepsilon - (m^{g.m.s.}_{B,\varepsilon})_\varepsilon (\bar{\psi}^{g.m.s.}_{B,\varepsilon})_\varepsilon (\psi^{g.m.s.}_{B,\varepsilon})_\varepsilon - \\
-(e^{g.m.s.}_{B,\varepsilon})_\varepsilon (\bar{\psi}^{g.m.s.}_{B,\varepsilon})_\varepsilon \left(\hat{A}^{g.m.s.}_{B,\varepsilon}\right)_\varepsilon (\psi^{g.m.s.}_{B,\varepsilon})_\varepsilon \; - \frac{1}{2\left(\xi^{g.m.s.}_{B,\varepsilon}\right)_\varepsilon} (\partial_\mu (A^{g.m.s.}_{\mu B,\varepsilon})_\varepsilon)^2,
\end{gathered} \qquad (4.4.24)$$

where

$$\begin{gathered}
(\psi^{g.m.s.}_{B,\varepsilon})_\varepsilon = \left((-Z_{2,\varepsilon})^{1/2}\right)_\varepsilon \psi, \;\; (A_{B,\varepsilon})_\varepsilon = \left((-Z_{3,\varepsilon})^{1/2}\right)_\varepsilon A, \\
(m_{B,\varepsilon})_\varepsilon = (Z_\varepsilon Z^{-1}_{2,\varepsilon})_\varepsilon m, \;\; (e_{B,\varepsilon})_\varepsilon = \left(Z_{1,\varepsilon} Z^{-1}_{2,\varepsilon} (-Z_{3,\varepsilon})^{-1/2}\right)_\varepsilon e, \\
\left(\xi_{B,\varepsilon}\right)_\varepsilon = (Z_{3,\varepsilon})_\varepsilon \xi.
\end{gathered} \qquad (4.4.25)$$

The gauge invariance connects the vertex Green function and the fermion propagator (the Ward identity), which leads to the identity $(Z_{1,\varepsilon})_\varepsilon = (Z_{2,\varepsilon})_\varepsilon$.

## 4.5.Quantum chromodynamics in a ghost sector via Colombeau generalized functions.

Consider now the non-Abelian gauge theories and, in particular, $QCD$. The Lagrangian of QCD has the form

$$\begin{gathered}
\mathcal{L}_{QÑD} = -\tfrac{1}{4}(F^a_{\mu\nu})^2 + \bar{\psi}(i\gamma^\mu \partial_\mu - m)\psi + g\bar{\psi}\gamma^\mu A^a_\mu T^a \psi - \frac{1}{2\xi}(\partial_\mu A^a_\mu)^2 \\
+\partial_\mu \bar{c}^a \partial_\mu c^2 + g f^{abc} \partial_\mu \bar{c}^a A^b_\mu c^c,
\end{gathered} \qquad (4.5.1)$$

where the stress tensor of the gauge field is now $F^a_{\mu\nu} = \partial_\mu A^a_\nu - \partial_\nu A^a_\mu + g f^{abc} A^b_\mu A^c_\nu$ and the last terms represent the Faddeev-Popov ghosts.

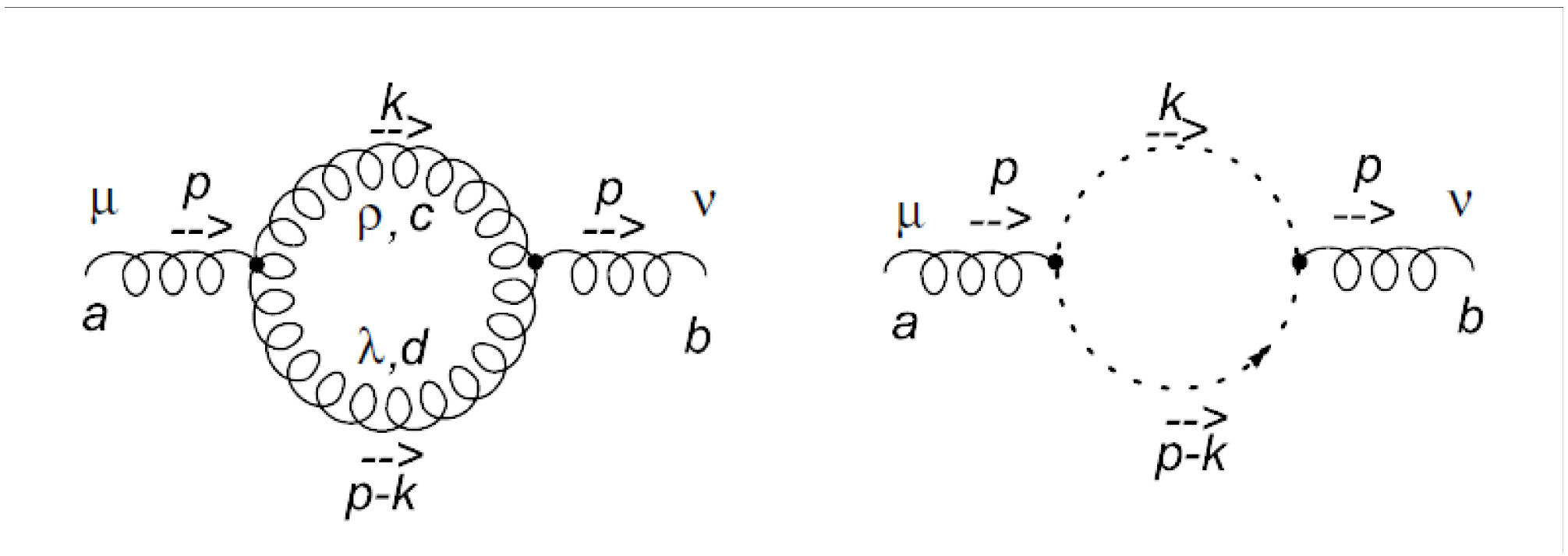

Fig.4.5.1.The vacuum polarization diagrams in the Yang-Mills theory

The complications which appear in non-Abelian theories are caused by the presence of many vertices with the same coupling as it follows from the gauge invariance. Hence, they have to renormalize the same way, i.e there appear new identities, called the Slavnov-Taylor identities. The full set of the counter-terms in $QCD$ are

$$\begin{aligned}\left(\Delta\mathcal{L}_\varepsilon^{QCD}\right)_\varepsilon = & -\frac{(Z_{3,\varepsilon})_\varepsilon - 1}{4}(\partial_\mu A_\nu^a - \partial_\nu A_\mu^a)^2 - g((Z_{1,\varepsilon})_\varepsilon - 1)f^{abc}A_\mu^a A_\nu^b \partial_\mu A_\nu^c - \\ & ((Z_{4,\varepsilon})_\varepsilon - 1)\frac{g^2}{4}f^{abc}f^{ade}A_\mu^b A_\nu^c A_\mu^d A_\nu^e + \\ & ((\tilde{Z}_{3,\varepsilon})_\varepsilon - 1)\partial_\mu \bar{c}^a \partial_\mu c^a + g((\tilde{Z}_{1,\varepsilon})_\varepsilon - 1)f^{abc}\partial_\mu \bar{c}^a A_\mu^b c^c \\ & +i((Z_{2,\varepsilon})_\varepsilon - 1)\bar{\psi}\hat{\partial}\psi - m((Z_\varepsilon)_\varepsilon - 1)\bar{\psi}\psi + g((Z_{1\psi,\varepsilon})_\varepsilon - 1)\bar{\psi}\hat{A}^a T^a \psi,\end{aligned} \tag{4.5.2}$$

**Remark 4.5.1**. We assume now that

$$\begin{aligned}&(Z_\varepsilon^{s.m.s.})_\varepsilon > 0, (Z_{1,\varepsilon}^{s.m.s.})_\varepsilon > 0, \left(\tilde{Z}_{1,\varepsilon}^{s.m.s.}\right)_\varepsilon > 0, (Z_{2,\varepsilon}^{s.m.s.})_\varepsilon > 0, \\ &(Z_{3,\varepsilon}^{s.m.s.})_\varepsilon > 0, \left(\tilde{Z}_{3,\varepsilon}^{s.m.s.}\right)_\varepsilon > 0, (Z_{4,\varepsilon}^{s.m.s.})_\varepsilon > 0.\end{aligned} \tag{4.5.3}$$

Where we abraviate $(Z_\varepsilon^{s.m.s.})_\varepsilon, (Z_{1,\varepsilon}^{s.m.s.})_\varepsilon$, etc., instead $(Z_\varepsilon)_\varepsilon, (Z_{1,\varepsilon})_\varepsilon$, etc., in standard matter sector.

Adding (4.5.2) to the initial Lagrangian (4.5.1) from (4.5.2) we obtain the $QCD_4$ bare Lagrangian in standard matter sector

$$
\begin{aligned}
&\left(\mathcal{L}^{Bare,s.m.s.}_{QCD,\varepsilon}\right)_\varepsilon = \mathcal{L}_{QCD} + \left(\Delta\mathcal{L}^{QCD}_\varepsilon\right)_\varepsilon = \\
&-\frac{(Z_{3,\varepsilon})_\varepsilon}{4}(\partial_\mu A^a_\nu - \partial_\nu A^a_\mu)^2 - g(Z_{1,\varepsilon})_\varepsilon f^{abc}A^a_\mu A^b_\nu \partial_\mu A^c_\nu \\
&-(Z_{4,\varepsilon})_\varepsilon \frac{g^2}{4} f^{abc}f^{ade}A^b_\mu A^c_\nu A^d_\mu A^e_\nu - (\tilde{Z}_{3,\varepsilon})_\varepsilon \partial_\mu \bar{c}^a \partial c^a - g(\tilde{Z}_{1,\varepsilon})_\varepsilon f^{abc}\partial_\mu \bar{c}^a A^b_\mu c^c \\
&+i(Z_{2,\varepsilon})_\varepsilon \bar{\psi}\hat{\partial}\psi - m(Z_\varepsilon)_\varepsilon \bar{\psi}\psi + g(Z_{1\psi,\varepsilon})_\varepsilon \bar{\psi}\hat{A}^a T^a \psi - \frac{1}{2\xi}(\partial_\mu A^a_\mu)^2 \\
&= -\frac{1}{4}(\partial_\mu (A^{a,s.m.s.}_{\nu B,\varepsilon})_\varepsilon - \partial_\nu (A^{a,s.m.s.}_{\mu B})_\varepsilon)^2 - \\
&g\left(Z_{1,\varepsilon}Z^{-3/2}_{3,\varepsilon}\right)_\varepsilon f^{abc}(A^{a,s.m.s.}_{\mu B,\varepsilon})_\varepsilon A^{b,s.m.s.}_{\nu B,\varepsilon}\partial_\mu A^{c,s.m.s.}_{\nu B,s.m.s.} \\
&-(Z_{4,\varepsilon}Z^{-2}_{3,\varepsilon})_\varepsilon \frac{g^2}{4} f^{abc}f^{ade}\left[\left(A^{b,s.m.s.}_{\mu B,\varepsilon}\right)_\varepsilon\right][(A^c_{\nu B,\varepsilon})_\varepsilon]\left[\left(A^{d,s.m.s.}_{\mu B,\varepsilon}\right)_\varepsilon\right][(A^{e,s.m.s.}_{\nu B,\varepsilon})_\varepsilon] + \\
&[(\partial_\mu \bar{c}^{a,s.m.s.}_{B,\varepsilon})_\varepsilon][(\partial_\mu c^{a,s.m.s.}_{B,\varepsilon})_\varepsilon] + \\
&g\left[\left(\tilde{Z}_{1,\varepsilon}\tilde{Z}^{-1}_{3,\varepsilon}Z^{-1/2}_{3,\varepsilon}\right)_\varepsilon\right] f^{abc}[(\partial_\mu \bar{c}^{a,s.m.s.}_{B,\varepsilon})_\varepsilon]\left[\left(A^{b,s.m.s.}_{\mu B,\varepsilon}\right)_\varepsilon\right][(c^{c,s.m.s.}_{B,\varepsilon})_\varepsilon] \\
&+\frac{(Z^{-1}_{3,\varepsilon})_\varepsilon}{2\xi}\left((\partial_\mu A^{a,s.m.s.}_{\mu B,\varepsilon})_\varepsilon\right)^2 + i[(\bar{\psi}^{s.m.s.}_{B,\varepsilon})_\varepsilon]\left[\left(\hat{\partial}\psi^{s.m.s.}_{B,\varepsilon}\right)_\varepsilon\right] - \\
&m(Z_\varepsilon Z^{-1}_{2,\varepsilon})_\varepsilon[(\bar{\psi}^{s.m.s.}_{B,\varepsilon})_\varepsilon][(\psi^{s.m.s.}_{B,\varepsilon})_\varepsilon] + \\
&g\left[\left(Z_{1\psi,\varepsilon}Z^{-1}_{2,\varepsilon}Z^{-1/2}_{3,\varepsilon}\right)_\varepsilon\right][(\bar{\psi}^{s.m.s.}_{B,\varepsilon})_\varepsilon]\left[\left(\hat{A}^{a,s.m.s.}_{B,\varepsilon}\right)_\varepsilon\right]T^a[(\psi^{s.m.s.}_{B,\varepsilon})_\varepsilon].
\end{aligned}
\tag{4.5.4}
$$

This results in the relations between the renormalized and the bare fields and couplings

$$
\begin{aligned}
&(\psi^{s.m.s.}_{B,\varepsilon})_\varepsilon = Z^{1/2}_2\psi,\ \left(A^{s.m.s.}_{B,\varepsilon}\right)_\varepsilon = Z^{1/2}_3 A,\ (c^{s.m.s.}_{B,\varepsilon})_\varepsilon = \tilde{Z}^{1/2}_3 c, \\
&(m^{s.m.s.}_{B,\varepsilon})_\varepsilon = ZZ^{-1}_2 m,\ (g^{s.m.s.}_{B,\varepsilon})_\varepsilon = Z_1 Z^{-3/2}_3 g,\ \left(\xi^{s.m.s.}_{B,\varepsilon}\right)_\varepsilon = Z_3\xi, \\
&Z_1Z^{-1}_3 = \tilde{Z}_1\tilde{Z}^{-1}_3,\ Z_4 = Z^2_1 Z^{-1}_3,\ Z_{1\psi}Z^{-1}_2 = Z_1Z^{-1}_3.
\end{aligned}
\tag{4.5.4}
$$

The last line of equalities follows from the requirement of identical renormalization of the coupling in various vertices and represents the Slavnov-Taylor identities for the singular parts.The explicit form of the renormalization constants in the lowest approximation follows from the one-loop diagrams of $QCD$. Aa usual, one has to take the singular part with the opposite sign. For instance one has in the $\overline{MS}$ scheme:

$$
\begin{aligned}
&(Z_{2,\varepsilon})_\varepsilon = 1 - \frac{g^2}{16\pi^2}\left(\frac{C_F}{\varepsilon}\right)_\varepsilon, (Z_{3,\varepsilon})_\varepsilon = 1 + \frac{g^2}{16\pi^2}\left(\frac{5}{3\varepsilon}C_A - \frac{4}{3\varepsilon}T_f n_f\right)_\varepsilon = \\
&1 + \frac{g^2}{16\pi^2}\left(\frac{1}{\varepsilon}\right)_\varepsilon\left(\frac{5}{3}C_A - \frac{4}{3}T_f n_f\right), \\
&(Z_\varepsilon)_\varepsilon = 1 - \frac{g^2}{16\pi^2}\left(\frac{4C_F}{\varepsilon}\right)_\varepsilon, (\tilde{Z}_{1,\varepsilon})_\varepsilon = 1 - \frac{g^2}{16\pi^2}\left(\frac{C_A}{2\varepsilon}\right)_\varepsilon, \\
&(\tilde{Z}_{2,\varepsilon})_\varepsilon = 1 + \frac{g^2}{16\pi^2}\left(\frac{C_A}{2\varepsilon}\right)_\varepsilon, (Z_{g,\varepsilon})_\varepsilon = \left(\tilde{Z}_{1,\varepsilon}\tilde{Z}^{-1}_{2,\varepsilon}Z^{-1/2}_{3,\varepsilon}\right)_\varepsilon = \\
&1 - \frac{g^2}{16\pi^2}\left(\frac{11}{6\varepsilon}C_A - \frac{4}{3\varepsilon}T_f n_f\right)_\varepsilon = 1 - \frac{g^2}{16\pi^2}\left(\frac{1}{\varepsilon}\right)_\varepsilon\left(\frac{11}{6}C_A - \frac{4}{3}T_f n_f\right).
\end{aligned}
\tag{4.5.5}
$$

where the following notation for the Casimir operators of the gauge group is used

$$
f^{abc}f^{dbc} = C_A\delta^{ad},\ (T^aT^a)_{ij} = C_F\delta_{ij},\ Tr(T^aT^b) = T_f\delta^{ab}. \tag{4.5.6}
$$

For the $SU(N)$ group and the fundamental representation of the fermion fields they are equal to

$$C_A = N,\ C_F = \frac{N^2-1}{2N},\ T_f = \frac{1}{2}. \tag{4.5.7}$$

**Remark 4.5.2**. We assume now that

$$\begin{gathered}(Z_\varepsilon^{g.m.s.})_\varepsilon < 0, (Z_{1,\varepsilon}^{g.m.s.})_\varepsilon < 0, (Z_{2,\varepsilon}^{g.m.s.})_\varepsilon < 0,\\ (Z_{3,\varepsilon}^{g.m.s.})_\varepsilon < 0, (Z_{4,\varepsilon}^{g.m.s.})_\varepsilon < 0.\end{gathered} \tag{4.5.8}$$

Where we abraviate $(Z_\varepsilon^{g.m.s.})_\varepsilon, (Z_{1,\varepsilon}^{g.m.s.})_\varepsilon$,etc., instead $(Z_\varepsilon)_\varepsilon, (Z_{1,\varepsilon})_\varepsilon$,etc., in ghost matter sector.

Adding (4.5.2) to the initial Lagrangian (4.5.1) from (4.5.5) we obtain the $QCD_4$ bare Lagrangian in a ghost matter sector

$$\begin{gathered}\left(\mathcal{L}_{QCD,\varepsilon}^{Bare,g.m.s.}\right)_\varepsilon = \mathcal{L}_{QCD} + \left(\Delta\mathcal{L}_\varepsilon^{QCD}\right)_\varepsilon =\\
-\frac{(Z_{3,\varepsilon}^{g.m.s.})_\varepsilon}{4}(\partial_\mu A_\nu^a - \partial_\nu A_\mu^a)^2 - g(Z_{1,\varepsilon}^{g.m.s.})_\varepsilon f^{abc} A_\mu^a A_\nu^b \partial_\mu A_\nu^c\\
-(Z_{4,\varepsilon})_\varepsilon \frac{g^2}{4} f^{abc} f^{ade} A_\mu^b A_\nu^c A_\mu^d A_\nu^e - (\tilde{Z}_{3,\varepsilon}^{g.m.s.})_\varepsilon \partial_\mu \bar{c}^a \partial c^a - g(\tilde{Z}_{1,\varepsilon}^{g.m.s.})_\varepsilon f^{abc} \partial_\mu \bar{c}^a A_\mu^b c^c\\
+i(Z_{2,\varepsilon}^{g.m.s.})_\varepsilon \bar{\psi}\hat{\partial}\psi - m(Z_\varepsilon^{g.m.s.})_\varepsilon \bar{\psi}\psi + g(Z_{1\psi,\varepsilon}^{g.m.s.})_\varepsilon \bar{\psi}\hat{A}^a T^a \psi - \frac{1}{2\xi}(\partial_\mu A_\mu^a)^2\\
= -\frac{1}{4}(\partial_\mu A_{\nu B,\varepsilon}^{a,g.m.s.} - \partial_\nu A_{\mu B}^{a,g.m.s.})^2 -\\
g\left[(Z_{1,\varepsilon}^{g.m.s.})_\varepsilon\right]\left[(Z_{3,\varepsilon}^{g.m.s.})_\varepsilon\right]^{-3/2} f^{abc}\left[(A_{\mu B,\varepsilon}^{a,g.m.s.})_\varepsilon\right]\left[\left(A_{\nu B,\varepsilon}^{b,g.m.s.}\right)_\varepsilon\right]\left[(\partial_\mu A_{\nu B,\varepsilon}^{c,g.m.s.})_\varepsilon\right]\\
-(Z_{4,\varepsilon} Z_{3,\varepsilon}^{-2})_\varepsilon \frac{g^2}{4} f^{abc} f^{ade}\left[\left(A_{\mu B,\varepsilon}^{b,g.m.s.}\right)_\varepsilon\right]\left[(A_{\nu B,\varepsilon}^{c,g.m.s.})_\varepsilon\right]\left[\left(A_{\mu B,\varepsilon}^{d,g.m.s.}\right)_\varepsilon\right]\left[(A_{\nu B,\varepsilon}^{e,g.m.s.})_\varepsilon\right] +\\
\left[(\partial_\mu \bar{c}_{B,\varepsilon}^{a,g.m.s.})_\varepsilon\right]\left[(\partial_\mu c_{B,\varepsilon}^{a,g.m.s.})_\varepsilon\right] +\\
g\left[\left(\tilde{Z}_{1,\varepsilon}\tilde{Z}_{3,\varepsilon}^{-1} Z_{3,\varepsilon}^{-1/2}\right)_\varepsilon\right] f^{abc}\left[(\partial_\mu \bar{c}_{B,\varepsilon}^{a,g.m.s.})_\varepsilon\right]\left[\left(A_{\mu B,\varepsilon}^{b,g.m.s.}\right)_\varepsilon\right]\left[(c_{B,\varepsilon}^{c,g.m.s.})_\varepsilon\right]\\
+\frac{\left[(Z_{3,\varepsilon}^{g.m.s.})_\varepsilon\right]^{-1}}{2\xi}\left((\partial_\mu A_{\mu B,\varepsilon}^{a,g.m.s.})_\varepsilon\right)^2 + i(\bar{\psi}_{B,\varepsilon}^{g.m.s.})_\varepsilon\left(\hat{\partial}\psi_{B,\varepsilon}^{g.m.s.}\right)_\varepsilon -\\
m(Z_\varepsilon Z_{2,\varepsilon}^{-1})_\varepsilon \bar{\psi}_{B,\varepsilon}^{g.m.s.} \psi_{B,\varepsilon}^{s.m.s.} + g\left[\left(Z_{1\psi,\varepsilon} Z_{2,\varepsilon}^{-1} Z_{3,\varepsilon}^{-1/2}\right)_\varepsilon\right]\bar{\psi}_{B,\varepsilon}^{g.m.s.} \hat{A}_{B,\varepsilon}^{a,g.m.s.} T^a \psi_{B,\varepsilon}^{g.m.s.}.\end{gathered} \tag{4.5.9}$$

The relations between the renormalized and the bare fields and couplings reads

$$\begin{gathered}(\psi_{B,\varepsilon}^{g.m.s.})_\varepsilon = \left[((-Z_{2,\varepsilon}^{g.m.s.}))_\varepsilon\right]^{1/2}\psi,\ (A_{B,\varepsilon}^{g.m.s.})_\varepsilon = \left[(-Z_{3,\varepsilon}^{g.m.s.})_\varepsilon\right]^{1/2} A,\\
\left(c_{B,\varepsilon}^{s.m.s.}\right)_\varepsilon = \left[(\tilde{Z}_{3,\varepsilon}^{g.m.s.})_\varepsilon\right]^{1/2} c, (m_{B,\varepsilon}^{g.m.s.})_\varepsilon = \left[(Z_\varepsilon^{g.m.s.})_\varepsilon\right]\left[(Z_{2,\varepsilon}^{g.m.s.})_\varepsilon\right]^{-1} m,\\
\left(g_{B,\varepsilon}^{g.m.s.}\right)_\varepsilon = Z_1\left[(-Z_{3,\varepsilon}^{g.m.s.})_\varepsilon\right]^{-3/2} g,\ (\xi_{B,\varepsilon}^{g.m.s.})_\varepsilon = \left[(Z_{3,\varepsilon}^{g.m.s.})_\varepsilon\right]\xi,\\
\left[(Z_{1,\varepsilon}^{g.m.s.})_\varepsilon\right]\left[(Z_{3,\varepsilon}^{g.m.s.})_\varepsilon\right]^{-1} = \tilde{Z}_1 \tilde{Z}_3^{-1},\ Z_4 = Z_1^2 Z_3^{-1},\ Z_{1\psi} Z_2^{-1} = Z_1 Z_3^{-1}.\end{gathered} \tag{4.5.10}$$

## 4.6.The general structure of the $\mathcal{R}$-operation via Colombeau generalized functions.

The structure of the counter-terms as functions of the field operators depends on the type of a theory. According to the canonical classification [8],[20],the QFT theories are divided into three classes: *superrenormalizable* (a finite number of divergent diagrams), *renormalizable* (a finite number of types of divergent diagrams) and *non-renormalizable* (a infinite number of types of divergent diagrams). Accordingly, in the first case one has a finite number of counter-terms; in the second case, a infinite number of counter-terms but they repeat the structure of the initial Lagrangian, and in the last case, one has an infinite number of structures with an increasing number of the fields and derivatives.

**Remark 5.6.1.** In the case of renormalizable and superrenormalizable theories, since the Colombeau counter-terms repeat the structure of the initial Lagrangian, the result of the introduction of counter-terms can be represented as

$$(\mathcal{L}_\varepsilon)_\varepsilon + (\Delta\mathcal{L}_\varepsilon)_\varepsilon = (\mathcal{L}_\varepsilon^{\mathbf{Bare}})_\varepsilon = \mathcal{L}((\phi_{B,\varepsilon})_\varepsilon, \{(g_{B,\varepsilon})_\varepsilon\}, \{(m_{B,\varepsilon})_\varepsilon\}), \qquad (4.6.1)$$

i.e., $(\mathcal{L}_\varepsilon^{Bare})_\varepsilon$ is the same Lagrangian $(\mathcal{L}_\varepsilon)_\varepsilon$ but with the fields, masses and coupling constants being the "bare" ones related to the renormalized quantities by the multiplicative equalities.

**Remark 5.6.2**. Note that: (i) for standard sector $(Z_{i,\varepsilon}(\{g\},(1/\varepsilon)_\varepsilon))_\varepsilon > 0$ and corresponding
equalities reads

$$\begin{aligned}
\left(\phi_{i,\varepsilon}^{Bare,s.m.s.}\right)_\varepsilon &= \left[(Z_{i,\varepsilon}(\{g_\varepsilon\},(1/\varepsilon))^{1/2})_\varepsilon\right]\phi,\\
\left(g_{i,\varepsilon}^{Bare,s.m.s.}\right)_\varepsilon &= \left[(Z_{g,\varepsilon}^{i,s.m.s.}(\{g_\varepsilon\},1/\varepsilon))_\varepsilon\right](g_{i,\varepsilon})_\varepsilon,\\
\left(m_{i,\varepsilon}^{Bare,s.m.s.}\right)_\varepsilon &= \left[(Z_{m,\varepsilon}^{i,s.m.s.}(\{g\},(1/\varepsilon)_\varepsilon))_\varepsilon\right]m_i,
\end{aligned} \qquad (4.6.2.a)$$

(ii) for a ghost sector $(Z_{i,\varepsilon}(\{g\},(1/\varepsilon)_\varepsilon))_\varepsilon < 0$ and corresponding equalities reads

$$\begin{aligned}
\left(\phi_{i,\varepsilon}^{Bare,g.m.s.}\right)_\varepsilon &= \left[\left(-Z_{i,\varepsilon}^{i,g.m.s.}(\{g_\varepsilon\},(1/\varepsilon))^{1/2}\right)_\varepsilon\right]\phi,\\
\left(g_{i,\varepsilon}^{Bare,g.m.s.}\right)_\varepsilon &= \left[\left(Z_{g,\varepsilon}^{i,g.m.s.}(\{g_\varepsilon\},1/\varepsilon)\right)_\varepsilon\right](g_{i,\varepsilon})_\varepsilon,\\
\left(m_{i,\varepsilon}^{Bare,g.m.s.}\right)_\varepsilon &= \left[\left(Z_{m,\varepsilon}^{i,g.m.s.}(\{g_\varepsilon\},(1/\varepsilon)_\varepsilon)\right)_\varepsilon\right]m_i,
\end{aligned} \qquad (4.6.2.b)$$

where the Colombeau renormalization constants $(Z_{i,\varepsilon})_\varepsilon$ depend on the renormalized parameters and the parameter of regularization, where for definiteness we have chosen $(1/\varepsilon)_\varepsilon \in \widetilde{\mathbb{R}}$. In some cases the renormalization can be nondiagonal and the renormalization constants become matrices.The renormalization constants are not unique and depend on the renormalization scheme. This arbitrariness, however, does not influence the observables expressed through the renormalized quantities. We will come back to this problem later when discussing the group of renormalization. In the gauge theories $(Z_{i,\varepsilon})_\varepsilon$ may depend on the choice of the gauge though in the minimal subtraction scheme the renormalizations of the masses and the couplings are gauge invariant.In the minimal schemes the renormalization constants do not depend on dimensional parameters like masses and do not depend on the arrangement of external momenta in the diagrams. This property allows one to simplify the calculation of the counter-terms putting the masses and some external momenta to zero, as it was exemplified above by calculation of the two-loop diagrams. In making this trick, however, one has to be careful not to create artificial infrared divergences. Since in dimensional regularization they also have the form of poles in $(\varepsilon)_\varepsilon \in \widetilde{\mathbb{R}}$, this may lead to the wrong answers.In renormalizable theory the finite Green function is obtained from the "bare" one, i.e., is calculated from the "bare" Lagrangian by multiplication on the corresponding Colombeau renormalization constant

$$(\Gamma_\varepsilon(\{p^2\},\mu^2,g_{\mu,\varepsilon}))_\varepsilon = [(Z_{\Gamma,\varepsilon}(1/\varepsilon,g_{\mu,\varepsilon}))_\varepsilon](\Gamma_\varepsilon^{\mathbf{Bare}}(\{p^2\},1/\varepsilon,g_\varepsilon^{\mathbf{Bare}}))_\varepsilon, \qquad (4.6.3)$$

where in the n-th order of perturbation theory the "bare" parameters in the RHS of the Eq.(4.6.3) have to be expressed in terms of the renormalized ones with the help of relations (4.6.2.a) or (4.6.2.b) taken in the $(n-1)$-th order. The remaining constant $(Z_{\Gamma,\varepsilon})_\varepsilon$ creates the counter-term of the n-th order of the form $(\Delta\mathcal{L}_\varepsilon)_\varepsilon = ((Z_{\Gamma,\varepsilon})_\varepsilon - 1)(O_{\Gamma,\varepsilon})_\varepsilon$, where the Colombeau generalized operator $(O_{\Gamma,\varepsilon})_\varepsilon$ reflects the corresponding Green function. If the Green function is finite by itself (for instance, has many legs), then one has to remove the divergences only in the subgraphs and the corresponding renormalization Colombeau constant $(Z_{\Gamma,\varepsilon})_\varepsilon = 1$.

Note that since the propagator is inverse to the operator quadratic in fields in the Lagrangian, the renormalization of the propagator is also inverse to the renormalization of the 1-particle irreducible two-point Green function:(i) for standard sector

$$D^{s.m.s.}(p^2,\mu^2,(g_{\mu,\varepsilon})_\varepsilon) = (D^{s.m.s.}_\varepsilon(p^2,\mu^2,g_{\mu,\varepsilon}))_\varepsilon = \left[\left(Z^{s.m.s.}_{2,\varepsilon}((1/\varepsilon)_\varepsilon,g_{\mu,\varepsilon})\right)_\varepsilon\right]^{-1} D^{\mathbf{Bare},s.m.s.}(p^2,(1/\varepsilon)_\varepsilon,(g^{\mathbf{Bare},s.m.s.}_\varepsilon)_\varepsilon). \quad (4.6.4.a)$$

and (ii) for a ghost sector

$$D^{g.m.s.}(p^2,\mu^2,(g_{\mu,\varepsilon})_\varepsilon) = (D^{g.m.s.}_\varepsilon(p^2,\mu^2,g_{\mu,\varepsilon}))_\varepsilon = \left[(Z^{g.m.s.}_{2,\varepsilon}((1/\varepsilon)_\varepsilon,g_{\mu,\varepsilon}))_\varepsilon\right]^{-1} D^{\mathbf{Bare},g.m.s.}(p^2,(1/\varepsilon)_\varepsilon,\left(g^{\mathbf{Bare},g.m.s.}_\varepsilon\right)_\varepsilon). \quad (4.6.4.b)$$

correspondingly. The propagator renormalization constant is also the renormalization constant of the corresponding field, but the fields themselves, contrary to the masses and couplings, do not enter into the expressions for observables.

We would like to stress once more that the $\mathcal{R}$-operation works independently on the fact renormalizable or non-renormalizable the theory is. In local theory the counter-terms are local anyway. But only in renormalizable theory the counter-terms are reduced to the multiplicative renormalization of the finite number of fields and parameters.

One can perform the $\mathcal{R}$-operation for each diagram separately. For this purpose one has first of all to subtract the divergences in the subgraphs and then subtract the divergence in the diagram itself which has to be local. This serves as a good test that the divergences in the subgraphs are subtracted correctly. In this case the $\mathcal{R}$-operation can be symbolically written in a factorized form

$$\mathcal{R}(G_\varepsilon)_\varepsilon = \prod\nolimits_{div.subgraphs}(1-(M_{\gamma,\varepsilon})_\varepsilon)(G_\varepsilon)_\varepsilon, \quad (4.6.5)$$

where $(G_\varepsilon)_\varepsilon$ is the initial diagram, $M_\varepsilon$ is the $\varepsilon$-subtraction operator (for instance, subtraction of the $\varepsilon$-singular part of the regularized diagram) and the product goes over all divergent subgraphs including the diagram itself. By a subgraph we mean here the 1-particle irreducible diagram consisting of the vertices and lines of the diagram which is UV divergent. The 1-particle irreducible is called the diagram which can not be made disconnected by deleting of one line.

We have demonstrated above the application of the $\mathcal{R}$-operation to the two–loop diagrams in a scalar theory. Consider some other examples of diagrams with larger number of loops shown in Fig.4.6.1. They appear in the $\phi^4_4$ theory in the three-loop approximation.

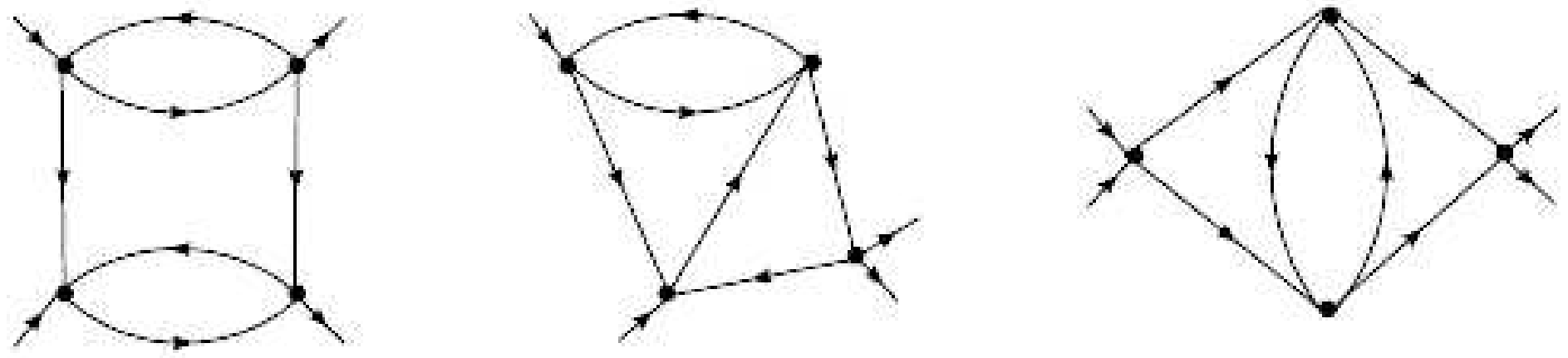

a b c

Fig.4.6.1.The multiloop diagrams in the φ[4] theory

In order to perform the $\mathcal{R}$-operation for these diagrams one first has to find out the divergent subgraphs. They are shown in Fig.4.6.2.

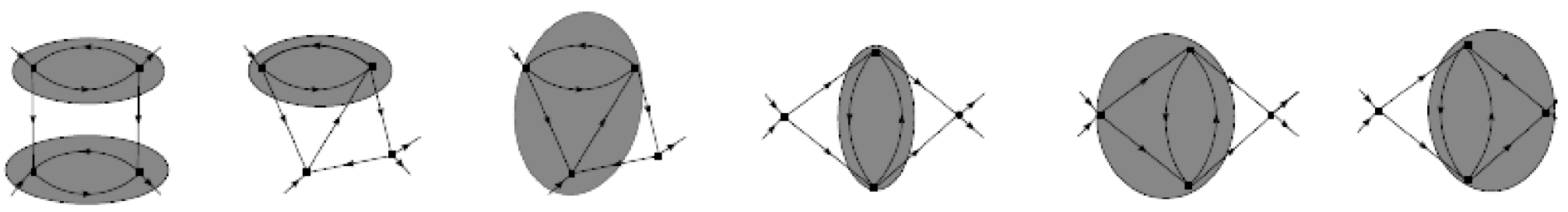

Fig.4.6.2.The divergent subgraphs in the diagrams of Fig.4.6.1.

Let us use the factorized representation of the $\mathcal{R}$-operation in the form of (4.6.5). For the three chosen diagrams one has, respectively,

$$(RG_{a,\varepsilon})_{\varepsilon} = (1-(M_{G_{\varepsilon}})_{\varepsilon})(1-M_{\gamma_1})(1-\left(M_{\gamma'_{1,\varepsilon}}\right)_{\varepsilon})(G_{a,\varepsilon})_{\varepsilon},$$

$$(RG_{b,\varepsilon})_{\varepsilon} = (1-(M_{G_{\varepsilon}})_{\varepsilon})(1-M_{\gamma_2})(1-(M_{\gamma_{1,\varepsilon}})_{\varepsilon})(G_{\acute{a},\varepsilon})_{\varepsilon}, \qquad (4.6.6)$$

$$(RG_{c,\varepsilon})_{\varepsilon} = (1-(M_{G_{\varepsilon}})_{\varepsilon})(1-(M_{\gamma_{2,\varepsilon}})_{\varepsilon})(1-\left(M_{\gamma'_{2,\varepsilon}}\right)_{\varepsilon})(1-(M_{\gamma_{1,\varepsilon}})_{\varepsilon})(G_{\hat{a},\varepsilon})_{\varepsilon},$$

where $\gamma_1$ and $\gamma_2$ are the one- and two-loop divergent subgraphs shown in Fig.4.6.2.The result of the application of the $\mathcal{R}$-operation without the last subtraction ( $\mathcal{R}'$-operation) for the diagrams of interest graphically is as follows:

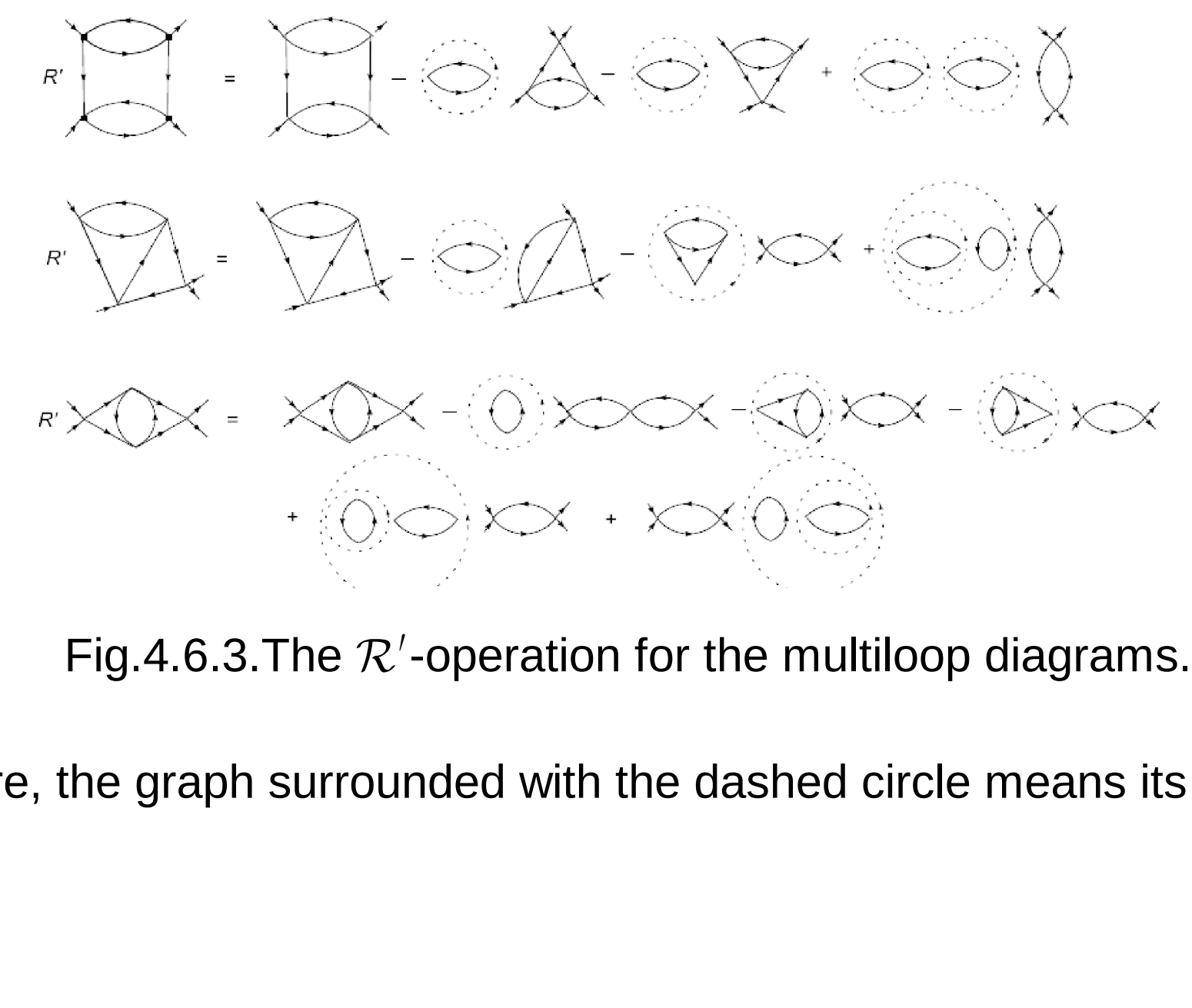

Fig.4.6.3.The $\mathcal{R}'$-operation for the multiloop diagrams.

Here, as before, the graph surrounded with the dashed circle means its singular part and the remaining graph is obtained by shrinking the singular subgraph to a point.

Let us demonstrate how the $\mathcal{R}'$-operation works for the diagram Fig.4.6.1a). Since the result of the $\mathcal{R}'$-operation does not depend on external momenta, we put two momenta on the diagonal to be equal to zero so that the integral takes the propagator form. Then we can use the method based on Fourier-transform, as it was explained above. One has

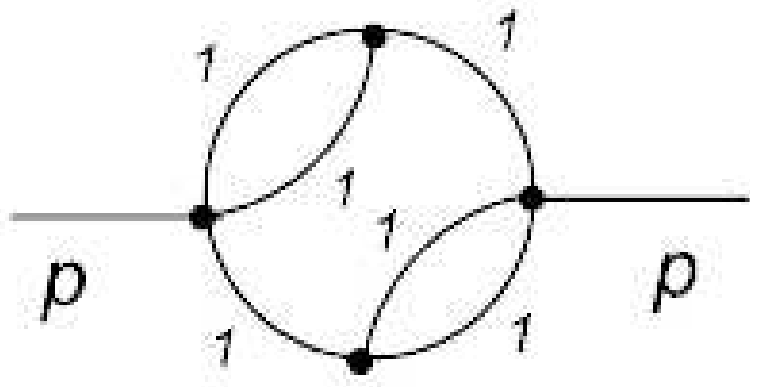

$$= \left[ \left( \Gamma(1-\varepsilon) \frac{\Gamma^2(1-\varepsilon)\Gamma(\varepsilon)}{\Gamma(2-2\varepsilon)} \right)_\varepsilon \right]^2$$

and

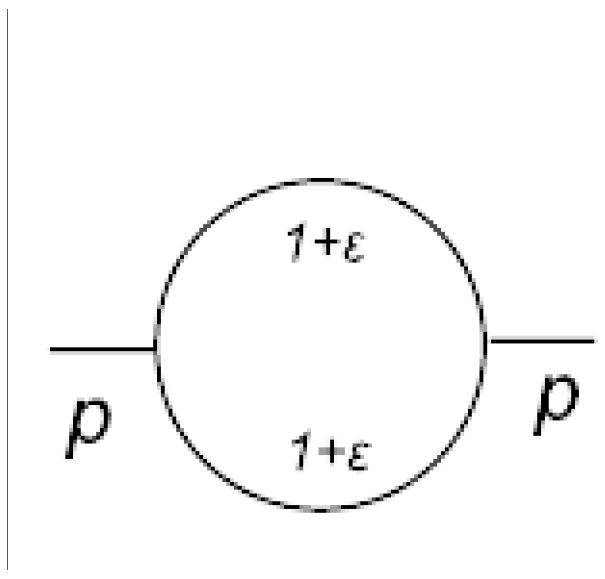

$$= \left[ \left( \Gamma(1-\varepsilon) \frac{\Gamma^2(1-\varepsilon)\Gamma(\varepsilon)}{\Gamma(2-2\varepsilon)} \right)_\varepsilon \right]^2 \left[ \left( \Gamma(1-\varepsilon) \frac{\Gamma^2(1-2\varepsilon)\Gamma(3\varepsilon)}{\Gamma^2(1+\varepsilon)\Gamma(2-4\varepsilon)} \right)_\varepsilon \right] \left( \left( \frac{\mu^2}{p^2} \right)^{3\varepsilon} \right)_\varepsilon \cong$$

$$\frac{1}{(\varepsilon^3(1-2\varepsilon)^2(1-4\varepsilon))_\varepsilon} \left( \left( \frac{\mu^2}{p^2} \right)^{3\varepsilon} \right)_\varepsilon .$$

Where we use the angular integration measure in the $4-2\varepsilon$ dimensional space accepted above, which results in the multiplication of the standard expression by $\Gamma(1-\varepsilon)$ in order to avoid the unwanted transcendental functions. Following the scheme shown in Fig.4.6.3 we get

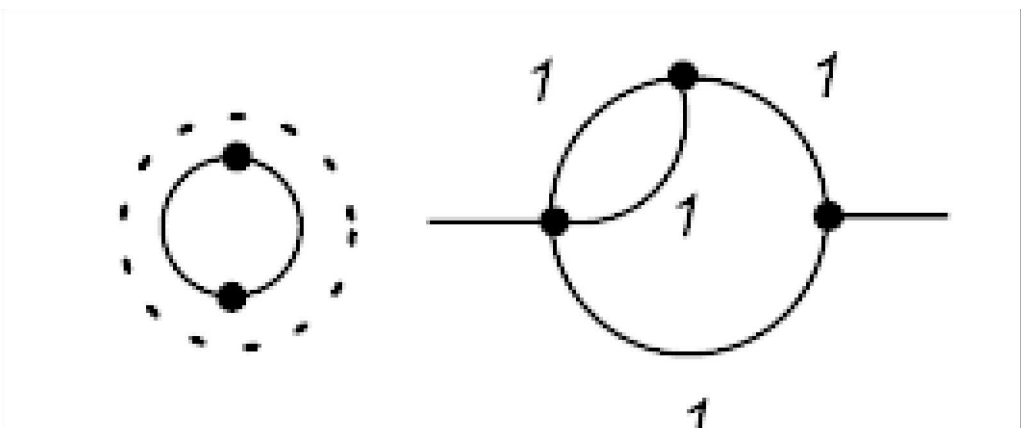

$$= \frac{1}{\varepsilon} \Gamma(1-\varepsilon) \frac{\Gamma^2(1-\varepsilon)\Gamma(\varepsilon)}{\Gamma(2-2\varepsilon)}$$

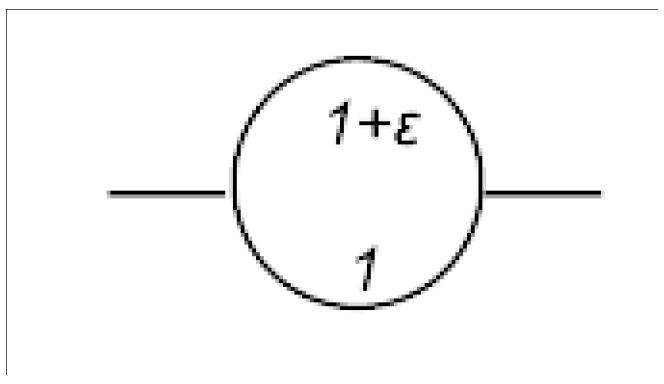

$$= \frac{1}{\varepsilon}\Gamma(1-\varepsilon)\frac{\Gamma^2(1-\varepsilon)\Gamma(\varepsilon)}{\Gamma(2-2\varepsilon)}\Gamma(1-\varepsilon)\frac{\Gamma(1-\varepsilon)\Gamma(1-2\varepsilon)\Gamma(2\varepsilon)}{\Gamma(1+\varepsilon)\Gamma(2-3\varepsilon)}\left(\frac{\mu^2}{p^2}\right)^{2\varepsilon} \cong \frac{1}{\varepsilon^3(1-2\varepsilon)(1-3\varepsilon)}\left(\frac{\mu^2}{p^2}\right)^{2\varepsilon}.$$

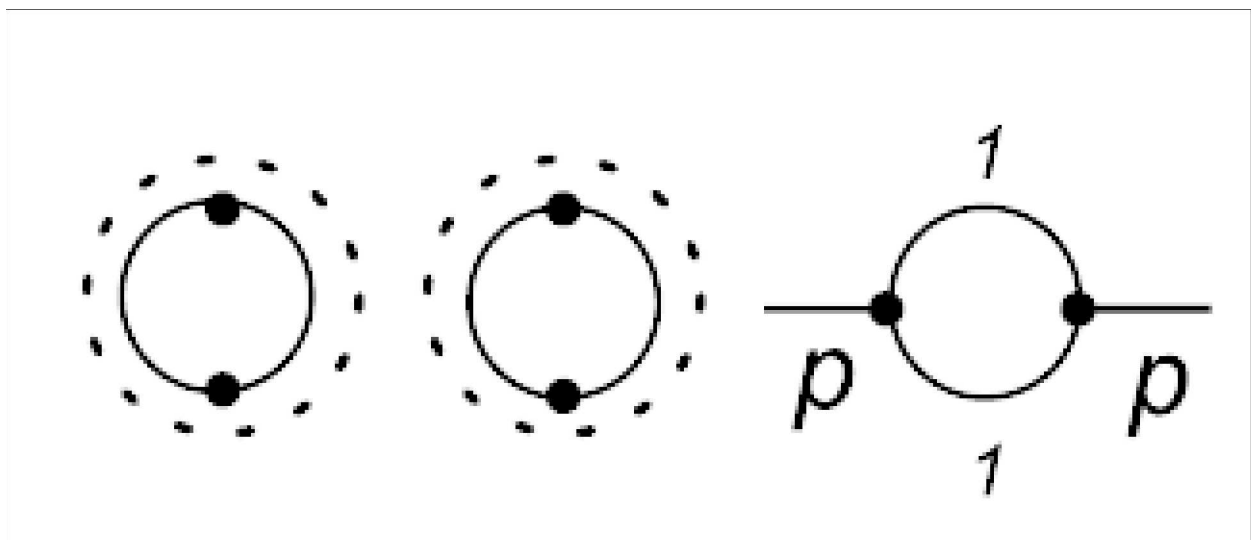

$$= \frac{1}{\varepsilon^2}\Gamma(1-\varepsilon)\frac{\Gamma^2(1-\varepsilon)\Gamma(\varepsilon)}{\Gamma(2-2\varepsilon)}\left(\frac{\mu^2}{p^2}\right)^{\varepsilon} \cong \frac{1}{\varepsilon^3(1-2\varepsilon)}\left(\frac{\mu^2}{p^2}\right)^{\varepsilon}.$$

Combining all together we get

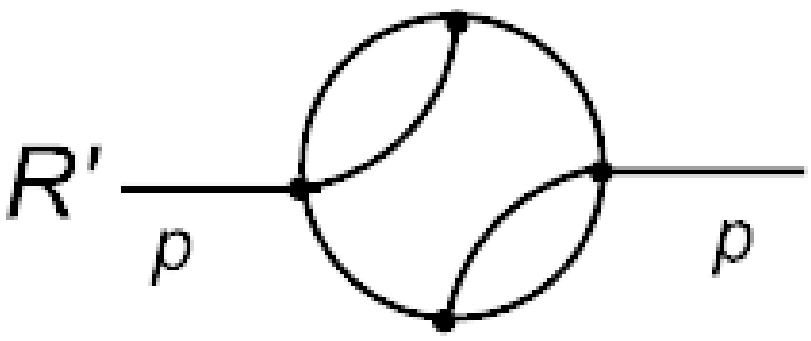

$$\cong \frac{1}{\varepsilon^3(1-2\varepsilon)^2(1-4\varepsilon)}\left(\frac{\mu^2}{p^2}\right)^{3\varepsilon} - 2\frac{1}{\varepsilon^3(1-2\varepsilon)}\left(\frac{\mu^2}{p^2}\right)^{\varepsilon} + \frac{1}{\varepsilon^3(1-2\varepsilon)}\left(\frac{\mu^2}{p^2}\right)^{\varepsilon}$$
$$= \frac{1-\varepsilon-\varepsilon^2}{\varepsilon^3}.$$

Note the cancellation of all nonlocal contributions. The singular part after the

$R'$-operation is always local.

The realization of the $\mathcal{R}'$-operation for each diagram $(G_\varepsilon)_\varepsilon$ allows one to find the contribution of a given diagram to the corresponding counter-term and, in the case of a renormalizable theory, to find the renormalization constant equal to

$$(Z_\varepsilon)_\varepsilon = 1 - \left(\mathcal{K}_\varepsilon\, \mathcal{R}' G_\varepsilon\right)_\varepsilon, \qquad (4.6.7)$$

where $\mathcal{K}_\varepsilon$ means the $\varepsilon$-extraction of the $\varepsilon$-singular part. Adding the contribution of various diagrams we get the resulting counter-term of a given order and, accordingly, the renormalization constant.

## 4.7.Renormalization Group in a ghost sector.

The procedure formulated above allows one to eliminate the ultraviolet divergences and get the finite expression for any Green function in any local quantum field theory. In renormalizable theories this procedure is reduced to the multiplicative renormalization of parameters (masses and couplings) and multiplication of the Green function by its own renormalization constant. This is true for any regularization and subtraction scheme. Thus, for example, in the canonical cutoff regularization and dimensional regularization via Colombeau generalized functions the relation between the "bare" and renormalized Green functions for standard matter sector looks like

$$\Gamma^{s.m.s}(\{p^2\},\mu^2,\{g_\mu\}) \;=Z_{\Gamma^{s.m.s}}(\Lambda^2/\mu^2,\{g_\mu\})\Gamma^{s.m.s}_{Bare}(\{p^2\},\Lambda,\{g^{s.m.s}_{Bare}\}) \qquad (4.7.1)$$

and

$$\Gamma^{s.m.s}(\{p^2\},\mu^2,\{(g_{\mu,\varepsilon})_\varepsilon\}) \;=Z_{\Gamma^{s.m.s}}((1/\varepsilon)_\varepsilon,\{(g_{\mu,\varepsilon})_\varepsilon\})\Gamma^{s.m.s}_{Bare}(\{p^2\},(1/\varepsilon)_\varepsilon,\{(g^{s.m.s}_{Bare,\varepsilon})_\varepsilon\}), \quad (4.7.2)$$

correspondingly and for a ghost sector looks like

$$\Gamma^{g.m.s}(\{p^2\},\mu^2,\{g_\mu\}) \;=Z_{\Gamma^{s.m.s}}(\Lambda^2/\mu^2,\{g_\mu\})\Gamma^{s.m.s}_{Bare}(\{p^2\},\Lambda,\{g^{s.m.s}_{Bare}\}) \qquad (4.7.3)$$

and

$$\Gamma^{g.m.s}(\{p^2\},\mu^2,\{(g^{s.m.s}_{\mu,\varepsilon})_\varepsilon\}) \;=Z_{\Gamma^{g.m.s}}((1/\varepsilon)_\varepsilon,\{(g^{s.m.s}_{\mu,\varepsilon})_\varepsilon\})\Gamma^{g.m.s}_{Bare}(\{p^2\},(1/\varepsilon)_\varepsilon,\{(g^{g.m.s}_{Bare,\varepsilon})_\varepsilon\}), \quad (4.7.4)$$

where $\{p^2\}$ is the set of external momenta, $\{g\}$ is the set of masses and couplings for standard matter sector

$$g^{s.m.s}_{Bare} = Z^{s.m.s}_g((\Lambda^2/\mu^2,\{g^{s.m.s}_\mu\})g, \qquad (4.7.5)$$

$$g^{s.m.s}_{Bare} = Z^{s.m.s}_g(((1/\varepsilon)_\varepsilon,\{(g^{s.m.s}_{\mu,\varepsilon})_\varepsilon\})g. \qquad (4.7.6)$$

and for a ghost sector

$$g^{g.m.s}_{Bare} = Z^{g.m.s.}_g((\Lambda^2/\mu^2,\{g^{g.m.s}_\mu\})g, \qquad (4.7.7)$$

$$g^{g.m.s}_{Bare} = Z^{g.m.s.}_g(((1/\varepsilon)_\varepsilon,\{(g^{g.m.s}_{\mu,\varepsilon})_\varepsilon\})g. \qquad (4.7.8)$$

correspondingly.In what follows we stick to dimensional regularization and rewrite relation (4.7.2) and (4.7.4) in the commun form

$$\Gamma^*_{Bare}(\{p^2\},(1/\varepsilon)_\varepsilon,\{g^*_{Bare}\}) \;=Z^{*-1}_{\Gamma^*}((1/\varepsilon)_\varepsilon,\{g_\mu\})\Gamma^*(\{p^2\},\mu^2,\{g_\mu\}) \qquad (4.7.9)$$

and

$$g^*_{Bare} = Z^*_g(((1/\varepsilon)_\varepsilon,\{(g^*_{\mu,\varepsilon})_\varepsilon\})(g^*_\varepsilon)_\varepsilon, \qquad (4.7.10)$$

where $\Gamma^*,\Gamma^*_{Bare}$ and $g^*_{Bare}$ stand for standard sector /ghost sector.

It is obvious that the LHS of this equation does not depend on the parameter of dimensional transmutation $\mu$ and, hence, the r.h.s. should not also depend on it. This allows us to write the functional equation for the renormalized Green function. Differentiating it with respect to the continuous parameter $\mu$ one can get the differential equation which has a practical value: solving this equation one can get the improved expression for the Green function which corresponds to summation of an infinite series of Feynman diagrams.

Consider an arbitrary Green function $(\Gamma^*_\varepsilon(\{p^2\},\mu^2,\{g^*_{\mu,\varepsilon}\}))_\varepsilon$ obeying equation (4.7.4) with the normalization condition

$$(\Gamma^*_\varepsilon(\{p^2\},\mu^2,0))_\varepsilon = 1. \qquad (4.7.11)$$

Differentiating the Eq.(4.7.4) with respect to $\mu^2$ we get:

$$\begin{aligned}\mu^2\frac{d}{d\mu^2}(\Gamma_\varepsilon^*)_\varepsilon &= \mu^2\frac{\partial}{\partial\mu^2}(\Gamma_\varepsilon^*)_\varepsilon + \mu^2\left(\frac{\partial g_\varepsilon^*}{\partial\mu^2}\right)_\varepsilon\left(\frac{\partial}{\partial g_\varepsilon^*}\Gamma_\varepsilon^*\right)_\varepsilon \\ &= \mu^2\frac{d\ln(Z_{\Gamma_\varepsilon^*})_\varepsilon}{d\mu^2}(Z_{\Gamma_\varepsilon^*}^*)_\varepsilon(\Gamma_{\varepsilon,Bare}^*)_\varepsilon,\end{aligned} \tag{4.7.12}$$

or

$$\begin{aligned}&\mu^2\frac{\partial}{\partial\mu^2}(\Gamma_\varepsilon^*(\{p^2\},\mu^2,g_{\mu,\varepsilon}))_\varepsilon + (\beta^*(g_\varepsilon))_\varepsilon\left(\frac{\partial}{\partial g_\varepsilon^*}\Gamma_\varepsilon^*(\{p^2\},\mu^2,g_{\mu,\varepsilon}^*)\right)_\varepsilon \\ &+(\gamma_{\Gamma_\varepsilon^*}^*\Gamma_\varepsilon^*(\{p^2\},\mu^2,g_{\mu,\varepsilon}^*))_\varepsilon = 0,\end{aligned} \tag{4.7.13}$$

where we have introduced the so-called beta function $(\beta(g_\varepsilon^*))_\varepsilon$ and the anomaly dimension of the Green function $(\gamma_{\Gamma_\varepsilon^*}^*(g_\varepsilon^*))_\varepsilon$ defined as

$$(\beta^*(g_\varepsilon^*))_\varepsilon = \mu^2\left(\frac{dg_\varepsilon^*}{d\mu^2}\right)_\varepsilon\Bigg|_{(g_{\varepsilon,bare}^*)_\varepsilon}, \tag{4.7.14}$$

and

$$(\gamma_{\Gamma_\varepsilon^*}^*(g_\varepsilon^*))_\varepsilon = -\mu^2\left(\frac{d\ln Z_{\Gamma_\varepsilon^*}^*}{d\mu^2}\right)_\varepsilon\Bigg|_{(g_{\varepsilon,bare}^*)_\varepsilon}. \tag{4.7.15}$$

Equation (4.7.13) is called the renormalization group equation in Colombeau generalized functions.The solution of the renormalization group equation can be written in terms of characteristics:

$$\left(\Gamma_\varepsilon^*\left(e^t\frac{\{p^2\}}{\mu^2},g_\varepsilon^*\right)\right)_\varepsilon = \left[\left(\Gamma_\varepsilon^*\left(\frac{\{p^2\}}{\mu^2},\bar{g}(t,g_\varepsilon^*)\right)\right)_\varepsilon\right]\exp\left(\int_0^t\left[(\gamma_{\Gamma_\varepsilon^*}^*(\bar{g}(t,g_\varepsilon^*)))_\varepsilon\right]dt\right), \tag{4.7.16}$$

where the characteristic equation is (for definiteness we restrict ourselves to a single coupling)

$$\frac{d}{dt}(\bar{g}(t,g_\varepsilon^*))_\varepsilon = (\beta^*(\bar{g}_\varepsilon^*))_\varepsilon,\ \ (\bar{g}_\varepsilon^*(0,g_\varepsilon^*))_\varepsilon = (g_\varepsilon^*)_\varepsilon. \tag{4.7.17}$$

Consider now the product

$$[(g_\varepsilon^*)_\varepsilon]\left(\Gamma_\varepsilon^*\left(\frac{\{p^2\}}{\mu^2},g_\varepsilon^*\right)\right)_\varepsilon. \tag{4.7.18}$$

If $(\Gamma_\varepsilon^*)_\varepsilon$ is the $n$-point function, then the renormalization of the coupling $(g_\varepsilon^*)_\varepsilon$ is given by

$$(g_{Bare,\varepsilon}^*)_\varepsilon = (Z_{\Gamma_\varepsilon^*}^*)_\varepsilon\left(Z_{2,\varepsilon}^{*-n/2}\right)_\varepsilon(g_\varepsilon^*)_\varepsilon, \tag{4.7.19}$$

and the product (4.7.18) is renormalized as

$$(g_\varepsilon^*)_\varepsilon(\Gamma_\varepsilon^*)_\varepsilon = \left(Z_{2,\varepsilon}^{*n/2}\right)_\varepsilon(g_{Bare,\varepsilon}^*)_\varepsilon(\Gamma_{Bare,\varepsilon}^*)_\varepsilon. \tag{4.7.20}$$

Hence, one has the same equation as (4.7.9) with solution (4.7.16) but with $(Z_{\Gamma_\varepsilon^*}^*)_\varepsilon = Z_{2,\varepsilon}^{*n/2}$ and $(\gamma_{\Gamma_\varepsilon^*}^*)_\varepsilon = -n/2(\gamma_{2,\varepsilon}^*)_\varepsilon$. (Recall that the anomalous dimension $(\gamma_{2,\varepsilon}^*)_\varepsilon$ is defined with respect to the Colombeau renormalization constant $(Z_{2,\varepsilon}^{*-1})_\varepsilon$.)

Furthermore, one can construct the so-called *invariant charge* by multiplying the product (4.7.18) by the corresponding propagators

$$(\xi^*_\varepsilon)_\varepsilon = (g^*_\varepsilon)_\varepsilon\left(\Gamma^*_\varepsilon\left(\frac{\{p^2\}}{\mu^2},g^*_\varepsilon\right)\right)_\varepsilon\left(\prod_{i=1}^n D^{*1/2}_\varepsilon\left(\frac{p_i^2}{\mu^2},g^*_\varepsilon\right)\right)_\varepsilon. \qquad (4.7.21)$$

The invariant charge $(\xi^*_\varepsilon)_\varepsilon$, being RG-invariant, obeys the RG equation without the anomalous dimension and plays an important role in the formulation of the renormalization group together with the effective charge. In some cases, for instance in the MOM subtraction scheme, the effective and invariant charges coincide.

The usefulness of solution (4.7.16) is that it allows one to sum up an infinite series of logs coming from the Feynman diagrams in the infrared ($t \to -\infty$) or ultraviolet ($t \to \infty$) regime and improve the usual perturbation theory expansions. This in its turn extends the applicability of perturbation theory and allows one to study the infrared or the ultraviolet asymptotics of the Green functions.

To demonstrate the power of the RG, let us consider the invariant charge in a theory with a single coupling and restrict ourselves to the massless case. Let the perturbative expansion be

$$\left(\xi^*_\varepsilon\left(\frac{p^2}{\mu^2},g^*_\varepsilon\right)\right)_\varepsilon = (g^*_\varepsilon)_\varepsilon(1 + b(g^*_\varepsilon)_\varepsilon \ln\frac{p^2}{\mu^2} + \ldots). \qquad (4.7.22)$$

The $\beta^*$ function in the one-loop approximation is given by

$$(\beta^*(g^*_\varepsilon))_\varepsilon = b(g^{*2}_\varepsilon)_\varepsilon. \qquad (4.7.23)$$

Notice that the coefficient $b$ of the logarithm in Eq.(4.7.22) coincides with that of the $\beta^*$ function. Alternatively the $\beta^*$ function can be defined as the derivative of the invariant charge with respect to logarithm of momentum

$$(\beta^*(g^*_\varepsilon))_\varepsilon = p^2\left(\frac{d}{dp^2}\xi^*_\varepsilon\left(\frac{p^2}{\mu^2},g^*_\varepsilon\right)\Bigg|_{p^2=\mu^2}\right)_\varepsilon \qquad (4.7.24)$$

This definition is useful in the MOM scheme where the mass is not considered as a coupling but as a parameter and the renormalization constants depend on it. We will come back to the discussion of this question below when considering different definitions of the mass. According to Eq.(4.7.16) (with vanishing anomalous dimension) the RG-improved expression for the invariant charge corresponding to the perturbative expression (4.7.22) is:

$$\left(\xi^*_{RG,\varepsilon}\left(\frac{p^2}{\mu^2},g^*_\varepsilon\right)\right)_\varepsilon = \left(\xi_{PT,\varepsilon}\left(1,\bar{g}^*_\varepsilon\left(\frac{p^2}{\mu^2},g^*_\varepsilon\right)\right)\right)_\varepsilon = \left(\bar{g}^*_\varepsilon\left(\frac{p^2}{\mu^2},g^*_\varepsilon\right)\right)_\varepsilon, \qquad (4.7.25)$$

where we have put in eq.(4.7.16) $p^2 = \mu^2$ and then replaced $t$ by $t = \ln p^2/\mu^2$. The effective coupling is a solution of the characteristic equation

$$\left(\frac{d}{dt}\bar{g}^*_\varepsilon(t,g^*_\varepsilon)\right)_\varepsilon = b(\bar{g}^{*2}_\varepsilon)_\varepsilon,\ \left(\bar{g}^*_\varepsilon(0,g^*_\varepsilon)\right)_\varepsilon = (g^*_\varepsilon)_\varepsilon,\ t = \ln\frac{p^2}{\mu^2}. \qquad (4.7.26)$$

The solution of this equation is

$$(\bar{g}^*_\varepsilon(t,g^*_\varepsilon))_\varepsilon = \frac{(g^*_\varepsilon)_\varepsilon}{1 - bt(g^*_\varepsilon)_\varepsilon}. \qquad (4.7.27)$$

Being expanded over $t$, the geometrical progression (4.7.27) reproduces the expansion (4.7.22); however, it sums the infinite series of terms of the form $[(g^{*n}_\varepsilon)_\varepsilon]t^n$. This is called the leading log approximation (LLA) in QFT. To get the correction to the LLA, one has to consider the next term in the expansion of the $\beta^*$ function. Then one can sum up the next series of terms of the form $[(g_\varepsilon^{*n})_\varepsilon]t^{n-1}$ which is called the next to leading log approximation (NLLA), etc. This procedure allows one to describe the leading asymptotics of the Green functions for $t \to \pm\infty$.Let us consider now the Green function with non-zero anomalous dimension. Let its perturbative expansion be

$$\left(\Gamma_\varepsilon^*\left(\frac{p^2}{\mu^2}, g_\varepsilon^*\right)\right)_\varepsilon = 1 + [(g_\varepsilon^*)_\varepsilon]c\ln\frac{p^2}{\mu^2} + \ldots \tag{4.7.28}$$

Then in the one-loop approximation the anomalous dimension is

$$(\gamma^*(g_\varepsilon^*))_\varepsilon = c(g_\varepsilon^*)_\varepsilon. \tag{4.7.29}$$

Again the coefficient of the logarithm coincides with that of the anomalous dimension. In analogy with Eq.(4.7.24) the anomalous dimension can be defined as a derivative with respect to the logarithm of momentum

$$(\gamma^*(g_\varepsilon^*))_\varepsilon = p^2\left(\frac{d}{dp^2}\ln\Gamma_\varepsilon^*\left(\frac{p^2}{\mu^2}, g_\varepsilon^*\right)\bigg|_{p^2=\mu^2}\right)_\varepsilon. \tag{4.7.30}$$

Substituting (4.7.29) into Eq.(4.7.16), one has in the exponent

$$\begin{gathered}\left(\int_0^t \gamma^*(\bar{g}_\varepsilon^*(t, g_\varepsilon^*)dt\right)_\varepsilon = \\ \left(\int_{g_\varepsilon^*}^{\bar{g}_\varepsilon^*}\frac{\gamma^*(g_\varepsilon^*)}{\beta^*(g_\varepsilon^*)}dg_\varepsilon^*\right)_\varepsilon = \frac{c}{b}\left(\int_{g_\varepsilon^*}^{\bar{g}_\varepsilon^*}\frac{g_\varepsilon^*}{g_\varepsilon^{*2}}dg_\varepsilon^*\right)_\varepsilon = \frac{c}{b}\left(\ln\frac{\bar{g}_\varepsilon^*}{g_\varepsilon^*}\right)_\varepsilon.\end{gathered} \tag{4.7.31}$$

This gives for the Green function the improved expression

$$(\Gamma_{RG,\varepsilon}^*)_\varepsilon = \left[\left(\frac{\bar{g}_\varepsilon^*}{g_\varepsilon^*}\right)_\varepsilon\right]^{-c/b} = \left[\frac{1}{1 - bt(g_\varepsilon^*)_\varepsilon}\right]^{c/b} \approx 1 + ct + \ldots \tag{5.7.32}$$

Thus, one again reproduces the perturbative expansion, but expression (4.7.32) again contains the whole infinite sum of the leading logs. To get the NLLA, one has to take into account the next term in eq.(4.7.29) together with the next term of expansion of the $\beta^*$ function. All the formulas can be easily generalized to the case of multiple couplings and masses.

## The effective coupling in a ghost sector

By virtue of the central role played by the effective coupling in RG formulas, consider it in more detail. The behaviour of the effective coupling is determined by the $\beta^*$ function. Qualitatively, the $\beta^*$ function can exhibit the behaviour shown in Fig.4.7.1. We restrict ourselves to the region of small couplings.In the first case, the $\beta^*$-function is positive. Hence, with increasing momentum the effective coupling unboundedly increases. This situation is typical of most of the models of QFT in standard matter sector in the one-loop approximation when $(\beta^*(g_\varepsilon^*))_\varepsilon = b(g_\varepsilon^{*2})_\varepsilon$ and $b > 0$. The solution of the RG equation for the effective coupling in this case has the form of a geometric progression (4.7.27).

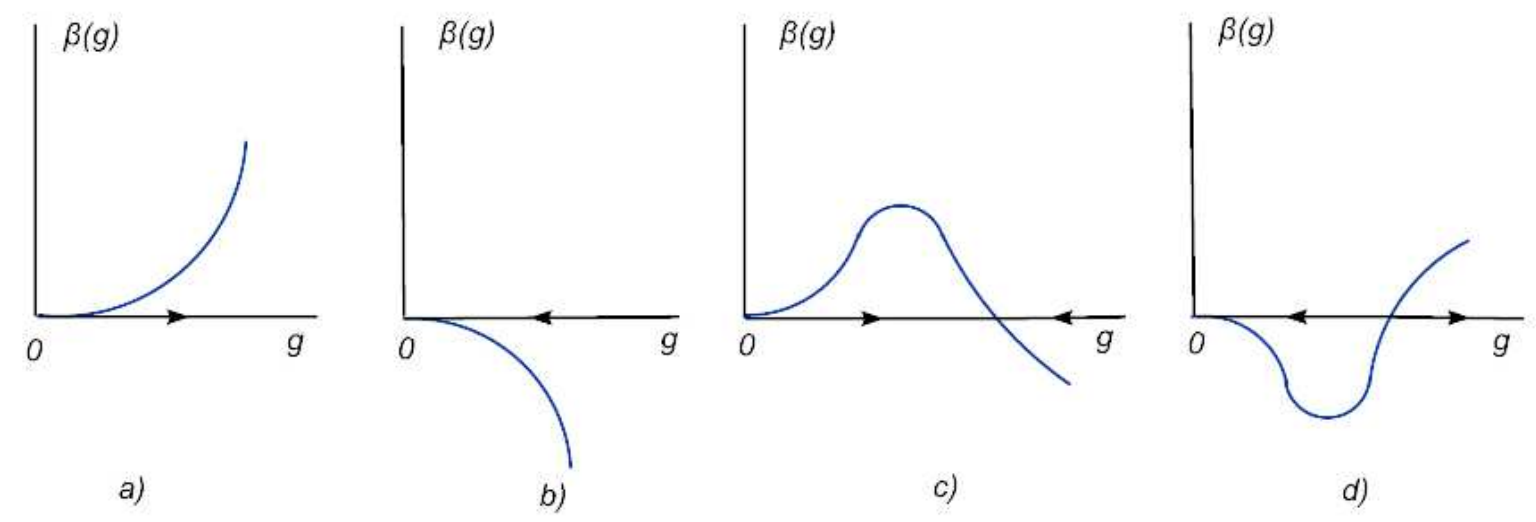

Fig.4.7.1.The possible form of the $\beta^*$-function.
The arrows show the behaviour of the
effective coupling in the UV regime: $t \to \infty$.

In the second case, the $\beta$-function is negative and, hence, the effective coupling decreases with increasing momentum. This situation appears in the one-loop approximation when $b < 0$, which takes place in the gauge theories. Here we also have a pole but in the infrared region.

In the third case, the $\beta$-function has zero: at first, it is positive and then is negative. This means that for small initial values the effective coupling increases; and for large ones, decreases. In both the cases, with increasing momentum it tends to the fixed value defined by the zero of the $\beta$-function. This is the so-called ultraviolet stable *fixed point*. It appears in some models in higher orders of perturbation theory.

## 4.8.Dimensional regularization and the $\overline{MS}$ scheme in a ghost sector

Consider now the calculation of the $(\beta^*(g^*_\varepsilon))_\varepsilon$ function and the anomalous dimensions in some particular models within the dimensional regularization and the minimal subtraction scheme. Note that in transition from dimension $4$ to $4-2\varepsilon$ the dimension of the coupling is changed and the "bare" coupling acquires the dimension $[g^*_{B,\varepsilon}] = 2\varepsilon$. That is why the relation between the "bare" and renormalized coupling contains the factor $(\mu^2)^\varepsilon$

$$(g^*_{B,\varepsilon})_\varepsilon = ((\mu^2)^\varepsilon Z_{g^*_\varepsilon} g^*_\varepsilon)_\varepsilon. \qquad (4.8.1)$$

Hence, even before the renormalization when $Z_g = 1$, in order to compensate this factor the dimensionless coupling $g$ should depend on $\mu$. Differentiating Eq.(4.8.1) with respect to $\mu^2$ one gets

$$0 = (\varepsilon Z_{g^*_\varepsilon} g^*_\varepsilon)_\varepsilon + \left(\frac{d\log Z_{g^*_\varepsilon}}{d\log\mu^2}\right)_\varepsilon (Z_{g^*_\varepsilon} g^*_\varepsilon)_\varepsilon + (Z_{g^*_\varepsilon})_\varepsilon \left(\frac{dg^*_\varepsilon}{d\log\mu^2}\right)_\varepsilon, \qquad (4.8.2)$$

i.e.,

$$(\beta_{4-2\varepsilon}(g^*_\varepsilon))_\varepsilon \equiv \left(\frac{dg^*_\varepsilon}{d\log\mu^2}\right)_\varepsilon = -\varepsilon(g^*_\varepsilon)_\varepsilon + (g^*_\varepsilon)_\varepsilon \left(\frac{d\log Z_{g^*_\varepsilon}}{d\log\mu^2}\right)_\varepsilon = -\varepsilon g^*_\varepsilon + \beta_4(g^*_\varepsilon)_\varepsilon. \qquad (4.8.3)$$

In the $\overline{MS}$ scheme the renormalization constants are given by the pole terms in $1/\varepsilon$ expansion and so does the bare coupling. They can be written as

$$(Z_{\Gamma^*_\varepsilon})_\varepsilon = 1 + \sum_{n=1}^{\infty} \frac{(c_n(g^*_\varepsilon))_\varepsilon}{(\varepsilon^n)_\varepsilon} = 1 + \sum_{n=1}^{\infty} \sum_{m=n}^{\infty} \frac{(c_{nm}(g^{*m}_\varepsilon))_\varepsilon}{(\varepsilon^n)_\varepsilon}. \qquad (4.8.4)$$

And similarly

$$(g^*_{Bare,\varepsilon})_\varepsilon = ((\mu^2)^\varepsilon)_\varepsilon\left[ (g^*_\varepsilon)_\varepsilon + \sum_{n=1}^\infty \frac{(a_n(g^*_\varepsilon))_\varepsilon}{(\varepsilon^n)_\varepsilon} \right]$$
$$= ((\mu^2)^\varepsilon)_\varepsilon\left[ (g^*_\varepsilon)_\varepsilon + \sum_{n=1}^\infty \sum_{m=n}^\infty \frac{(a_{nm}(g_\varepsilon^{*m+1}))_\varepsilon}{(\varepsilon^n)_\varepsilon} \right]. \quad (4.8.5)$$

Differentiating eq.(4.8.4) with respect to $\ln\mu^2$ and having in mind the definitions (4.7.14) and (4.7.15), one has:

$$-\left[ 1 + \sum_{n=1}^\infty \frac{(c_n(g^*_\varepsilon))_\varepsilon}{(\varepsilon^n)_\varepsilon} \right](\gamma_{\Gamma^*_\varepsilon}(g^*_\varepsilon))_\varepsilon$$
$$= [-(\varepsilon g^*_\varepsilon)_\varepsilon + (\beta(g^*_\varepsilon))_\varepsilon]\left( \frac{d}{dg^*_\varepsilon} \sum_{n=1}^\infty \frac{c_n(g^*_\varepsilon)}{\varepsilon^n} \right)_\varepsilon . \quad (4.8.6)$$

Equalizing the coefficients of equal powers of $\varepsilon$, one finds

$$(\gamma_{\Gamma^*_\varepsilon}(g^*_\varepsilon))_\varepsilon = (g^*_\varepsilon)_\varepsilon\left( \frac{d}{dg^*_\varepsilon} c_1(g^*_\varepsilon) \right)_\varepsilon ,$$
$$(g^*_\varepsilon)_\varepsilon\left( \frac{d}{dg^*_\varepsilon} c_n(g) \right)_\varepsilon = \quad (4.8.7)$$
$$(\gamma_{\Gamma^*_\varepsilon}(g^*_\varepsilon))_\varepsilon (c_{n-1}(g^*_\varepsilon))_\varepsilon + (\beta(g^*_\varepsilon))_\varepsilon\left( \frac{d}{dg^*_\varepsilon} c_{n-1}(g^*_\varepsilon) \right)_\varepsilon , \quad n \geq 2.$$

One sees that the coefficients of higher poles $(c_n(g^*_\varepsilon))_\varepsilon$, $n \geq 2$ are completely defined by that of the lowest pole $(c_1(g^*_\varepsilon))_\varepsilon$ and the $\beta^*$ function. In its turn the $\beta^*$-function is also defined by the lowest pole. To see this, consider Eq.(4.7.20). Differentiating it with respect to $\ln\mu^2$ one has

$$(\varepsilon)_\varepsilon\left[ (g^*_\varepsilon)_\varepsilon + \sum_{n=1}^\infty \frac{(a_n(g^*_\varepsilon))_\varepsilon}{(\varepsilon^n)_\varepsilon} \right] +$$
$$[-(\varepsilon g^*_\varepsilon)_\varepsilon + (\beta(g^*_\varepsilon))_\varepsilon]\left[ 1 + \left( \frac{d}{dg^*_\varepsilon} \sum_{n=1}^\infty \frac{a_n(g^*_\varepsilon)}{\varepsilon^n} \right)_\varepsilon \right] = 0. \quad (4.8.8)$$

Equalizing the coefficients of equal powers of $\varepsilon$, one finds

$$(\beta^*(g^*_\varepsilon))_\varepsilon = (g^*_\varepsilon)_\varepsilon\left( \frac{d}{dg^*_\varepsilon} a_1(g^*_\varepsilon) \right)_\varepsilon - (a_1(g^*_\varepsilon))_\varepsilon, \quad (4.8.9)$$

and

$$(g^*_\varepsilon)_\varepsilon\left( \frac{d}{dg^*_\varepsilon} a_n(g^*_\varepsilon) \right)_\varepsilon - (a_n(g^*_\varepsilon))_\varepsilon = (\beta(g^*_\varepsilon))_\varepsilon\left( \frac{d}{dg^*_\varepsilon} a_{n-1}(g^*_\varepsilon) \right)_\varepsilon , \quad n \geq 2. \quad (4.8.10)$$

Thus, knowing the coefficients of the lower poles one can reproduce all the higher order divergences. This means that they are not independent, all the information about them is connected in the lowest pole. In particular, substituting in (4.8.10) the perturbative expansion given by Eq.(4.8.5) one can solve the recurrent equation and find for the highest pole term

$$(a_{nn}(g^*_\varepsilon))_\varepsilon = (a^n_{11}(g^*_\varepsilon))_\varepsilon, \quad (4.8.11)$$

i.e. in the leading order one has the geometric progression

$$(g_{Bare,\varepsilon})_\varepsilon = \frac{(g_\varepsilon^*)_\varepsilon(\mu^{2\varepsilon})_\varepsilon}{1-(\varepsilon^{-1})_\varepsilon(g_\varepsilon)_\varepsilon(a_{11}(g_\varepsilon^*))_\varepsilon}, \tag{4.8.12}$$

which reflects the fact that the effective coupling in the leading log approximation (LLA) is also given by a geometric progression (4.8.12).The pole equations are easily generalized for the multiple couplings case, the higher poles are also expressed through the lower ones though the solutions of the RG equations are more complicated.

Consider now some particular models and calculate the corresponding $\beta^*$-functions and the anomalous dimensions.

**The $\phi_4^4$ theory**

**Standard matter sector**

We remind that standard matter sector of the $\phi^4$ theory defined by the inequality

$$1-\frac{1}{24}\left(\frac{g_\varepsilon}{\varepsilon}\right)_\varepsilon > 0. \tag{4.8.13}$$

i.e.,

$$\frac{1}{24}\left(\frac{g_\varepsilon}{\varepsilon}\right)_\varepsilon < 1. \tag{4.8.14}$$

**Remark.4.8.1.**Note that

$$(Z_{2,\varepsilon}^{-1})_\varepsilon = \left[1-\frac{1}{24}\left(\frac{g_\varepsilon}{\varepsilon}\right)_\varepsilon\right]^{-1} = 1+\frac{1}{24}\left(\frac{g_\varepsilon}{\varepsilon}\right)_\varepsilon+\ldots \tag{4.8.15}$$

The renormalization constants in the $\overline{MS}$ scheme up to two loops are given by Eq.(4.3.14) -Eq.(4.3.16), where $(g_\varepsilon)_\varepsilon = \dfrac{(\lambda_\varepsilon)_\varepsilon}{16\pi^2}$ :

$$\begin{aligned}
(Z_{4,\varepsilon})_\varepsilon &= 1+\frac{3}{2}\left(\frac{g_\varepsilon}{\varepsilon}\right)_\varepsilon+\frac{9}{4}\left(\frac{g_\varepsilon^2}{\varepsilon^2}\right)_\varepsilon-\frac{3}{2}\left(\frac{1}{\varepsilon}\right)_\varepsilon,\\
(Z_{2,\varepsilon}^{-1})_\varepsilon &= \left[1-\frac{1}{24}\left(\frac{g_\varepsilon}{\varepsilon}\right)_\varepsilon\right]^{-1} \simeq 1+\frac{1}{24}\left(\frac{g_\varepsilon}{\varepsilon}\right)_\varepsilon,\\
(Z_{g,\varepsilon})_\varepsilon &= (Z_{4,\varepsilon}Z_{2,\varepsilon}^{-2})_\varepsilon = 1+\frac{3}{2}\left(\frac{g_\varepsilon}{\varepsilon}\right)_\varepsilon+\left(\frac{9}{4}\left(\frac{g^2}{\varepsilon^2}\right)_\varepsilon-\frac{17}{12}\left(\frac{1}{\varepsilon}\right)_\varepsilon\right).
\end{aligned} \tag{4.8.16}$$

Notice that the higher pole coefficient $a_{22} = 9/4$ in the last expression is the square of the lowest pole one $a_{11} = 3/2$ in accordance with Eq.(4.8.11).Applying now Eq.(4.8.7) and Eq.(4.8.9) we get

$$\begin{aligned}
\gamma_4(g) &= \frac{3}{2}g-3g^2,\\
\gamma_2(g) &= \frac{1}{12}g^2,\\
\beta(g) = g(\gamma_4+2\gamma_2) &= \frac{3}{2}g^2-\frac{17}{6}g^2.
\end{aligned} \tag{4.8.17}$$

One can see from Eqs.(4.8.17) that the first coefficient of the $\beta$-function is $3/2$, i.e., the $\phi^4$ theory in standard sector belongs to the type of theories shown in Fig.4.7.1a). In the leading log approximation (LLA) one has a Landau pole behaviour. In the two-loop approximation (NLLA) the $\beta$-function gets a non-trivial zero and the effective coupling possesses an UV fixed point like the one shown in Fig.4.7.1). However, this fixed point is unstable with respect to higher orders and is not reliable. Here we encounter the problem of divergence of perturbation series in quantum field theory, they are the so-called asymptotic series which have a zero radius of convergence.

**The $\phi_4^4$ theory**

**Ghost matter sector**

We remind that ghost matter sector of the $\phi^4$ theory defined by the inequality

$$1 - \frac{1}{24}\left(\frac{g_\varepsilon}{\varepsilon}\right)_\varepsilon < 0, \tag{4.8.17}$$

i.e.,

$$\frac{1}{24}\left(\frac{g_\varepsilon}{\varepsilon}\right)_\varepsilon > 1. \tag{4.8.18}$$

**Remark.4.8.2.** Note that

$$\begin{gathered}(Z_{2,\varepsilon}^{-1})_\varepsilon = \left[1 - \frac{1}{24}\left(\frac{g_\varepsilon}{\varepsilon}\right)_\varepsilon\right]^{-1} = \\ -\left(\frac{\varepsilon}{g_\varepsilon}\right)_\varepsilon\left(1 + \frac{1}{24}\left(\frac{\varepsilon}{g_\varepsilon}\right)_\varepsilon + \ldots\right) = -\left(\frac{\varepsilon}{g_\varepsilon}\right)_\varepsilon - \frac{1}{24}\left(\frac{\varepsilon^2}{g_\varepsilon^2}\right)_\varepsilon - \ldots\end{gathered} \tag{4.8.19}$$

The renormalization constants in the $\overline{MS}$ scheme up to two loops are given by Eq.(4.3.14)- Eq.(4.3.16), where $(g_\varepsilon)_\varepsilon \equiv \dfrac{(\lambda_\varepsilon)_\varepsilon}{16\pi^2}$ :

$$\begin{gathered}(Z_{4,\varepsilon})_\varepsilon = 1 + \frac{3}{2}\left(\frac{g_\varepsilon}{\varepsilon}\right)_\varepsilon + \frac{9}{4}\left(\frac{g_\varepsilon^2}{\varepsilon^2}\right)_\varepsilon - \frac{3}{2}\left(\frac{1}{\varepsilon}\right)_\varepsilon, \\ (Z_{2,\varepsilon}^{-1})_\varepsilon = \left[1 - \frac{1}{24}\left(\frac{g_\varepsilon}{\varepsilon}\right)_\varepsilon\right]^{-1} \simeq -\left(\frac{\varepsilon}{g_\varepsilon}\right)_\varepsilon - \frac{1}{24}\left(\frac{\varepsilon^2}{g_\varepsilon^2}\right)_\varepsilon, \\ (Z_{g,\varepsilon})_\varepsilon = (Z_{4,\varepsilon}Z_{2,\varepsilon}^{-2})_\varepsilon = -\left[1 + \frac{3}{2}\left(\frac{g_\varepsilon}{\varepsilon}\right)_\varepsilon + \left(\frac{9}{4}\left(\frac{g_\varepsilon^2}{\varepsilon^2}\right)_\varepsilon - \frac{17}{12}\left(\frac{1}{\varepsilon}\right)_\varepsilon\right)\right] \times \\ \left[\left(\frac{\varepsilon}{g_\varepsilon}\right)_\varepsilon + \frac{1}{24}\left(\frac{\varepsilon^2}{g_\varepsilon^2}\right)_\varepsilon\right] \simeq -\frac{9}{4}\left(\frac{g_\varepsilon}{\varepsilon}\right)_\varepsilon + const.\end{gathered} \tag{4.8.20}$$

Applying now Eq.(4.8.7) and Eq.(4.8.9) we get

$$(\beta(g_\varepsilon))_\varepsilon \simeq -\frac{9}{4}(g_\varepsilon^2)_\varepsilon. \tag{5.8.21}$$

One can see from Eqs.(4.8.21) that the first coefficient of the $\beta$-function is $-9/4$, i.e., the $\phi^4$ theory in a ghost sector belongs to the type of theories shown in Fig.4.7.1b).

# 5.Renormalizability-of-Higher-Derivative-Quantum-Gravity.

## 5.1.The Higher-Derivative Theories of Gravitation.Green's functions.

Adding quadratic products of the curvature tensor to the gravitational action leads to field equations in which some terms involve four derivatives. While it is not the purpose of this paper to investigate the novel consequences of these classical field equations, a brief summary of some of the salient features is in order to give a grounding to the following discussion of renormalization.

Gravitational actions which include terms quadratic in the curvature tensor are renormalizable. The necessary Slavnov identities are derived from Becchi-Rouet-Stora (BRS) transformations of the gravitational and Faddeev-Popov ghost fields. In general, non-gauge-invariant divergences do arise, but they may be absorbed by nonlinear renormalizations of the gravitational and ghost fields and of the BRS transformations

[13].The generic expression of the action reads

$$I_{sym} = -\int d^4x\sqrt{-g}\,(\alpha R_{\mu\nu}R^{\mu\nu} - \beta R^2 + 2\kappa^{-2}R), \tag{5.1.1}$$

where the curvature tensor and the Ricci is defined by $R^{\lambda}_{\mu\alpha\nu} = \partial_\nu \Gamma^{\lambda}_{\mu\alpha}$ and $R_{\mu\nu} = R^{\lambda}_{\mu\lambda\nu}$ correspondingly, $\kappa^2 = 32\pi G,$ we used the signature $(-+++)$. The convenient definition of the gravitational field variable in terms of the contravariant metric density reads

$$\kappa h^{\mu\nu} = g^{\mu\nu}\sqrt{-g} - \eta^{\mu\nu}. \tag{5.1.2}$$

Analysis of the linearized radiation shows that there are eight dynamical degrees of freedom in the field. Two of these excitations correspond to the familiar massless spin-2 graviton. Five more correspond to a massive spin-2 particle with mass $m_2$. The eighth corresponds to a massive scalar particle with mass $m_0$. Although the linearized field energy of the massless spin-2 and massive scalar excitations is positive definite, the linearized energy of the massive spin-2 excitations is negative definite. This feature is characteristic of higher-derivative models, and poses the major obstacle to their physical interpretation.

In the quantum theory, there is an alternative problem which may be substituted for the negative energy. It is possible to recast the theory so that the massive spin-2 eigenstates of the free-fieid Hamiltonian have positive-definite energy, but also negative norm in the state vector space.

These negative-norm states cannot be excluded from the physical sector of the vector space without destroying the unitarity of the S matrix. The requirement that the graviton propagator behave like $p^{-4}$ for large momenta makes it necessary to choose the indefinite-metric vector space over the negative-energy states.

The presence of massive quantum states of negative norm which cancel some of the divergences due to the massless states is analogous to the Pauli-Villars regularization of other field theories. For quantum gravity, however, the resulting improvement in the ultraviolet behavior of the theory is sufficient only to make it renormalizable, but not finite.

The gauge choice which we adopt in order to defining the quantum theory is the canonical harmonic gauge: $\partial_\nu h^{\mu\nu} = 0$. Corresponding Green's functions are then given by a generating functional

$$\begin{aligned} Z(T_{\mu\nu}) = N\int\Big[\prod\nolimits_{\mu\leq\nu} dh^{\mu\nu}\Big][dC^{\sigma}][d\overline{C}_{\tau}]\delta^4(F^{\tau}) \\ \exp\Big[\,i\Big(I_{sym} + \int d^4x\overline{C}_{\tau}\vec{F}^{\tau}_{\mu\nu}D^{\mu\nu}_{\alpha}C^{\alpha} + \kappa\int d^4xT_{\mu\nu}h^{\mu\nu}\Big)\Big]. \end{aligned} \tag{5.1.3}$$

Here $F^{\tau} = \vec{F}^{\tau}_{\mu\nu}h^{\mu\nu}, \vec{F}^{\tau}_{\mu\nu} = \delta^{\tau}_{\mu}\vec{\partial}_\nu$ and the arrow indicates the direction in which the derivative acts. $N$ is an normalization constant. $C^{\sigma}$ is the Faddeev-Popov ghost field, and $\overline{C}_{\tau}$ is the antighost field. Notice that both $C^{\sigma}$ and $\overline{C}_{\tau}$ are anticommuting quantities. $D^{\mu\nu}_{\alpha}$ is the operator which generates gauge transformations in $h^{\mu\nu}$, given an arbitrary spacetime-dependent vector $\xi^{\alpha}(x)$ corresponding to $x^{\mu'} = x^{\mu} + \kappa\xi^{\mu}$ and where

$$D^{\mu\nu}_{\alpha}\xi^{\alpha}(x) = \partial^{\mu}\xi^{\nu} + \partial^{\nu}\xi^{\mu} - \eta^{\mu\nu}\partial_{\alpha}\xi^{\alpha} + \kappa(\partial_{\alpha}\xi^{\mu}h^{\alpha\nu} + \partial_{\alpha}\xi^{\nu}h^{\alpha\mu} - \xi^{\alpha}\partial_{\alpha}h^{\mu\nu} - \partial_{\alpha}\xi^{\alpha}h^{\mu\nu}.) \tag{5.1.4}$$

In the functional integral (5.1.3), we have written the metric for the gravitational field as $\Big[\prod_{\mu\leq\nu} dh^{\mu\nu}\Big]$ without any local factors of $g = \det(g_{\mu\nu})$. Such factors do not contribute to the Feynman rules because their effect is to introduce terms proportional to $\delta^4(0)\int d^4x\ln(-g)$ into the effective action and $\delta^4(0)$ is set equal to zero in dimensional regularization.

In calculating the generating functional (5.1.3.) by using the loop expansion, one may represent the $\delta$-function which fixes the gauge as the limit of a Gaussian, discarding an infinite normalization constant

$$\delta^4(F^\tau) \sim \lim_{\Delta\to 0} \exp\Big[ i\Big(\tfrac{1}{2}\Delta^{-1}\int d^4xF_\tau F^\tau\Big)\Big]. \qquad (5.1.5)$$

In this expression, the index $\tau$ has been lowered using the flat-space metric tensor $\eta_{\mu\nu}$. For the remainder of this paper, we shall adopt the standard approach to the covariant quantization of gravity, in which only Lorentz tensors occur, and all raising and lowering of indices is done with respect to flat space. The graviton propagator may be calculated from $I_{sym} + \frac{1}{2}\Delta^{-1}\int d^4xF_\tau F^\tau$ in the usual fashion, letting $\Delta \to 0$ after inverting. The expression $\frac{1}{2}\Delta^{-1}\int d^4xF_\tau F^\tau$ contains only two derivatives. Consequently, there are parts of the graviton propagator which behave like $p^{-2}$ for large momenta. Specifically, the $p^{-2}$ terms consist of everything but those parts of the propagator which are transverse in all indices. These terms give rise to unpleasant infinities already at the one-loop order. For example, the graviton self-energy diagram shown in Fig.5.1.1 has a divergent part with the general structure $(\partial^4 h)^2$. Such divergences do cancel when they are connected to tree diagrams whose outermost lines are on the mass shell, as they must if the S matrix is to be made finite without introducing counterterms for them. However, they greatly complicate the renormalization of Green's functions.

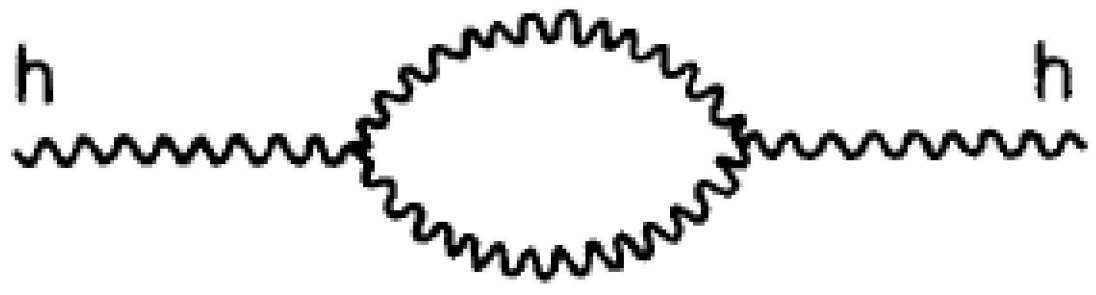

Fig.5.1.1.The one-loop graviton self-energy diagram.

We may attempt to extricate ourselves from the situation described in the last paragraph by picking a different weighting functional. Keeping in mind that we want no part of the graviton propagator to fall off slower than $p^{-4}$ for large momenta, we now choose the weighting functional [12]

$$\omega_4(e^\tau) = \exp\Big[ i\Big(\tfrac{1}{2}\Delta^{-1}\int d^4xe_\tau\Box^2e^\tau\Big)\Big], \qquad (5.1.6)$$

where $e^\tau$ is any four-vector function.The corresponding gauge-fixing term in the effective action is

$$-\tfrac{1}{2}\kappa^2\Delta^{-1}\int d^4xF_\tau\Box^2F^\tau. \qquad (5.1.7)$$

The graviton propagator resulting from the gauge-fixing term (5.1.7) is derived in [13]. For most values of the parameters $\alpha$ and $\beta$ in $I_{sym}$ it satisfies the requirement that all its leading parts fall off like $p^{-4}$ for large momenta. There are, however, specific choices of these parameters which must be avoided. If $\alpha = 0,$ the massive spin-2 excitations disappear, and inspection of the graviton propagator shows that some terms then behave like $k^{-2}$. Likewise, if $3\beta - \alpha = 0,$ the massive scalar excitation disappears, and there are again terms in the propagator which behave like $p^{-2}$. However, even if we avoid the special cases $\alpha = 0$ and $3\beta - \alpha = 0,$ and if we use the propagator derived from (5.1.7), we still do not obtain a clean renormalization of the Green's functions. We now turn to the implications of gauge invariance.Before we write down the BRS transformations for gravity, let us first establish the commutation relation for gravitational gauge transformations, which reveals the group structure of the theory. Take the gauge transformation (5.1.4) of $h^{\mu\nu},$ generated by $\xi^{\mu}$ and perform a second gauge transformation, generated by $\eta^{\mu},$ on the $h^{\mu\nu}$ fields appearing there. Then antisymmetrize in $\xi^{\mu}$ and $\eta^{\mu}$.The result is

$$\frac{\delta D^{\mu\nu}_{\alpha}}{\delta h^{\rho\sigma}} D^{\rho\sigma}_{\beta}(\xi^{\alpha}\eta^{\beta} - \eta^{\alpha}\xi^{\beta}) = \kappa D^{\mu\nu}_{\lambda}(\partial_{\alpha}\xi^{\lambda}\eta^{\alpha} - \partial_{\alpha}\xi^{\alpha}\eta^{\lambda}), \tag{5.1.8}$$

where the repeated indices denote both summation over the discrete values of the indices and integration over the spacetime arguments of the functions or operators indexed.

The BRS transformations for gravity appropriate for the gauge-fixing term (5.1.6) are [13]

$$\begin{gathered}(a)\ \delta_{\mathbf{BRS}} h^{\mu\nu} = \kappa D^{\mu\nu}_{\alpha} C^{\alpha}\delta\lambda, (b)\ \delta_{\mathbf{BRS}} C^{\alpha} = -\kappa^{2}\partial_{\beta}C^{\alpha}C^{\beta}\delta\lambda,\\ (c)\ \delta_{\mathbf{BRS}}\overline{C}_{\tau} = -\kappa^{3}\Delta^{-1}\Box^{2}F_{\tau}\delta\lambda,\end{gathered} \tag{5.1.9}$$

where $\delta\lambda$ is an infinitesimal anticommuting constant parameter.The importance of these transformations resides in the quantities which they leave invariant. Note that

$$\delta_{\mathbf{BRS}}(\partial_{\beta}C^{\sigma}C^{\beta}) = 0 \tag{5.1.10}$$

and

$$\delta_{\mathbf{BRS}}(D^{\mu\nu}_{\alpha}C^{\alpha}) = 0. \tag{5.1.11}$$

As a result of Eq. (5.1.11), the only part of the ghost action which varies under the BRS transformations is the antighost $\overline{C}_{\tau}$. Accordingly, the transformation (2.2.9c) has been chosen to make the variation of the ghost action just cancel the variation of the gauge-fixing term. Therefore, the entire effective action is BRS invariant:

$$\delta_{\mathbf{BRS}}\left(I_{sim} - \tfrac{1}{2}\kappa^{2}\Delta^{-1}F_{\tau}\Box^{2}F^{\tau} + \overline{C}_{\tau}F^{\tau}_{\mu\nu}D^{\mu\nu}_{\alpha}C^{\alpha}\right) = 0. \tag{5.1.12}$$

Equations (5.1.9), (5.1.10), and (5.1.12) now enable us to write the Slavnov identities in an economical way. In order to carry out the renormalization program, we will need to have Slavnov identities for the proper vertices.

## 5.1.1. Slavnov identities for Green's functions

First consider the Slavnov identities for Green's functions.

$$\begin{gathered}Z(T_{\mu\nu},\overline{\beta}_{\sigma},\beta^{\tau},K_{\mu\nu},L_{\sigma}) = N\int\left[\prod\nolimits_{\mu\leq\nu} dh^{\mu\nu}\right][dC^{\sigma}][d\overline{C}_{\tau}]\\ \exp\left[i\tilde{\Sigma}(h^{\mu\nu},C^{\sigma},\overline{C}_{\tau},K_{\mu\nu},L_{\sigma},\overline{\beta}_{\sigma}C^{\sigma}) + \overline{\beta}_{\sigma}C^{\sigma} + \overline{C}_{\tau}\beta^{\tau} + \kappa T_{\mu\nu}h^{\mu\nu}\right].\end{gathered} \tag{5.1.13}$$

Anticommuting sources have been included for the ghost and antighost fields, and the effective action $\tilde{\Sigma}$ has been enlarged by the inclusion of BRS invariant couplings of the ghosts and gravitons to some external fields $K_{\mu\nu}$ (anticommuting) and $L_\sigma$ (commuting),

$$\tilde{\tilde{\Sigma}} = I_{sim} - \frac{1}{2}\kappa^2\Delta^{-1}F_\tau\Box^2 F^\tau + \overline{C}_\tau \vec{F}^\tau_{\mu\nu} D^{\mu\nu}_\alpha C^\alpha + \kappa K_{\mu\nu} D^{\mu\nu}_\alpha + \kappa^2 L_\sigma \partial_\beta C^\sigma C^\beta. \qquad (5.1.14)$$

$\tilde{\tilde{\Sigma}}$ is BRS invariant by virtue of Eq.(5.1.9), Eq.(5.1.10), and Eq.(5.1.12). We may use the new couplings to write this invariance as

$$\frac{\delta\tilde{\tilde{\Sigma}}}{\delta K_{\mu\nu}}\frac{\delta\tilde{\tilde{\Sigma}}}{\delta h^{\mu\nu}} + \frac{\delta\tilde{\tilde{\Sigma}}}{\delta L_\sigma}\frac{\delta\tilde{\tilde{\Sigma}}}{\delta C^\sigma} + \kappa^3\Delta^{-1}\Box^2 F_\tau \frac{\delta\tilde{\tilde{\Sigma}}}{\delta\overline{C}_\tau}. \qquad (5.1.15)$$

In this equation, and throughout this subsection, we use left variational derivatives with respect to anticommuting quantities: $\delta f(C^\sigma) = \delta C^\tau \delta f/\delta C^\tau$. Equation (2.2.15) may be simplified by rewriting it in terms of a reduced effective action,

$$\Sigma = \tilde{\tilde{\Sigma}} + \frac{1}{2}\kappa^2\Delta^{-1}F_\tau\Box^2 F^\tau. \qquad (5.1.16)$$

Substitution of (5.1.16) into (5.1.15) gives

$$\frac{\delta\Sigma}{\delta K_{\mu\nu}}\frac{\delta\Sigma}{\delta h^{\mu\nu}} + \frac{\delta\Sigma}{\delta L_\sigma}\frac{\delta\Sigma}{\delta C^\sigma} = 0, \qquad (5.1.17)$$

where we have used the relation

$$\kappa^{-1}\vec{F}^\tau_{\mu\nu}\frac{\delta\Sigma}{\delta K_{\mu\nu}} - \frac{\delta\Sigma}{\delta\overline{C}_\tau} = 0. \qquad (5.1.18)$$

Note that a measure

$$\left[\prod\nolimits_{\mu\leq\nu} dh^{\mu\nu}\right][dC^\sigma][d\overline{C}_\tau] \qquad (5.1.19)$$

is BRS invariant since for infinitesimal transformations, the Jacobian is $1$, because of the trace relations

$$(a)\ \frac{\delta^2\tilde{\tilde{\Sigma}}}{\delta K_{(\mu\nu)}\delta h^{(\mu\nu)}} = 0,$$
$$(b)\ \frac{\delta^2\tilde{\tilde{\Sigma}}}{\delta C^\sigma \delta L_\sigma} = 0, \qquad (5.1.20)$$

both of which follow from $\int d^4x \partial_\alpha C^\alpha = 0$. The parentheses surrounding the indices in (5.1.20a) indicate that the summation is to be carried out only for $\mu \leq \nu$.

**Remark 5.1.1.** Note that the Slavnov identity for the generating functional of Green's functions is obtained by performing the BRS transformations (5.1.9) on the integration variables in the generating functional (5.1.13). This transformation does not change the value of the generating functional and therefore we obtain

$$N\int\left[\prod\nolimits_{\mu\leq\nu} dh^{\mu\nu}\right][dC^\sigma][d\overline{C}_\tau] \times$$
$$\left(\kappa^2 T_{\mu\nu}D^{\mu\nu}_\alpha - \kappa^2\overline{\beta}_\sigma\partial_\beta C^\sigma C^\beta + \kappa^3\Delta^{-1}\beta^\tau\Box^2\vec{F}_{\tau\mu\nu}h^{\mu\nu}\right) \times \qquad (5.1.21)$$
$$\exp\left[i\left(\tilde{\tilde{\Sigma}} + \kappa T_{\mu\nu}h^{\mu\nu} + \overline{\beta}_\sigma C^\sigma + \overline{C}_\tau\beta^\tau\right)\right] = 0.$$

Another identity which we shall need is the ghost equation of motion. To derive this equation, we shift the antighost integration variable $\overline{C}_\tau$ to $\overline{C}_\tau + \delta\overline{C}_\tau$, again with no resulting change in the value of the generating functional:

$$N\int\left[\prod\nolimits_{\mu\leq\nu} dh^{\mu\nu}\right][dC^\sigma][d\overline{C}_\tau]\left(\frac{\delta\tilde{\tilde{\Sigma}}}{\delta C^\sigma} + \beta^\tau\right)\exp\left[i\left(\tilde{\tilde{\Sigma}} + \kappa T_{\mu\nu}h^{\mu\nu} + \overline{\beta}_\sigma C^\sigma + \overline{C}_\tau\beta^\tau\right)\right] \qquad (5.1.22)$$

We define now the generating functional of connected Green's functions as the logarithm of the functional (5.1.13),

$$W[T_{\mu\nu},\overline{\beta}_{\sigma},\beta^{\tau},K_{\mu\nu},L_{\sigma}] = -i\ln Z[T_{\mu\nu},\overline{\beta}_{\sigma},\beta^{\tau},K_{\mu\nu},L_{\sigma}]. \tag{5.1.23}$$

and make use of the couplings to the external fields $K_{\mu\nu}$ and $L_{\sigma}$ to rewrite (5.1.22) in terms of $W$

$$\kappa T_{\mu\nu}\frac{\delta W}{\delta K_{\mu\nu}} - \overline{\beta}_{\sigma}\frac{\delta W}{\delta L_{\sigma}} + \kappa^{2}\Delta^{-1}\beta^{\tau}\Box^{2}\overrightarrow{F}_{\tau\mu\nu}\frac{\delta W}{\delta T_{\mu\nu}} = 0. \tag{5.1.24}$$

Similarly, we get the ghost equation of motion:

$$\kappa^{-1}\overrightarrow{F}^{\tau}_{\mu\nu}\frac{\delta W}{\delta K_{\mu\nu}} + \beta^{\tau} = 0. \tag{5.1.25}$$

## 5.1.2. Proper vertices

A Legendre transformation takes us from the generating functional of connected Green's functions (2.2.23) to the generating functional of proper vertices. First, we define the expectation values of the gravitational, ghost, and antighost fields in the presence of the sources $T_{\mu\nu},\overline{\beta}_{\sigma}$, and $\beta^{\tau}$ and the external fields $K_{\mu\nu}$ and $L_{\sigma}$

$$(a)\ h^{\mu\nu}(x) = \frac{\delta W}{\kappa\delta T_{\mu\nu}(x)},(b)\ C^{\sigma}(x) = \frac{\delta W}{\delta\overline{\beta}_{\sigma}(x)},(c)\ \overline{C}_{\tau}(x) = \frac{\delta W}{\delta\beta^{\tau}(x)}. \tag{5.1.26}$$

We have chosen to denote the expectation values of the fields by the same symbols which were used for the fields in the effective action (5.1.14).

The Legendre transformation can now be performed, giving us the generating functional of proper vertices as a functional of the new variables (5.1.26) and the external fields $K_{\mu\nu}$ and $L_{\sigma}$

$$\widetilde{\Gamma}[h^{\mu\nu},C^{\sigma},\overline{C}_{\tau},K_{\mu\nu},L_{\sigma}] = W[T_{\mu\nu},\overline{\beta}_{\sigma},\beta^{\tau},K_{\mu\nu},L_{\sigma}] - \kappa T_{\mu\nu}h^{\mu\nu} - \overline{\beta}_{\sigma}C^{\sigma} - \overline{C}_{\tau}\beta^{\tau}. \tag{5.1.27}$$

In this equation, the quantities $T_{\mu\nu}$, $\overline{\beta}_{\sigma}$, and $\beta^{\tau}$ are given implicitly in terms of $h^{\mu\nu},C^{\sigma},\overline{C}_{\tau},K_{\mu\nu},$ and $L_{\sigma}$ by Eq.(5.1.26).The relations dual to (5.1.26) are

$$(a)\ \kappa T_{\mu\nu}(x) = -\frac{\delta\widetilde{\Gamma}}{\delta h^{\mu\nu}(x)},(b)\ \overline{\beta}_{\sigma}(x) = \frac{\delta\widetilde{\Gamma}}{\delta C^{\sigma}(x)},(c)\ \beta_{\tau}(x) = -\frac{\delta\widetilde{\Gamma}}{\delta\overline{C}_{\tau}(x)}. \tag{5.1.28}$$

Since the external fields $K_{\mu\nu}$ and $L_{\sigma}$ do not participate in the Legendre transformation (5.1.26), for them we have the relations

$$(a)\ \frac{\delta\widetilde{\Gamma}}{\delta K_{\mu\nu}(x)} = \frac{\delta W}{\delta K_{\mu\nu}(x)},(b)\ \frac{\delta\widetilde{\Gamma}}{\delta L_{\sigma}(x)} = \frac{\delta W}{\delta L_{\sigma}(x)}. \tag{5.1.29}$$

Finally, the Slavnov identity for the generating functional of proper vertices is obtained by transcribing (5.1.24) using the relations (5.1.26), (5.1.28), and (5.1.29)

$$\frac{\delta\widetilde{\Gamma}}{\delta K_{\mu\nu}}\frac{\delta\widetilde{\Gamma}}{\delta h^{\mu\nu}} + \frac{\delta\widetilde{\Gamma}}{\delta L_{\sigma}}\frac{\delta\widetilde{\Gamma}}{\delta C^{\sigma}} + \kappa^{3}\Delta^{-1}\Box^{2}\overrightarrow{F}_{\tau\mu\nu}h^{\mu\nu}\frac{\delta\widetilde{\Gamma}}{\delta C^{\sigma}} = 0. \tag{5.1.30}$$

We also have the ghost equation of motion,

$$\kappa^{-1}\overrightarrow{F}^{\tau}_{\mu\nu}\frac{\delta\widetilde{\Gamma}}{\delta K_{\mu\nu}} - \frac{\delta\widetilde{\Gamma}}{\delta C^{\sigma}} = 0. \tag{5.1.31}$$

Since Eq. (5.1.30) has exactly the same form as (5.1.15), we follow the example set by (5.1.16) and define a reduced generating functional of the proper vertices,

$$\Gamma = \widetilde{\Gamma} + \tfrac{1}{2}\kappa^{2}\Delta^{-1}\left(\overrightarrow{F}_{\tau\mu\nu}h^{\mu\nu}\right)\Box^{2}\left(\overrightarrow{F}^{\tau}_{\rho\sigma}h^{\rho\sigma}\right). \tag{5.1.32}$$

Substituting this into (5.1.30) and (5.1.31), the Slavnov identity becomes

$$\frac{\delta\Gamma}{\delta K_{\mu\nu}}\frac{\delta\Gamma}{\delta h^{\mu\nu}} + \frac{\delta\Gamma}{\delta L_\sigma}\frac{\delta\Gamma}{\delta C^\sigma} = 0. \qquad (5.1.33)$$

and the ghost equation of motion becomes

$$\kappa^{-1}\vec{F}^{\tau}_{\mu\nu}\frac{\delta\Gamma}{\delta K_{\mu\nu}} - \frac{\delta\Gamma}{\delta\overline{C}_\tau} = 0. \qquad (5.1.34)$$

Equations (5.1.33) and (5.1.34) are of exactly the same form as (5.5) and (5.6). This is as it should be, since at the zero-loop order

$$\Gamma^{(0)} = \Sigma. \qquad (5.1.35)$$

# 5.1.3. Structure of the divergences and renormalization equation.

The Slavnov identity (5.1.33) is quadratic in the functional $\Gamma$. This nonlinearity is reflected in the fact that the renormalization of the effective action generally also involves the renormalization of the BRS transformations which must leave the effective action invariant.

The canonical approach uses the Slavnov identity for the generating functional of proper vertices to derive a linear equation for the divergent parts of the proper vertices. This equation is then solved to display the structure of the divergences. From this structure, it can be seen how to renormalize the effective action so that it remains invariant under a renormalized set of BRS transformations [13].

Suppose that we have successfully renormalized the reduced effective action up to $n-1$ loop order; that is, suppose we have constructed a quantum extension of $\Sigma$ which satisfies Eqs. (5.1.17) and (5.1.18) exactly, and which leads to finite proper vertices when calculated up to order $n-1$. We will denote this renormalized quantity by $\Sigma^{(n-1)}$. In general, it contains terms of many different orders in the loop expansion, including orders greater than $n-1$. The $n-1$ loop part of the reduced generating functional of proper vertices will be denoted by $\Gamma^{(n-1)}$.

When we proceed to calculate $\Gamma^{(n)}$, we find that it contains divergences. Some of these come from $n$-loop Feynman integrals. Since all the subintegrals of an $n$-loop Feynman integral contain less than w loops, they are finite by assumption. Therefore, the divergences which arise from w-loop Feynman integrals come only from the overall divergences of the integrals, so the corresponding parts of $\Gamma^{(n)}$ are local in structure. In the dimensional regularization procedure, these divergences are of order $\epsilon^{-1} = (d-4)^{-1}$, where $d$ is the dimensionality of spacetime in the Feynman integrals.

There may also be divergent parts of $\Gamma^{(n)}$ which do not arise from loop integrals, and which contain higher-order poles in the regulating parameter $\epsilon$. Such divergences comes from $n$-loop order parts of $\Sigma^{(n-1)}$ which are necessary to ensure that (5.1.17) is satisfied. Consequently, they too have a local structure. We may separate the divergent and finite parts of $\Gamma^{(n)}$:

$$\Gamma^{(n)} = \Gamma^{(n)}_{\mathbf{div}} + \Gamma^{(n)}_{\mathbf{finite}}. \qquad (5.1.36)$$

If we insert this breakup into Eq. (5.1.20), and keep only the terms of the equation which are of $n$-loop order, we get

$$
\begin{gathered}
\frac{\delta\Gamma^{(n)}_{\mathbf{div}}}{\delta K_{\mu\nu}}\frac{\delta\Gamma^{(0)}}{\delta h^{\mu\nu}}+\frac{\delta\Gamma^{(0)}}{\delta K_{\mu\nu}}\frac{\delta\Gamma^{(n)}_{\mathbf{div}}}{\delta h^{\mu\nu}}+\frac{\delta\Gamma^{(n)}_{\mathbf{div}}}{\delta L_{\sigma}}\frac{\delta\Gamma^{(0)}}{\delta C^{\sigma}}+\frac{\delta\Gamma^{(0)}}{\delta L_{\sigma}}\frac{\delta\Gamma^{(n)}_{\mathbf{div}}}{\delta C^{\sigma}}= \\
-\sum\nolimits_{i=0}^{n}\left[\frac{\delta\Gamma^{(n-i)}_{\mathbf{finite}}}{\delta K_{\mu\nu}}\frac{\delta\Gamma^{(i)}_{\mathbf{finite}}}{\delta h^{\mu\nu}}+\frac{\delta\Gamma^{(n-i)}_{\mathbf{finite}}}{\delta L_{\sigma}}\frac{\delta\Gamma^{(i)}_{\mathbf{finite}}}{\delta C^{\sigma}}\right].
\end{gathered}
\tag{5.1.37}
$$

Since each term on the right-hand side of (5.1.37) remains finite as $\epsilon \to 0$, while each term on the left-hand side contains a factor with at least a simple pole in e, each side of the equation must vanish separately. Remembering the Eq.(5.1.35), we can write the following equation, called the renormalization equation:

$$
\Re\Gamma^{(n)}_{\mathbf{div}} = 0, \tag{5.1.38}
$$

where

$$
\Re = \frac{\delta\Sigma}{\delta h^{\mu\nu}}\frac{\delta}{\delta K_{\mu\nu}}+\frac{\delta\Sigma}{\delta C^{\sigma}}\frac{\delta}{\delta L_{\sigma}}+\frac{\delta\Sigma}{\delta K_{\mu\nu}}\frac{\delta}{\delta h^{\mu\nu}}+\frac{\delta\Sigma}{\delta L_{\sigma}}\frac{\delta}{\delta C^{\sigma}}. \tag{5.1.39}
$$

Similarly by collecting the $n$-loop order divergences in the ghost equation of motion (5.1.34) we get

$$
\kappa^{-1}\vec{F}^{\tau}_{\mu\nu}\frac{\delta\Gamma^{(n)}_{\mathbf{div}}}{\delta K_{\mu\nu}}-\frac{\delta\Gamma^{(n)}_{\mathbf{div}}}{\delta\overline{C}_{\tau}} = 0. \tag{5.1.40}
$$

In order to construct local solutions to Eqs. (5.1.38) and (5.1.40) remind that the operator $\Re$ defined in (5.1.39) is nilpotent [13]:

$$
\Re^{2} = 0. \tag{5.1.41}
$$

Equation (5.1.41) gives us the local solutions to Eq.(5.1.38) of the form

$$
\Gamma^{(n)}_{\mathbf{div}} = \Im(h^{\mu\nu})+\Re[X(h^{\mu\nu},C^{\sigma},\overline{C}_{\tau},K_{\mu\nu},L_{\sigma})], \tag{5.1.42}
$$

where $\Im$ is an arbitrary gauge-invariant local functional of $h^{\mu\nu}$ and its derivatives, and $X$ is an arbitrary local functional of $h^{\mu\nu},C_{\sigma},\overline{C}_{\tau},K_{\mu\nu}$ and $L^{\sigma}$ and their derivatives. In order to satisfy the ghost equation of motion (5.1.40) we require that

$$
\Gamma^{(n)}_{\mathbf{div}} = \Gamma^{(n)}_{\mathbf{div}}\left(h^{\mu\nu},C^{\sigma},K_{\mu\nu}-\kappa^{-1}\overline{C}_{\tau}\vec{F}^{\tau}_{\mu\nu},L_{\sigma}\right). \tag{5.1.43}
$$

## 5.1.4. Ghost number and power counting

Structure of the effective action (5.1.14) shows that we may define the following conserved quantity, called ghost number [13]:

$$
\begin{gathered}
N_{\mathbf{ghost}}[h^{\mu\nu}] = 0, N_{\mathbf{ghost}}[C_{\sigma}] = +1, N_{\mathbf{ghost}}[\overline{C}_{\tau}] = -1, \\
N_{\mathbf{ghost}}[K_{\mu\nu}] = -1, N_{\mathbf{ghost}}[L_{\sigma}] = -2.
\end{gathered}
\tag{5.1.44}
$$

From Eqs.(2.2.44) follows that

$$
N_{\mathbf{ghost}}[\Sigma] = N_{\mathbf{ghost}}[\Gamma] = 0. \tag{5.1.45}
$$

Since

$$
N_{\mathbf{ghost}}[\Re] = +1, \tag{5.1.46}
$$

we require of the functional $X(\cdot)$ that

$$
N_{\mathbf{ghost}}[X] = -1. \tag{5.1.47}
$$

In order to complete analysis of the structure of $\Gamma^{(n)}_{\mathbf{div}}$, we must supplement the symmetry equations (5.1.42), (5.1.43), and (5.1.47) with the constraints on the divergences which arise from power counting. Accordingly, we introduce the following notations:

$n_E$ = number of graviton vertices with two derivatives,

$n_G$ = number of antighost-graviton-ghost vertices,

$n_K$ = number of K-graviton-ghost vertices,

$n_L$ = number of L-ghost-ghost vertices,

$I_G$ = number of internal-ghost propagators,

$E_C$=number of external ghosts,

$E_{\overline{C}}$=number of external antighosts.

Since graviton propagators behave like $p^{-4}$, and ghost propagators like $p^{-2}$, we are led by standard power counting to the degree of divergence of an arbitrary diagram,

$$D = 4 - 2n_E + 2I_G - 2n_G - 3n_L - 3n_K - E_{\overline{C}}. \qquad (5.1.48)$$

The last term in (5.1.48) arises because each external antighost line carries with it a factor of external momentum. We can make use of the topological relation

$$2I_G - 2n_G = 2n_L + n_K - E_C - E_{\overline{C}} \qquad (5.1.49)$$

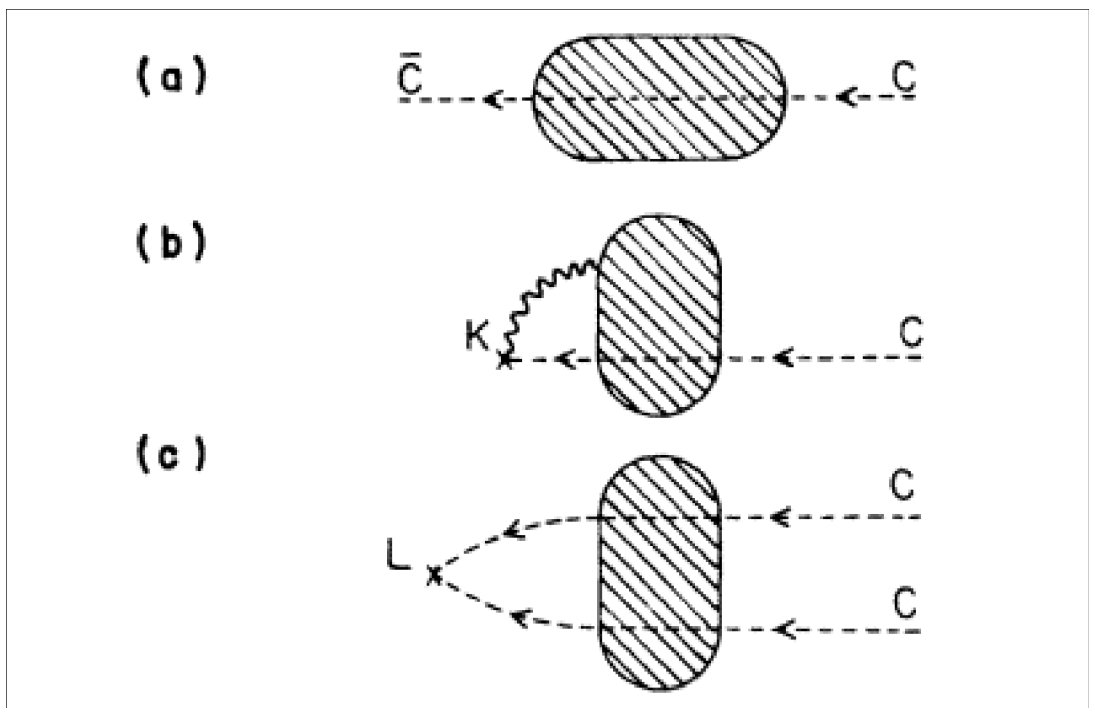

Fig.5.1.2.The three types of divergent diagram which involve external ghost lines. Arbitrarily many gravitons may emerge from each of the central regions,(a) Ghost action type,(b) $K$ type, (c) $L$ type.

to write the degree of divergence as

$$D = 4 - 2n_E - n_L - 2n_K - E_C - 2E_{\overline{C}}. \qquad (5.1.50)$$

Together with conservation of ghost number,Eq. (5.1.50) enables us to catalog three different types of divergent structures involving ghosts. These are illustrated in Fig.5.1.2. Each of the three types has degree of divergence $D = 1 - 2n_E$. Consequently, all the divergences which involve ghosts have $n_E = 0$. Since the degree of divergence is then $1$, the associated divergent structures in $\Gamma^{(n)}_{\mathbf{div}}$ have an extra derivative appearing on one of the fields. Diagrams whose external lines are all gravitons have degree of divergence $D = 4 - 2n_E$. Combining (5.1.50) with (5.1.47), (5.2.43), and (5.1.42), we can finally write the most general expression for $\Gamma^{(n)}_{\mathbf{div}}$ which satisfies all the constraints of symmetries and power counting:

$$\Gamma^{(n)}_{\mathbf{div}} = \Im(h^{\mu\nu}) + \Re\Big[ \Big(K_{\mu\nu} - \kappa^{-1}\overline{C}_\tau \vec{F}^{\tau}_{\mu\nu}\Big) P^{\mu\nu}(h^{\alpha\beta}) + L_\sigma Q^{\sigma}_{\tau}(h^{\alpha\beta}) C^\tau \Big], \qquad (5.1.51)$$

where $P^{\mu\nu}(h^{\alpha\beta})$ and $Q^{\sigma}_{\tau}(h^{\alpha\beta})$ are arbitrary Lorentz-covariant functions of the gravitational field $h^{\mu\nu}$, but not of its derivatives, at a single spacetime point. $\Im(h^{\mu\nu})$ is a local gauge-invariant functional of $h^{\mu\nu}$ containing terms with four, two, and zero derivatives.Expanding (5.1.51), we obtain an array of possible divergent structures:

$$\Gamma^{(n)}_{\mathbf{div}} = \Im(h^{\mu\nu}) + \frac{\delta I_{sym}}{\delta h^{\mu\nu}} P^{\mu\nu} + \left(\kappa K_{\rho\sigma} - \overline{C}_\tau \vec{F}^\tau_{\rho\sigma}\right)\left(\frac{\delta D^{\rho\sigma}_\alpha}{\delta h^{\mu\nu}} C^\alpha\right) P^{\mu\nu} -$$
$$-\left(\kappa K_{\rho\sigma} - \overline{C}_\tau \vec{F}^\tau_{\rho\sigma}\right) \frac{\delta P^{\rho\sigma}}{\delta h^{\mu\nu}} D^{\mu\nu}_\alpha C^\alpha - \left(\kappa K_{\mu\nu} - \overline{C}_\tau \vec{F}^\tau_{\mu\nu}\right) D^{\mu\nu}_\sigma (Q^\sigma_\epsilon C^\epsilon) - \kappa^2 L_\sigma \partial_\beta (Q^\sigma_\tau C^\tau) C^\beta \quad (5.1.52)$$
$$-\kappa^2 L_\sigma \partial_\beta C^\sigma Q^\beta_\tau C^\tau - \kappa L_\sigma \frac{\delta Q^\sigma_\tau}{\delta h^{\mu\nu}} C^\tau D^{\mu\nu}_\alpha C^\alpha + \kappa^2 L_\sigma Q^\sigma_\tau \partial_\beta C^\tau C^\beta.$$

The breakup between the gauge-invariant divergences S and the rest of (5.1.52) is determined only up to a term of the form [13]

$$\int d^4x (\eta^{\mu\nu} + \kappa h^{\mu\nu}) \frac{\delta I_{sym}}{\kappa \delta h^{\mu\nu}}, \quad (5.1.53)$$

which can be generated by adding to $P^{\mu\nu}$ a term proportional to $\eta^{\mu\nu} + \kappa h^{\mu\nu} = \sqrt{g}\, g^{\mu\nu}$.The profusion of divergences allowed by (5.1.52) appears to make the task of renormalizing the effective action rather complicated. Although the many divergent structures do pose a considerable nuisance for practical calculations, the situation is still reminiscent in principle of the renormalization of Yang-Mills theories. There, the non-gauge-invariant divergences may be eliminated by a number of field renormalizations. We shall find the same to be true here, but because the gravitational field $h^{\mu\nu}$ carries no weight in the power counting, there is nothing to prevent the field renormalizations from being nonlinear, or from mixing the gravitational and ghost fields. The corresponding renormalizations procedure considered in [13].

**Remark 5.1.2**.We assume now that:

(i) The local Poincaré group of momentum space is deformed at some fundamental high-energy cutoff $\Lambda_*$ [9],[10].

(ii) The canonical quadratic invariant $\|p\|^2 = \eta^{ab} p_a p_b$ collapses at high-energy cutoff $\Lambda_*$ and being replaced by the non-quadratic invariant:

$$\|p\|^2 = \frac{\eta^{ab} p_a p_b}{(1 + l_{\Lambda_*} p_0)}. \quad (5.1.54)$$

(iii) The canonical concept of Minkowski space-time collapses at a small distances $l_{\Lambda_*} = \Lambda_*^{-1}$ to fractal space-time with Hausdorff-Colombeau negative dimension and therefore the canonical Lebesgue measure $d^4x$ being replaced by the Colombeau-Stieltjes measure (see section III)

$$(d\eta(x,\varepsilon))_\varepsilon = (v_\varepsilon(s(x)) d^4x)_\varepsilon, \quad (5.1.55)$$

where

$$(v_\varepsilon(s(x)))_\varepsilon = (v_\varepsilon(x))_\varepsilon = \left(\left(|s(x)|^{|D^-|} + \varepsilon\right)^{-1}\right)_\varepsilon,$$
$$s(x) = \sqrt{x_\mu x^\mu}, D^- < 0, \quad (5.1.56)$$

see subsection IV.2.

(iv) The canonical concept of local momentum space collapses at fundamental high-energy cutoff $\Lambda_*$ to fractal momentum space with Hausdorff-Colombeau negative dimension and therefore the canonical Lebesgue measure $d^3\mathbf{k}$,where $\mathbf{k} = (k_x, k_y, k_z)$ being replaced by the Hausdorff-Colombeau measure

$$d^{D^+,D^-}\mathbf{k} \triangleq \frac{\Delta(D^-)d^{D^+}\mathbf{k}}{\left(|\mathbf{k}|^{|D^-|}+\varepsilon\right)_\varepsilon} = \frac{\Delta(D^+)\Delta(D^-)p^{D^+-1}dp}{(p^{|D^-|}+\varepsilon)_\varepsilon}, \tag{5.1.57}$$

see subsection 3.3-3.4. Note that the integral over measure $d^{D^+,D^-}\mathbf{k}$ is given by formula (3.3.16).

**Remark 5.1.3.**Note that by assumption (iii) mentioned above, the generic expression (2.2.1) of the action becomes to the following form

$$\begin{aligned}(I_{sym}(\varepsilon))_\varepsilon &= -\left(\int d^4x v_\varepsilon(x)\sqrt{-g}\,(\alpha R_{\mu\nu}R^{\mu\nu} - \beta R^2 + 2\kappa^{-2}R)\right)_\varepsilon \\ &= -\int d^4x(v_\varepsilon(x))_\varepsilon\sqrt{-g}\,(\alpha R_{\mu\nu}R^{\mu\nu} - \beta R^2 + 2\kappa^{-2}R)_\varepsilon.\end{aligned} \tag{5.1.58}$$

Corresponding Green's functions are then given by a generating functional

$$\begin{aligned}&(Z_\varepsilon(T_{\mu\nu}))_\varepsilon = N\int\left[\prod\nolimits_{\mu\leq\nu} dh^{\mu\nu}\right][dC^\sigma][d\bar{C}_\tau]\delta^4_\nu(F^\tau) \\ &\exp\left[i\left((I_{sym}(\varepsilon))_\varepsilon + \int((v_\varepsilon(x))_\varepsilon)d^4x\bar{C}_\tau\vec{F}^\tau_{\mu\nu}D^{\mu\nu}_\alpha C^\alpha + \kappa\int((v_\varepsilon(x))_\varepsilon)d^4xT_{\mu\nu}h^{\mu\nu}\right)\right].\end{aligned} \tag{5.1.59}$$

**Remark 5.1.4.(I)**The renormalizable models which we have considered in this section many years regarded only as constructs for a study of the ultraviolet problem of quantum gravity. The difficulties with unitarity appear to preclude their direct acceptability as canonical physical theories in locally Minkowski space-time. In canonical case they do have only some promise as phenomenological models.

(**II**) However, for their unphysical behavior may be restricted to *arbitrarily large energy scales* $\Lambda_*$ mentioned above by an appropriate limitation on the renormalized masses $m_2$ and $m_0$. Actually, it is only the massive spin-two excitations of the field which give the trouble with unitarity and thus require a very large mass. The limit on the mass $m_0$ is determined only by the observational constraints on the static field.

## 5.2.Cleaner methods

The renormalization procedure described in the last section is sufficiently complicated to make practical calculations unappealing. We now turn to other choices of the gauge-fixing term which greatly simplify matters by eliminating the need for the field and transformation renormalizations.

**A. Unweighted gauge condition**

Explicit calculations of samples of the nongauge-invariant divergences allowed by (6.17) reveal that they depend upon the gauge-fixing parameter $\Delta$ which was introduced into the effective action by the weighting functional (3.6). This suggests that if we take the limit $\Delta \to 0$, all the field and transformation renormalizations may disappear. This limit as $\Delta \to 0$ returns us to the unweighted gauge condition

$$\partial_\nu h^{\mu\nu} = 0 \tag{5.2.1}$$

with the same Feynman rules as those obtained using the simple Gaussian representation (3.4) of the gauge-fixing $\delta$-function.The graviton propagator in the limit $\Delta \to 0$ maybe calculated as suggested above, setting $\Delta = 0$ in the propagator calculated for finite $\Delta$ (cf. sec.6.4), or by substituting the gauge condition (8.1) into the linearized classical field equations and then inverting. The resulting propagator is constructed entirely from projectors which are transverse in all their indices:

$$D^{\Delta=0}_{\mu\nu\rho\sigma} = \frac{1}{(2\pi)^4 i}\left(\frac{2P^{(2)}_{\mu\nu\rho\sigma}(k)}{k^2(\alpha\kappa^2k^2+\gamma)} - \frac{2P^{(0-s)}_{\mu\nu\rho\sigma}(k)}{k^2[(3\beta-\alpha)]\kappa^2k^2+\frac{1}{2}\gamma}\right). \tag{5.2.2}$$

The definitions of the projectors $P^{(2)}$ and $P^{(0-s)}$ are given in sec.6.4.The antighost-graviton- ghost interaction is

$$V_{\overline{C}hC} = \partial_\sigma \overline{C}_\alpha \partial_\rho h^{\rho\sigma} C^\alpha + \partial_\sigma \overline{C}_\alpha \partial_\rho h^{\rho\sigma} C^\alpha + \partial^2_{\rho\sigma} \overline{C}_\alpha h^{\rho\sigma} C^\alpha. \tag{5.2.3}$$

The first two terms in this expression contain the gauge condition (5.2.1), and consequently do not connect to the graviton propagator (5.2.2). Similarly, integration by parts in the remaining term can be used to move the derivatives onto the ghost field $C^\alpha$. When these derivatives fall on $h^{\rho\sigma}$ they form the gauge condition (5.2.1) again,so we have effectively

$$V_{\overline{C}hC} \approx \partial^2_\sigma \overline{C}_\alpha \partial_\rho h^{\rho\sigma} C^\alpha \approx \overline{C}_\alpha \partial_\rho h^{\rho\sigma} \partial^2_\sigma C^\alpha. \tag{5.2.4}$$

The symbol ≈ is used to indicate that terms containing $\partial_\rho h^{\rho\sigma}$ hpo or $\partial^2_{\rho\sigma} h^{\rho\sigma}$ have been dropped, since they do not connect to the graviton propagator.

The power-counting rule given in Sec. 5.1 must be modified as a results of (5.2.4). In one-particle-irreducible (1PI) diagrams, there is a separate vertex $V_{\overline{C}hC}$ for each external ghost and antighost line. Consequently, each of these lines carries with it two factors of external momentum. The resulting degree of divergence of an arbitrary 1PI diagram is

$$D^{1PI}_{(\Delta=0)} = 4 - 2n_E - n_L - 2n_K - 3E_C - 3E_C. \tag{5.2.5}$$

This result would hold even if we had not chosen (2.2) as our definition of the gravitational field variable. However, the simple relation (5.2.4) is dependent upon that choice, which accords with the harmonic gauge condition (5.2.1). Otherwise there would be a complicated cancellation between vertices.

From the power-counting rule (5.2.5), we see that each of the three types of diagrams shown in Fig. 5.1.2 is now convergent: The ghost action type has $D^{1PI}_{(\Delta=0)} = -2 - 2n_E$, the $K$ type has $D^{1PI}_{(\Delta=0)} = -1 - 2n_E$ and the $L$ type has $D^{1PI}_{(\Delta=0)} = -3 - 2n_E$ Therefore, there are no parts of $\Gamma^{(n)}_{\text{div}(\Delta=0)}$ which depend upon ghosts:

$$(a)\ \frac{\delta\Gamma^{(n)}_{\text{div}}(\Delta=0)}{\delta C^\sigma} = 0, (b)\frac{\delta\Gamma^{(n)}_{\text{div}}(\Delta=0)}{\delta K_{\mu\nu}}, (c)\frac{\delta\Gamma^{(n)}_{\text{div}}(\Delta=0)}{\delta L_\sigma}. \tag{5.2.6}$$

Insertion of Eqs. (6.2.6) into the renormalization Eq. (5.2.3) yields

$$\frac{\delta\Sigma}{\delta K_{\mu\nu}}\frac{\delta\Gamma^{(n)}_{\text{div}}(\Delta=0)}{\delta h^{\mu\nu}} = 0. \tag{5.2.7}$$

Together with (6.2.6a), this implies that $\Gamma^{(n)}_{\text{div}}(\Delta=0)$ is gauge invariant. All the divergences may therefore be eliminated by renormalizations of the parameters $\alpha$, $\beta$ and $\gamma$ in $I_{sym}$ and by the addition of a cosmological counterterm. The field variables and the BRS transformations do not need to be renormalized.The contrast between the complicated renormalization procedure which one must use when the quantum theory is defined with the gauge -fixing term (3.7) and the much simpler procedure for the unweighted gauge condition is reminiscent of the situation in the axial gauge in Yang-Mills theory. There, the ghosts decouple entirely from the Yang-Mills fields if one uses the unweighted axial gauge condition. However, if one smears the axial gauge with a weighting functional, the resulting propagator does connect to the ghosts, and then there arise non-gauge-invariant divergences. These Yang-Mills divergences are similar to those we would have obtained in the gravitational theory had we kept the two-derivative gauge-fixing term derived from (3.5). In both cases, the part of the propagator which depends upon the gauge-fixing parameter has a bad asymptotic behavior for large momenta, leading to nongauge-invariant divergences of progressively higher order as the calculation proceeds in perturbation theory.

Taking the limit $\Delta \to 0$ is necessary for the axial gauge quantization of Yang-Mills theory to avoid these artifactual divergences. However, this limit is less useful in other gauges: Although one obtains an improvement in the power counting just as we have found for gravitation, the improvement is not sufficient to eliminate all the nongauge-invariant divergences, and one must still renormalize the Yang-Mills gauge transformation. Thus, although taking the limit $\Delta \to 0$ is perfectly acceptable in Yang-Mills theory, it is generally of no particular advantage, and has not been much used in the literature.

**B. Third-derivative gauge**

Since we are dealing with theories in which the classical field equations involve fourth derivatives, the Cauchy data which must be initially specified to determine the classical evolution of the field include the values of the field and up to its third derivatives on some spacelike hyper surface. Accordingly, we should also be prepared to use gauge conditions which involve up to third derivatives. A gauge condition of this type which has the same structure as the harmonic gauge condition (5.2.1) is

$$\kappa^2\Box^2\partial_\nu h^{\mu\nu} = 0. \tag{5.2.8}$$

If we weight the gauge condition (6.2.8) with the Gaussian functional (3.5), we get the gauge-fixing term

$$\tfrac{1}{2}\kappa^4\Delta^{-1}(\Box^2\partial_\nu h^{\mu\nu})(\Box^2\partial_\lambda h^\lambda_\mu). \tag{5.2.9}$$

Another way to arrive at (5.2.9) is to start from the usual harmonic gauge condition (5.2.1) and to weight it with the functional

$$\omega_6(e_\tau) = \exp\Big[ i\Big(\tfrac{1}{2}\kappa^4\Delta^{-1}\int(\Box^2 e_\tau)(\Box^2 e_\tau)\Big)\Big]. \tag{5.2.10}$$

When we obtain (5.2.9) this second way, it is clear that the ghost action which we must use is exactly the same that we had before in the generating functional (3,2). This also follows from the first method of arriving at (5.2.9), because we may always redefine the antighost field: $\Box^2\overline{C}_\tau \to \overline{C}_\tau$.

The gauge-fixing term (5.2.9) requires us to change the BRS transformation of the antighost field $\overline{C}_\tau$. The new transformation is

$$\delta_{(6)BRS}\overline{C}_\tau = \kappa^5\Delta^{-1}\Box^2\Box^2 F_\tau\delta\lambda. \tag{5.2.11}$$

The Slavnov identities for the generating functionals of Green's functions and of proper vertices must be changed too, but the identity for the reduced generating functional of proper vertices,

$$\Gamma_{(6)} = \widetilde{\Gamma}_{(6)} - \tfrac{1}{2}\kappa^4\Delta^{-1}(\Box^2 F_{\tau\mu\nu}h^{\mu\nu})(\Box^2 F^\tau_{\rho\sigma}h^{\rho\sigma}) \tag{5.2.12}$$

remains the same as (5.20). Consequently, the renormalization equation is the same as (6.3).The Feynman rules which we obtain using (5.2.9) differ from those obtained using (3.7) only in the replacement of the factors of $\Delta\kappa^{-2}k^{-4}$ in the graviton propagator by $\Delta\kappa^{-4}k^{-6}$. This change brings about a reduction in the degree of divergence of those parts of diagrams which depend on the parameter $\Delta$. The degree of divergence is reduced by 2 for each factor of $\Delta$, so that once again all three types of diagram involving ghosts shown in Fig. 6.1.2 are convergent. The renormalization equation then implies that all the divergences in $\Gamma^{(n)}_{(6)\,\mathrm{div}}$ are gauge invariant.

## 5.3.Coupling to fields of standard matter and to fields of physical ghost matter.

Now that we know how to carry out the renormalization procedure for a purely gravitational model, it is straightforward to include coupling to other renormalizable fields. As an example, we discuss a massive scalar field in interaction with the gravitational field alone, adding to the action (5.1.1) the additional term

$$\mathcal{L}_\varphi = \int d^4x\Big(-\frac{1}{2}\partial_\mu\varphi\partial_\nu g^{\mu\nu} - \frac{1}{2}m^2\varphi^2\Big)\sqrt{-g}\,. \qquad (5.3.1)$$

The BRS transformations must now include a transport term for the scalar field,

$$\delta_{BRS}\varphi = -\kappa^2\partial_\mu\varphi C^\mu\delta\lambda. \qquad (5.3.2)$$

This transformation is nilpotent:

$$\delta_{BRS}(\partial_\mu\varphi C^\mu) = 0. \qquad (5.3.3)$$

In order to write the Slavnov identities, we make use of (5.3.3) by adding a term coupling the scalar and ghost fields to a new anticommuting external field $B(x)$:

$$\Sigma = I_{sym} + I_\varphi + \Big(\kappa K_{\mu\nu} - \overline{C}_\tau \overleftarrow{F}^{\tau}_{\mu\nu}\Big)D^{\mu\nu}_\alpha C^\alpha + \kappa^2 L_\sigma\partial_\beta C^\sigma C^\beta - \kappa^2 B\partial_\mu\varphi C^\mu. \qquad (5.3.4)$$

In the generating functional of Green's functions, the scalar field is coupled to a source $J(x)$; the Legendre transformation then trades this dependence on $J(x)$ for a dependence on $\varphi(x) = \delta W/\delta J(x)$ in the generating functional of proper vertices. The Slavnov identity for the reduced generating functional of proper vertices reads

$$\frac{\delta\Gamma_\varphi}{\delta B}\frac{\delta\Gamma_\varphi}{\delta\varphi} + \frac{\delta\Gamma_\varphi}{\delta K_{\mu\nu}}\frac{\delta\Gamma_\varphi}{\delta h^{\mu\nu}} + \frac{\delta\Gamma_\varphi}{\delta L_\sigma}\frac{\delta\Gamma_\varphi}{\delta C^\sigma} = 0. \qquad (5.3.5)$$

As before, this identity leads to the renormalization equation for $\Gamma^{(n)}_{\varphi\,\mathrm{div}}$.Power counting, using the unweighted gauge condition, gives the degree of divergence of an arbitrary 1PI diagram,

$$D^{1PI}_\varphi(\Delta = 0) = 4 - 2n_E - n_L - 2n_K - n_B - 3E_C - 3E_{\overline{C}} - 2E_S, \qquad (5.3.6)$$

where $n_B$ is the number of $B$-scalar-ghost vertices and $E_S$ is the number of external scalar lines. The external scalar lines are counted twice in (5.3.6) because of the linkage of scalar fields and derivatives in the interaction between scalars and gravitons (the mass term is super-renormalizable and is not included in the power counting). This linkage is similar to the linkage of ghosts and derivatives which we have already encountered.The power- counting rule (5.3.6), together with the conservation of ghost number, shows that all 1PI diagrams with external ghost lines are convergent, so that

$$\frac{\delta\Sigma_\varphi}{\delta B}\frac{\delta\Gamma^{(n)}_{\varphi\,\mathrm{div}}}{\delta\varphi} + \frac{\delta\Sigma_\varphi}{\delta K_{\mu\nu}}\frac{\delta\Gamma^{(n)}_{\varphi\,\mathrm{div}}}{\delta h^{\mu\nu}} = 0 \qquad (5.3.7)$$

Consequently,$\Gamma^{(n)}_{\varphi\,\mathrm{div}}$ is gauge invariant. The only gauge-invariant structures consistent with (5.3.6) are

$$\Gamma^{(n)}_{\varphi\,\mathrm{div}} = \alpha^{(n)}\int R_{\mu\nu}R^{\mu\nu}\sqrt{-g} - \beta^{(n)}\int R^2\sqrt{-g} + \gamma^{(n)}\kappa^{-2}\int R\sqrt{-g} - \lambda^{(n)}\kappa^{-4}\int\sqrt{-g}$$
$$+\frac{1}{2}f^{(n)}\int\partial_\mu\varphi\partial_\nu g^{\mu\nu}\sqrt{-g} + \frac{1}{2}(\delta m^2)^{(n)}\int\varphi^2\sqrt{-g}. \quad (5.3.8)$$

These divergences may be eliminated by renormalizations of the appropriate coefficients in $I_{sym}$ and $I_\varphi$, and by the addition of a cosmological counterterm. It should be noted that the absense of a term like $\int R\varphi^2\sqrt{-g}$ in (6.3.8) is due to the linkage of scalars and derivatives. If this linkage were broken by the inclusion in (6.3.1) of a scalar self-interaction $\int\varphi^4\sqrt{-g}$, then it would be nec-essary to include as well the nonminimal gravitational-scalar interaction.

The scalar field example shows that once renormalizability has been established for a purely gravitational model, the inclusion of couplings to other renormalizable fields poses no further problems (except possibly the necessity for a nonminimal gravitational-scalar interaction). In particular, the Faddeev-Popov ghost machinery remains unrenormalized just as it did in the purely gravitational case.The allowed divergences may be summarized by assigning a power-counting weight to each field, and then requiring that divergent structures be gauge invariant and of power-counting weight four or less. It is necessary to take into account any linkages of fields and derivatives in the interactions by augmenting the weight of a field by the number of derivatives linked to it. The weight of the gravitational field is zero, and before linkages with derivatives are taken into account, the weights of other fields are simply given by their canonical dimensions.

# 5.4.The graviton propagator.

The inversion of the gravitational kinetic matrix which is necessary to calculate the graviton propagator involves a substantial amount of Lorentz algebra on symmetric rank-two tensors. To organize the calculation, it is convenient to use a set of orthogonal projectors in momentum space. We choose a set of projectors which emphasises transversality,16 since this is important in Sec. 5.3.These projectors are constructed using the transverse and longitudinal projectors for vector quantities,

$$(a)\ \theta_{\mu\nu} = \eta_{\mu\nu} - \frac{k_\mu k_\nu}{k^2}, (b)\ \omega_{\mu\nu} = \frac{k_\mu k_\nu}{k^2}. \quad (5.4.1)$$

The four projectors for symmetric rank-two tensors then reads

$$(a)\ P^{(2)}_{\mu\nu\rho\sigma} = \frac{1}{2}(\theta_{\mu\rho}\theta_{\nu\sigma} + \theta_{\mu\sigma}\theta_{\nu\rho}) - \frac{1}{3}\theta_{\mu\nu}\theta_{\rho\sigma},$$
$$(b)\ P^{(1)}_{\mu\nu\rho\sigma} = \frac{1}{2}(\theta_{\mu\rho}\omega_{\nu\sigma} + \theta_{\mu\sigma}\omega_{\nu\rho} + \theta_{\nu\rho}\omega_{\mu\sigma} + \theta_{\nu\sigma}\omega_{\mu\rho}), \quad (5.4.2)$$
$$(c)\ P^{(0-s)}_{\mu\nu\rho\sigma} = \frac{1}{3}\theta_{\mu\nu}\theta_{\rho\sigma}, (d)\ P^{(0-w)}_{\mu\nu\rho\sigma} = \omega_{\mu\nu}\omega_{\rho\sigma}.$$

For a massive tensor field in the rest frame, the projectors (5.2.2.a)-(5.2.2.d) select out the spin-two, spin-one, and two spin-zero parts of the field.However, the projectors (5.2.2) do not span the operator space of the gravitational field equations. In order to have a complete basis, we must also include the two spin-zero transfer operators,

$$(a)\ P^{(0-sw)}_{\mu\nu\rho\sigma} = \frac{1}{\sqrt{3}}\theta_{\mu\nu}\omega_{\rho\sigma}, (b)\ P^{(0-sw)}_{\mu\nu\rho\sigma} = \frac{1}{\sqrt{3}}\omega_{\mu\nu}\theta_{\rho\sigma} \quad (5.4.3)$$

The orthogonality relations of the projectors (5.2.2) and the transfer operators (5.2.3) are

$$(a)\ P^{(i-a)}P^{(j-b)} = \delta^{ij}\delta^{ab}P^{(j-b)},(b)\ P^{(i-ab)}P^{(j-cd)} = \delta^{ij}\delta^{bc}P^{(j-a)},$$
$$(c)\ P^{(i-a)}P^{(j-bc)} = \delta^{ij}\delta^{ab}P^{(j-ac)},(d)\ P^{(i-ab)}P^{(j-c)} = \delta^{ij}\delta^{bc}P^{(j-ac)}. \tag{5.4.4}$$

where $i$ and $j$ run from $0$ to $2$, and $a$ and $b$ take on the values $w$ and $s$.In order to calculate the graviton propagator, we must first write out the part of the effective action (5.1.14) which is purely quadratic in the gravitational field $h^{\mu\nu}$. Going over to momentum space and using (5.2.2) and (5.2.3), we get

$$\frac{1}{4}\int d^4k h^{\mu\nu}(-k)\Big(-(\alpha\kappa^2k^2+\gamma)k^2P^{(2)}_{\mu\nu\rho\sigma}(k)+\Delta^{-1}\kappa^2k^4P^{(1)}_{\mu\nu\rho\sigma}(k)+$$
$$\Big\{3k^2\Big[(3\beta-\alpha)\kappa^2k^2+\tfrac{1}{2}\gamma\Big]+2\Delta^{-1}\kappa^2k^4\Big\}P^{(0-w)}_{\mu\nu\rho\sigma}(k)+ \tag{5.4.5}$$
$$k^2\Big[(3\beta-\alpha)\kappa^2k^2+\tfrac{1}{2}\gamma\Big]\Big\{P^{(0-s)}_{\mu\nu\rho\sigma}(k)-\sqrt{3}\Big[P^{(0-ws)}_{\mu\nu\rho\sigma}(k)+P^{(0-sw)}_{\mu\nu\rho\sigma}(k)\Big]\Big\}\Big)h^{\rho\sigma}(k).$$

The combination of parameters $(3\beta-\alpha)$ which occurs throughout this expression is an echo of the conformally invariant action
$\int d^4x\sqrt{-g}\,(R_{\mu\nu}R^{\mu\nu}-1/3R^2) = 1/2\int d^4x\sqrt{-g}\,C_{\mu\nu\alpha\beta}C^{\mu\nu\alpha\beta}$ where $C_{\mu\nu\alpha\beta}$ is the Weyl tensor. The orthogonality relations (5.4.4) may now be used in inverting the kinetic matrix shown in (5.4.5) to obtain the graviton propagator:

$$D_{\mu\nu\rho\sigma}(k) = \frac{1}{(2\pi)^4i}\left(\frac{2P^{(2)}_{\mu\nu\rho\sigma}(k)}{k^2(\alpha\kappa^2k^2+\gamma)}-\frac{2P^{(0-s)}_{\mu\nu\rho\sigma}(k)}{k^2\big([3\beta-\alpha]\kappa^2k^2+\frac{1}{2}\gamma\big)}-\frac{2\Delta P^{(1)}_{\mu\nu\rho\sigma}(k)}{\kappa^2k^4}\right)$$
$$\frac{\Delta\Big(3P^{(0-s)}_{\mu\nu\rho\sigma}(k)-\sqrt{3}\Big[P^{(0-sw)}_{\mu\nu\rho\sigma}(k)+P^{(0-ws)}_{\mu\nu\rho\sigma}(k)\Big]+P^{(0-w)}_{\mu\nu\rho\sigma}(k)\Big)}{\kappa^2k^4} \tag{5.4.6}$$

To determine the propagator (5.4.6) completely, we must specify how the $k_0$ integration contour is to skirt the poles in calculating Feynman integrals. We do this in the customary way by including ie terms in the denominators of the individual poles, which must first be obtained by separating (5.4.6) into partial fractions. Ignoring for the moment the terms proportional to $\Delta$, we find

$$D_{\mu\nu\rho\sigma}(k) = \frac{1}{(2\pi)^4i}\times$$
$$\left(\frac{2\Big[P^{(2)}_{\mu\nu\rho\sigma}(k)-2P^{(0-s)}_{\mu\nu\rho\sigma}(k)\Big]}{\gamma k^2}-\frac{2P^{(2)}_{\mu\nu\rho\sigma}(k)}{\gamma\big(k^2+\gamma[\alpha\kappa^2]^{-1}\big)}+\frac{4P^{(0-s)}_{\mu\nu\rho\sigma}(k)}{\gamma\big(k^2+\gamma[2\kappa^2(3\beta-\alpha)]^{-1}\big)}\right). \tag{5.4.7}$$

Normally, one requires that quantum states have positive-definite norm and energy. Such states give rise to poles in the propagator with positive residues. Since both the massless pole and the pole at $k^2=-\gamma[2\kappa^2(3\beta-\alpha)]^{-1}$ in (5.4.7) do have positive residues, we shift them in the standard fashion, replacing the denominators respectively by

$$(k^2-i\epsilon) \tag{5.4.8}$$

and by

$$\big\{k^2+\gamma[2(3\beta-\alpha)\kappa^2]^{-1}-i\epsilon\big\}. \tag{5.4.9}$$

On the other hand, the negative residue of the massive spin-two pole at $k^2 = -\gamma[\alpha\kappa]^{-1}$ faces us with a choice between two unfortunate alternatives: to give up either the positive definiteness of the norm or of the energy of the corresponding quantum states.

Both choices give the required negative residue, but they differ in the way the pole must be shifted.If the massive spin-two states are taken to have negative norm, the situation is analogous to a Pauli-Villars regularized theory. We recall that in the usual derivation of the propagator, one starts from $\langle 0|T[h_{\mu\nu}(x)h_{\rho\sigma}(x')]|0\rangle$, transforms to momentum space, and sums over a complete set of momentum eigenstates inserted between the two field operators. The only difference in the present case is that the negative-norm states must be accompanied by a vector space metric factor of $(-1)$ in the sum over states. This gives rise to a negative residue for the massive spin-two pole, but does not affect the location of the pole, whose denominator is consequently given by

$$\left(k^2+\gamma[\alpha\kappa]^{-1}-i\epsilon\right). \tag{5.4.10}$$

As the Pauli-Villars analogy leads us to expect, the choice (5.4.10), together with (5.4.8) and (5.4.9), gives a high-energy behavior of the total propagator which is like $k^{-4}$. To see this, one may, for example, perform a Wick rotation into Euclidean space and then drop the re terms. This is allowedbecause (5.4.10), (5.4.8), and (5.4.9) all shift the poles in the same way. If the massive spin-two states are taken to have negative energy, the pole in the propagator acquires a negative residue for a different reason.In this case, there are no vector space metric factors in the sum over states, but the expansion of the field operators into creation and annihilation operators involves normalization factors $(2|k_0|)^{-1/2}=(-2k_0)^{-1/2}$. These contribute an overall minus sign to the massive spin-two part of the propagator. In addition, the sign of the energy flow for a given time ordering is opposite to that for a positive-energy field, so the denominator of the pole is now given by

$$\left(k^2+\gamma[\alpha\kappa]^{-1}+i\epsilon\right). \tag{5.4.11}$$

The difference between the poles given by (5.4.10) and (5.4.11) is a term proportional to $\delta\left(k^2+\gamma[\alpha\kappa]^{-1}\right)$. While the choice of (5.4.10) leads to the desired behavior, this additional term effectively spoils the high-energy behavior of (5.4.11). Thus, our power-counting requirements lead us to adopt an indefinite-metric state vector space, following the analogy to Pauli-Villars regularization.The pure $k^4$ terms in (5.4.6), proportional to $\Delta$, may be handled by confluence, replacing them by $\zeta^{-1}\left[(k^2-i\epsilon)^{-1}-(k^2+\zeta-i\epsilon)^{-1}\right]$, and then letting $\zeta\to 0$ at the end of the calculation.

# 6. Hausdorff-Colombeau measure and associated negative Hausdorff-Colombeau dimensions.

## 6.1.Fractional Integration in negative dimensions.

Let $\mu_H^{D^+}$ be a Hausdorff measure [33-34] and $X\subset\mathbb{R}^n, D^+<n$ is measurable set. Let $s(x)$ be a function $s:X\to\mathbb{R}$ such that is symmetric with respect to some centre $x_0\in X,$ i.e. $s(x)=$ constant for all $x$ satisfying $d(x,x_0)=r$ for arbitrary values of $r.$ Then the

integral in respect to Hausdorff measure over $n$-dimensional metric space $X$ is then given by [33]:

$$\int_X s(x)d\mu_H^{D^+}=\frac{2\pi^{D^+/2}}{\Gamma(D^+/2)}\int_0^\infty s(r)r^{D^+-1}dr. \tag{6.1.1}$$

The integral in RHS of the Eq.(3.1.1) is known in the theory of the Weyl fractional calculus where, the Weyl fractional integral $W^Df(x),$ is given by

$$W^{D^+}f(x) = \frac{1}{\Gamma(D^+)}\int_0^\infty (t-x)^{D^+-1}f(t)dt. \quad (6.1.2)$$

**Remark 6.1.1**. In order to extend the Weyl fractional integral (6.1.1) in negative dimensions we apply the Colombeau generalized functions [21-25] and Colombeau generalized numbers [23].Recall that Colombeau algebras $\mathcal{G}(\Omega)$ of the Colombeau generalized functions defined as follows [21-22].

Let $\Omega$ be an open subset of $\mathbb{R}^n$. Throughout this paper, for elements of the space $C^\infty(\Omega)^{(0,1]}$ of sequences of smooth functions indexed by $\varepsilon \in (0,1]$ we shall use the canonical notation $(u_\varepsilon)_\varepsilon$ so $u_\varepsilon \in C^\infty(\Omega),\ \varepsilon \in (0,1]$.

**Definition 6.1.1**.We set $\mathcal{G}(\Omega) = \mathcal{E}_M(\Omega)/\mathcal{N}(\Omega)$, where

$$\begin{gathered}\mathcal{E}_M(\Omega) = \left\{(u_\varepsilon)_\varepsilon \in C^\infty(\Omega)^{(0,1]} \middle| \forall K \subset\subset \Omega, \forall \alpha \in \mathbb{N}^n \exists p \in \mathbb{N} \text{ with} \right.\\ \left. \sup_{x\in K}|u_\varepsilon(x)| = O(\varepsilon^{-p}) \text{ as } \varepsilon \to 0\right\},\\ \mathcal{N}(\Omega) = \left\{(u_\varepsilon)_\varepsilon \in C^\infty(\Omega)^{(0,1]} \middle| \forall K \subset\subset \Omega, \forall \alpha \in \mathbb{N}^n \forall q \in \mathbb{N} \right.\\ \left. \sup_{x\in K}|u_\varepsilon(x)| = O(\varepsilon^{q}) \text{ as } \varepsilon \to 0\right\}.\end{gathered} \quad (6.1.3)$$

Notice that $\mathcal{G}(\Omega)$ is a differential algebra.Equivalence classes of sequences $(u_\varepsilon)_\varepsilon$ will be denoted by $\mathbf{cl}[(u_\varepsilon)_\varepsilon]$.is a differential algebra containing $D'(\Omega)$ as a linear subspace and $C^\infty(\Omega)$ as subalgebra.

**Definition 6.1.2**. Weyl fractional integral $(W_\varepsilon^{D^{--}}f(x))_\varepsilon$ in negative dimensions $D^- < 0$, $D^- \neq 0,-1,\ldots,-n,\ldots,n \in \mathbb{N}$ is given by

$$\begin{gathered}W^{D^-}f(x) = \frac{1}{\Gamma(D^-)}\left(\int_\varepsilon^\infty (t-x)^{D^--1}f(t)dt\right)_\varepsilon\\ \text{or}\\ \left(W_\varepsilon^{D^-}f(x)\right)_\varepsilon = \frac{1}{\Gamma(D^-)}\left(\int_0^\infty \frac{1}{\varepsilon+(t-x)^{|D^-|+1}}f(t)dt\right)_\varepsilon,\end{gathered} \quad (6.1.4)$$

where $\varepsilon \in (0,1]$ and $\int_0^\infty |f(t)dt| < \infty$.Note that $(W_\varepsilon^{D^{--}}f(x))_\varepsilon \in \mathcal{G}(\mathbb{R})$.Thus in order to obtain apropriate extension of the Weyl fractional integral $W^{D^+}f(x)$ on the negative dimensions $D^- < 0$ the notion of the Colombeau generalized functions is essentially importent.

**Remark 6.1.2**.Thus in negative dimensions from Eq.(6.1.1) we formally obtain

$$\left(\int_X s(x)d\mu_{HC,\varepsilon}^{D^{--}}\right)_\varepsilon = \frac{2\pi^{D^{--}/2}}{\Gamma(D^-/2)}\left(\int_0^\infty \frac{s(r)}{\varepsilon+r^{|D^-|+1}}dr\right)_\varepsilon = \left(I_\varepsilon^{D^-}\right)_\varepsilon, \quad (6.1.5)$$

where $\varepsilon \in (0,1]$ and $D^- \neq 0,-2,\ldots,-2n,\ldots,n \in \mathbb{N}$ and where $(\mu_{HC,\varepsilon}^{D^-})_\varepsilon$ is apropriate generalized Colombeau outer measure.Namely Hausdorff-Colombeau outer measure.

**Remark 6.1.3**. Note that: if $s(0) \neq 0$ the quantity $(I_\varepsilon^{D^+,D^{--}})_\varepsilon$ takes infinite large value in sense of Colombeau generalized numbers ,i.e., $(I_\varepsilon^{D^+,D^{--}})_\varepsilon =_{\widetilde{\mathbb{R}}} \widetilde{\infty}$ , see Definition 3.3.2 and Definition 3.3.3.

**Remark 6.1.4**.We apply throught this paper more general definition then definition (6.1.2):

$$\left(\int_X s(x)d\mu_{HC,\varepsilon}^{D^+,D^{--}}\right)_\varepsilon = \frac{4\pi^{D^+/2}\pi^{D^{--}/2}}{\Gamma(D^+/2)\Gamma(D^-/2)}\left(\int_0^\infty \frac{r^{D^+-1}s(r)}{\varepsilon+r^{|D^{--}|+1}}dr\right)_\varepsilon = (I_\varepsilon^{D^+,D_-})_\varepsilon, \quad (6.1.6)$$

where $\varepsilon \in (0,1]$ and $D^+ \geq 1,\ D^- \neq 0,-2,\ldots,-2n,\ldots,n \in \mathbb{N}$ and where $\left(\mu_{HC,\varepsilon}^{D^+,D^{---}}\right)_\varepsilon$ is apropriate generalized Colombeau outer measure.Namely Hausdorff-Colombeau outer measure. In subsection 3.3 we pointed out that there exist Colombeau generalized measure $\left(d\mu_{HC,\varepsilon}^{D^+,D^{--}}\right)_\varepsilon$ and therefore Eq.(6.1.4) gives apropriate extension of the Eq.(6.1.1) on the negative Hausdorff-Colombeau dimensions.

.

# 6.2.Hausdorff measure and associated positive Hausdorff dimension.Colombeau-Feynman path integral in $D^+ = 4$ from dimensional regularization.

Recall that the classical Hausdorff measure [33],[22] originate in Caratheodory's construction, which is defined as follows: for each metric space $X$, each set $F = \{E_i\}_{i\in\mathbb{N}}$ of subsets $E_i$ of $X$, and each positive function $\zeta^+(E)$, such that $0 \leq \zeta^+(E_i) \leq \infty$ whenever $E_i \in F$, a preliminary measure $\phi_\delta^+(E)$ can be constructed corresponding to $0 < \delta \leq +\infty$ , and then a final measure $\mu^+(E)$, as follows: for every subset $E \subset X$, the preliminary measure $\phi_\delta^+(E)$ is defined by

$$\phi_\delta^+(E) = \inf_{\{E_i\}_{i\in\mathbb{N}}} \left\{\sum_{i\in\mathbb{N}} \zeta^+(E_i) | E \subset \bigcup_{i\in\mathbb{N}} E_i, diam(E_i) \leq \delta\right\}. \qquad (6.2.1)$$

Since $\phi_{\delta_1}^+(E) \geq \phi_{\delta_2}^+(E)$ for $0 < \delta_1 < \delta_2 \leq +\infty$, the limit

$$\mu^+(E) = \lim_{\delta\to 0_+} \phi_\delta^+(E) = \sup_{\delta>0} \phi_\delta^+(E) \qquad (6.2.2)$$

exists for all $E \subset X$. In this context, $\mu^+(E)$ can be called the result of Caratheodory's construction from $\zeta^+(E)$ on $F$. $\phi_\delta^+(E)$ can be referred to as the size $\delta$ approximating positive measure. Let $\zeta^+(E_i,d^+)$ be for example

$$\zeta^+(E_i,d^+) = \Theta(d^+)[diam(E_i)]^{d^+}, 0 < \Theta(d^+), \qquad (6.2.3)$$

for non-empty subsets $E_i, i \in \mathbb{N}$ of $X$. Where $\Theta(d^+)$ is some geometrical factor,depends on the geometry of the sets $E_i$, used for covering. When $F$ is the set of all non-empty subsets of $X$, the resulting measure $\mu_H^+(E,d^+)$ is called the $d^+$-dimensional Hausdorff measure over $X$; in particular, when $F$ is the set of all (closed or open) balls in $X$,

$$\Theta(d^+) \triangleq \Omega(d^+) = \pi^{\frac{d^+}{2}}(2^{-d^+})\Gamma\left(1 + \frac{d^+}{2}\right). \qquad (6.2.4)$$

Consider a measurable metric space $(X,\mu_H(d)), X \subseteqq \mathbb{R}^n, d \in (-\infty,+\infty)$. The elements of $X$ are denoted by $x,y,z,\ldots,$ and represented by $n$-tuples of real numbers $x = (x_1,x_2,\ldots,x_n)$

The metric $d(x,y)$ is a function $d : X \times X \to R_+$ is defined in $n$ dimensions by

$$d(x,y) = \sum_{ij}[\delta_{ij}(y_i - x_i)(y_j - x_j)]^{1/2} \qquad (6.2.5)$$

and the diameter of a subset $E \subset X$ is defined by

$$diam(E) = \sup\{d(x,y)|x,y \in E\}. \qquad (6.2.6)$$

**Definition 6.2.1**. The Hausdorff measure $\mu_H^+(E,D^+)$ of a subset $E \subset X$ with the associated Hausdorff positive dimension $D^+ \in \mathbb{R}_+$ is defined by canonical way

$$\mu_H^+(E,D^+) = \lim_{\delta\to 0}\left[\inf_{\{E_i\}_{i\in\mathbb{N}}} \left\{\sum_{i\in\mathbb{N}} \zeta^+(E_i,D^+) | E \subset \bigcup_i E_i, \forall i(diam(E_i) < \delta)\right\}\right], \qquad (6.2.7)$$

$$D^+(E) = \sup\{d^+ \in \mathbb{R}_+ | d^+ > 0, \mu_H^+(E,d^+) = +\infty\}.$$

**Definition 6.2.2**. Remind that a function $f : X \to \mathbb{R}$ defined in a measurable space $(X,\Sigma,\mu)$, is called a simple function if there is a finite disjoint set of sets $\{E_1,,\ldots,E_n\}$ of measurable sets and a finite set $\{\alpha_1,,\ldots,\alpha_n\}$ of real numbers such that $f(x) = \alpha_i$ if $x \in E_i$ and $f(x) = 0$ if $x \notin E_i$. Thus $f(x) = \sum_{i=1}^{n} \alpha_i \chi_{E_i}(x)$, where $\chi_{E_i}(x)$ is the characteristic function of $E_i$. A simple function $f$ on a measurable space $(X,\Sigma,\mu)$ is integrable if $\mu(E_i) < \infty$ for every index $i$ for which $\alpha_i \neq 0$. The Lebesgue-Stieltjes integral of $f$ is defined by

$$\int f d\mu = \sum_{i=1}^{n} \alpha_i \mu(E_i). \tag{6.2.8}$$

A continuous function is a function for which $\lim_{x\to y} f(x) = f(y)$ whenever $\lim_{x\to y} d(x,y) = 0$.

The Lebesgue-Stieltjes integral over continuous functions can be defined as the limit of infinitesimal covering diameter: when $\{E_i\}_{i\in\mathbb{N}}$ is a disjoined covering and $x_i \in E_i$ by definition (6.2.2) one obtains

$$\int_X f(x) d\mu_H^+(x,D^+) = \lim_{diam(E_i)\to 0} \left[ \sum_{\cup E_i = X} f(x_i) \inf_{\{E_{ij}\} \text{ with } \cup_j E_{ij} \supset E_i} \sum_j \zeta^+(E_{ij}, D^+(E_{ij})) \right]. \tag{6.2.9}$$

From now on, $X$ is assumed metrically unbounded, i.e. for every $x \in X$ and $r > 0$ there exists a point $y$ such that $d(x,y) > r$. The standard assumption that $D^+$ is uniquely defined in all subsets $E$ of $X$ requires $X$ to be regular (homogeneous, uniform) with respect to the measure,i.e. $\mu_H^+(B_r(x),D^+) = \mu_H^+(B_r(y),D^+)$ for all elements $x,y \in X$ and (convex) balls $B_r(x)$ and $B_r(y)$ of the form $B_{r>0}(x) = \{z|d(x,z) \leq r; x,z \in X\}$. In the limit $diam(E_i) \to 0$, the infimum is satisfied by the requirement that the variation over all coverings $\{E_{ij}\}_{ij\in\mathbb{N}}$ is replaced by one single covering $E_i$, such that $\cup_j E_{ij} \to E_i \ni x_i$. Hence

$$\int_X f(x) d\mu_H^+(x,D^+) = \lim_{diam(E_i)\to 0} \sum_{\cup E_i = X} f(x_i) \zeta^+(E_i, D^+). \tag{6.2.10}$$

The range of integration $X$ may be parametrised by polar coordinates with $r = d(x,0)$ and angle $\Omega$. $\{E_{r_i,\Omega_i}\}_{i\in\mathbb{N}}$ can be thought of as spherically symmetric covering around a centre at the origin. In the limit, the function $\zeta^+(E_{r,\Omega},D^+)$ defined by Eq.(3.2.7) is given by

$$d\mu_H^+(x,D^+) = \lim_{diam(E_{r,\omega})\to 0} \zeta^+(E_{r,\Omega},D^+) = d\Omega^{D^+-1} r^{D^+-1} dr. \tag{6.2.11}$$

Let us assume now for simplification that $f(x) = f(\|x\|) = f(r)$ and $\lim_{r\to\infty} f(r) = 0$. The integral over a $D^+$-dimensional metric space $X$ is then given by

$$\int_X f(x) d\mu_H^+(x,D^+) = \int_X f(x) d^{D^+} x = \frac{2\pi^{\frac{D^+}{2}}}{\Gamma\left(1 + \frac{D^+}{2}\right)} \int_0^\infty f(r) r^{D^+-1} dr. \tag{6.2.12}$$

The integral defined in (6.2.12) satisfies the following conditions.

(i) Linearity:

$$\int_X [a_1 f_1(x) + a_2 f_2(x)] d\mu_H^+(x,D^+) = a_1 \int_X f_1(x) d\mu_H^+(x,D^+) + a_2 \int_X f_2(x) d\mu_H^+(x,D^+). \tag{6.2.13}$$

(ii) Translational invariance:

$$\int_X f(x + x_0) d\mu_H^+(x,D^+) = \int_X f(x) d\mu_H^+(x,D^+) \tag{6.2.14}$$

since $d\mu_H^+(x-x_0,D^+) = d\mu_H^+(x,D^+)$.

(iii) Scaling property:

$$\int_X f(ax)d\mu_H^+(x,D^+) = a^{-D^+}\int_X f(x)d\mu_H^+(x,D^+) \qquad (6.2.15)$$

since $d\mu_H^+(x/a,D^+) = a^{-D^+}d\mu_H^+(x,D^+)$.

(iv) The generalised $\delta^{D^+}(x)$ function:

The generalised $\delta^{D^+}(x)$ function for sets with non-integer Hausdorff dimension exists and can be defined by formula

$$\int_X f(x)\delta^{D^+}(x-x_0)d\mu_H^+(x,D^+) = f(x_0). \qquad (6.2.16)$$

(v) The Fourier-Stieltjes transform is given by the following definition

$$\begin{gathered} g(x) = (2\pi)^{-D^+}\int_{\hat{X}\subset\mathbb{R}_k^n}\hat{g}(k)\exp(ikx)d\mu_H^+(k,D^+),\\ \hat{g}(k) = \int_{X\subset\mathbb{R}_x^n} g(x)\exp(-ikx)d\mu_H^+(x,D^+). \end{gathered} \qquad (6.2.17)$$

(vi) The following equality holds

$$\int_{\hat{X}}\exp(ikx)d\mu_H^+(k,D^+) = (2\pi)^{D^+}\delta^{D^+}(x). \qquad (6.2.18)$$

## 6.2.1.Colombeau-Feynman path integral in $D^+ = 4$ from dimensional regularization via fractal spacetime.What is wrong with dimensional regularization via fractal spacetime.

In oreder to obtain Colombeau-Feynman path integral related to dimensional regularization let us consider a free scalar field with action in Hausdorff dimensions $D_\varepsilon^+ = D-\varepsilon,\ \varepsilon\in(0,1], D\in\mathbb{N}$ :

$$(\mathbf{S}_{0,\varepsilon})_\varepsilon = -\frac{1}{2}\left(\int_{X_\varepsilon} d\mu_H^+(x,D_\varepsilon^+)\,\phi_\varepsilon(x)\,P(\Box)\,\phi_\varepsilon(x)\right)_\varepsilon, \qquad (6.2.19)$$

where $\cup_{\varepsilon\in(0,1]}X_\varepsilon = \mathbb{R}^D, [(\phi_\varepsilon(x))_\varepsilon]\in\mathcal{G}(\mathbb{R}^D)$ and $P(\cdot)$ is a polinomial.

**Definition 6.2.3**. Assume that (i) $(\phi_\varepsilon(x))_\varepsilon\in\mathcal{G}(\mathbb{R}_x^D)$ and

(ii) there exist Colombeau generalized function $\left(\tilde{\phi}_\varepsilon(k)\right)_\varepsilon\in\mathcal{G}(\mathbb{R}_k^D)$ such that

$$\begin{gathered} (\phi_\varepsilon(x))_\varepsilon = \left((2\pi)^{-D_\varepsilon^+}\int_{\hat{X}_\varepsilon} d\mu_H^+(k,D_\varepsilon^+)\,\tilde{\phi}_\varepsilon(k)\,e^{ik\cdot x}\right)_\varepsilon =\\ ((2\pi)^{-D_\varepsilon^+})_\varepsilon\int_{\mathbb{R}_k^D}(d\mu_H^+(k,D_\varepsilon^+))_\varepsilon\left(\tilde{\phi}_\varepsilon(k)\right)_\varepsilon e^{ik\cdot x}, \end{gathered} \qquad (6.2.20)$$

and

$$\left(\tilde{\phi}_\varepsilon(k)\right)_\varepsilon = \left(\int_{X_\varepsilon} d\mu_H^+(x,D_\varepsilon^+)\,\phi_\varepsilon(x)e^{-ik\cdot x}\right)_\varepsilon = \int_{\mathbb{R}_x^D}(d\mu_H^+(x,D_\varepsilon^+)\,)_\varepsilon(\phi_\varepsilon(x))_\varepsilon e^{-ik\cdot x}. \qquad (6.2.21)$$

Then we will say that: (1) $\left(\tilde{\phi}_\varepsilon(k)\right)_\varepsilon$ is Colombeau Fourier–Stieltjes transform of the field

$(\phi_\varepsilon(x))_\varepsilon$ and abraviate

$$\left(\tilde{\phi}_\varepsilon(k)\right)_\varepsilon = (\mathcal{F}S)[(\phi_\varepsilon(x))_\varepsilon](k), \tag{6.2.22}$$

(2) $(\phi_\varepsilon(x))_\varepsilon$ is inverse Colombeau Fourier–Stieltjes transform of the field $\left(\tilde{\phi}_\varepsilon(k)\right)_\varepsilon$ and abraviate

$$(\phi_\varepsilon(x))_\varepsilon = (\mathcal{F}S)^{-1}\left[\left(\tilde{\phi}_\varepsilon(k)\right)_\varepsilon\right](x). \tag{6.2.23}$$

**Definition 6.2.4**. We will denote:

(i) the set of the Colombeau generalized functions $\left(\tilde{\phi}_\varepsilon(k)\right)_\varepsilon \in \mathcal{G}(\mathbb{R}^D_k)$ which is Colombeau Fourier–Stieltjes transform by $(\mathcal{F}S)\left[\mathcal{G}\left(\widehat{X}_\varepsilon\right)\right]$ or by $\mathcal{G}^{(\mathcal{F}S)}\left(\widehat{X}_\varepsilon\right)$

(ii) the set of the Colombeau generalized functions $(\phi_\varepsilon(x))_\varepsilon \in \mathcal{G}(\mathbb{R}^n_x)$ which is inverse Colombeau-Fourier–Stieltjes transform by $(\mathcal{F}S)^{-1}[\mathcal{G}_x(\mathbb{R}^D)]$ or $\mathcal{G}^{(\mathcal{F}S)^{-1}}(\mathbb{R}^D_x)$.

(iii) Note that we assume that in both cases (i) and (ii) the Eqs.(6.2.20)-(6.2.21) are satisfies.

**Remark 6.2.1**.Note that $\mathcal{G}^{(\mathcal{F}S)}(\mathbb{R}^D_k) \times \mathcal{G}^{(\mathcal{F}S)^{-1}}(\mathbb{R}^D_x)$ is the linear topological subspace of Colombeau algebra $\mathcal{G}(\mathbb{R}^D_k) \times \mathcal{G}(\mathbb{R}^D_x) \supsetneqq \mathcal{G}^{(\mathcal{F}S)}(\mathbb{R}^D_k) \times \mathcal{G}^{(\mathcal{F}S)^{-1}}(\mathbb{R}^D_x)$.

The free partition function $(\mathbf{Z}_{0,\varepsilon}[J_\varepsilon])_\varepsilon$ in the presence of a local source $(J_\varepsilon(x))_\varepsilon \in \mathcal{G}(\mathbb{R}^D_x)$ is

$$\begin{gathered}(\mathbf{Z}_{0,\varepsilon}[J_\varepsilon])_\varepsilon = \\ \left(\int_{(\phi_\varepsilon)_\varepsilon \in \mathcal{G}(\mathbb{R}^D_x)} [\mathcal{D}(\phi_\varepsilon)] \exp\left\{i[\mathbf{S}_{0,\varepsilon} + \int_{X_\varepsilon} d\mu^+_H(x, D^+_\varepsilon)\, J_\varepsilon(x)\phi_\varepsilon(x)]\right\}\right)_\varepsilon \triangleq \\ \int_{(\phi_\varepsilon)_\varepsilon \in \mathcal{G}^{(\mathcal{F}S)^{-1}}(\mathbb{R}^D_x)} [\mathcal{D}(\phi_\varepsilon)_\varepsilon] \exp\left\{[(\mathbf{S}_{0,\varepsilon})_\varepsilon + \int_{X_\varepsilon} \varepsilon(d\mu^+_H(x, D^+_\varepsilon))_\varepsilon (J_\varepsilon(x))_\varepsilon (\phi_\varepsilon(x))_\varepsilon]\right\} \triangleq \\ \int_{(\phi_\varepsilon)_\varepsilon \in \mathcal{G}^{(\mathcal{F}S)^{-1}}(\mathbb{R}^D_x)} [\mathcal{D}(\phi_\varepsilon)_\varepsilon]\, e^{i(\mathbf{S}_{J_\varepsilon})_\varepsilon}.\end{gathered} \tag{6.2.24}$$

From Eq.(6.2.16) we obtain

$$\left(\int_{X_\varepsilon} f_\varepsilon(x)\delta^{D^+_\varepsilon}(x - x_0)d\mu^+_H(x, D^+_\varepsilon)\right)_\varepsilon = (f_\varepsilon(x_0))_\varepsilon. \tag{6.2.25}$$

From Eq.(6.2.18) we obtain

$$\left(\int_{\widehat{X}_\varepsilon} \exp(ikx)d\mu^+_H(k, D^+_\varepsilon)\right)_\varepsilon = \left((2\pi)^{D^+_\varepsilon}\delta^{D^+_\varepsilon}(x)\right)_\varepsilon. \tag{6.2.26}$$

From Eq.(6.2.24) and Eqs.(6.2.22)-(6.2.23) and Eqs.(6.2.25)-(6.2.26) we obtain

$$
\begin{aligned}
&(\mathbf{S}_{J_\varepsilon})_\varepsilon = \\
&\frac{1}{2}\int_{X_\varepsilon}(d\mu_H^+(x,D_\varepsilon^+))_\varepsilon \times \\
&\left\{\int_{\widehat{X}_\varepsilon}\frac{(d\mu_H^+(k_1,D_\varepsilon^+))_\varepsilon}{(2\pi)^{D_\varepsilon^+}}\int_{\widehat{X}_\varepsilon}\frac{(d\mu_H^+(k_1,D_\varepsilon^+))_\varepsilon}{(2\pi)^{D_\varepsilon^+}}e^{i(k_1+k_2)\cdot x}\times\right. \\
&\left[-\left(\left(\widetilde{\phi}_\varepsilon(k_1)\right)_\varepsilon\right)\left((f_\varepsilon(-k_2^2))_\varepsilon\right)\left(\left(\widetilde{\phi}_\varepsilon(k_2)\right)_\varepsilon\right)\right. \\
&\left.\left.+\left(\widetilde{J}_\varepsilon(k_1)\right)_\varepsilon\left(\widetilde{\phi}_\varepsilon(k_2)\right)_\varepsilon+\left(\widetilde{J}_\varepsilon(k_2)\right)_\varepsilon\left(\widetilde{\phi}_\varepsilon(k_1)\right)_\varepsilon\right]\right\} = \\
&=\frac{1}{2}\int_{\widehat{X}_\varepsilon}\frac{(d\mu_H^+(-k,D_\varepsilon^+))_\varepsilon}{(2\pi)^{D_\varepsilon^+}}\left[-\left(\left(\widetilde{\phi}_\varepsilon(-k)\right)_\varepsilon\right)((f_\varepsilon(-k^2))_\varepsilon)\left(\left(\widetilde{\phi}_\varepsilon(k)\right)_\varepsilon\right)+\right. \\
&\left.\left(\left(\widetilde{J}_\varepsilon(-k)\right)_\varepsilon\right)\left(\widetilde{\phi}_\varepsilon(k)\right)_\varepsilon+\left(\widetilde{J}_\varepsilon(k)\right)_\varepsilon\left(\left(\widetilde{\phi}_\varepsilon(-k)\right)_\varepsilon\right)_\varepsilon\right] \\
&=\frac{1}{2}\int_{\widehat{X}_\varepsilon}\frac{(d\mu_H^+(-k,D_\varepsilon^+))_\varepsilon}{(2\pi)^{D_\varepsilon^+}}\times \\
&\left[-\left((\widetilde{\varphi}_\varepsilon(-k))_\varepsilon\right)\left(\left(f(-k^2)\right)_\varepsilon\right)\left((\widetilde{\varphi}_\varepsilon(k))_\varepsilon\right)+\frac{\left(\left(\widetilde{J}_\varepsilon(-k)\right)_\varepsilon\right)\left(\left(\widetilde{J}_\varepsilon(k)\right)_\varepsilon\right)}{(f_\varepsilon(-k^2))_\varepsilon}\right], \\
&(\widetilde{\varphi}_\varepsilon)_\varepsilon=\left(\widetilde{\phi}_\varepsilon(k)\right)_\varepsilon-\frac{\left(\widetilde{J}_\varepsilon(-k)\right)_\varepsilon}{(f_\varepsilon(-k^2))_\varepsilon},(f_\varepsilon(-k_2^2))_\varepsilon=(P(-k_2^2)+i\epsilon)_\varepsilon.
\end{aligned}
\tag{6.2.27}
$$

Thus Eq.(6.2.24) becomes

$$
\begin{aligned}
&(Z_{0,\varepsilon}[J_\varepsilon])_\varepsilon = \\
&\left\{\int_{(\phi_\varepsilon)_\varepsilon\in\mathcal{G}^{(\mathcal{FS})^{-1}_\varrho}(\mathbb{R}^D_x)}[\mathcal{D}(\phi_\varepsilon)_\varepsilon]\times\right. \\
&\left.\exp\left[-\frac{i}{2}\int_{\widehat{X}_\varepsilon}\frac{(d\mu_H^+(k,D_\varepsilon^+))_\varepsilon}{(2\pi)^{D_\varepsilon^+}}\left(\left(\widetilde{\phi}_\varepsilon(-k)\right)_\varepsilon\right)((f_\varepsilon(-k^2))_\varepsilon)\left(\left(\widetilde{\phi}_\varepsilon(k)\right)_\varepsilon\right)\right]\right\}\times \\
&\exp\left[\frac{i}{2}\int_{\widehat{X}_\varepsilon}\frac{(d\mu_H^+(-k,D_\varepsilon^+))_\varepsilon}{(2\pi)^{D_\varepsilon^+}}\frac{\left(\left(\widetilde{J}_\varepsilon(-k)\right)_\varepsilon\right)\left(\left(\widetilde{J}_\varepsilon(k)\right)_\varepsilon\right)}{(f_\varepsilon(-k^2))_\varepsilon}\right] \\
&=Z_0[0]\exp\left[\frac{i}{2}\int_{\widehat{X}_\varepsilon}\frac{(d\mu_H^+(-k,D_\varepsilon^+))_\varepsilon}{(2\pi)^{D_\varepsilon^+}}\frac{\left(\left(\widetilde{J}_\varepsilon(-k)\right)_\varepsilon\right)\left(\left(\widetilde{J}_\varepsilon(k)\right)_\varepsilon\right)}{(f_\varepsilon(-k^2))_\varepsilon}\right].
\end{aligned}
\tag{6.2.28}
$$

Therefore the exponent in Eq.(6.2.28) can be written as

$$
\begin{aligned}
&\int_{\widehat{X}_\varepsilon}\frac{(d\mu_H^+(-k,D_\varepsilon^+))_\varepsilon}{(2\pi)^{D_\varepsilon^+}}\frac{\left(\left(\widetilde{J}_\varepsilon(-k)\right)_\varepsilon\right)\left(\left(\widetilde{J}_\varepsilon(k)\right)_\varepsilon\right)}{(f_\varepsilon(-k^2))_\varepsilon} = \\
&\int_{\widehat{X}_\varepsilon}\frac{(d\mu_H^+(-k,D_\varepsilon^+))_\varepsilon}{(2\pi)^{D_\varepsilon^+}}\times \\
&\int_{X_\varepsilon}((d\mu_H^+(x,D_\varepsilon^+))_\varepsilon)\int_{X_\varepsilon}(d\mu_H^+(y,D_\varepsilon^+)_\varepsilon)e^{ik\cdot(x-y)}\frac{\left(\left(\widetilde{J}_\varepsilon(-k)\right)_\varepsilon\right)\left(\left(\widetilde{J}_\varepsilon(k)\right)_\varepsilon\right)}{(f_\varepsilon(-k^2))_\varepsilon},
\end{aligned}
\tag{6.2.29}
$$

so that, if $(\varrho_\varepsilon(-k))_\varepsilon = (d\varrho_\varepsilon(k))_\varepsilon$, the free partition function reads

$$(Z_{0,\varepsilon}[J])_\varepsilon =$$
$$((Z_{0,\varepsilon}[0])_\varepsilon) \times \tag{6.2.30}$$
$$\exp\left[\frac{i}{2}\int_{\hat{X}_\varepsilon}[(d\mu_H^+(x,D_\varepsilon^+))_\varepsilon]\int_{\hat{X}_\varepsilon}[(d\mu_H^+(x,D_\varepsilon^+))_\varepsilon](J_\varepsilon(x)((G(x-y;\varepsilon))_\varepsilon)J_\varepsilon(y))_\varepsilon\right]$$

where

$$(G(x-y;\varepsilon))_\varepsilon = \frac{1}{(2\pi)^{D_\varepsilon^+}}\int_{\hat{X}_\varepsilon}[(d\mu_H^+(-k,D_\varepsilon^+))_\varepsilon]\frac{\exp[ik\cdot(x-y)]}{(f_\varepsilon(-k^2))_\varepsilon}. \tag{6.2.31}$$

Thus, we have recovered the usual definition of the propagator as the solution of the Green equation in Hausdorff dimensions $D_\varepsilon^+ = D-\varepsilon, \varepsilon \in (0,1]$ in the sense of Colombeau generalized functions

$$(P(\Box)\,G(x-y;\varepsilon))_\varepsilon = (\delta^{D_\varepsilon^+}(x-y))_\varepsilon. \tag{6.2.32}$$

Defining the perturbative quantum field theory in momentum space makes a derivation of covariant Feynman rules in $D_\varepsilon^+$-dimensional momentum space straightforward [20].The only difference to conventional Feynman rules is the substitution of the measure in the momentum integral

$$(2\pi)^{-4}d^4k \to (2\pi)^{-D_\varepsilon^+}d\mu_H^+(k,D_\varepsilon^+). \tag{6.2.33}$$

For symmetric kernels, a representation of $d\mu_H^+(k,D^+)$ in terms of spherical coordinates is

$$\left((2\pi)^{-D_\varepsilon^+}d\mu_H^+(k,D_\varepsilon^+)\right)_\varepsilon = \left((2\pi)^{-D_\varepsilon^+}d\Omega^{D_\varepsilon^+-1}k^{D_\varepsilon^+-1}dk\right)_\varepsilon. \tag{6.2.34}$$

In the following the electron self-energy $(\Sigma_\varepsilon)_\varepsilon$, the vacuum polarisation $(\Pi_\varepsilon)_\varepsilon$ and the vertex correction $\Lambda_\varepsilon$ for each $\varepsilon \in (0,1]$ are enumerated as a function of the Hausdorff dimension $D_\varepsilon^+ = D-\varepsilon$. The lowest-order contribution to the vacuum polarisation

$$(\Pi_{\mu\nu,\varepsilon}(q))_\varepsilon = -e^2\left(\mathbf{Tr}\int\frac{d^{D_\varepsilon^+}k}{(2\pi)^{D_\varepsilon^+}}\left(\gamma_\mu\frac{i}{\gamma k-m+i\epsilon}\gamma_\nu\frac{i}{\gamma k-\gamma q-m+i\epsilon}\right)\right)_\varepsilon \tag{6.2.35}$$

can be written in the following form

$$\begin{gathered}\Pi_{\mu\nu,\varepsilon}(q) = (q_\mu q_\nu - q^2 g_{\mu\nu})(\Pi_\varepsilon(q))_\varepsilon,\\ (\Pi_\varepsilon(q))_\varepsilon = \alpha(2^{4-D_\varepsilon^+}\pi^{1-D_\varepsilon^+/2}\Gamma(2-D_\varepsilon^+/2)m^{D_\varepsilon^+-4}F(2-D_\varepsilon^+/2,2,5/2;-q^2/4m^2))_\varepsilon\\ F(a,b,c;z) = 1+\frac{ab}{c}\frac{z}{1!}+\ldots+\frac{z^{n-1}}{(n-1)!}\prod_{i=0}^{n-1}\frac{(a+i)(b+i)}{c+i}+\ldots.\end{gathered} \tag{6.2.36}$$

The lowest-order contribution to the electron self-energy

$$(\Sigma_\varepsilon(p))_\varepsilon = -ie^2\left(\int\frac{d^{D_\varepsilon^+}k}{(2\pi)^{D_\varepsilon^+}}\left(\frac{-i}{k^2-\lambda^2+i\epsilon}\gamma^\mu\frac{i}{\gamma p-\gamma k-m+i\epsilon}\gamma_\mu\right)\right)_\varepsilon \tag{6.2.37}$$

can be written in the following form

$$(\Sigma_\varepsilon(p))_\varepsilon = (A_\varepsilon)_\varepsilon - (\gamma p - m)(B_\varepsilon)_\varepsilon + (\gamma p-m)^2(\sigma_\varepsilon(p))_\varepsilon, \tag{6.2.38}$$

where $(A_\varepsilon)_\varepsilon, (B_\varepsilon)_\varepsilon \in \widetilde{\mathbb{R}}$ and $(\sigma_\varepsilon(p))_\varepsilon \in \mathcal{G}(\mathbb{R}_p^4)$ are given by

$$(A_\varepsilon)_\varepsilon = -3\alpha(2^{2-D^+_\varepsilon}\pi^{1-D^+_\varepsilon/2}m^{D^+_\varepsilon-3}\Gamma(2-D^+_\varepsilon/2))_\varepsilon,$$
$$(B_\varepsilon)_\varepsilon = -(A_\varepsilon)_\varepsilon m^{-1}, \tag{6.2.39}$$
$$(\sigma_\varepsilon(p))_\varepsilon \sim -\alpha(p^2+m^2)^2(2^{2-D^+_\varepsilon}\pi^{1-D^+_\varepsilon/2}\Gamma(2-D^+_\varepsilon/2))_\varepsilon.$$

The lowest-order contribution to the vertex term $(\Lambda_\varepsilon)_\varepsilon$, with the photon momentum $q$ and two outgoing electron momenta $p$ and $p'$, is given by

$$(\Lambda_\varepsilon(q,p,p'))_\varepsilon = -e^2\left(\int\frac{d^{D^+_\varepsilon}k}{(2\pi)^{D^+_\varepsilon}}\left(\frac{-i}{k^2-\lambda^2+i\epsilon}\gamma_\nu\frac{i}{\gamma p-\gamma k-m+i\epsilon}\gamma_\mu\frac{i}{\gamma p'-\gamma k-m+i\epsilon}\gamma^\nu\right)\right)_\varepsilon \tag{6.2.40}$$

can for $q = p'-p,$ be written in the following form

$$(\Lambda_\varepsilon(q,p,p'))_\varepsilon = [(B_\varepsilon)_\varepsilon + (g_\varepsilon(q))_\varepsilon]\gamma_\mu + \alpha\frac{i\sigma_{\mu\nu}q^\nu}{2m}(2^{3-D^+_\varepsilon}\pi^{1-D^+_\varepsilon/2D^+_\varepsilon-3}\Gamma(3-D^+_\varepsilon/2))_\varepsilon, \tag{6.2.41}$$

where $(g_\varepsilon(q))_\varepsilon$ is a function proportional to $(\Gamma(3-D^+_\varepsilon/2))_\varepsilon$ vanishing for $q^2\to 0,$ which will not be enumerated here. The term proportional to $\sigma_{\mu\nu}q^\nu$yields contributions to the anomalous magnetic moment and to the $l\neq 0$ Lamb shift.

**Remark 6.2.2.**Note that in this subsection the dimensional $\varepsilon$-regularization $\varepsilon\in(0,1]$ of quantum electrodynamics is considered as QFT with a fractal support of the quantum fields. These quantum fields are well defined as Colombeau quantum fields [22] for Hausdorff dimensions $D^+_\varepsilon$ arbitrarily close but smaller than $D^+ = 4.$

Let us consider the dimensional renormalisation of the bare two-point Green function $(S^{Bare}_\varepsilon)_\varepsilon$ of the electron. The full propagator $(S_\varepsilon)_\varepsilon$ can be formally written as the analytic continuation of a sum over self-energy diagrams $(S_\varepsilon(p))_\varepsilon = (\gamma p-(m^{Bare}_\varepsilon)_\varepsilon-(\Sigma_\varepsilon(p))_\varepsilon+i\epsilon)^{-1}$, where $(m^{Bare}_\varepsilon)_\varepsilon$ is the bare electron mass and $(\Sigma_\varepsilon(p))_\varepsilon$ is the proper self-energy. Substituting for $(\Sigma_\varepsilon(p))_\varepsilon$ its lowest-order contribution Eq.(6.2.37), and recalling Eq.(6.2.38), yields

$$(S_\varepsilon(p))_\varepsilon = \frac{(Z_{2,\varepsilon})_\varepsilon}{(\gamma p-m+i\epsilon)}[1+(Z_{2,\varepsilon})_\varepsilon(\gamma p-m)(\sigma_\varepsilon(p))_\varepsilon]^{-1}, \tag{6.2.42}$$

where the physical mass $m\in\mathbb{R}_+$ and the renormalisation constant $(Z_{2,\varepsilon})_\varepsilon\in\widetilde{\mathbb{R}}$ are defined by

$$\begin{gathered}m = (m^{Bare}_\varepsilon)_\varepsilon-(A_\varepsilon)_\varepsilon = \\ (m^{Bare}_\varepsilon)_\varepsilon+3\alpha(2^{2-D^+_\varepsilon}\pi^{1-D^+_\varepsilon/2}m^{D^+_\varepsilon-3}\Gamma(2-D^+_\varepsilon/2))_\varepsilon = \\ (m^{Bare}_\varepsilon)_\varepsilon+3\alpha\left(2^{-2+\varepsilon}\pi^{-1+\varepsilon/2}m^{1-\varepsilon}\Gamma\left(\frac{\varepsilon}{2}\right)\right)_\varepsilon \simeq \\ (m^{Bare}_\varepsilon)_\varepsilon+\frac{3\alpha m}{4\pi}\left(\left(\frac{2}{\varepsilon}\right)_\varepsilon-0.577\right)\end{gathered} \tag{6.2.43}$$

and

$$(Z_{2,\varepsilon})_\varepsilon = 1-(A_\varepsilon)_\varepsilon m^{-1}. \tag{6.2.44}$$

From Eq.(6.2.43) one obtains

$$(m_\varepsilon^{Bare})_\varepsilon = m - \frac{3\alpha m}{4\pi}\left(\frac{2}{\varepsilon} - 0.577\right) =$$
$$m\left(1 - \frac{3\alpha}{2\pi}\frac{1}{\varepsilon} + \frac{3\alpha}{4\pi}0.577\right) \simeq m\left(1 - \frac{3\alpha}{2\pi}\left(\frac{1}{\varepsilon}\right)_\varepsilon\right). \qquad (6.2.45)$$

Thus

$$(m_\varepsilon^{Bare})_\varepsilon = m(Z_\varepsilon)_\varepsilon = m\left(1 - \frac{3\alpha}{2\pi}\left(\frac{1}{\varepsilon}\right)_\varepsilon\right), \qquad (6.2.46)$$ where

$$(Z_\varepsilon)_\varepsilon = 1 - \frac{3\alpha}{2\pi}\left(\frac{1}{\varepsilon}\right)_\varepsilon. \qquad (6.2.47)$$

**Remark 6.2.3**. Note that standard sector of $QED_4$ defined by the condition

$$(Z_\varepsilon(\alpha(\varepsilon)))_\varepsilon > 0, \qquad (6.2.48)$$

see Sec. IV.3. In particular for a given $\alpha$ standard sector of $QED_4$ defined by the following
condition

$$\frac{3\alpha}{2\pi\varepsilon} < 1. \qquad (6.2.49)$$

The bare photon propagator is renormalised by the formal summation of vacuum polarisation diagrams, whose lowest-order contribution given by Eq.(6.2.36).Again, $(\Pi_\varepsilon(q))_\varepsilon$ can be expanded around the mass shell $q^2 = 0$, yielding

$$(\Pi_\varepsilon(q))_\varepsilon = (P_\varepsilon)_\varepsilon + q^2(\pi_\varepsilon(q^2))_\varepsilon. \qquad (6.2.50)$$

The full photon propagator can be written as $(q_\mu q_\nu - q^2 g_{\mu\nu})(\Delta_\varepsilon(q^2))_\varepsilon$, with

$$(\Delta_\varepsilon(q^2))_\varepsilon = [1 - (P_\varepsilon)_\varepsilon - q^2(\pi_\varepsilon(q^2))_\varepsilon]^{-1} \qquad (6.2.51)$$

The term in brackets contributes to the renormalisation of the bare charge $(e_\varepsilon^{Bare})_\varepsilon$ , which relates to the renormalised charge $e$ by

$$(e_\varepsilon^{Bare})_\varepsilon = e[1 - (P_\varepsilon)_\varepsilon - q^2(\pi_\varepsilon(q^2))_\varepsilon] = e(Z_{3,\varepsilon})_\varepsilon$$
$$(Z_{3,\varepsilon})_\varepsilon = 1 - (P_\varepsilon)_\varepsilon - q^2(\pi_\varepsilon(q^2))_\varepsilon \qquad (6.2.52)$$

For zero momentum transfer $q^2 \to 0$ and for $\varepsilon \ll 1$ $(Z_{3,\varepsilon})_\varepsilon = 1 - (P_\varepsilon)_\varepsilon$ reduces to

$$(Z_{3,\varepsilon})_\varepsilon = 1 - (\Pi_\varepsilon(q^2 = 0))_\varepsilon =$$
$$1 - \alpha(2^{4-D_\varepsilon^\pm}\pi^{1-D_\varepsilon^\pm/2}\Gamma(2 - D_\varepsilon^\pm/2)m^{D_\varepsilon^\pm-4})_\varepsilon \qquad (6.2.53)$$

yielding

$$(\alpha_\varepsilon^{Bare})_\varepsilon = \frac{[(e_\varepsilon^{Bare})_\varepsilon]^2}{4\pi} = (\alpha(\varepsilon)Z_{3,\varepsilon}^{-2}(\alpha(\varepsilon)))_\varepsilon. \qquad (6.2.54)$$

**Remark 6.2.4**.Note that standard sector of $QED_4$ except the condition (6.2.48),defined by
the condition $(Z_{3,\varepsilon}^{-2}(\alpha(\varepsilon)))_\varepsilon > 0$,see Sec. IV.3. In particular for a given $\alpha$ standard sector of
$QED_4$ except the condition (6.2.49),defined by the following condition

$$Z_{3,\varepsilon}^{-2}(\alpha) > 0. \tag{6.2.55}$$

All other contributions to the renormalisation of the electric charge cancel, as can be explicitly seen by a summation of the lowest-order radiative corrections to the charge

$$(1 - (B_\varepsilon)_\varepsilon + (L_\varepsilon)_\varepsilon - (P_\varepsilon)_\varepsilon)(e_\varepsilon^{Bare})_\varepsilon. \tag{6.2.56}$$

As can be shown from (6.2.38) and (6.2.41), $(L_\varepsilon)_\varepsilon$ equals $(B_\varepsilon)_\varepsilon$ and only the $(Z_{3,\varepsilon})_\varepsilon$, factor remains for the renormalisation of the electric charge.

Let us consider corrections to the magnetic moment due to vertex corrections as (6.2.40).In particular, the term proportional to $\sigma_{\mu\nu}q^\nu$ remains finite for Hausdorff dimensions smaller than six. It gives rise to low-order contributions to the anomalous magnetic moment as well as the $l \neq 0$ splitting of energy levels in atoms (Lamb shift).Utilising the expansion of the gamma function into a polynomial $\Gamma(1+z) = \sum_{i=0}^{\infty} c_i z^i$ with coefficients $c_0 = 1, c_{n+1} = (n+1)^{-1}\sum_{i=0}^{n}(-1)^{i+1}s_{i+1}s_{n-i}$ and $s_1 = 0.577, s_n = \zeta(n)$ for $n \geq 2, \mathrm{Re}\,z > 0,$ where $\zeta(n)$ is the Riemann zeta function (e.g. $s_2 = \pi^2/6$) ,one obtains for small $\varepsilon = 4 - D^+ \ll 1$

$$\begin{aligned} \alpha_e(D^+ = 4) - \alpha_e(D^+) &= \alpha[4\pi^{-1} - 2^{2-D^+}\pi^{1-D^+/2}\Gamma(3 - D^+/2)] \simeq \\ &\frac{\alpha}{8\pi}(0.577 + \log\pi)(4 - D^+). \end{aligned} \tag{6.2.57}$$

Here, $\alpha_e$ is the form factor of the electromagnetic current proportional to $\sigma_{\mu\nu}q^\nu$. Presently the difference between experimental and theoretical values of $\alpha_e$ suggests

$$D^+ \geq 4 - (5.3 \pm 2.7) \times 10^{-7}, \tag{6.2.58}$$

i.e., $\varepsilon \leq (5.3 \pm 2.7) \times 10^{-7}$ at the scale of the Compton wavelength of the electron [37].

**Remark 6.2.5**. Note that for $\alpha = 1/137$ from (6.2.58) one obtains

$$\frac{3\alpha}{2\pi\varepsilon} = \frac{3}{2\pi \times 137 \times 3 \times 10^{-7}} = 11617 \gg 1 \tag{6.2.59}$$

Obviously this inequality in a contradiction with inequality (6.2.49) and therefore dimensional renormalization breaks down for $\alpha = 1/137$.

Similarly, corrections to the $l \neq 0$ levels for a hydrogen-like atom can be derived for $\delta E = \Delta E(D^+ = 4) - \Delta E(D^+)$ :

$$\begin{aligned} \delta E\left(2p\tfrac{1}{2}\right) &= \frac{\alpha^3}{2\pi n^3}[0.57722 + \log\pi](4 - D^+)Ry_\infty \times \\ &\begin{cases} [(l + 1/2)] & j = l + 1/2 \\ [-l(l + 1/2)]^{-1} & j = l - 1/2 \end{cases} \end{aligned} \tag{6.2.60}$$

This correction yields

$$\delta E_{l\neq 0} = \frac{\alpha^3 Ry_\infty}{24\pi}[0.57722 + \log\pi](4 - D^+) \leq 0.03 \pm 0.01 MHz \tag{6.2.61}$$

and the lower bound

$$D^+ \geq 4 - (1 \pm 0.3) \times 10^{-3}, \tag{6.2.62}$$

i.e., $\varepsilon \leq (1 \pm 0.3) \times 10^{-3}$.

**Remark 6.2.6**.Note that for $\alpha = 1/137$ from (6.2.62) one obtains

$$\frac{3\alpha}{2\pi\varepsilon} = \frac{3}{2\pi \times 137 \times 3 \times 10^{-3}} = 1.1617 > 1. \tag{6.2.63}$$

Obviously this inequality in a contradiction with inequality (6.2.49) and therefore dimensional renormalization breaks down for $\alpha = 1/137$.

# 6.3.Hausdorff-Colombeau measure and associated negative Hausdorff-Colombeau dimensions.

During last 20 years the notion of negative dimension in geometry was many developed,

see for example [15],[35]-[41].

Remind that canonical difenitions of noninteger positive dimension alwais equipped with an finite measure. For instance Hausdorff–Besicovich dimension equipped with Hausdorff measure $d\mu_H^+(x, D^+)$ [35]-[38].

**Remark 6.3.1.**Note that in the case where a finite measure $\mu$ on a metric space $(X,d)$ is exact-dimensional, i.e. there exists $\alpha \geq 0$ such that

$$\mu\left(X \setminus \left\{ x \,\middle|\, \lim_{r\to 0} \frac{\ln \mu(B(x,r))}{\ln r} = \alpha \right\}\right) \tag{6.3.1}$$

many different definitions of dimensions of $\mu$ collapse to the value $\alpha$.

Given a metric space $(X,d)$, let $B(X)$ denote the Borel $\sigma$-algebra in $X$, and let $BM(X)$ be the class of non-null finite Borel measures defined on $X$.

**Definition 6.3.1**. A measure $\mu \in BM(\mathbb{R}^D)$ is called:

(i) of lower exact dimension $\underline{\alpha}$ if

$$\mu\left(\mathbb{R}^D \setminus \left\{ x \,\middle|\, \liminf_{r\to 0} \frac{\ln \mu(B(x,r))}{\ln r} = \underline{\alpha} \right\}\right), \tag{6.3.2}$$

(ii) of apperer exact dimension $\overline{\alpha}$ if

$$\mu\left(\mathbb{R}^D \setminus \left\{ x \,\middle|\, \limsup_{r\to 0} \frac{\ln \mu(B(x,r))}{\ln r} = \overline{\alpha} \right\}\right). \tag{6.3.3}$$

**Remark 6.3.2.**It is well known that $\mu$ is of lower exact dimension $\underline{\alpha}$ if and only if $\dim_H(\mu) = \dim_H^*(\mu) = \underline{\alpha}$, and $\mu$ is of apperer exact dimension $\overline{\alpha}$ if and only if $\dim_P(\mu) = \dim_P^*(\mu) = \overline{\alpha}$. Here

$$\begin{aligned} \dim_H^*(\mu) &\triangleq \inf\{\dim_H(A) | \mu(A) = \|\mu\|, A \in B(\mathbb{R}^D)\}, \\ \dim_P^*(\mu) &\triangleq \inf\{\dim_P(A) | \mu(A) = \|\mu\|, A \in B(\mathbb{R}^D)\}, \end{aligned} \tag{6.3.4}$$

where $\|\mu\|$ is the total mass of $\mu$.

Let us consider now example of a space of noninteger positive dimension equipped with the Haar measure.On the closed interval $0 \leq x \leq 1$ there is a scale $0 \leq \sigma \leq 1$ of Cantor dust with the Haar measure equal to $x^\sigma$ for any interval $(0,x)$ similar to the entire given set of the Cantor dust. The direct product of this scale by the Euclidean cube of dimension $k-1$ gives the entire scale $k+\sigma$, where $k \in \mathbb{Z}$ and $\sigma \in (0,1)$ [36].

Let $\mu_D$ be a Lebesgue measure on $\mathbb{R}^D$ and let $\mathbf{S}_D$ be an $D$-dimensional ball of radius $R$. In the spherical coordinates, the volume $\mu_D(\mathbf{S}_D)$ of the ball is equal to

$$\frac{2^D}{\Gamma\left(1+\frac{D}{2}\right)} \cdot \int_0^R r^{D-1} dr = \frac{2^D}{\Gamma\left(1+\frac{D}{2}\right)} \cdot R^D. \tag{6.3.5}$$

Here $r^{n-1}$ stands for the *classical regular* density [36].

Let $\delta(x - x_0)dx$ be a Dirac measure with a support $\{x_0\}$ and total mass 1 on $\mathbb{R}$, i.e. for any mesurable set $A \subseteq \mathbb{R}$ $\delta_{x_0}(A) = \int_{A\subseteq\mathbb{R}} \delta(x - x_0)dx = 1$ iff $x_0 \in A$ and $\delta_{x_0}(A) = \int_{A\subseteq\mathbb{R}} \delta(x - x_0)dx = 0$ iff $x_0 \notin A$.

**Remark 6.3.3.** Note that: (i) Dirac measure $\delta(x - x_0)dx$ has a representation as Colombeau generalized measure

$$\delta(x - x_0)dx \approx_{\widetilde{\mathbb{R}}} \frac{1}{\pi}\left(\frac{\varepsilon dx}{(x - x_0)^2 + \varepsilon^2}\right)_{\varepsilon}, \tag{6.3.6}$$

$\varepsilon \in (0,1]$, i.e.

$$\int_{\mathbb{R}} \delta(x - x_0)f(x)dx \approx_{\widetilde{\mathbb{R}}} \left(\int_{\mathbb{R}} \delta(x - x_0;\varepsilon)f(x)dx\right)_{\varepsilon} = \frac{1}{\pi}\left(\int_{\mathbb{R}} \frac{\varepsilon f(x)dx}{(x - x_0)^2 + \varepsilon^2}\right)_{\varepsilon}; \tag{6.3.7}$$

(ii) here Colombeau generalized function $(\delta(x - x_0;\varepsilon))_{\varepsilon}$

$$(\delta(x - x_0;\varepsilon))_{\varepsilon} \triangleq \left(\frac{1}{\pi}\frac{\varepsilon}{(x - x_0)^2 + \varepsilon^2}\right)_{\varepsilon} \tag{6.3.8}$$

stands for the *nonclassical nonregular* density.

**Remark 6.3.4.** Note that: (i) for any $x \neq x_0$

$$\left(\frac{\varepsilon}{(x - x_0)^2 + \varepsilon^2}\right)_{\varepsilon} \approx_{\widetilde{\mathbb{R}}} 0_{\widetilde{\mathbb{R}}}, \tag{6.3.9}$$

(ii) for $x = x_0$

$$(\delta(x - x_0;\varepsilon))_{\varepsilon}|_{x=x_0} = \left(\frac{1}{\varepsilon}\right)_{\varepsilon}. \tag{6.3.10}$$

Given a metric space $(X,d)$, let $B(X)$ denote the Borel $\sigma$-algebra in $X$, and let $CBM(X)$ be the class of non-null hyperfinite Colombeau-Borel measures defined on $X$.

**Definition 6.3.2.** Let $(\mu_\varepsilon)_\varepsilon \in CBM(\mathbb{R}^D)$ be Colombeau-Borel measure and let $(I_\varepsilon(d^-,x,r))_\varepsilon$ be

$$(I_\varepsilon(d^-,x,r))_\varepsilon = \int_{B(x,r)} \|x - y\|^{|d^-|}(d\mu_\varepsilon(y))_\varepsilon. \tag{6.3.11}$$

Colombeau-Borel measure $(\mu_\varepsilon)_\varepsilon \in CBM(\mathbb{R}^D)$ is called exact-negative-dimensional iff there exists $D^-$ such that $-\infty < D^- < 0$ and $\forall x \in \mathbb{R}^D$ :

$$|D^-| = \inf\{|d_x^-| | \lim_{r\to 0} \lim_{\varepsilon\to 0} I_\varepsilon(d_x^-,x,r) = 0\}. \tag{6.3.12}$$

For instance Colombeau-Borel measure $(\mu_\varepsilon)_\varepsilon \in CBM(\mathbb{R})$ :

$$(\mu_\varepsilon)_\varepsilon = \left(\frac{dx}{(x - x_0)^k + \varepsilon^2}\right)_\varepsilon \tag{6.3.13}$$

is exact-negative-dimensional measure with $D^- = -k$.

**Definition 6.3.3.** Colombeau-Borel measure $(\mu_\varepsilon)_\varepsilon \in CBM(\mathbb{R}^D)$ is called:

(i) of lower negative dimension $\underline{D}^- < 0$ if

$$|\underline{D}^-| = \inf_{x\in\mathbb{R}^D} \inf\{|d_x^-| | \lim_{r\to 0} \lim_{\varepsilon\to 0} I_\varepsilon(d_x^-,x,r) = 0\}. \tag{6.3.14}$$

(ii) of apperer negative dimension $\overline{D}^- < 0$ if

$$|\overline{D}^-| = \sup_{x\in\mathbb{R}^D} \inf\{|d_x^-| | \lim_{r\to 0} \lim_{\varepsilon\to 0} I_\varepsilon(d_x^-,x,r) = 0\}. \tag{6.3.15}$$

For instance Colombeau-Borel measure $(\mu_\varepsilon)_\varepsilon \in CBM(\mathbb{R})$ :

$$(\mu_\varepsilon)_\varepsilon = \left(\frac{dx}{(x-x_1)^{k_1}+\varepsilon^2}\right)_\varepsilon + \left(\frac{dx}{(x-x_2)^{k_2}+\varepsilon^2}\right)_\varepsilon, \tag{6.3.16}$$

where $x_1 \neq x_2$ and $k_1 > k_2$, is lower negative dimensional with $\underline{D}^- = k_2$ and apperer negative dimensional with $\overline{D}^- = k_1$.

**Remark 6.3.5.** In this subsection we define generalized Hausdorff-Colombeau measure.In subsection VI.4 we will prove that negative dimensions of fractals equipped with the Hausdorff- Colombeau measure in natural way.

Let $\Omega$ be an open subset of $\mathbb{R}^n$, let $X$ be metric space $X \subsetneqq \mathbb{R}^n$ and let $F$ be a set $F = \{E_i\}_{i\in\mathbb{N}}$ of subsets $E_i$ of $X$. Let $\zeta(E,x,\check{x})$ be a function $\zeta : F\times\Omega\times\Omega \to \mathbb{R}$. Let $C_F^\infty(\Omega)$ be a set of the all functions $\zeta(E,x,\check{x})$ such that $\zeta(E,x,\check{x}) \in C^\infty(\Omega\times\Omega)$ whenever $E \in F$. Throughout this paper, for elements of the space $C_F^\infty(\Omega\times\Omega)^{(0,1]}$ of sequences of smooth functions indexed by $\varepsilon \in (0,1]$ we shall use the canonical notations $(\zeta_\varepsilon(E,x,\check{x}))_\varepsilon$ and $(\zeta_\varepsilon)_\varepsilon$ so $\zeta_\varepsilon \in C_F^\infty(\Omega\times\Omega)$, $\varepsilon \in (0,1]$.

**Definition 6.3.4.** We set $\mathcal{G}_F(\Omega) = \mathcal{E}_M(F,\Omega)/\mathcal{N}(F,\Omega)$, where

$$\begin{aligned}\mathcal{E}_M(F,\Omega) = \Big\{&(\zeta_\varepsilon(E,x,\check{x}))_\varepsilon \in C_F^\infty(\Omega\times\Omega)^{(0,1]} \,\Big|\, \forall K \subset\subset \Omega, \forall\alpha \in \mathbb{N}^n \exists p \in \mathbb{N} \text{ with} \\ &\sup_{E\in F; x\in K}|\zeta_\varepsilon(E,x,\check{x})| = O(\varepsilon^{-p}) \text{ as } \varepsilon \to 0\Big\}, \\ \mathcal{N}(F,\Omega) = \Big\{&(\zeta_\varepsilon(E,x,\check{x}))_\varepsilon \in C_F^\infty(\Omega\times\Omega)^{(0,1]} \,\Big|\, \forall K \subset\subset \Omega, \forall\alpha \in \mathbb{N}^n \forall q \in \mathbb{N} \\ &\sup_{E\in F; x\in K}|\zeta_\varepsilon(E,x,\check{x})| = O(\varepsilon^{q}) \text{ as } \varepsilon \to 0\Big\}.\end{aligned} \tag{6.3.17}$$

Notice that $\mathcal{G}_F(\Omega)$ is a differential algebra.Equivalence classes of sequences $(\zeta_\varepsilon)_\varepsilon = (\zeta_\varepsilon(E,x,\check{x}))_\varepsilon$ will be denoted by $\mathbf{cl}[(\zeta_\varepsilon)_\varepsilon]$ or simply $[(\zeta_\varepsilon)_\varepsilon]$.

**Definition 6.3.5.** We denote by $\widetilde{\mathbb{R}}$ the ring of real,Colombeau generalized numbers. Recall that by definition $\widetilde{\mathbb{R}} = \mathcal{E}_M(\mathbb{R})/\mathcal{N}(\mathbb{R})$ [21-23], where

$$\begin{aligned}\mathcal{E}_M(\mathbb{R}) &= \{(x_\varepsilon)_\varepsilon \in \mathbb{R}^{(0,1]} | (\exists\alpha \in \mathbb{R})(\exists\varepsilon_0 \in (0,1])\forall\varepsilon \leq \varepsilon_0[|x_\varepsilon| \leq \varepsilon^\alpha]\}, \\ \mathcal{N}(\mathbb{R}) &= \{(x_\varepsilon)_\varepsilon \in \mathbb{R}^{(0,1]} | (\forall\alpha \in \mathbb{R})(\exists\varepsilon_0 \in (0,1])\forall\varepsilon \leq \varepsilon_0[|x_\varepsilon| \leq \varepsilon^\alpha]\}.\end{aligned} \tag{6.3.18}$$

Notice that the ring $\widetilde{\mathbb{R}}$ arises naturally as the ring of constants of the Colombeau algebras $\mathcal{G}(\Omega)$. Recall that there exists natural embedding $\widetilde{r} : \mathbb{R} \hookrightarrow \widetilde{\mathbb{R}}$ such that for all $r \in \mathbb{R}, \widetilde{r} = (r_\varepsilon)_\varepsilon$ where $r_\varepsilon \equiv r$ for all $\varepsilon \in (0,1]$. We say that $r$ is standard number and abbreviate $r \in \mathbb{R}$ for short. The ring $\widetilde{\mathbb{R}}$ can be endowed with the structure of a partially ordered ring: for $r,s \in \widetilde{\mathbb{R}}$ $\eta \in \mathbb{R}_+, \eta \leq 1$ we abbreviate $r \leq_{\widetilde{\mathbb{R}},\eta} s$ or simply $r \leq_{\widetilde{\mathbb{R}}} s$ if and only if there are representatives $(r_\varepsilon)_\varepsilon$ and $(s_\varepsilon)_\varepsilon$ with $r_\varepsilon \leq s_\varepsilon$ for all $\varepsilon \in (0,\eta]$. Colombeau generalized number $r \in \widetilde{\mathbb{R}}$ with representative $(r_\varepsilon)_\varepsilon$ we abbreviate $\mathbf{cl}[(r_\varepsilon)_\varepsilon]$.

**Definition 6.3.6.** (i) Let $\delta \in \widetilde{\mathbb{R}}$. We say that $\delta$ is infinite small Colombeau generalized number and abbreviate $\delta \approx_{\widetilde{\mathbb{R}}} \widetilde{0}$ if there exists representative $(\delta_\varepsilon)_\varepsilon$ and some $q \in \mathbb{N}$ such that $|\delta_\varepsilon| = O(\varepsilon^q)$ as $\varepsilon \to 0$. (ii) Let $\Delta \in \widetilde{\mathbb{R}}$. We say that $\Delta$ is infinite large Colombeau generalized number and abbreviate $\Delta =_{\widetilde{\mathbb{R}}} \widetilde{\infty}$ if $\Delta_{\widetilde{\mathbb{R}}}^{-1} \approx_{\widetilde{\mathbb{R}}} \widetilde{0}$. (iii) Let $\mathbb{R}_\infty$ be $\mathbb{R} \cup \{\infty\}$ We say that $\Theta \in \widetilde{\mathbb{R}}_\infty$ is infinite Colombeau generalized number and abbreviate $\Theta =_{\widetilde{\mathbb{R}}} \infty_{\widetilde{\mathbb{R}}}$ if there exists representative $(\Theta_\varepsilon)_\varepsilon$ where $\Theta_\varepsilon = \infty$ for all $\varepsilon \in (0,1]$. Here we set $\mathcal{E}_M(\mathbb{R}_\infty) = \mathcal{E}_M(\mathbb{R}) \cup \{(\Theta_\varepsilon)_\varepsilon\}$, $\mathcal{N}(\mathbb{R}_\infty) = \mathcal{N}(\mathbb{R})$ and $\widetilde{\mathbb{R}}_\infty = \mathcal{E}_M(\mathbb{R}_\infty)/\mathcal{N}(\mathbb{R}_\infty)$.

**Definition 6.3.7**.The singular Hausdorff-Colombeau measure originate in Colombeau generalization of canonical Caratheodory's construction, which is defined as follows: for each metric space $X$, each set $F = \{E_i\}_{i\in\mathbb{N}}$ of subsets $E_i$ of $X$, and each Colombeau generalized function $(\zeta_\varepsilon(E,x,\breve{x}))_\varepsilon$, such that: (i) $0 \leq (\zeta_\varepsilon(E,x,\breve{x}))_\varepsilon$ , (ii) $(\zeta_\varepsilon(E,\breve{x},\breve{x}))_\varepsilon =_{\widetilde{\mathbb{R}}} \widetilde{\infty}$, whenever $E \in F$, a preliminary Colombeau measure $(\phi_\delta(E,x,\breve{x},\varepsilon))_\varepsilon$ can be constructed corresponding to $0 < \delta \leq +\infty$ , and then a final Colombeau measure $(\mu_\varepsilon(E,x,\breve{x}))_\varepsilon$, as follows: for every subset $E \subset X$, the preliminary Colombeau measure $(\phi_\delta(E,x,\breve{x},\varepsilon))_\varepsilon$ is defined by

$$\phi_\delta(E,x,\breve{x},\varepsilon) = \sup_{\{E_i\}_{i\in\mathbb{N}}} \left\{\sum_{i\in\mathbb{N}} \zeta_\varepsilon(E_i,x,\breve{x}) | E \subset \bigcup_{i\in\mathbb{N}} E_i, diam(E_i) \leq \delta\right\}. \quad (6.3.19)$$

Since for all $\varepsilon \in (0,1]$ : $\phi_{\bar{\delta}_1}(E,x,\breve{x},\varepsilon) \geq \phi_{\bar{\delta}_2}(E,x,\breve{x},\varepsilon)$ for $0 < \delta_1 < \delta_2 \leq +\infty$, the limit

$$(\mu(E,x,\breve{x},\varepsilon))_\varepsilon = \left(\lim_{\delta\to 0_+} \phi_\delta(E,x,\breve{x},\varepsilon)\right)_\varepsilon = \left(\inf_{\delta>0} \phi_\delta(E,x,\breve{x},\varepsilon)\right)_\varepsilon \quad (6.3.20)$$

exists for all $E \subset X$. In this context, $(\mu(E,x,\breve{x},\varepsilon))_\varepsilon$ can be called the result of Caratheodory's construction from $(\zeta_\varepsilon(E,x,\breve{x}))_\varepsilon$ on $F$ and $(\phi_\delta(E,x,\breve{x},\varepsilon))_\varepsilon$ can be referred to as the size $\delta$ approximating Colombeau measure.

**Definition 6.3.8**. Let $(\zeta_\varepsilon(E_i,D^+,D^-,x,\breve{x}))_\varepsilon$ be

$$(\zeta_\varepsilon(E_i,D^+,D^-,x,\breve{x}))_\varepsilon = \begin{cases} \left(\dfrac{\Omega_1(D^+)\Omega_2(D^-)[diam(E_i)]^{D^+}}{[d(x,\breve{x})]^{|D^-|} + \varepsilon}\right)_\varepsilon & \text{if } x \in E_i \\ 0 & \text{if } x \notin E_i \end{cases} \quad (6.3.21)$$

where $\varepsilon \in (0,1], \Omega_1(D^+), |\Omega_2(D^-)| > 0$. In particular, when $F$ is the set of all (closed or open) balls in $X$,

$$\Omega_1(D^+) = \frac{2^{-D^+}\Gamma\left(\frac{1}{2}\right)^{D^+}}{\Gamma\left(1+\frac{D^+}{2}\right)} = \frac{2^{-D^+}\pi^{\frac{D^+}{2}}}{\Gamma\left(1+\frac{D^+}{2}\right)} \quad (6.3.22)$$

and

$$\Omega_2(D^-) = \frac{2^{-D^-}\pi^{\frac{D^-}{2}}}{\Gamma\left(1+\frac{D^-}{2}\right)}, \quad (6.3.23)$$

$$D^- \neq -2,-4,-6,\ldots,-2(n+1),\ldots$$

**Definition 6.3.9**. The Hausdorff-Colombeau singular measure $(\mu_{HC}(E,D^+,D^-,x,\breve{x},\varepsilon))_\varepsilon$ of a subset $E \subset X$ with the associated Hausdorff-Colombeau dimension $\breve{D}^+(D^-) \in \mathbb{R}_+, D^- \in \mathbb{R}_+$, is defined by

$$(\mu_{HC}(E,\breve{D}^+,D^-,x,\breve{x},\varepsilon))_\varepsilon =$$

$$\left(\lim_{\delta\to 0}\left[\sup_{\{E_i\}_{i\in\mathbb{N}}} \left\{\sum_{i\in\mathbb{N}} (\zeta_\varepsilon(E_i,\breve{D}^+,D^-,x,\breve{x}))_\varepsilon | E \subset \bigcup_i E_i, \forall i(diam(E_i) < \delta)\right\}\right]\right)_\varepsilon, \quad (6.3.24)$$

$$\breve{D}^+ = \sup\left\{D^+ > 0 | (\mu_{HC}(E,D^+,D^-,x,\breve{x},\varepsilon))_\varepsilon = \infty_{\widetilde{\mathbb{R}}}\right\},$$

The Colombeau-Lebesgue-Stieltjes integral over continuous functions $f : X \to \mathbb{R}$ can be evaluated similarly as in subsection III.3,(but using the limit in sense of Colombeau generalized functions) of infinitesimal covering diameter when $\{E_i\}_{i\in\mathbb{N}}$ is a disjoined covering and $x_i \in E_i$ :

$$\left(\int_X f(x)d\mu_{HC}(E,D^+,D^-,x,\breve{x},\varepsilon)\right)_\varepsilon = \left(\lim_{diam(E_i)\to 0}\left[\sum_{\cup E_i=X} f(x_i) \inf_{\{E_{ij}\}\text{ with }\cup_j E_{ij}\supset E_i} \sum_j \zeta_\varepsilon(E_i,D^+,D^-,x_i,\breve{x})\right]\right)_\varepsilon . \tag{6.3.25}$$

We assume now that $X$ is metrically unbounded, i.e. for every $x \in X$ and $r > 0$ there exists a point $y$ such that $d(x,y) > r$. The standard assumption that $\breve{D}^+$ and $\breve{D}^-$ is uniquely defined in all subsets $E$ of $X$ requires $X$ to be regular (homogeneous, uniform) with respect to the measure,i.e. $(\mu_{HC}^-(B_r(\breve{x}),\breve{D}^+,\breve{D}^-,x,\breve{x},\varepsilon))_\varepsilon = (\mu_{HC}^-(B_r(\breve{y}),\breve{D}^+,\breve{D}^-,x',\breve{y},\varepsilon))_\varepsilon$, where $d(x,\breve{x}) = d(x',\breve{y})$ for all elements $\breve{x},\breve{y},x,x' \in X$ and convex balls $B_r(\breve{x})$ and $B_r(\breve{y})$ of the form $B_r(\breve{x}) = \{z|d(\breve{x},z) \le r;\breve{x},z \in X\}$ and $B_r(\breve{y}) = \{z|d(\breve{y},z) \le r;\breve{y},z \in X\}$. In the limit $diam(E_i) \to 0$, the infimum is satisfied by the requirement that the variation over all coverings $\{E_{ij}\}_{ij\in\mathbb{N}}$ is replaced by one single covering $E_i$, such that $\cup_j E_{ij} \to E_i \ni x_i$. Therefore

$$\left(\int_X f(x)d\mu_{HC}(E,\breve{D}^+,\breve{D}^-,x,\breve{x},\varepsilon)\right)_\varepsilon = \left(\lim_{diam(E_i)\to 0} \sum_{\cup E_i=X} f(x_i)\zeta_\varepsilon(E_i,\breve{D}^+,\breve{D}^-,x_i,\breve{x})\right)_\varepsilon . \tag{6.3.26}$$

Assume that $f(x) = f(r), r = \|r\|$. The range of integration $X$ may be parametrised by polar coordinates with $r = d(x,0)$ and angle $\omega$. $\{E_{r_i,\omega_i}\}$ can be thought of as spherically symmetric covering around a centre at the origin. Thus

$$\left(\int_X f(r)d\mu_{HC}(E_x,\breve{D}^+,\breve{D}^-,x,\breve{x},\varepsilon)\right)_\varepsilon = \left(\lim_{diam(E_i)\to 0} \sum_{\cup E_i=X} f(r_i)\zeta_\varepsilon(E_i,\breve{D}^+,\breve{D}^-,x_i,\breve{x})\right)_\varepsilon . \tag{6.3.27}$$

Notice that the metric set $X \subset \mathbb{R}^n$ can be tesselated into regular polyhedra; in particular it is always possible to divide $\mathbb{R}^n$ into parallelepipeds of the form

$$\Pi_{i_1,\ldots,i_n} = \{x = (x_1,\ldots,x_n) \in X|x_j = (i_j-1)\Delta x_j + \gamma_j, 0 \le \gamma_j \le \Delta x_j, j = 1,\ldots,n\}. \tag{6.3.28}$$

For $n = 2$ the polyhedra $\Pi_{i_1,i_2}$ is shown in figure 6.3.1.Since $X$ is uniform

$$(d\mu_{HC}(E,\breve{D}^+,\breve{D}^-,x,\breve{x},\varepsilon))_\varepsilon = \left(\lim_{diam(\Pi_{i_1,\ldots,i_n})} \zeta_\varepsilon(\Pi_{i_1,\ldots,i_n},\breve{D}^+,\breve{D}^-,x,\breve{x})\right)_\varepsilon = \left(\lim_{diam(\Pi_{i_1,\ldots,i_n})} \prod_{j=1}^n \left(\frac{\Delta x_j}{|x_j-\breve{x}_j|^{|\breve{D}^-|}+\varepsilon}\right)^{\frac{\breve{D}^+}{n}}\right)_\varepsilon \triangleq \tag{6.3.29}$$

$$\triangleq \left(\prod_{j=1}^n \frac{d^{\frac{\breve{D}^+}{n}}x_j}{\left(|x_j-\breve{x}_j|^{|\breve{D}^-|}+\varepsilon\right)^{\frac{\breve{D}^+}{n}}}\right)_\varepsilon .$$

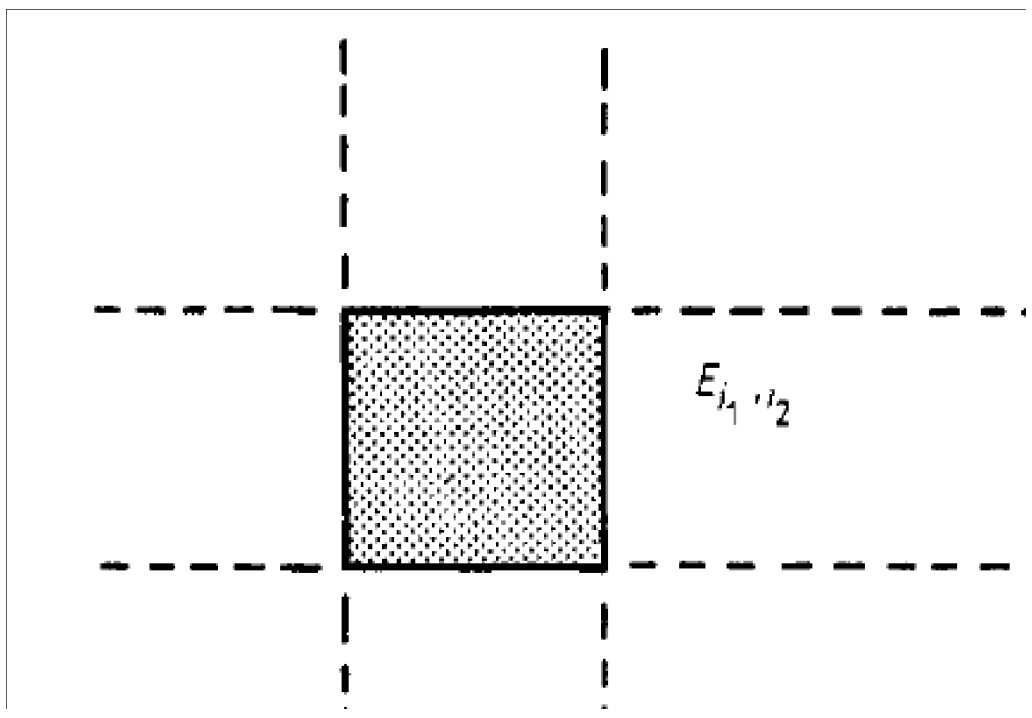

Fig.3.3.1. The polyhedra covering for $n = 2$.

Notice that the range of integration $X$ may also be parametrised by polar coordinates with $r = d(x,0)$ and angle $\Omega$. $E_{r,\Omega}$ can be thought of as spherically symmetric covering around a centre at the origin (see figure 3.3.2 for the two-dimensional case). In the limit, the Colombeau generaliza function $(\zeta_\varepsilon(E_{r,\Omega},\breve{D}^{+},\breve{D}^{-},r,\breve{r}))_\varepsilon$ is given by

$$(d\mu_{HC}(E_{r,\Omega},\breve{D}^{+},\breve{D}^{-},r,\breve{r},\Omega,\varepsilon))_\varepsilon =$$

$$\left(\lim_{diam(\Pi_{i_1,\dots,i_n})}\zeta_\varepsilon(E_{r,\Omega},\breve{D}^{+},\breve{D}^{-},r,\breve{r},\Omega)\right)_\varepsilon \triangleq \frac{d\Omega^{\breve{D}^{+}-1}r^{\breve{D}^{+}-1}dr}{\left((r-\breve{r})^{|\breve{D}^{-}|}+\varepsilon\right)_\varepsilon}. \qquad (6.3.30)$$

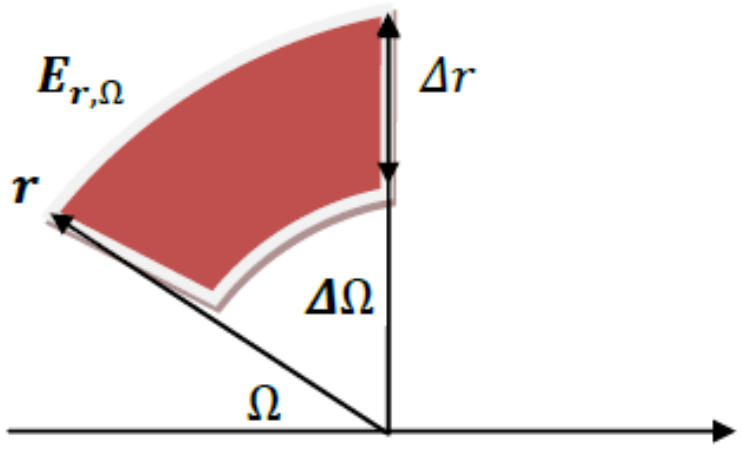

Fig.6.3.2.The spherical covering $E_{r,\Omega}$.

When $f(x)$ is symmetric with respect to some centre $\breve{x} \in X$, i.e. $f(x) =$ constant for all $x$ satisfying $d(x,\breve{x}) = r$ for arbitrary values of $r$, then chainge of the variable

$$x \to z = x - \breve{x} \qquad (6.3.31)$$

can be performed to shift the centre of symmetry to the origin (since $X$ is not a linear space, (6.3.15) need not be a map of $X$ onto itself and (6.3.15) is measure preseming). The integral over metric space $X$ is then given by formula

$$\left(\int_X f(x)d\mu_{HC}(E_x,\breve{D}^{+},\breve{D}^{-},x,\breve{x},\varepsilon)\right)_\varepsilon = \frac{4\pi^{D^{+}/2}\pi^{D^{-}/2}}{\Gamma(D^{+}/2)\Gamma(D^{-}/2)}\left(\int_0^\infty \frac{r^{D^{+}-1}f(r)}{r^{|D^{-}|}+\varepsilon}dr\right)_\varepsilon. \qquad (6.3.32)$$

The Colombeau integral defined in (6.3.27) satisfies the following conditions.

(i) Linearity:

$$\left(\int_X [a_1 f_1(x) + a_2 f_2(x)] d\mu_{HC}(E_x, \breve{D}^+, \breve{D}^-, x, \breve{x}, \varepsilon)\right)_\varepsilon =$$
$$a_1 \left(\int_X f_1(x) d\mu_{HC}(E_x, \breve{D}^+, \breve{D}^-, x, \breve{x}, \varepsilon)\right)_\varepsilon + a_2 \left(\int_X f_2(x) d\mu_{HC}(E_x, \breve{D}^+, \breve{D}^-, x, \breve{x}, \varepsilon)\right)_\varepsilon . \quad (6.3.33)$$

(ii) Translational invariance:

$$\left(\int_X f(x + x_0) d\mu_{HC}(E_x, \breve{D}^+, \breve{D}^-, x, \breve{x}, \varepsilon)\right)_\varepsilon = \left(\int_X f(x) d\mu_{HC}(E_x, \breve{D}^+, \breve{D}^-, x, \breve{x}, \varepsilon)\right)_\varepsilon \quad (6.3.34)$$

since

$$(d\mu_{HC}(E_{x-x_0}, \breve{D}^+, \breve{D}^-, x - x_0, \breve{x} - x_0, \varepsilon))_\varepsilon = (d\mu_{HC}(E_x, \breve{D}^+, \breve{D}^-, x, \breve{x}, \varepsilon))_\varepsilon \quad (6.3.35) .$$

(iii) Scaling property:

$$\left(\int_X f(a \cdot x) d\mu_{HC}(E_x, \breve{D}^+, \breve{D}^-, x, \breve{x}, \varepsilon)\right)_\varepsilon = a^{-D^+ + D^-} \left(\int_X f(x) d\mu_{HC}(E_x, \breve{D}^+, \breve{D}^-, x, \breve{x}, \varepsilon)\right)_\varepsilon \quad (6.3.36)$$

## 6.4.Main properties of the Hausdorff-Colombeau metric measures with associated negative Hausdorff-Colombeau dimensions.

**Definition 6.4.1**. An outer Colombeau metric measure on a set $X \subseteqq \mathbb{R}^n$ is a Colombeau
generalized function $[(\phi_\varepsilon(E))_\varepsilon] \in \mathcal{G}_F(\Omega)$ (see Definition 6.3.1) defined on all subsets of $X$ satisfies the following properties:.

(i) Null empty set:The empty set has zero Colombeau outer measure

$$[(\phi_\varepsilon(\varnothing))_\varepsilon] = 0. \quad (6.4.1)$$

(ii) Monotonicity: For any two subsets $A$ and $B$ of $X$

$$A \subseteqq B \Rightarrow [(\phi_\varepsilon(A))_\varepsilon] \leq_{\widetilde{\mathbb{R}}} [(\phi_\varepsilon(B))_\varepsilon]. \quad (6.4.2)$$

(iii) Countable subadditivity: For any sequence $\{A_j\}$ of subsets of $X$ pairwise disjoint or not

$$\left[(\phi_\varepsilon(\cup_{j=1}^\infty A_j))_\varepsilon\right] \leq_{\widetilde{\mathbb{R}}} \left[\left(\sum_{j=1}^\infty \phi_\varepsilon(A_j)\right)_\varepsilon\right]. \quad (6.4.3)$$

(iv) Whenever $d(A,B) = \inf\{d_n(x,y) | x \in A, y \in B\} > 0$

$$[(\phi_\varepsilon(A \cup B))_\varepsilon] = [(\phi_\varepsilon(A))_\varepsilon] + [(\phi_\varepsilon(B))_\varepsilon], \quad (6.4.4)$$

where $d_n(x,y)$ is the usual Euclidean metric:$d_n(x,y) = \sqrt{\sum (x_i - y_i)^2}$.

**Definition 6.4.2**. We say that outer Colombeau metric measure $(\mu_\varepsilon)_\varepsilon, \varepsilon \in (0,1]$ is a Colombeau measure on $\sigma$-algebra of subests of $X \subseteqq \mathbb{R}^n$ if $(\mu_\varepsilon)_\varepsilon$ satisfies the following property:

$$\left[(\phi_\varepsilon(\cup_{j=1}^\infty A_j))_\varepsilon\right] = \left[\left(\sum_{j=1}^\infty \phi_\varepsilon(A_j)\right)_\varepsilon\right]. \quad (6.4.5)$$

**Definition 6.4.3**.If $U$ is any non-empty subset of $n$-dimensional Euclidean space, $\mathbb{R}^n$, the
diamater $|U|$ of $U$ is defined as

$$|U| = \sup\{d(x,y) | x,y \in U\} \quad (6.4.6)$$

If $F \subseteq \mathbb{R}^n$, and a collection $\{U_i\}_{i \in \mathbb{N}}$ satisfies the following conditions:

(i) $|U_i| < \delta$ for all $i \in \mathbb{N}$, (ii) $F \subseteq \cup_{i\in\mathbb{N}} U_i$, then we say the collection $\{U_i\}_{i\in\mathbb{N}}$ is a $\delta$-cover of $F$.

**Definition 6.4.4.** If $F \subseteq \mathbb{R}^n$ and $s, q, \delta > 0$, we define Hausdorff-Colombeau content:

$$(H_\delta^{s,q}(F,\varepsilon))_\varepsilon = \left(\inf\left\{\sum_{i=1}^{\infty} \frac{|U_i|^s}{\|x_i\|^q + \varepsilon}\right\}\right)_\varepsilon \tag{6.4.7}$$

where the infimum is taken over all $\delta$-covers of $F$ and where $x_i = (x_{i,1}, \ldots x_{i,n}) \in U_i$ for all $i \in \mathbb{N}$ and $\|x\|$ is the usual Euclidean norm: $\|x\| = \sqrt{\sum_{j=1}^{n} x_j^2}$.

Note that for $0 < \delta_1 < \delta_2 < 1, \varepsilon \in (0,1]$ we have

$$H_{\delta_1}^{s,q}(F,\varepsilon) \geq H_{\delta_2}^{s,q}(F,\varepsilon) \tag{6.4.8}$$

since any $\delta_1$ cover of $F$ is also a $\delta_2$ cover of $F$, i.e. $H_{\delta_1}^{s,q}(F,\varepsilon)$ is increasing as $\delta$ decreases.

**Definition 6.4.5.** We define the $(s,q)$-dimensional Hausdorff-Colombeau (outer) measure as:

$$(H^{s,q}(F,\varepsilon))_\varepsilon = \left(\lim_{\delta\to 0} H_\delta^{s,q}(F,\varepsilon)\right)_\varepsilon. \tag{6.4.9}$$

**Theorem 6.4.1.** For any $\delta$-cover, $\{U_i\}_{i\in\mathbb{N}}$ of $F$, and for any $\varepsilon \in (0,1]$, $t > s$ :

$$H^{t,q}(F,\varepsilon) \leq \delta^{t-s} H^{s,q}(F,\varepsilon). \tag{6.4.10}$$

**Proof.** Consider any $\delta$-cover $\{U_i\}_{i\in\mathbb{N}}$ of $F$. Then each $|U_i|^{t-s} \leq \delta^{t-s}$ since $|U_i| \leq \delta$, so:

$$|U_i|^t = |U_i|^{t-s}|U_i|^s \leq \delta^{t-s}|U_i|^s. \tag{6.4.11}$$

From (6.4.11) follows that

$$\frac{|U_i|^t}{\|x_i\|^q + \varepsilon} \leq \frac{\delta^{t-s}|U_i|^s}{\|x_i\|^q + \varepsilon} \tag{6.4.12}$$

and summing (6.4.11) over all $i \in \mathbb{N}$ we obtain

$$\sum_{i=1}^{\infty} \frac{|U_i|^t}{\|x_i\|^q + \varepsilon} \leq \delta^{t-s} \sum_{i=1}^{\infty} \frac{|U_i|^s}{\|x_i\|^q + \varepsilon}. \tag{6.4.13}$$

Thus (3.4.10) follows by taking the infimum.

**Theorem 6.4.2.** (1) If $(H^{s,q}(F,\varepsilon))_\varepsilon <_{\widetilde{\mathbb{R}}} \infty_{\widetilde{\mathbb{R}}}$, and if $t > s$, then $(H^{t,q}(F,\varepsilon))_\varepsilon = 0_{\widetilde{\mathbb{R}}}$.

(2) If $0_{\widetilde{\mathbb{R}}} <_{\widetilde{\mathbb{R}}} (H^{s,q}(F,\varepsilon))_\varepsilon$, and if $t < s$, then $(H^{t,q}(F,\varepsilon))_\varepsilon = \infty_{\widetilde{\mathbb{R}}}$.

**Proof.** (1) The result follows from (6.4.10) after taking limits, since $\forall\varepsilon \in (0,1]$ by definitions follows that $H^{s,q}(F,\varepsilon) < \infty$.

(2) From (6.4.10) $\forall\varepsilon \in (0,1], \forall\delta > 0$ follows that

$$\frac{1}{\delta^{s-t}} H^{s,q}(F,\varepsilon) \leq H^{t,q}(F,\varepsilon). \tag{6.4.14}$$

After taking limit $\delta \to 0$, we obtain $H^{t,q}(F,\varepsilon) = \infty$, since $\lim_{\delta\to 0}(\delta^{s-t})^{-1} = \infty$ and $\lim_{\delta\to 0} H_\delta^{s,q}(F,\varepsilon) = H^{s,q}(F,\varepsilon) > 0$.

**Definition 6.4.6.** We define now the Hausdorff-Colombeau dimension $\dim_{HC}(F,q)$ of a set $F$ (relative to $q > 0$) as

$$\dim_{HC}(F,q) = \sup\left\{s = s(q) | (H^{s,q}(F,\varepsilon))_\varepsilon = \infty_{\widetilde{\mathbb{R}}}\right\} = \inf\left\{s = s(q) | (H^{s,q}(F,\varepsilon))_\varepsilon = 0_{\widetilde{\mathbb{R}}}\right\}. \tag{6.4.15}$$

**Remark 6.4.1.** From theorem 3.4.2 follows that for any fixed $q = \check{q}$ :
$(H^{s,\check{q}}(F,\varepsilon))_{\varepsilon} = 0_{\widetilde{\mathbb{R}}}$ or $(H^{s,\check{q}}(F,\varepsilon))_{\varepsilon} = \infty_{\widetilde{\mathbb{R}}}$ everywhere except at a unique value $s$, where this
value may be finite. As a function of $s, H^{s,\check{q}}(F,\varepsilon)$ is decreasing function. Therefore, the graph of $H^{s,\check{q}}(F,\varepsilon)$ will have a unique value where it jumps from $\infty$ to $0$.

**Remark 6.4.2.** Note that the graph of $(H^{s,\check{q}}(F,\varepsilon))_{\varepsilon}$ for a fixed $q = \check{q}$ is

$$(H^{s,\check{q}}(F,\varepsilon))_{\varepsilon} = \begin{cases} \infty_{\widetilde{\mathbb{R}}} & \text{if} \quad s < \dim_{HC}(F,\check{q}) \\ 0_{\widetilde{\mathbb{R}}} & \text{if} \quad s > \dim_{HC}(F,\check{q}) \\ \text{undetermined} & \text{if} \quad s = \dim_{HC}(F,\check{q}) \end{cases} \tag{6.4.16}$$

**Definition 6.4.7.** We say that fractal $\mathcal{F} \subset \mathbb{R}^n$ has a negative dimension relative to $\check{q} > 0$ iff $\dim_{HC}(F,q) - \check{q} < 0$.

# 7. Scalar quantum field theory in spacetime with Hausdorff-Colombeau negative dimensions.

## 7.1. Equation of motion and Hamiltonian.

Scalar quantum field theory and quantum gravity in spacetime with noninteger positive Hausdorff dimensions developed in papers [42]-[45]. Quantum mechanics in negative dimensions developed in papers [40],[41] Scalar quantum field theory and quantum gravity in spacetime with Hausdorff-Colombeau negative dimensions originally developed in paper [15]. In this section only free scalar quantum field in spacetime with negative dimensions briefly is considered.

A negative-dimensional spacetime structures is a desirable feature of superrenormalizable spacetime models of quantum gravity, and the most simply way to obtain it is to let the effective dimensionality of the multifractal universe to change at different scales. A simple realization of this feature is via suitable extended fractional calculus and the definition of a fractional action. Note that below we use canonical isotropic scaling such that:

$$[x^{\mu}] = -1, \mu = 0,1,\ldots,D_{\mathbf{t}} - 1 \tag{7.1.1}$$

while replacing the standard measure with a nontrivial Colombeau-Stieltjes measure in negative Hausdorff-Colombeau dimension $D^{-}$,

$$\begin{gathered} d^{D_{\mathbf{t}}}x \to d^{D^{-}}x = (d\eta(x,D^{-},\varepsilon))_{\varepsilon}, \\ [\eta] = D_{\mathbf{t}} \cdot \alpha, \alpha \in [1,-\infty). \end{gathered} \tag{7.1.2}$$

Here $D_{\mathbf{t}}$ is the topological (positive integer) dimension of embedding spacetime and $\alpha$ is a parameter. Any Colombeau integrals on net multifractals can be approximated by the left-sided Colombeau-Riemann–Liouville complex milti-fractional integral of a function $\mathcal{L}(t)$ :

$$\left(\int_0^{\bar{t}} d\eta(x,\varepsilon)\mathcal{L}(t)\right)_\varepsilon \propto \left(I_{\bar{t},\varepsilon}^{\{z_i(\bar{t})\}}\right)_\varepsilon \triangleq \left(\sum_{i=1}^m \int_\varepsilon^{\bar{t}} \frac{[(\bar{t}-t)+i\varepsilon]^{z_i(\bar{t})-1}}{\Gamma(z_i(\bar{t}))}\mathcal{L}(t)dt\right)_\varepsilon,$$
$$(\eta(t,\varepsilon))_\varepsilon = \left(\frac{\bar{t}^{\,z_i(\bar{t})} - [(\bar{t}-t)+i\varepsilon]^{z_i(\bar{t})}}{\Gamma(z_i(\bar{t})+1)}\right)_\varepsilon, \tag{7.1.3}$$

where $\varepsilon \in (0,1]$, $\bar{t}$ is fixed and the order $z(\bar{t})$ is (related to) the complex Hausdorff-Colombeau dimensions of the set. In particular if $z_i \in \mathbb{C}, i = 1,2,\ldots,m$ is a complex parameters an integrals on net multifractals can be approximated by finite sum of the left-sided Colombeau-Riemann–Liouville complex fractional integral of a function $\mathcal{L}(t)$

$$\left(\int_0^{\bar{t}} d\eta(x,\varepsilon)\mathcal{L}(t)\right)_\varepsilon \propto \left(I_{\bar{t},\varepsilon}^{\{z_i\}_{i=1}^m}\right)_\varepsilon =$$
$$\sum_{i=1}^m \left(I_{\bar{t},\varepsilon}^{z_i}\right)_\varepsilon \triangleq \sum_{i=1}^m \left(\frac{1}{\Gamma(z_i)}\int_\varepsilon^{\bar{t}} d[(\bar{t}-t)+i\varepsilon]^{z_i-1}\mathcal{L}(t)\right)_\varepsilon, \tag{7.1.4}$$
$$(\eta(t,\varepsilon))_\varepsilon = \sum_{i=1}^m \left(\frac{\bar{t}^{\,z_i} - [(\bar{t}-t)+i\varepsilon]^{z_i}}{\Gamma(z_i+1)}\right)_\varepsilon.$$

Note that a change of variables $t \to \bar{t}-t$ trasforms Eq. (7.1.4) into the form

$$\left(\int_0^{\bar{t}} d\eta(x,\varepsilon)\mathcal{L}(t)\right)_\varepsilon = \sum_{i=1}^m \left(\int_0^{\bar{t}} dt \frac{[t+i\varepsilon]^{z_i-1}}{\Gamma(z(\bar{t}))}\mathcal{L}(\bar{t}-t)\right)_\varepsilon. \tag{4.1.5}$$

The Colombeau-Riemann–Liouville multifractional integral (7.1.5) can be mapped onto a Colombeau-Weyl multifractional integral in the formal limit $\bar{t} \to +\infty$. We assume otherwise,so that there exists $\lim_{\bar{t}\to+\infty} z(\bar{t})$ and $\lim_{\bar{t}\to+\infty} \mathcal{L}(\bar{t}-t) = \mathcal{L}[q(t),\dot{q}(t)]$. In particular if

$z \in \mathbb{C}$ is a complex parameter a change of variables $t \to \bar{t}-t$ trasforms eq. (7.1.5) into the form

$$\sum_{i=1}^m \left(I_{\bar{t},\varepsilon}^{z_i}\right)_\varepsilon = \sum_{i=1}^m \left(\int_\varepsilon^{\bar{t}} dt \frac{[t+i\varepsilon]^{z_i-1}}{\Gamma(z_i)}\mathcal{L}[q(t),\dot{q}(t)]\right)_\varepsilon. \tag{7.1.6}$$

This form will be the most convenient for defining a Colombeau-Stieltjes field theory action.In $D_{\mathbf{t}}$ dimensions, we consider now the action

$$(\mathbf{S}_\varepsilon)_\varepsilon = \left(\int_M d\eta(x,\varepsilon)\mathcal{L}[\varphi_\varepsilon(x),\partial_\mu\varphi_\varepsilon(x)]\right)_\varepsilon, \tag{7.1.7}$$

where $\mathcal{L}[\varphi,\partial_\mu\varphi]$ is the Lagrangian density of the scalar field $(\varphi_\varepsilon(x))_\varepsilon$ and where

$$(d\eta(x,\varepsilon))_\varepsilon = \sum_{i=1}^m \prod_{\mu=0}^{D_{\mathbf{t}}-1}\left(f_{\mu,i}(x,\varepsilon)\right)_\varepsilon dx^\mu, \left(f_{\mu,i}(x,\varepsilon)\right)_\varepsilon : M \to \widetilde{\mathbb{R}}, \tag{7.1.8}$$

is some Colombeau–Stieltjes measure. We denote with pair $(M,(d\eta(x,\varepsilon))_\varepsilon)$ the metric spacetime $M$ equipped with Colombeau-Stieltjes measure $(d\eta(x,\varepsilon))_\varepsilon$. The former can be taken to be the canonical scalar field Lagrangian,

$$(\mathcal{L}[\varphi_\varepsilon(x),\partial_\mu\varphi_\varepsilon(x)])_\varepsilon = -\frac{1}{2}(\partial_\mu\varphi_\varepsilon\partial^\mu\varphi_\varepsilon)_\varepsilon - (V_\varepsilon(\varphi_\varepsilon))_\varepsilon, \tag{7.1.9}$$

where $V_\varepsilon(\varphi_\varepsilon)$ is a potential and contraction of Lorentz indices is given via the Minkowski metric $\eta_{\mu\nu} = (-,+,\ldots+)_{\mu\nu}$. As for the Colombeau-Stieltjes measure, we make the multifractal spacetime isotropic choice

$$\left(f_{(\mu,i),\varepsilon}\right)_{\varepsilon} = (f_{i,\varepsilon})_{\varepsilon}, \mu = 0,1,\ldots,D_{\mathbf{t}}-1; i = 1,\ldots,m. \quad (7.1.10)$$

Hence the scalar field action (4.1.7) reads

$$\begin{gathered}(\mathbf{S}_{\varepsilon})_{\varepsilon} = \left(\int_{M} d\eta(x,\varepsilon)\mathcal{L}[\varphi_{\varepsilon}(x),\partial_{\mu}\varphi_{\varepsilon}(x)]\right)_{\varepsilon} = \\ \sum_{j=1}^{m}\left(\int d^{D_{t}}x v_{\varepsilon,j}(x)\left[\frac{1}{2}\partial_{\mu}\varphi_{\varepsilon}\partial^{\mu}\varphi_{\varepsilon} + V(\varphi_{\varepsilon})\right]\right)_{\varepsilon} = \\ \sum_{j=1}^{m}\left(\int d^{D_{t}}x v_{\varepsilon,j}(x,D_{,j}^{-})\left[\frac{1}{2}\partial_{\mu}\varphi_{\varepsilon}\partial^{\mu}\varphi_{\varepsilon} + V(\varphi_{\varepsilon})\right]\right)_{\varepsilon},\end{gathered} \quad (7.1.11)$$

where $(v_{\varepsilon,j}(x))_{\varepsilon} = (v_{\varepsilon,j}(x,D_{,j}^{-}))_{\varepsilon}$ is a coordinate-dependent Lorentz scalar

$$(v_{\varepsilon,j}(x,D_{j}^{-}))_{\varepsilon} = \left(\frac{1}{[s_{j}(x)]^{D_{\mathbf{t}}(|\alpha_{j}-1|)} + \varepsilon}\right)_{\varepsilon}. \quad (7.1.12)$$

We define now the Dirac distribution in negative dimensions $D_{j}^{-}, j = 1,\ldots,m$ as Colombeau generalized function by equation

$$\sum_{j=1}^{m}\left(\int dv_{j}(x,D_{,j}^{-},\varepsilon)\delta_{\{v_{j}\}}^{(D_{j}^{-})}(x,\varepsilon)\right)_{\varepsilon} = m. \quad (7.1.13)$$

In particular for the case $m = 1$

$$\left(\int dv(x,D_{,j}^{-},\varepsilon)\delta_{\{v\}}^{(D^{-})}(x,\varepsilon)\right)_{\varepsilon} = 1. \quad (7.1.14)$$

Invariance of the action under the infinitesimal shift $\varphi(x) \to \varphi(x) + \delta\varphi(x)$ gives the equation of motion for a generic weight $(v_{i,\varepsilon})_{\varepsilon}, i = 1,\ldots,m$ :

$$\left(\frac{\partial\mathcal{L}}{\partial\varphi_{\varepsilon}}\right)_{\varepsilon} - \sum_{i=1}^{m}\left(\left[\left(\frac{\partial_{\mu}v_{i,\varepsilon}}{v_{i,\varepsilon}}\right) + \frac{d}{dx^{\mu}}\right]\frac{\partial\mathcal{L}}{\partial(\partial_{\mu}\varphi_{\varepsilon})}\right)_{\varepsilon} = 0. \quad (7.1.15)$$

In particular for for the case $m = 1$ we obtain

$$\left(\frac{\partial\mathcal{L}}{\partial\varphi_{\varepsilon}}\right)_{\varepsilon} - \left(\left[\frac{\partial_{\mu}v_{\varepsilon}}{v_{\varepsilon}} + \frac{d}{dx^{\mu}}\right]\frac{\partial\mathcal{L}}{\partial(\partial_{\mu}\varphi_{\varepsilon})}\right)_{\varepsilon} = 0. \quad (7.1.16)$$

From Eq.(7.1.11) and Eq.(7.1.15) we obtain

$$(\Box\varphi_{\varepsilon})_{\varepsilon} + \sum_{i=1}^{m}\left\{\left(\left[\frac{\partial_{\mu}v_{i,\varepsilon}}{v_{i,\varepsilon}}\right]\partial^{\mu}\varphi_{\varepsilon}\right)_{\varepsilon}\right\} - \left(\frac{d}{d\varphi_{\varepsilon}}V(\varphi_{\varepsilon})\right)_{\varepsilon} = 0. \quad (7.1.17)$$

where $\Box = \partial_{\mu}\partial^{\mu}$. In particular for for the case $m = 1$ we obtain

$$(\Box\varphi_{\varepsilon})_{\varepsilon} + \left(\left[\frac{\partial_{\mu}v_{\varepsilon}}{v_{\varepsilon}}\right]\partial^{\mu}\varphi_{\varepsilon}\right)_{\varepsilon} - \left(\frac{d}{d\varphi_{\varepsilon}}V(\varphi_{\varepsilon})\right)_{\varepsilon} = 0. \quad (7.1.18)$$

Note that the Hamiltonian is no longer an integral of motion. Let us define now the momentum

$$(\pi_{\varphi,\varepsilon})_{\varepsilon} \triangleq \left(\frac{\delta\mathbf{S}_{\varepsilon}}{\delta\phi_{\varepsilon}}\right)_{\varepsilon} = (\dot{\varphi}_{\varepsilon})_{\varepsilon}, \quad (7.1.19)$$

where dots indicate (total) derivatives with respect to time and we have taken Eq. (7.1.13) into account. Defining the Lagrangian in negative dimensions $D_{i}^{-}, i = 1,\ldots,m$

$$(L_{\varepsilon})_{\varepsilon} = \sum_{i=1}^{m}\left(\int d^{D_{\mathbf{t}}-1}\mathbf{x}v_{i,\varepsilon}(x,D_{i}^{-})\mathcal{L}\right)_{\varepsilon}, \quad (7.1.20)$$

where $\mathbf{x} = (x_1, \ldots, x_{D_t - 1}), d^{D_t - 1}\mathbf{x} = dx_1 \times \ldots \times dx_{D_t - 1}$ the Hamiltonian is

$$\begin{gathered}(\mathbf{H}_\varepsilon)_\varepsilon = \sum_{i=1}^{m}\left(\int d^{D_i^-}\mathbf{x}v_{i,\varepsilon}(x, D_i^-)(\pi_{\varphi,\varepsilon}\dot{\varphi}_\varepsilon)\right)_\varepsilon - (L_\varepsilon)_\varepsilon = \\ \sum_{i=1}^{m}\left(\int d^{D_i^-}\mathbf{x}v_{i,\varepsilon}(x, D_i^-)\left(\frac{1}{2}\pi_{\varphi,\varepsilon}^2 + \frac{1}{2}\partial_i\varphi\partial^i\varphi + V_\varepsilon\right)\right)_\varepsilon = \\ \sum_{i=1}^{m}\left(\int d^{D_i^-}\mathbf{x}v_{i,\varepsilon}(x, D_i^-)\widetilde{H}_\varepsilon\right)_\varepsilon, \\ \left(\widetilde{H}_\varepsilon\right)_\varepsilon = \left(\frac{1}{2}\pi_{\varphi,\varepsilon}^2 + \frac{1}{2}\partial_i\varphi\partial^i\varphi + V_\varepsilon\right)_\varepsilon.\end{gathered} \tag{7.1.21}$$

The definition of the equal-time Poisson brackets is

$$\begin{gathered}\left(\{A_\varepsilon(\mathbf{x}), B_\varepsilon(\mathbf{x}')\}_{\{v_{i,\varepsilon}\}}\right)_\varepsilon \triangleq \\ \sum_{i=1}^{m}\left(\int d^{D_t - 1}\mathbf{y}v_{i,\varepsilon}(\mathbf{y})\left[\frac{\delta A_\varepsilon(\mathbf{x}_i)}{\delta\varphi_\varepsilon(\mathbf{y}_i)}\frac{\delta B_\varepsilon(\mathbf{x}_i')}{\delta\pi_{\varphi_\varepsilon}(\mathbf{y}_i)} - \frac{\delta A_\varepsilon(\mathbf{x}_i)}{\delta\pi_{\varphi_\varepsilon}(\mathbf{y}_i)}\frac{\delta B_\varepsilon(\mathbf{x}_i')}{\delta\varphi_\varepsilon(\mathbf{y}_i)}\right]\right)_\varepsilon,\end{gathered} \tag{7.1.22}$$

gives the Hamilton equations

$$\begin{gathered}(\dot{\varphi}_\varepsilon)_\varepsilon = (\{\varphi_\varepsilon, \mathbf{H}_\varepsilon\})_\varepsilon, \\ (\dot{\pi}_{\varphi_\varepsilon})_\varepsilon = (\{\pi_{\varphi_\varepsilon}, \mathbf{H}_\varepsilon\})_\varepsilon - \sum_{i=1}^{m}\left(\frac{\dot{v}_{i,\varepsilon}}{v_{i,\varepsilon}}\pi_{\varphi_\varepsilon}\right)_\varepsilon\end{gathered} \tag{7.1.22}$$

equivalent to eqs. (7.1.19) and (7.1.14), respectively. Therefore, time evolution of an observable $(O_\varepsilon)_\varepsilon = (O(\varphi_\varepsilon, \pi_{\varphi_\varepsilon}, x, \varepsilon))_\varepsilon$ is

$$\left(\frac{dO_\varepsilon}{dt}\right)_\varepsilon = (\partial_t O_\varepsilon)_\varepsilon + \left(\{O_\varepsilon, \mathbf{H}_\varepsilon\}_{\{v_{i,\varepsilon}\}}\right)_\varepsilon - \left(\frac{\dot{v}_{i,\varepsilon}}{v_{i,\varepsilon}}\pi_{\varphi_\varepsilon}\frac{\partial O_\varepsilon}{\partial\pi_{\varphi_\varepsilon}}\right)_\varepsilon. \tag{7.1.23}$$

# 7.2.Propagator in configuration space with Hausdorff-Colombeau negative dimensions.

In this subsection we define the generalized vacuum-to-vacuum amplitude in Hausdorff-Colombeau negative dimensions by

$$\begin{gathered}(Z[J_\varepsilon])_\varepsilon = \left(\int_{(\varphi_\varepsilon)_\varepsilon \in \mathcal{G}(\mathbb{R}_k^D)}\mathbf{D}[\varphi_\varepsilon]\exp\left[i\sum_{j=1}^{m}\int_{\mathbb{R}_x^D}d\eta_{j,\varepsilon}(\mathcal{L}_\varepsilon(\varphi_\varepsilon) + \varphi_\varepsilon J_\varepsilon)\right]\right)_\varepsilon = \\ \int_{(\varphi_\varepsilon)_\varepsilon \in \mathcal{G}(\mathbb{R}_k^D)}\mathbf{D}[(\varphi_\varepsilon)_\varepsilon]\exp\left[i\sum_{j=1}^{m}\left(\int_{\mathbb{R}_x^D}d\eta_{j,\varepsilon}(\mathcal{L}_\varepsilon(\varphi_\varepsilon) + \varphi_\varepsilon J_\varepsilon)\right)_\varepsilon\right] = \\ \int_{(\varphi_\varepsilon)_\varepsilon \in \mathcal{G}(\mathbb{R}_k^D)}\mathbf{D}[(\varphi_\varepsilon)_\varepsilon]\exp\left[i\sum_{j=1}^{m}\int_{\mathbb{R}_x^D}(d\eta_{j,\varepsilon})_\varepsilon((\mathcal{L}_\varepsilon(\varphi_\varepsilon))_\varepsilon + (\varphi_\varepsilon)_\varepsilon(J_\varepsilon)_\varepsilon)\right]\end{gathered} \tag{7.2.1}$$

where $(\varphi_\varepsilon)_\varepsilon \in \mathcal{G}(\mathbb{R}_k^D)$ and $(J_\varepsilon)_\varepsilon \in \mathcal{G}(\mathbb{R}_k^D)$ is a source. Integration by parts in the exponent leads to the Lagrangian density for a free field as

$$(\mathcal{L}_\varepsilon)_\varepsilon = \frac{1}{2}\left(\varphi_\varepsilon\left(\Box + \sum_{j=1}^{m}\frac{\partial_\mu v_{j,\varepsilon}}{v_{j,\varepsilon}}\partial^\mu - m^2\right)\varphi_\varepsilon\right)_\varepsilon = \frac{1}{2}(\varphi_\varepsilon\mathcal{K}_\varepsilon\varphi_\varepsilon)_\varepsilon, \tag{7.2.2}$$

where

$$\mathcal{K}_\varepsilon = \Box + \sum_{j=1}^{m}\frac{\partial_\mu v_{j,\varepsilon}}{v_{j,\varepsilon}}\partial^\mu - m^2; j = 1, \ldots, m. \tag{7.2.3}$$

In particular for for the case $m = 1$ we obtain

$$\mathcal{K}_\varepsilon = \Box + \frac{\partial_\mu v_\varepsilon}{v_\varepsilon}\partial^\mu - m^2. \quad (7.2.4)$$

The propagator is the Green function $(G_\varepsilon(x))_\varepsilon$ solving the equation

$$(\mathcal{K}_\varepsilon G_\varepsilon(x))_\varepsilon = (\delta_v^{D^-}(x,\varepsilon))_\varepsilon, \quad (7.2.5)$$

where $D^- = D_\mathbf{t}(\alpha - 1) < 0$. By virtue of Lorentz covariance, the Green function $(G_\varepsilon(x))_\varepsilon$ must depend only on the Lorentz interval $s^2 = x_\mu x^\mu = x_i x^i - t^2$, where $x^0 = t$ and $i = 1,\ldots,D_\mathbf{t} - 1$. In particular, $(v_\varepsilon(x))_\varepsilon = (v_\varepsilon(s(x)))_\varepsilon$ with the correct scaling property is

$$(v_\varepsilon(s(x)))_\varepsilon = \left(\left(|s(x)|^{|D^-|} + \varepsilon\right)^{-1}\right)_\varepsilon, s(x) = \sqrt{x_\mu x^\mu}. \quad (7.2.6)$$

Note that

$$\partial_\mu = \frac{x_\mu}{(s(x)+\varepsilon)_\varepsilon}\partial_s, \Box = \partial_s^2 + \frac{D_\mathbf{t} - 1}{(s(x)+\varepsilon)_\varepsilon}\partial_s. \quad (7.2.7)$$

Hence the inhomogeneous equation (7.2.5) reads

$$\left(\partial_s^2 + \frac{D_\mathbf{t}\alpha - 1}{(s(x)+\varepsilon)_\varepsilon}\partial_s - m^2\right)(G_\varepsilon(x))_\varepsilon = (\delta_v^{D^-}(x,\varepsilon))_\varepsilon. \quad (7.2.8)$$

We first consider the Euclidean propagator and denote with $r = \sqrt{x_i x^i + t^2}$ the Wick-rotated Lorentz invariant. In the massless case, the solution of the homogeneous equations for any $\alpha < 0$ is

$$(G_\varepsilon(r))_\varepsilon = \mathcal{C}r^{2\beta}, \beta = \frac{2 + D|\alpha|}{2}, \quad (7.2.9)$$

where $D = D_\mathbf{t}$ and where $\mathcal{C}$ is a normalization constant. The right-hand side of Eq.(7.2.8) is not the usual $\delta$-function defined in radial coordinates. To find the latter note that

$$\begin{aligned} 1 &= \int_{\mathbb{R}_x^D} (d^D x v_\varepsilon(x))_\varepsilon \delta^{(D^-)}_{(v_\varepsilon)_\varepsilon}(x) = \left(\int_{\mathbb{R}_x^D} d^D x v_\varepsilon(x)\delta^{(D^-)}_{v_\varepsilon}(x)\right)_\varepsilon \\ &= \Omega_D\Omega_{D^-}\int_{\mathbb{R}_x^D} dr v_\varepsilon(r)\, r^{D-1}\delta^{(D^-)}_{(v_\varepsilon)_\varepsilon}(x) = \int_{\mathbb{R}_x^D} dr\,\delta(r), \end{aligned} \quad (7.2.10)$$

where $\delta^{(D^-)}_{(v_\varepsilon)_\varepsilon}(x) \triangleq \left(\delta^{(D^-)}_{v_\varepsilon}(x)\right)_\varepsilon$, $\Omega_D = 2\pi^{D/2}/\Gamma(D/2)$ and $\Omega_{D^-} = 2\pi^{D^-/2}/\Gamma(D^-/2)$, (see

Definition 6.3.8.) Therefore,

$$\begin{aligned} &\int_{\mathbb{R}_x^D} (d\varrho_\varepsilon(x))_\varepsilon \delta^{(D^-)}_{(v_\varepsilon)_\varepsilon}(x) = \left(\int_{\mathbb{R}_x^D} d\varrho_\varepsilon(x)\delta^{(D^-)}_{v_\varepsilon}(x)\right)_\varepsilon \\ &\int_{\mathbb{R}_x^D}(d\varrho_\varepsilon(x))_\varepsilon\left[\frac{r^{1-D}}{\Omega_D\Omega_{D^-}(v_\varepsilon(r))_\varepsilon}\delta(r)\right] \end{aligned} \quad (7.2.11)$$

In order to find the propagator also for $r = 0$ we can take some test function $\phi$ and compute

$$\Omega_D\Omega_{D^-}\langle\mathcal{C}_\varepsilon G, \phi\rangle = \left(\lim_{\epsilon\to 0}\Omega_D\Omega_{D^-}\int_\epsilon^{+\infty} dr\,\mathcal{C}_\varepsilon G(r)\phi(r)\right)_\varepsilon, \quad (7.2.12)$$

where

$$\begin{aligned} (\mathcal{C}_\varepsilon)_\varepsilon &= (v_\varepsilon(r)r^{D-1}\mathcal{K}_\varepsilon)_\varepsilon = \partial_r(r^{-(D|\alpha|+1)}\partial_r) - r^{-(D|\alpha|+1)}m^2 \\ &= \partial_r(r^{-(D|\alpha|+1)}\partial_r) - r^{-(D|\alpha|+1)}m^2 = \partial_r(r^{D^-}\partial_r) - r^{D^-}m^2. \end{aligned} \quad (7.2.13)$$

where $D^- = -(D|\alpha| + 1)$. Thus

$$\Omega_D\Omega_{D^-}\langle\mathcal{K}_0 G,\phi\rangle = \lim_{\epsilon\to 0}\Omega_D\Omega_{D^-}\int_\epsilon^{+\infty} dr\,G(r)\partial_r(r^{-(D\alpha+1)}\partial_r\phi(r))$$
$$= C'\Omega_D\Omega_{D^-}(2+D|\alpha|)\phi(0)\,, \tag{7.2.14}$$

where $C'$ is an constant and where we have used Eq.(7.2.9) and integrated by parts once (since the boundary terms vanish). The last line must be equal to $\langle\delta,\phi\rangle$, thus fixing $C'$. Then, the Green function for $m=0$ finally reads

$$G(r) = \frac{1}{\Omega_D\Omega_{D^-}(2+D|\alpha|)}(r^2)^{\frac{2+D|\alpha|}{2}}\,. \tag{7.2.15}$$

Let us consider now the massive case.The solution of the homogeneous equation $(\mathcal{K}_\varepsilon G_\varepsilon(r))_\varepsilon = 0$ for any $\alpha < 0$ is

$$(G_\varepsilon(r))_\varepsilon = \left(\frac{r}{m}\right)^{\frac{2+D\cdot|\alpha|}{2}}\left[\mathcal{C}_1 K_{-\frac{2+D\cdot|\alpha|}{2}}(mr) + \mathcal{C}_2 I_{-\frac{2+D\cdot|\alpha|}{2}}(mr)\right], \tag{7.2.16}$$

where $D = D_{\mathbf{t}}$, $\mathcal{C}_1,\mathcal{C}_2$ are constants and $K_\nu$ and $I_\nu$ are the modified Bessel functions. The function $I_\nu(z)$ is

$$I_\nu(z) = \sum_{k=0}^\infty \frac{(z/2)^{\nu+2k}}{k!\Gamma(\nu+k+1)}. \tag{7.2.17}$$

Formula (7.2.17) is valid providing $\nu \neq -1,-2,-3,\ldots.$

$$I_{-|\nu|}(z) = \sum_{k=0}^\infty \frac{(z/2)^{-|\nu|+2k}}{k!\Gamma(-|\nu|+k+1)} \tag{7.2.18}$$

Formula (7.2.18) is obtained by replacing $\nu$ in (7.2.17) with a $-\nu$. Note that

$$K_\nu(z) = \frac{\pi}{2\sin\nu\pi}[I_{-\nu}(z) - I_\nu(z)]. \tag{7.2.19}$$

It follows

$$K_{-|\nu|}(z) = -\frac{\pi}{2\sin|\nu|\pi}[I_{|\nu|}(z) - I_{-|\nu|}(z)]. \tag{7.2.20}$$

The modified Bessel functions $I_{|\nu|}(z)$ and $I_{-|\nu|}(z)$ have the following asymptotic forms for $z \to 0$ :

$$I_{|\nu|}(z) \simeq \frac{1}{\Gamma(|\nu|+1)}\left(\frac{z}{2}\right)^{|\nu|}, I_{-|\nu|}(z) \simeq \frac{1}{\Gamma(-|\nu|+1)}\left(\frac{z}{2}\right)^{-|\nu|},$$
$$\nu \neq -1,-2,-3,\ldots. \tag{7.2.21}$$

From (7.2.16) we obtain

$$(G_\varepsilon(r))_\varepsilon = \left(\frac{r}{m}\right)^{\frac{2+D\cdot|\alpha|}{2}}\left[\mathcal{C}_1 K_{-\frac{2+D\cdot|\alpha|}{2}}(mr) + \mathcal{C}_2 I_{-\frac{2+D\cdot|\alpha|}{2}}(mr)\right] =$$
$$\left(\frac{r}{m}\right)^{\frac{2+D\cdot|\alpha|}{2}}\left\{-\frac{\mathcal{C}_1\pi}{2\sin\frac{2+D\cdot|\alpha|}{2}\pi}\left[I_{\frac{2+D\cdot|\alpha|}{2}}(z) - I_{-\frac{2+D\cdot|\alpha|}{2}}(mr)\right] + \mathcal{C}_2 I_{-\frac{2+D\cdot|\alpha|}{2}}(mr)\right\} =$$
$$\left(\frac{r}{m}\right)^{\frac{2+D\cdot|\alpha|}{2}}\left\{-\frac{\mathcal{C}_1\pi}{2\sin\frac{2+D\cdot|\alpha|}{2}\pi}I_{\frac{2+D\cdot|\alpha|}{2}}(mr) + I_{-\frac{2+D\cdot|\alpha|}{2}}(mr)\left[\frac{\mathcal{C}_1\pi}{2\sin\frac{2+D\cdot|\alpha|}{2}\pi} + \mathcal{C}_2\right]\right\}. \tag{7.2.22}$$

Since for small $m \simeq 0$ the solution must agree with the massless case (4.2.15), we set

$$\frac{\mathcal{C}_1\pi}{2\sin\frac{2+D\cdot|\alpha|}{2}\pi}+\mathcal{C}_2=0. \qquad (4.2.23)$$

Therefore

$$(G_\varepsilon(r))_\varepsilon=-\mathcal{C}_1\left(\frac{r}{m}\right)^{\frac{2+D\cdot|\alpha|}{2}}\frac{\pi}{2\sin\frac{2+D\cdot|\alpha|}{2}\pi}I_{\frac{2+D\cdot|\alpha|}{2}}(mr). \qquad (4.2.24)$$

To find the solution of the inhomogeneous equation, one exploits the fact that the mass term does not contribute near the origin. Expanding Eq. (7.2.24) at $mr\simeq 0$, we find

$$(G_\varepsilon(r))_\varepsilon=-\mathcal{C}_1\left(\frac{r}{m}\right)^{\frac{2+D\cdot|\alpha|}{2}}\left(\frac{mr}{2}\right)^{\frac{2+D\cdot|\alpha|}{2}}\frac{\pi}{2\Gamma\left(\frac{2+D\cdot|\alpha|}{2}+1\right)\sin\frac{2+D\cdot|\alpha|}{2}\pi} \qquad (7.2.25)$$

which must coincide with Eq.(7.2.15). This gives the coefficient $C_1$ and the propagator reads

$$\begin{aligned}G(r)&=\frac{1}{2\pi^{\frac{D}{2}}\pi^{\frac{D^-}{2}}}\frac{\Gamma\left(\frac{D}{2}\right)\Gamma\left(\frac{D\cdot|\alpha|}{2}\right)}{\Omega_{D^-}}\left(\frac{r}{m}\right)^{\frac{2+D\cdot|\alpha|}{2}}I_{\frac{2+D|\alpha|}{2}}(mr)=\\&=\Upsilon(D,D^-)\left(\frac{r}{m}\right)^{\frac{2+D\cdot|\alpha|}{2}}I_{\frac{2+D|\alpha|}{2}}(mr),\end{aligned} \qquad (7.2.26)$$

where we let

$$\Upsilon(D,D^-)=\frac{1}{2\pi^{\frac{D}{2}}\pi^{\frac{D^-}{2}}}\frac{\Gamma\left(\frac{D}{2}\right)\Gamma\left(\frac{D\cdot|\alpha|}{2}\right)}{\Omega_{D^-}}. \qquad (7.2.27)$$

We can analytically continue the Helmholtz propagator (7.2.26) to the Klein–Gordon propagator according to the prescriptions: (i) multiply $G$ times the imaginary unit $i$, due to Wick rotation of the time direction; (ii) replace $r^2$ with $s^2+i\varepsilon$, where the positive sign of the extra infinitesimal term corresponds to the causal Feynman propagator. Thus Feynman propagator reads

$$G(s)=-i\Upsilon(D,D^-)\left(\frac{m^2}{s^2+i\varepsilon}\right)^{-\frac{D|\alpha|}{4}-\frac{1}{2}}I_{\frac{D|\alpha|}{2}+1}\left(m\sqrt{s^2+i\varepsilon}\right),\qquad s>0, \quad (7.2.28)$$

The massless propagator is

$$G(s)=\frac{i}{\Omega_D\Omega_{D^-}(2+D\alpha)}(s^2+i\varepsilon)^{1+\frac{D\cdot|\alpha|}{2}}, \qquad (7.2.29)$$

In order to calculate the propagator in momentum space, we can start from the Euclidean one and then analytically continue the result as usual. In the Lorentzian propagators the substitution $k^2\to|\mathbf{k}|^2-(k^0)^2-i\varepsilon$ is understood.Massless propagator in momentum space is given by Fourier transform (in the sense of generalized functions, see [32]) of the Eq.(7.2.15)

$$\begin{aligned}\widehat{G}(k)&=\frac{1}{\Omega_D\Omega_{D^-}(D|\alpha|+2)}\mathcal{F}[r^{2+D|\alpha|}]=\\&\frac{(2\pi)^{\frac{D}{2}}2^{\frac{2+D(|\alpha|+1)}{2}}}{\Omega_D\Omega_{D^-}\Gamma\left(-\frac{2+D|\alpha|}{2}\right)(D|\alpha|+2)}\frac{1}{k^{2+D(|\alpha|+1)}},\end{aligned} \qquad (7.2.30)$$

where $D^-=-D(|\alpha|+1)$. The Fourier transform of the massive propagator in radial coordinates is

$$\widehat{G}(k,m^2) = \Omega_D \int dr\, r^{D-1}\, G(r,m^2)\, e^{-ik\bullet x}. \tag{7.2.31}$$

The integrand in Eq.(7.3.21) is not radial but we can choose a frame where $k_\mu x^\mu = -kr\cos\theta,\ k = |k^\mu|$, and the angular integral is

$$\begin{gathered}\int d\Omega_D\, e^{-ik\bullet x} = \Omega_{D-1}\int_0^\pi d\theta(\sin\theta)^{D-2} e^{ikr\cos\theta} = \\ \Omega_{D-1}\sqrt{\pi}\,\Gamma\left(\frac{D-1}{2}\right)\left(\frac{2}{kr}\right)^{\frac{D}{2}-1} J_{\frac{D}{2}-1}(kr),\end{gathered} \tag{7.2.32}$$

where $J_\gamma(\beta) = \int_0^\pi e^{i\beta\cos u}(\sin u)^{2\gamma} du$. Thus we obtain

$$\widehat{G}(k,m^2) = \Omega_D \Gamma\left(\frac{D}{2}\right)\left(\frac{2}{k}\right)^{\frac{D}{2}-1}\int_0^{+\infty} dr\, r^{D-1}[G(r,m^2)] J_{\frac{D}{2}-1}(kr). \tag{7.2.33}$$

Now we take the massive propagator $G(r,m^2)$ (7.2.26)

$$\begin{gathered}\widehat{G}(k,m^2) = \Upsilon(D,D^-)\Gamma\left(\tfrac{D}{2}\right)\left(\frac{2}{k}\right)^{\frac{D}{2}-1}\left(\frac{m}{2}\right)^{-\frac{D|\alpha|}{2}-1}\times \\ \int_0^{+\infty} dr\, r^{\frac{2+D(|\alpha|+1)}{2}+1} I_{\frac{D|\alpha|}{2}+1}(mr) J_{\frac{D}{2}-1}(kr).\end{gathered} \tag{7.2.34}$$

# 7.3.Green’s functions in spacetime with Hausdorff-Colombeau negative dimensions.

We consider now a self-interecting scalar field $(\varphi_\varepsilon)_\varepsilon \in \mathcal{G}(\mathbb{R}^D_k)$ describing by the action

$$(S_\varepsilon)_\varepsilon = \int_{\mathbb{R}^D_x} dv_\varepsilon\left[\frac{1}{2}(\partial_\mu\varphi_\varepsilon\partial^\mu\varphi_\varepsilon)_\varepsilon - \frac{1}{2}m^2(\varphi^2_\varepsilon)_\varepsilon + (V_\varepsilon(\varphi_\varepsilon))_\varepsilon\right], \tag{7.3.1}$$

where

$$(V_\varepsilon(x))_\varepsilon \in \mathcal{G}(\mathbb{R}^D). \tag{7.3.2}$$

Corresponding generalized vacuum-to-vacuum amplitude in Hausdorff-Colombeau negative dimensions reads

$$\begin{gathered}(Z_M[J_\varepsilon])_\varepsilon = \\ \int_{(\varphi_\varepsilon)_\varepsilon\in\tilde{\mathcal{G}}(\mathbb{R}^D_x)} \mathbf{D}[(\varphi_\varepsilon)_\varepsilon]\times \\ \exp\left[i\left(\int_{\mathbb{R}^D_x} dv_\varepsilon\left(\frac{1}{2}(\partial_\mu\varphi_\varepsilon\partial^\mu\varphi_\varepsilon)_\varepsilon - \frac{1}{2}m^2(\varphi^2_\varepsilon)_\varepsilon + (V_\varepsilon(\varphi_\varepsilon))_\varepsilon + \varphi_\varepsilon J_\varepsilon\right)\right)_\varepsilon\right] = \\ \int_{(\varphi_\varepsilon)_\varepsilon\in\tilde{\mathcal{G}}(\mathbb{R}^D_x)} \mathbf{D}[(\varphi_\varepsilon)_\varepsilon]\times \\ \exp\left[i\int_{\mathbb{R}^D_x}(dv_\varepsilon)_\varepsilon\left(\frac{1}{2}(\partial_\mu\varphi_\varepsilon\partial^\mu\varphi_\varepsilon)_\varepsilon - \frac{1}{2}m^2(\varphi^2_\varepsilon)_\varepsilon + (V(\varphi_\varepsilon))_\varepsilon + (\varphi_\varepsilon)_\varepsilon(J_\varepsilon)_\varepsilon\right)\right]\end{gathered} \tag{7.3.3}$$

where $\tilde{\mathcal{G}}(\mathbb{R}^D_x) \subsetneqq \mathcal{G}(\mathbb{R}^D_x)$ is an topological linear subspace of Colombeau algebra $\mathcal{G}(\mathbb{R}^D_k)$, $(\varphi_\varepsilon)_\varepsilon \in \tilde{\mathcal{G}}(\mathbb{R}^D_x)$ and $(J_\varepsilon)_\varepsilon \in \tilde{\mathcal{G}}(\mathbb{R}^D_x)$ is a source.

**Remark 7.3.1.**Note that in (7.3.3) we integrate over an topological linear subspace $\tilde{\mathcal{G}}(\mathbb{R}^D_x) \subsetneqq \mathcal{G}(\mathbb{R}^D_x)$ of Colombeau algebra $\mathcal{G}(\mathbb{R}^D_x)$ but not over full Colombeau algebra $\mathcal{G}(\mathbb{R}^D_x)$.

We will be write for short the expression $(Z[J_\varepsilon])_\varepsilon$ in the following form

$$
\begin{aligned}
&(Z_M[J_\varepsilon])_\varepsilon = \mathbf{N}_M \int_{\widetilde{\mathcal{G}}(\mathbb{R}^D_x)} \mathbf{D}[(\varphi_\varepsilon)_\varepsilon] \times \\
&\exp\Big[ i\Big\langle \frac{1}{2}(\partial_\mu\varphi_\varepsilon\partial^\mu\varphi_\varepsilon)_\varepsilon - \frac{1}{2}m^2(\varphi_\varepsilon^2)_\varepsilon + (V_\varepsilon(\varphi_\varepsilon))_\varepsilon + (\varphi_\varepsilon)_\varepsilon(J_\varepsilon)\Big\rangle \Big],
\end{aligned}
\tag{7.3.4}
$$

where $\mathbf{N}_M$ is a normalizing constant, the $\langle\ldots\rangle_v$ now means integration with nontrivial measure $(d^D v_\varepsilon(x))_\varepsilon$ over spacetime, and $(J_\varepsilon)_\varepsilon \in \widetilde{\mathcal{G}}(\mathbb{R}^D_x)$ is a source. The integrand in (7.3.4) is oscillatory and even path integrals are not well defined. There are two canonical ways to resolve this problem:

(i) put in a convergence factor $\exp\big[-\frac{1}{2}\langle(\varphi_\varepsilon^2)_\varepsilon\rangle\big]$ with $\epsilon > 0$, or

(ii) define $(Z[J_\varepsilon])_\varepsilon$ in Euclidean space by setting $x_0 = i\bar{x}_0, d^D x = -id^D\bar{x}$, $(\partial_\mu\varphi_\varepsilon\partial^\mu\varphi_\varepsilon)_\varepsilon = -(\bar{\partial}_\mu\varphi_\varepsilon\bar{\partial}^\mu\varphi_\varepsilon)_\varepsilon$, where the bar denotes Euclidean space variables, $\bar{\partial}_\mu = \partial/\partial\bar{x}_\mu$ .

Then Eq.(7.3.4) becomes

$$
\begin{aligned}
&(Z_E[J_\varepsilon])_\varepsilon = \mathbf{N}_E \int_{\widetilde{\mathcal{G}}(\mathbb{R}^D_x)} \mathbf{D}[(\varphi_\varepsilon)_\varepsilon] \times \\
&\exp\Big[ -\Big\langle \frac{1}{2}(\bar{\partial}_\mu\varphi_\varepsilon\bar{\partial}^\mu\varphi_\varepsilon)_\varepsilon + \frac{1}{2}m^2(\varphi_\varepsilon^2)_\varepsilon - (V_\varepsilon(\varphi_\varepsilon))_\varepsilon - (\varphi_\varepsilon)_\varepsilon(J_\varepsilon)_\varepsilon\Big\rangle \Big]
\end{aligned}
, \tag{7.3.5}
$$

where for instance

$$
(V_\varepsilon(\varphi_\varepsilon))_\varepsilon = \sum_{k=3}^m c_k(\varphi_\varepsilon^k)_\varepsilon. \tag{7.3.6}
$$

The exponent of the integrand is now negative definite for positive $m$ and $V$.

In either case, the generating functional is used to manufacture the Green's functions which are the coefficients of the functional expansion

$$
\begin{aligned}
&(Z[J_\varepsilon])_\varepsilon = \sum_{N=0}^\infty \frac{i^N}{N!}\big\langle (J_{1,\varepsilon}J_{2,\varepsilon}\ldots J_{N,\varepsilon}G^{(N)}(1,2,\ldots,N;\varepsilon))_\varepsilon\big\rangle_{\rho_1,\rho_2,\ldots,\rho_N}. = \\
&\sum_{N=0}^\infty \frac{i^N}{N!}\big\langle (J_{1,\varepsilon})_\varepsilon(J_{2,\varepsilon})_\varepsilon\ldots(J_{N,\varepsilon})_\varepsilon(G^{(N)}(1,2,\ldots,N;\varepsilon))_\varepsilon\big\rangle_{\rho_1,\rho_2,\ldots,\rho_N}.
\end{aligned}
\tag{7.3.7}
$$

Thus

$$
\begin{aligned}
&(G^{(N)}(1,2,\ldots,N;\varepsilon))_\varepsilon = \\
&\frac{1}{i^N}\left(\frac{\delta}{\delta(J_{1,\varepsilon}\varrho_{1,\varepsilon})}\frac{\delta}{\delta(J_{2,\varepsilon}\varrho_{2,\varepsilon})}\cdots\frac{\delta}{\delta(J_{N,\varepsilon}\varrho_{N,\varepsilon})}Z[J_\varepsilon]\right)_\varepsilon\Bigg|_{(\mathbf{J}_\varepsilon)_\varepsilon=0} = \\
&\frac{1}{i^N}\frac{\delta}{\delta(J_{1,\varepsilon}\varrho_{1,\varepsilon})_\varepsilon}\frac{\delta}{\delta(J_{2,\varepsilon}\varrho_{2,\varepsilon})_\varepsilon}\cdots\frac{\delta}{\delta(J_{N,\varepsilon}\varrho_{N,\varepsilon})_\varepsilon}(Z[J_\varepsilon])_\varepsilon\Bigg|_{(\mathbf{J}_\varepsilon)_\varepsilon=0} = \\
&\frac{1}{i^N}\frac{\delta}{\delta(J^\rho_{1,\varepsilon})_\varepsilon}\frac{\delta}{\delta(J^\rho_{2,\varepsilon})_\varepsilon}\cdots\frac{\delta}{\delta(J^\rho_{N,\varepsilon})_\varepsilon}(Z[J_\varepsilon])_\varepsilon\Bigg|_{(\mathbf{J}^\rho_\varepsilon)_\varepsilon=0}
\end{aligned}
\tag{7.3.8}
$$

where $(J_{i,\varepsilon})_\varepsilon = (J_{i,\varepsilon}(x_i))_\varepsilon, (\varrho_{i,\varepsilon})_\varepsilon = (\varrho_{i,\varepsilon}(x_i))_\varepsilon, (J^\rho_{i,\varepsilon})_\varepsilon = (J_{i,\varepsilon}(x_i)\varrho_{i,\varepsilon}(x_i))_\varepsilon, i = 1,2,\ldots,N$, $(\mathbf{J}_\varepsilon)_\varepsilon = ((J_{1,\varepsilon})_\varepsilon,(J_{2,\varepsilon})_\varepsilon,\ldots,(J_{N,\varepsilon})_\varepsilon)$ , $(\mathbf{J}^\rho_\varepsilon)_\varepsilon = \big((J^\rho_{1,\varepsilon})_\varepsilon,(J^\rho_{2,\varepsilon})_\varepsilon,\ldots,(J^\rho_{N,\varepsilon})_\varepsilon\big)$ and $\langle\ldots\rangle_{\rho_1,\rho_2,\ldots,\rho_N}$. means integration with nontrivial Colombeau–Stieltjes measure $d^D\varrho_\varepsilon(x_1)\times d^D\varrho_\varepsilon(x_2)\times\ldots\times d^D\varrho_\varepsilon(x_N)$. We evaluate now $(Z_M[J_\varepsilon])_\varepsilon$ when $(V_\varepsilon(\varphi_\varepsilon))_\varepsilon = 0$. We choose to do it in Minkowski. Let $(Z^M_{0,\varepsilon}[J_\varepsilon])_\varepsilon$ be

$$(Z^{M}_{0,\varepsilon}[J_\varepsilon])_\varepsilon =$$

$$\mathbf{N}_M \int_{\mathcal{G}(\mathbb{R}^D_x)} \mathbf{D}[(\varphi_\varepsilon)_\varepsilon] \exp\left[ i\left\langle \frac{1}{2}(\partial_\mu\varphi_\varepsilon\partial^\mu\varphi_\varepsilon)_\varepsilon - \frac{1}{2}(m^2 - i\epsilon)(\varphi^2_\varepsilon)_\varepsilon + (\varphi_\varepsilon)_\varepsilon (J_\varepsilon)\right\rangle_\rho \right]. \tag{7.3.9}$$

We assume now that $(\phi_\varepsilon(x))_\varepsilon, (J_\varepsilon(x))_\varepsilon \in \widetilde{\mathcal{G}}(\mathbb{R}^D), (\varrho_\varepsilon(x))_\varepsilon \in \mathcal{G}(\mathbb{R}^D)$ and introduce the $D$-dimensional Colombeau Fourier–Stieltjes transform $(\widetilde{\varphi}_\varepsilon(k))_\varepsilon \triangleq (\mathcal{F}S)_\varrho[(\phi_\varepsilon(x))_\varepsilon](k)$, of the field $(\varphi_\varepsilon(x))_\varepsilon$ with weight $\left\{\left(\varrho^{(1)}_\varepsilon(x)\right)_\varepsilon, \left(\varrho^{(2)}_\varepsilon(k)\right)_\varepsilon\right\}$ using the following formal definitions: where the Colombeau-Fourier– Stieltjes transform $F_v$ of a function $(G_\varepsilon(x))_\varepsilon \in \widetilde{\mathcal{G}}(\mathbb{R}^D)$ and its inverse are defined as

$$\left(\widetilde{G}_\varepsilon(k)\right)_\varepsilon = \left(\int d\varrho_\varepsilon(x)\, G_\varepsilon(x)\, e^{-ik\cdot x}\right)_\varepsilon = \int (d\varrho_\varepsilon(x))_\varepsilon (G_\varepsilon(x))_\varepsilon e^{-ik\cdot x} \triangleq (F_{v_\varepsilon}[G_\varepsilon(x)])_\varepsilon \tag{7.3.10}$$

and

$$(G_\varepsilon(x))_\varepsilon = \frac{1}{(2\pi)^D}\left(\int d\varrho_\varepsilon(k)\, \widetilde{G}_\varepsilon(k)\, e^{ik\cdot x}\right)_\varepsilon = \frac{1}{(2\pi)^D}\int (d\varrho_\varepsilon(k))_\varepsilon \left(\widetilde{G}_\varepsilon(k)\right)_\varepsilon e^{ik\cdot x}. \tag{7.3.11}$$

correspondingly,where $k\cdot x = k^0x^0 - \vec{k}\cdot\vec{x}$. Using now the definition of the $D$-dimensional Colombeau–Dirac distribution with nontrivial Colombeau–Stieltjes measure $(d\varrho_\varepsilon(k))_\varepsilon$

$$(\delta_{\varrho_\varepsilon}(k))_\varepsilon = \left(\int_{\mathbb{R}^D_k} d\varrho_\varepsilon(x)\, e^{-ik\cdot x}\right)_\varepsilon \triangleq \int_{\mathbb{R}^D_k} (d\varrho_\varepsilon(x))_\varepsilon e^{-ik\cdot x} \tag{7.3.12}$$

the exponent of the integrand in Eq.(7.3.9) is easily expressed in terms of the Colombeau- -Fourier–Stieltjes transforms of $(\varphi_\varepsilon(x))_\varepsilon$ and $(J_\varepsilon(x))_\varepsilon$. Finally it reads

$$\frac{i}{2}\left(\int d\varrho_\varepsilon(-k)\left[\widetilde{\varphi}'_\varepsilon(k)(k^2 - m^2 + i\epsilon)\widetilde{\varphi}'(-k) - \widetilde{J}_\varepsilon(k)(k^2 - m^2 + i\epsilon)^{-1}\widetilde{J}_\varepsilon(-k)\right]\right)_\varepsilon, \tag{7.3.13}$$

where

$$\left(\widetilde{\varphi}'_\varepsilon(k)\right)_\varepsilon = (\widetilde{\varphi}_\varepsilon(k))_\varepsilon + (k^2 - m^2 + i\epsilon)\left(\widetilde{J}_\varepsilon(k)\right)_\varepsilon. \tag{7.3.14}$$

**Definition 7.3.1**. Let $\widetilde{\mathcal{G}}(\mathbb{R}^D_x) \subsetneq \mathcal{G}(\mathbb{R}^D_x)$ be a maximal topological linear subspace of the Colombeau algebra $\mathcal{G}(\mathbb{R}^D_k)$ such that for any $(\varphi_\varepsilon)_\varepsilon \in \widetilde{\mathcal{G}}(\mathbb{R}^D_x)$ bouth Eqs.(7.3.10)-(7.3.11) are satisfied.

**Remark 7.3.2**. Note that we willin to choose a source $\left(\widetilde{J}_\varepsilon(k)\right)_\varepsilon \in \widetilde{\mathcal{G}}(\mathbb{R}^D_k)$ such that $\left(\widetilde{\varphi}'_\varepsilon(k)\right)_\varepsilon \in \widetilde{\mathcal{G}}(\mathbb{R}^D_k), (\varphi'_\varepsilon(x))_\varepsilon \in \widetilde{\mathcal{G}}(\mathbb{R}^D_x)$ and therefore

$$\int_{\mathcal{G}(\mathbb{R}^D_k)} \mathbf{D}\left[(\widetilde{\varphi}_\varepsilon)_\varepsilon\right]\exp[\ldots] = \int_{\mathcal{G}(\mathbb{R}^D_k)} \mathbf{D}\left[\left(\widetilde{\varphi}'_\varepsilon\right)_\varepsilon\right]\exp[\ldots]. \tag{7.3.15}$$

Putting it all together, we obtain

$$\begin{gathered}(Z^M_{0,\varepsilon}[J_\varepsilon])_\varepsilon = \mathbf{N}_M \exp\left[-\frac{i}{2}\left(\int d\varrho_\varepsilon(-k)\frac{\widetilde{J}_\varepsilon(k)\widetilde{J}_\varepsilon(-k)}{k^2 - m^2 + i\epsilon}\right)_\varepsilon\right]\times \\ \int_{(\varphi'_\varepsilon)_\varepsilon \in \widetilde{\mathcal{G}}(\mathbb{R}^D_x)} \mathbf{D}[(\varphi'_\varepsilon)_\varepsilon]\exp\left[i\left\langle\frac{1}{2}(\partial_\mu\varphi'_\varepsilon\partial^\mu\varphi'_\varepsilon)_\varepsilon - \frac{1}{2}(m^2 - i\epsilon)(\varphi'^2_\varepsilon)_\varepsilon\right\rangle_\rho\right],\end{gathered} \tag{7.3.16}$$

where we observed that the $(\varphi'_\varepsilon)_\varepsilon$-dependent term was just the same as the $(\varphi_\varepsilon)_\varepsilon$ -dependent term in (4.3.9) with $(J_\varepsilon(x))_\varepsilon = 0$. Thus

$$(Z^M_{0,\varepsilon}[J_\varepsilon])_\varepsilon = (Z^M_{0,\varepsilon}[0])_\varepsilon \exp\left[-\frac{i}{2}\left(\int d\varrho_\varepsilon(-k)\frac{\widetilde{J}_\varepsilon(k)\widetilde{J}_\varepsilon(-k)}{k^2 - m^2 + i\epsilon}\right)_\varepsilon\right]. \tag{7.3.17}$$

By adjusting $\mathbf{N}_M$, we can take $(Z^M_{0,\varepsilon}[0])_\varepsilon = 1$. The important thing is that we have succeeded in finding the explicit dependence of $(Z^M_{0,\varepsilon}[J_\varepsilon])_\varepsilon$ on $(J_\varepsilon(x))_\varepsilon$. The use of the $(J_\varepsilon(x))_\varepsilon$ inverse Colombeau-Fourier–Stieltjes transform (7.3.11) yields

$$\left(\int_{\mathbb{R}^D_k} d\varrho_\varepsilon(-k)\frac{\widetilde{J}_\varepsilon(k)\widetilde{J}_\varepsilon(-k)}{k^2-m^2+i\epsilon}\right)_\varepsilon = \int_{\mathbb{R}^D_k}\frac{(d\varrho_\varepsilon(-k))_\varepsilon}{(2\pi)^D}\int_{\mathbb{R}^D_x}((d\varrho_\varepsilon(x))_\varepsilon)\int_{\mathbb{R}^D_y}((d\varrho(y))_\varepsilon)e^{ik\cdot(x-y)}\frac{\left(\left(\widetilde{J}_\varepsilon(-k)\right)_\varepsilon\right)\left(\left(\widetilde{J}_\varepsilon(k)\right)_\varepsilon\right)}{k^2-m^2+i\epsilon} \tag{7.3.18}$$

and so that, if $(\varrho_\varepsilon(-k))_\varepsilon = (\varrho_\varepsilon(k))_\varepsilon$, the free partition function reads

$$(Z_{0,\varepsilon}[J_\varepsilon])_\varepsilon = ((Z_{0,\varepsilon}[0])_\varepsilon)\times \exp\left[\frac{i}{2}\int_{\mathbb{R}^D_x}((d\varrho_\varepsilon(x))_\varepsilon)\int_{\mathbb{R}^D_y}((d\varrho_\varepsilon(y))_\varepsilon)(J_\varepsilon(x)((\Delta_F(x-y;\varepsilon))_\varepsilon)J_\varepsilon(y))_\varepsilon\right] = \exp\left[\frac{i}{2}\left\langle(J_\varepsilon(x)((\Delta_F(x-y;\varepsilon))_\varepsilon)J_\varepsilon(y))_\varepsilon\right\rangle_{\rho(x),\rho(y)}\right] \tag{7.3.19}$$

where

$$(\Delta_F(x-y;\varepsilon))_\varepsilon = \frac{1}{(2\pi)^D}\int_{\mathbb{R}^D_k}((d\varrho_\varepsilon(k))_\varepsilon)\frac{e^{ik\cdot(x-y)}}{k^2-m^2+i\epsilon}. \tag{7.3.20}$$

Thus, we have recovered the usual definition of the propagator as the solution of the Green equation in negative dimensions $D^- = D(\alpha-1), \alpha < 0$ :

$$(\Box + m^2)\,\Delta_F(x-y;\varepsilon))_\varepsilon = -(\delta_{\varrho_\varepsilon}(x-y)\,)_\varepsilon. \tag{7.3.21}$$

It is the Feynman propagator. We now interpret the Green's functions obtained from $(Z_{0,\varepsilon}[J])_\varepsilon$. From (7.3.8) we find

$$\begin{gathered}\left(G^{(2)}_0(x_1,x_2;\varepsilon)\right)_\varepsilon = (\Delta_F(x_1-x_2;\varepsilon))_\varepsilon\\ \left(G^{(4)}_0(x_1,x_2,x_3,x_4;\varepsilon)\right)_\varepsilon = \\ -[(\Delta_F(x_1-x_2;\varepsilon)\Delta_F(x_3-x_4;\varepsilon))_\varepsilon + (\Delta_F(x_1-x_3;\varepsilon)\Delta_F(x_2-x_4;\varepsilon))_\varepsilon + \\ +(\Delta_F(x_1-x_4;\varepsilon)\Delta_F(x_2-x_3;\varepsilon))_\varepsilon]\\ \text{etc. ...}\end{gathered} \tag{7.3.22}$$

together with the vanishing of the Green's functions with odd number of variables. This fact is easy to understand since $(Z_{0,\varepsilon}[J_\varepsilon])_\varepsilon$ depends only on $(J_\varepsilon)_\varepsilon$. In passing, note that all Green's functions are functions of only the difference of coordinates, reflecting the translation invariance of the theory. Another imprtent notes is that the higher Green's functions can all be represented only in terms of $\left(G^{(2)}_0(x_1,x_2;\varepsilon)\right)_\varepsilon$. Hence it would apper more convenient to set

$$(Z_{0,\varepsilon}[J_\varepsilon])_\varepsilon = (\exp[i\check{Z}_{0,\varepsilon}[J_\varepsilon]])_\varepsilon \tag{7.3.24}$$

and define new Green's functions in terms of $(\check{Z}_{0,\varepsilon}[J_\varepsilon])_\varepsilon$ :

$$(i\check{Z}[J_\varepsilon])_\varepsilon = \sum_{N=0}^{\infty}\frac{i^N}{N!}\left\langle\left(J_{1,\varepsilon}J_{2,\varepsilon}\ldots J_{N,\varepsilon}G^{(N)}_c(1,2,\ldots,N;\varepsilon)\right)_\varepsilon\right\rangle_{\rho_1,\rho_2,\ldots,\rho_N} .= \sum_{N=0}^{\infty}\frac{i^N}{N!}\left\langle(J_{1,\varepsilon})_\varepsilon(J_{2,\varepsilon})_\varepsilon\ldots(J_{N,\varepsilon})_\varepsilon\left(G^{(N)}_c(1,2,\ldots,N;\varepsilon)\right)_\varepsilon\right\rangle_{\rho_1,\rho_2,\ldots,\rho_N}. \tag{7.3.25}$$

Indeed since $(G^{(N)}(x_1,x_2,\ldots,x_N;\varepsilon))_\varepsilon$ depends only on differences of coordinates the Colombeau-Fourier transform

$$\int d^D x_1 \ldots d^D x_N \exp\left(-i\sum_{j=1}^N k_j x_j\right)(G^{(N)}(x_1,x_2,\ldots,x_N;\varepsilon))_\varepsilon \tag{7.3.26}$$

necessarily contains the usual $\delta$-function of $\sum_{j=1}^N k_j x_j$. So we can set

$$\begin{gathered}\left(\widehat{G}^{(N)}(p_1,p_2,\ldots,p_N;\varepsilon)\right)_\varepsilon = \\ (2\pi)^D \delta\left(\sum_{j=1}^N k_j\right) \int d^D x_1 \ldots d^D x_N \exp\left(-i\sum_{j=1}^N k_j x_j\right)(G^{(N)}(x_1,x_2,\ldots,x_N;\varepsilon))_\varepsilon\end{gathered} \tag{7.3.27}$$

with $\left(\widehat{G}^{(N)}(p_1,p_2,\ldots,p_N;\varepsilon)\right)_\varepsilon$ defined only when $p_1+p_2+\ldots+p_N=0$. For instance

$$\left(G_0^{(2)}(p,-p;\varepsilon)\right)_\varepsilon = \frac{(\varrho_\varepsilon(k))_\varepsilon}{k^2-m^2+i\epsilon} \tag{7.3.28}$$

gives the amplitude that a free particle of momentum $k$ and mass $m^2$ propagates in fractal spacetime. .

## 7.3.1.The Effective Action

Out of the generating functional we can construct local quantities which lend themselves to familiar interpretations. For instance,

$$\frac{\delta(Z_{0,\varepsilon}[J_\varepsilon])_\varepsilon}{\delta(J_\varepsilon^\rho(x))_\varepsilon} = -i\frac{\delta\langle (J_{1,\varepsilon}\Delta_c(x_1-x_2;\varepsilon)J_{2,\varepsilon})_\varepsilon\rangle_{\rho_1,\rho_2}}{\delta(J_\varepsilon^\rho(x))_\varepsilon}(Z_{0,\varepsilon}[J_\varepsilon])_\varepsilon \tag{7.3.29}$$

so that

$$\left(\varphi_{cl}^{(0)}(x;\varepsilon)\right)_\varepsilon \equiv -i\frac{\delta(\ln Z_{0,\varepsilon}[J_\varepsilon])_\varepsilon}{\delta(J_\varepsilon^\rho(x))_\varepsilon} = \frac{\delta(\breve{Z}_{0,\varepsilon}[J_\varepsilon])_\varepsilon}{\delta(J_\varepsilon^\rho(x))_\varepsilon} \tag{7.3.30}$$

(by using Eq.(7.3.21)) satisfies the classical equation of motion

$$(\Box+m^2)\left(\varphi_{cl}^{(0)}(x;\varepsilon)\right)_\varepsilon = (J_\varepsilon(x))_\varepsilon. \tag{7.3.31}$$

In fact, we can use Eq.(7.3.30) in order to replace $(J_\varepsilon(x))_\varepsilon$ in terms of $\left(\varphi_{cl}^{(0)}(x;\varepsilon)\right)_\varepsilon$. Formally it comes down to performing a functional Legendre transformation. Introducing

$$\left(\Gamma_{0,\varepsilon}[\varphi_{cl}^{(0)}(x;\varepsilon)]\right)_\varepsilon = (\breve{Z}_{0,\varepsilon}[J_\varepsilon])_\varepsilon - \left\langle\left(J_\varepsilon(x)\varphi_{cl}^{(0)}(x;\varepsilon)\right)_\varepsilon\right\rangle_\rho \tag{7.3.32}$$

we see by using Eq.(7.3.30) that is independent of $(J_\varepsilon(x))_\varepsilon$. In this case it is easy to find the explicit form of $\left(\Gamma_{0,\varepsilon}[\varphi_{cl}^{(0)}(x;\varepsilon)]\right)_\varepsilon$ by replacing $(J_\varepsilon(x))_\varepsilon$ in terms of $\left(\varphi_{cl}^{(0)}(x;\varepsilon)\right)_\varepsilon$. One finds using Eq.(7.3.21),Eq.(7.3.31) and integrating by parts we obtain

$$\begin{gathered}\left(\Gamma_{0,\varepsilon}[\varphi_{cl}^{(0)}(x;\varepsilon)]\right)_\varepsilon = \\ -\tfrac{1}{2}\left\langle\left[(\Box+m^2)\left(\varphi_{cl}^{(0)}(x_1;\varepsilon)\right)_\varepsilon\right](\Delta_F(x_1-x_2;\varepsilon))_\varepsilon\left[(\Box+m^2)\left(\varphi_{cl}^{(0)}(x_2;\varepsilon)\right)_\varepsilon\right]\right\rangle_{\rho_1,\rho_2} \\ -\left\langle\left(\varphi_{cl}^{(0)}(x;\varepsilon)\right)_\varepsilon(\Box+m^2)\left(\varphi_{cl}^{(0)}(x;\varepsilon)\right)_\varepsilon\right\rangle_\rho = \\ -\tfrac{1}{2}\left\langle\left(\varphi_{cl}^{(0)}(x;\varepsilon)\right)_\varepsilon(\Box+m^2)\left(\varphi_{cl}^{(0)}(x;\varepsilon)\right)_\varepsilon\right\rangle_\rho\end{gathered} \tag{7.3.33}$$

Integration by parts gives the final form

$$\left(\Gamma_{0,\varepsilon}[\varphi_{cl}^{(0)}(x;\varepsilon)]\right)_\varepsilon =$$
$$\frac{1}{2}\int (d\varrho_\varepsilon(x))_\varepsilon \left[\left(\partial_\mu \varphi_{cl}^{(0)}(x;\varepsilon)\partial^\mu \varphi_{cl}^{(0)}(x;\varepsilon) - m^2\left(\varphi_{cl}^{(0)}(x;\varepsilon)\right)^2\right)_\varepsilon\right]. \tag{4.3.34}$$

A very similar procedure can be carried out in the general case $V \neq 0$. We form now

$$(\varphi_{cl}(x;\varepsilon))_\varepsilon \equiv -i\frac{\delta(\ln Z_\varepsilon[J_\varepsilon])_\varepsilon}{\delta(J_\varepsilon^\rho(x))_\varepsilon} = \frac{\delta(\breve{Z}_\varepsilon[J_\varepsilon])_\varepsilon}{\delta(J_\varepsilon^\rho(x))_\varepsilon} \tag{7.3.34}$$

and try to compute the effective action in general case

$$(\Gamma_\varepsilon[\varphi_{cl}(x;\varepsilon)])_\varepsilon = (\breve{Z}_\varepsilon[J_\varepsilon])_\varepsilon - \langle (J_\varepsilon(x)\varphi_{cl}(x;\varepsilon))_\varepsilon\rangle_\rho \tag{7.3.35}$$

with now

$$(J_\varepsilon(x))_\varepsilon = \frac{\delta(\Gamma_\varepsilon[\varphi_{cl}(x;\varepsilon)])_\varepsilon}{\delta(\varphi_{cl}(x;\varepsilon))_\varepsilon} \tag{7.3.36}$$

as seen by differentiating Eq.(7.3.35) with respect to $(\varphi_{cl}(x;\varepsilon))_\varepsilon$.By the way, we observe that since $(\Gamma_\varepsilon[\varphi_{cl}(x;\varepsilon)])_\varepsilon$ is an effective action, (7.3.36) is proportional to its equation of motion comes, from extremizing $(\Gamma_\varepsilon[\varphi_{cl}(x;\varepsilon)])_\varepsilon$. In the $V = 0$ case this is obvious from Eq.(7.3.31).In order to derive an equation of motion for $(\varphi_{cl}(x;\varepsilon))_\varepsilon$, we have to write $(Z_\varepsilon^M[J_\varepsilon])_\varepsilon$ in a manageable form. We write

$$\begin{aligned}(Z_\varepsilon^M[J_\varepsilon])_\varepsilon &= \mathbf{N}_M \int_{\mathcal{G}(\mathbb{R}_x^D)} \mathbf{D}[(\varphi_\varepsilon)_\varepsilon] \times \\ \exp&\left[i\left\langle \frac{1}{2}(\partial_\mu\varphi_\varepsilon\partial^\mu\varphi_\varepsilon)_\varepsilon - \frac{1}{2}(m^2 - i\epsilon)(\varphi_\varepsilon^2)_\varepsilon + V((\varphi_\varepsilon)_\varepsilon) + (\varphi_\varepsilon)_\varepsilon(J_\varepsilon)_\varepsilon\right\rangle_\rho\right] \\ &= \mathbf{N}_M \int_{\mathcal{G}(\mathbb{R}_x^D)} \mathbf{D}[(\varphi_\varepsilon)_\varepsilon]\exp\left[-i\langle V((\varphi_\varepsilon)_\varepsilon)\rangle_\rho\right] \times \\ \exp&\left[i\left\langle \frac{1}{2}(\partial_\mu\varphi_\varepsilon\partial^\mu\varphi_\varepsilon)_\varepsilon - \frac{1}{2}(m^2 - i\epsilon)(\varphi_\varepsilon^2)_\varepsilon + (\varphi_\varepsilon)_\varepsilon(J_\varepsilon)_\varepsilon\right\rangle_\rho\right].\end{aligned} \tag{7.3.37}$$

Now immediately comes the canonical trick: observe that

$$\frac{1}{i}\frac{\delta}{\delta(J_\varepsilon^\rho(x))_\varepsilon}\exp\left[i\langle(\varphi_\varepsilon)_\varepsilon(J_\varepsilon)_\varepsilon\rangle_\rho\right] = (\varphi_\varepsilon(x))_\varepsilon\exp\left[i\langle(\varphi_\varepsilon)_\varepsilon(J_\varepsilon)_\varepsilon\rangle_\rho\right] \tag{7.3.38}$$

and since $(J_\varepsilon)_\varepsilon$ and $(\varphi_\varepsilon)_\varepsilon$ are independent variables, the same will be true for any function $(V(\varphi_\varepsilon))_\varepsilon$ of $(\varphi_\varepsilon)_\varepsilon$. In particular

$$\begin{aligned}&\exp\left[-i\langle(\varphi_\varepsilon)_\varepsilon(J_\varepsilon)_\varepsilon\rangle_\rho\right]\exp\left[i\langle(\varphi_\varepsilon)_\varepsilon(J_\varepsilon)_\varepsilon\rangle_\rho\right] = \\ &\exp\left[-i\left\langle V\left(\frac{1}{i}\frac{\delta}{\delta(J_\varepsilon^\rho(x))_\varepsilon}\right)\right\rangle_\rho\right]\exp\left[i\langle(\varphi_\varepsilon)_\varepsilon(J_\varepsilon)_\varepsilon\rangle_\rho\right].\end{aligned} \tag{7.3.39}$$

This allows us to take the $V$ dependent term out of the path integral

$$(Z_\varepsilon^M[J_\varepsilon])_\varepsilon = \exp\left[-i\left\langle V\left(\frac{1}{i}\frac{\delta}{\delta(J_\varepsilon^\rho(x))_\varepsilon}\right)\right\rangle_\rho\right]\times$$
$$\mathbf{N}_M\int_{\tilde{\mathcal{G}}(\mathbb{R}_x^D)}\mathbf{D}[(\varphi_\varepsilon)_\varepsilon]\exp\left[i\left\langle\frac{1}{2}(\partial_\mu\varphi_\varepsilon\partial^\mu\varphi_\varepsilon)_\varepsilon - \frac{1}{2}(m^2-i\epsilon)(\varphi_\varepsilon^2)_\varepsilon + (\varphi_\varepsilon)_\varepsilon(J_\varepsilon)_\varepsilon\right\rangle_\rho\right] \quad (7.3.40)$$
$$= \exp\left[-i\left\langle V\left(\frac{1}{i}\frac{\delta}{\delta(J_\varepsilon^\rho(x))_\varepsilon}\right)\right\rangle_\rho\right](Z_{0,\varepsilon}[J_\varepsilon])_\varepsilon$$

or

$$\exp\left[i(\check{Z}_{0,\varepsilon}[J_\varepsilon])_\varepsilon\right] = (Z_\varepsilon^M[J_\varepsilon])_\varepsilon = \exp\left[-i\left\langle V\left(\frac{1}{i}\frac{\delta}{\delta(J_\varepsilon^\rho(x))_\varepsilon}\right)\right\rangle_\rho\right] \quad (7.3.41)$$
$$\times\exp\left[\frac{i}{2}\left\langle(J_\varepsilon(x_1)((\Delta_F(x_1-x_2;\varepsilon))_\varepsilon)J_\varepsilon(x_2))_\varepsilon\right\rangle_{\rho_1,\rho_2}\right]$$

The equation Eq.(7.3.41) will be the starting point of the perturbative evaluation of $(Z_\varepsilon^M[J_\varepsilon])_\varepsilon$. For the moment, we use it to derive an equation for $(\varphi_{cl}(x;\varepsilon))_\varepsilon$ . From Eq.(7.3.41)

$$\frac{\delta(Z_\varepsilon^M[J_\varepsilon])_\varepsilon}{\delta(J_\varepsilon^\rho(x))_\varepsilon} = -i\exp\left[i\left\langle V\left(-i\frac{\delta}{\delta(J_\varepsilon^\rho)_\varepsilon}\right)\right\rangle_\rho\right]\langle(\Delta_F(x-x_1;\varepsilon)J_{1,\varepsilon})_\varepsilon\rangle_{\rho_1}\times$$
$$(Z_{0,\varepsilon}^M[J_\varepsilon])_\varepsilon = -i\exp\left[-i\left\langle V\left(-i\frac{\delta}{\delta(J_\varepsilon^\rho)_\varepsilon}\right)\right\rangle_\rho\right]\langle(\Delta_F(x-x_1;\varepsilon)J_{1,\varepsilon})_\varepsilon\rangle_{\rho_1}\times \quad (7.3.42)$$
$$\exp\left[i\left\langle V\left(-i\frac{\delta}{\delta(J_\varepsilon^\rho)_\varepsilon}\right)\right\rangle_\rho\right](Z_\varepsilon^M[J_\varepsilon])_\varepsilon.$$

From Eq.(4.3.42) follows that

$$(\Box_x+m^2)\frac{\delta(Z_\varepsilon^M[J_\varepsilon])_\varepsilon}{\delta(J_\varepsilon^\rho(x))_\varepsilon} = iO_x(Z_\varepsilon^M[J_\varepsilon])_\varepsilon \quad (7.3.43)$$

where

$$O_x = \exp\left[-i\left\langle V\left(-i\frac{\delta}{\delta(J_\varepsilon^\rho)_\varepsilon}\right)\right\rangle_\rho\right](J_\varepsilon(x))_\varepsilon\exp\left[i\left\langle V\left(-i\frac{\delta}{\delta(J_\varepsilon^\rho)_\varepsilon}\right)\right\rangle_\rho\right]. \quad (7.3.44)$$

We can evaluate $O_x$ by means of yet another trick. We set now

$$O_x(\lambda) =$$
$$\exp\left[-i\lambda\left\langle V\left(-i\frac{\delta}{\delta(J_\varepsilon^\rho)_\varepsilon}\right)\right\rangle_\rho\right](J_\varepsilon(x))_\varepsilon\exp\left[-i\lambda\left\langle V\left(-i\frac{\delta}{\delta(J_\varepsilon^\rho)_\varepsilon}\right)\right\rangle_\rho\right], \quad (7.3.45)$$

where $\lambda \in \mathbb{R}$ is a real parameter. Clearly

$$\frac{dO_x(\lambda)}{d\lambda} =$$

$$\exp\left[-i\lambda\left\langle V\left(-i\frac{\delta}{\delta(J_\varepsilon^\rho)_\varepsilon}\right)\right\rangle_\rho\right]\left[-iV\left(\left\langle -i\frac{\delta}{\delta(J_\varepsilon^\rho)_\varepsilon}\right\rangle_\rho\right),(J_\varepsilon(x))_\varepsilon\right]\times \quad (7.3.46)$$

But

$$\left[i\lambda\left\langle V\left(-i\frac{\delta}{\delta(J_\varepsilon^\rho)_\varepsilon}\right)\right\rangle_\rho\right].$$

$$\left[-iV\left(\left\langle -i\frac{\delta}{\delta(J_\varepsilon^\rho)_\varepsilon}\right\rangle_\rho\right),(J_\varepsilon(x))_\varepsilon\right] = -iV'\left(\left\langle -i\frac{\delta}{\delta(J_\varepsilon^\rho(y))_\varepsilon}\right\rangle_\rho\right)(\delta_{\rho_\varepsilon}^D(x-y))_\varepsilon, \quad (7.3.47)$$

where $V'$ is the derivative of $V$ with respect to its argument. Integrating Eq.(4.3.46) over $y,$ we find

$$\frac{dO_x(\lambda)}{d\lambda} = -V'\left(-i\frac{\delta}{\delta(J_\varepsilon^\rho(x))_\varepsilon}\right). \quad (7.3.48)$$

The equation Eq.(7.3.48) is now integrated over $\lambda$ to yield

$$O_x = O_x(1) = (J_\varepsilon(x))_\varepsilon - V'\left(-i\frac{\delta}{\delta(J_\varepsilon^\rho(x))_\varepsilon}\right). \quad (7.3.49)$$

Hence

$$(\Box_x + m^2)\frac{\delta(Z_\varepsilon^M[J_\varepsilon])_\varepsilon}{\delta(J_\varepsilon^\rho(x))_\varepsilon} = i\left((J_\varepsilon(x))_\varepsilon - V'\left(-i\frac{\delta}{\delta(J_\varepsilon^\rho(x))_\varepsilon}\right)\right)(Z_\varepsilon^M[J_\varepsilon])_\varepsilon \quad (7.3.50)$$

or

$$(\Box_x + m^2)(\varphi_{cl}(x;\varepsilon))_\varepsilon = (J_\varepsilon(x))_\varepsilon - \frac{1}{(Z_\varepsilon^M[J_\varepsilon])_\varepsilon}\left(V'\left(-i\frac{\delta}{\delta(J_\varepsilon^\rho(x))_\varepsilon}\right)\right)(Z_\varepsilon^M[J_\varepsilon])_\varepsilon. \quad (7.3.51)$$

The last term clearly resembles a force. For instance, take

$$V = \frac{\lambda}{4!}\varphi^4, \quad (7.3.52)$$

where $\lambda$ is dimensionless. Then

$$\frac{1}{(Z_\varepsilon^M[J_\varepsilon])_\varepsilon}\left(V'\left(-i\frac{\delta}{\delta(J_\varepsilon^\rho(x))_\varepsilon}\right)\right)(Z_\varepsilon^M[J_\varepsilon])_\varepsilon =$$

$$-\frac{\lambda}{3!}i^3\frac{1}{(Z_\varepsilon^M[J_\varepsilon])_\varepsilon}\frac{\delta^3(Z_\varepsilon^M[J_\varepsilon])_\varepsilon}{\delta\left((J_\varepsilon^\rho(x))_\varepsilon\right)^3} = \quad (7.3.53)$$

$$\frac{\lambda}{3!}\left[(\varphi_{cl}^3(x;\varepsilon))_\varepsilon - \frac{\delta^2(\varphi_{cl}(x;\varepsilon))_\varepsilon}{\delta\left((J_\varepsilon^\rho(x))_\varepsilon\right)^2} - 3i(\varphi_{cl}(x;\varepsilon))_\varepsilon\frac{\delta(\varphi_{cl}(x;\varepsilon))_\varepsilon}{\delta(J_\varepsilon^\rho(x))_\varepsilon}\right],$$

and finally we get

$$(\square_x + m^2)(\varphi_{cl}(x;\varepsilon))_\varepsilon =$$

$$(J_\varepsilon(x))_\varepsilon - \frac{\lambda}{3!}(\varphi_{cl}^3(x;\varepsilon))_\varepsilon + \frac{\lambda}{3!}\frac{\delta^2(\varphi_{cl}(x;\varepsilon))_\varepsilon}{\delta\left((J_\varepsilon^\rho(x))_\varepsilon\right)^2} + \frac{i\lambda}{4}\frac{\delta(\varphi_{cl}^2(x;\varepsilon))_\varepsilon}{\delta(J_\varepsilon^\rho(x))_\varepsilon}. \qquad (7.3.54)$$

The first two terms on the right hand side of the Eq.(7.3.54) as in classical case, give the classical equation of motion modified by the last two terms, which must amount to corrections from the quantum theory.

In the case $V \neq 0$, the explicit form of the effective action is, of course, not known but we can expand it functionally in terms of $(\varphi_{cl}(x;\varepsilon))_\varepsilon$ as

$$(\Gamma_\varepsilon[\varphi_{cl}(x;\varepsilon)])_\varepsilon = \int (d\rho_\varepsilon(x))_\varepsilon[-(V_\varepsilon^{\mathbf{e}}(\varphi_{cl}(x;\varepsilon)))_\varepsilon +$$
$$\frac{1}{2}(F_\varepsilon(\varphi_{cl}(x;\varepsilon))\partial_\mu\varphi_{cl}(x;\varepsilon)\partial^\mu\varphi_{cl}(x;\varepsilon))_\varepsilon + \text{higher order derivatives}\Big], \qquad (7.3.55)$$

where we take into account now local effects by including arbitrarily high higher order derivatives of $(\varphi_{cl}(x;\varepsilon))_\varepsilon$. We have arbitrary Colombeau generalized functions $(V_\varepsilon^{\mathbf{e}}(\varphi_{cl}(x;\varepsilon)))_\varepsilon, F_\varepsilon(\varphi_{cl}(x;\varepsilon))$ etc., to be determined. $(V_\varepsilon^{\mathbf{e}}(\bullet))_\varepsilon$ is clearly an effective potential. By expressing $(J_\varepsilon(x))_\varepsilon$ in terms of $(\varphi_{cl}(x;\varepsilon))_\varepsilon$ using (4.3.54) and integrating (4.3.36), we obtain that

$$(V_\varepsilon^{\mathbf{e}}(\varphi_{cl}(x;\varepsilon)))_\varepsilon = \frac{\lambda}{4!}(\varphi_{cl}^4(x;\varepsilon))_\varepsilon + \frac{m^2}{2}(\varphi_{cl}^2(x;\varepsilon))_\varepsilon + O(\hbar) \qquad (7.3.56)$$

and

$$(F_\varepsilon(\varphi_{cl}(x;\varepsilon)))_\varepsilon = 1 + \text{corrections}. \qquad (7.3.57)$$

Alternatively, we can expand the effective action in terms of $(\varphi_{cl}(x;\varepsilon))_\varepsilon$ by nonlocal way:

$$(\Gamma_\varepsilon[\varphi_{cl}(x;\varepsilon)])_\varepsilon = \left\langle \left(\Gamma_\varepsilon^{(N)}(1,2,\ldots,N)\varphi_{cl}(1;\varepsilon)\varphi_{cl}(2;\varepsilon)\ldots\varphi_{cl}(N;\varepsilon)\right)_\varepsilon \right\rangle_{\rho_1,\rho_2,\ldots,\rho_N} \qquad (7.3.58)$$

where the coefficients $\Gamma_\varepsilon^{(N)}(x_1,x_2,\ldots,x_N)$ are called the proper vertices. They depend only on the differences $x_i - x_j$ because of translation invariance so that their Colombeau-Fourier transforms are introduced by

$$\left(\widehat{\Gamma}_\varepsilon^{(N)}(p_1,p_2,\ldots,p_N)\right)_\varepsilon (2\pi)^D\delta(p_1+p_2+\ldots+p_N) =$$
$$\left(\int d^Dx_1 d^Dx_2\ldots d^Dx_N \exp[i(p_1x_1+p_2x_2+\ldots+p_Nx_N)]\widehat{\Gamma}_\varepsilon^{(N)}(p_1,p_2,\ldots,p_N)\right)_\varepsilon, \qquad (7.3.59)$$

where $\left(\widehat{\Gamma}_\varepsilon^{(N)}(p_1,p_2,\ldots,p_N)\right)_\varepsilon$ being defined only when the sum of its arguments vanishes.

## 7.4.Saddle-Point Evaluation of the Path Integral in negative dimensions.

We start from the Euclidean space definition of the generating functional in negative dimensions

$$(Z_E[J_\varepsilon])_\varepsilon = \mathbf{N}_E \int_{\mathcal{G}(\mathbb{R}_x^D)} \mathbf{D}[(\varphi_\varepsilon)_\varepsilon] \exp\{-\mathbf{S}_E[(\varphi_\varepsilon)_\varepsilon, (J_\varepsilon)_\varepsilon]\} = \mathbf{N}_E \int_{\mathcal{G}(\mathbb{R}_x^D)} \mathbf{D}[(\varphi_\varepsilon)_\varepsilon] \times$$
$$\exp\left[-\left\langle \frac{1}{2}(\bar{\partial}_\mu \varphi_\varepsilon \bar{\partial}_\mu \varphi_\varepsilon)_\varepsilon + \frac{1}{2} m^2 (\varphi_\varepsilon^2)_\varepsilon + (V_\varepsilon(\varphi_\varepsilon))_\varepsilon - (\varphi_\varepsilon)_\varepsilon (J_\varepsilon)_\varepsilon \right\rangle_\rho \right], \quad (7.4.1)$$

where

$$(\mathbf{S}_E[(\varphi_\varepsilon)_\varepsilon, (J_\varepsilon)_\varepsilon])_\varepsilon = \left\langle \frac{1}{2}(\bar{\partial}_\mu \varphi_\varepsilon \bar{\partial}_\mu \varphi_\varepsilon)_\varepsilon + \frac{1}{2} m^2 (\varphi_\varepsilon^2)_\varepsilon + (V_\varepsilon(\varphi_\varepsilon))_\varepsilon - (\varphi_\varepsilon)_\varepsilon (J_\varepsilon)_\varepsilon \right\rangle_\rho \quad (7.4.2)$$

**Remark 7.4.1**. Note that in Eq.(7.4.1) in contrast with Eq.(7.3.4) and Eq.(7.3.5), we integrate over full Colombeau algebra $\mathcal{G}(\mathbb{R}_x^D)$.

We then expand the action around a field configuration $(\varphi_{0,\varepsilon}(\bar{x}))_\varepsilon$

$$(\mathbf{S}_E[(\varphi_\varepsilon)_\varepsilon, (J_\varepsilon)_\varepsilon])_\varepsilon =$$
$$(\mathbf{S}_E[(\varphi_{0,\varepsilon})_\varepsilon, (J_\varepsilon)_\varepsilon])_\varepsilon + \left\langle \frac{\delta \mathbf{S}_E}{\delta(\varphi_\varepsilon(\bar{x}_1))_\varepsilon} (\varphi_\varepsilon(\bar{x}_1) - \varphi_{0,\varepsilon}(\bar{x}_1))_\varepsilon \right\rangle_\rho + \quad (7.4.3)$$
$$\frac{1}{2}\left\langle \frac{\delta^2 \mathbf{S}_E}{\delta(\varphi_\varepsilon(\bar{x}_1))_\varepsilon \delta(\varphi_\varepsilon(\bar{x}_2))_\varepsilon} (\varphi_\varepsilon(\bar{x}_1) - \varphi_{0,\varepsilon}(\bar{x}_1))_\varepsilon (\varphi_\varepsilon(\bar{x}_2) - \varphi_{0,\varepsilon}(\bar{x}_2))_\varepsilon \right\rangle_{\rho_1,\rho_2}$$

with the functional derivative evaluated at $(\varphi_{0,\varepsilon}(x))_\varepsilon$. We take $(\mathbf{S}_E[(\varphi_\varepsilon)_\varepsilon, (J_\varepsilon)_\varepsilon])_\varepsilon$ to be stationary at $(\varphi_{0,\varepsilon}(\bar{x}))_\varepsilon$, which means that $(\varphi_{0,\varepsilon}(\bar{x}))_\varepsilon$ obeys the classical equations of motion with the source term

$$\left.\frac{\delta \mathbf{S}_E}{\delta(\varphi_\varepsilon(\bar{x}))_\varepsilon}\right|_{(\varphi_\varepsilon)_\varepsilon = (\varphi_{0,\varepsilon})_\varepsilon} = -(\bar{\partial}_\mu \bar{\partial}_\mu \varphi_\varepsilon(\bar{x}))_\varepsilon -$$
$$\left( \frac{\bar{\partial}_\mu \rho_\varepsilon(\bar{x})}{\rho_\varepsilon(\bar{x})} \bar{\partial}_\mu \varphi_\varepsilon(\bar{x}) \right)_\varepsilon + m^2 (\varphi_\varepsilon^2(\bar{x}))_\varepsilon + (V_\varepsilon'(\varphi_\varepsilon(\bar{x})))_\varepsilon - (J_\varepsilon(\bar{x}))_\varepsilon = 0 \quad (7.4.4)$$

It follows that after integration by parts

$$(\mathbf{S}_E[(\varphi_{0,\varepsilon})_\varepsilon, (J_\varepsilon)_\varepsilon])_\varepsilon = \frac{1}{2} \int d^D \bar{x} (\rho_\varepsilon(\bar{x}))_\varepsilon \left( \left( 2 - \varphi_{0,\varepsilon} \frac{d}{d\varphi_{0,\varepsilon}} \right) (J_\varepsilon \varphi_{0,\varepsilon} + V_\varepsilon(\varphi_{0,\varepsilon})) \right)_\varepsilon \quad (7.4.5)$$

while

$$\frac{\delta^2 \mathbf{S}_E}{\delta(\varphi_{1,\varepsilon})_\varepsilon \delta(\varphi_{2,\varepsilon})_\varepsilon} = \left[ -\bar{\partial}_\mu \bar{\partial}_\mu - \left( \frac{\bar{\partial}_\mu \rho_\varepsilon}{\rho_\varepsilon} \right)_\varepsilon \bar{\partial}_\mu + m^2 + (V_\varepsilon''(\varphi_{1,\varepsilon}))_\varepsilon \right] \delta^D(x_1 - x_2) \quad (7.4.6)$$

is an operator. By using formal approuch of the saddle point evaluation, the generatin functional now becomes

$$(Z_E[J_\varepsilon])_\varepsilon \simeq \mathbf{N}_E \exp\{-\mathbf{S}_E[(\varphi_{0,\varepsilon})_\varepsilon, (J_\varepsilon)_\varepsilon]\} \int_{\mathcal{G}(\mathbb{R}_x^D)} \mathbf{D}[(\varphi_\varepsilon)_\varepsilon] \times$$
$$\exp\left[ -\frac{1}{2} \left\langle (\varphi_{1,\varepsilon})_\varepsilon \frac{\delta^2 \mathbf{S}_E}{\delta(\varphi_{1,\varepsilon})_\varepsilon \delta(\varphi_{2,\varepsilon})_\varepsilon} (\varphi_{2,\varepsilon})_\varepsilon \right\rangle_{\rho_1,\rho_2} \right] =$$
$$\mathbf{N}_E' \exp\{-\mathbf{S}_E[(\varphi_{0,\varepsilon})_\varepsilon, (J_\varepsilon)_\varepsilon]\} \times \quad (7.4.7)$$
$$\left( \left\{ \det\left[ -\bar{\partial}_\mu \bar{\partial}_\mu - \frac{\bar{\partial}_\mu \rho_\varepsilon}{\rho_\varepsilon} \bar{\partial}_\mu + m^2 + V_\varepsilon''(\varphi_{0,\varepsilon}) \right] \delta_{1,2}(\varepsilon) \right\}^{-\frac{1}{2}} \right)_\varepsilon$$

Where the Gaussian integral in Eq.(7.4.7) with the formal result is

$$\int_{\mathcal{G}(\mathbb{R}_x^D)} \mathbf{D}[(\varphi_\varepsilon)_\varepsilon]\exp\left[-\frac{1}{2}\left\langle(\varphi_{1,\varepsilon})_\varepsilon\frac{\delta^2\mathbf{S}_E}{\delta(\varphi_{1,\varepsilon})_\varepsilon\delta(\varphi_{2,\varepsilon})_\varepsilon}(\varphi_{2,\varepsilon})_\varepsilon\right\rangle_{\rho_1,\rho_2}\right] =$$

$$\int_{\mathcal{G}(\mathbb{R}_x^D)} \mathbf{D}[(\varphi_\varepsilon)_\varepsilon]\exp\left[-\bar{\partial}_\mu\bar{\partial}_\mu - \left(\frac{\bar{\partial}_\mu\rho_\varepsilon}{\rho_\varepsilon}\right)_\varepsilon\bar{\partial}_\mu + m^2 + V''_\varepsilon(\varphi_{1,\varepsilon})\right] = \tag{7.4.8}$$

$$\left(\left\{\det\left[-\bar{\partial}_\mu\bar{\partial}_\mu - \frac{\bar{\partial}_\mu\rho_\varepsilon}{\rho_\varepsilon}\bar{\partial}_\mu + m^2 + V''_\varepsilon(\varphi_{1,\varepsilon})\right]\delta_{1,2}(\varepsilon)\right\}^{-\frac{1}{2}}\right)_\varepsilon.$$

Clearly this expression needs some getting used to. We can rewrite it in a slightly more suggestive form by using the identity

$$(\det M_\varepsilon)_\varepsilon = (\exp[\mathbf{Tr}\ln M_\varepsilon])_\varepsilon. \tag{7.4.9}$$

Therefore

$$(Z_E[J_\varepsilon])_\varepsilon = \mathbf{N}'_E\exp\{-\mathbf{S}_E[(\varphi_{0,\varepsilon})_\varepsilon,(J_\varepsilon)_\varepsilon] -$$
$$\frac{1}{2}\left(\mathbf{Tr}\ln\left[\left\{-\bar{\partial}_\mu\bar{\partial}_\mu - \frac{\bar{\partial}_\mu\rho_\varepsilon}{\rho_\varepsilon}\bar{\partial}_\mu + m^2 + V''_\varepsilon(\varphi_{0,\varepsilon})\right\}\delta_{1,2}(\varepsilon)\right]\right)_\varepsilon\Bigg\}, \tag{7.4.10}$$

which clearly Indicates we are computing corrections to $(\check{Z}_E[J_\varepsilon])_\varepsilon$.The first term $\mathbf{S}_E[(\varphi_{0,\varepsilon})_\varepsilon,(J_\varepsilon)_\varepsilon]$ gives the classical conribution to the Green's functions (remember Dirac's identification). The next term, of $O(\hbar)$, gives the first quantum correction to the Green's functions.The determinant of an operator is understood to mean the product of its igenvalues. We start by computing the classical contributions to $(Z_E[J_\varepsilon])_\varepsilon$.First recoll that $(\varphi_{0,\varepsilon})_\varepsilon$, being the solution of (7.4.4), is a functional of $(J_\varepsilon)_\varepsilon$. The procedure is therefore very simple: (i) calculate the functional dependence of $(\varphi_{0,\varepsilon})_\varepsilon$ on $(J_\varepsilon)_\varepsilon$, (ii) insert it in (7.4.5) and, (iii) by comparing the resulting expression with the expansion (7.3.25), extract the Green's functions $G^{(N)}(1,..,N)$. The best case do this using perturbation theory. Specifically, take the potential and expand around $\lambda = 0$.We set

$$(\varphi_{0,\varepsilon}(x))_\varepsilon = \left(\varphi_\varepsilon^{(0)}(x)\right)_\varepsilon + \lambda\left(\varphi_\varepsilon^{(1)}(x)\right)_\varepsilon + \lambda^2\left(\varphi_\varepsilon^{(2)}(x)\right)_\varepsilon + \ldots \tag{7.4.11}$$

thus

$$(\mathbf{S}_E[(J_\varepsilon)_\varepsilon])_\varepsilon =$$
$$-\frac{1}{2}\int(d^D\rho_\varepsilon(x))_\varepsilon\Big[(J_\varepsilon)_\varepsilon\Big(\left(\varphi_\varepsilon^{(0)}(x)\right)_\varepsilon + \lambda\left(\varphi_\varepsilon^{(1)}(x)\right)_\varepsilon + \lambda^2\left(\varphi_\varepsilon^{(2)}(x)\right)_\varepsilon + \ldots\Big)$$
$$\frac{\lambda}{12}\Big(\left(\varphi_\varepsilon^{(0)}(x)\right)_\varepsilon + \lambda\left(\varphi_\varepsilon^{(1)}(x)\right)_\varepsilon + \lambda^2\left(\varphi_\varepsilon^{(2)}(x)\right)_\varepsilon + \ldots\Big)^4\Big] = \tag{7.4.12}$$
$$-\frac{1}{2}\int(d^D\rho_\varepsilon(x))_\varepsilon(J_\varepsilon)_\varepsilon\left(\varphi_\varepsilon^{(0)}(x)\right)_\varepsilon -$$
$$\frac{\lambda}{2}\int(d^D\rho_\varepsilon(x))_\varepsilon\left[(J_\varepsilon)_\varepsilon\left(\varphi_\varepsilon^{(1)}(x)\right)_\varepsilon + \frac{1}{12}\left(\left(\varphi_\varepsilon^{(0)}(x)\right)_\varepsilon\right)^4\right] + O(\lambda^2).$$

We define now the Euclidean Greent's function in obvious notation:

$$\left(\bar{\partial}_\mu\bar{\partial}_\mu - \frac{\bar{\partial}_\mu\rho_\varepsilon}{\rho_\varepsilon} - m^2\right)(G_{\varepsilon,xy})_\varepsilon = -(\delta^D_{\rho_\varepsilon,xy})_\varepsilon \tag{7.4.13}$$

it follows that

$$\left(\varphi_\varepsilon^{(0)}(x)\right)_\varepsilon = \langle (G_{\varepsilon,xy}J_{\varepsilon,a})_\varepsilon \rangle_{\rho_a} = \langle (G_{\varepsilon,xy})_\varepsilon (J_{\varepsilon,a})_\varepsilon \rangle_{\rho_a}$$

$$\left(\varphi_\varepsilon^{(1)}(x)\right)_\varepsilon = -\frac{1}{6}\langle (G_{\varepsilon,xy})_\varepsilon (G_{\varepsilon,ya})_\varepsilon (G_{\varepsilon,yb})_\varepsilon (G_{\varepsilon,yc})_\varepsilon (J_{\varepsilon,a})_\varepsilon (J_{\varepsilon,b})_\varepsilon (J_{\varepsilon,c})_\varepsilon \rangle_{\rho_a\rho_b\rho_c\rho_y}, \quad (7.4.14)$$

etc...

Thus

$$(\mathbf{S}_E[(J_\varepsilon)_\varepsilon])_\varepsilon = -\frac{1}{2}\langle (J_{\varepsilon,a})_\varepsilon (G_{\varepsilon,ab})_\varepsilon (J_{\varepsilon,b})_\varepsilon \rangle_{\rho_a\rho_b} +$$
$$\frac{\lambda}{4!}\langle (G_{\varepsilon,xa})_\varepsilon (G_{\varepsilon,xb})_\varepsilon (G_{\varepsilon,xc})_\varepsilon (G_{\varepsilon,xd})_\varepsilon (J_{\varepsilon,a})_\varepsilon (J_{\varepsilon,b})_\varepsilon (J_{\varepsilon,c})_\varepsilon (J_{\varepsilon,d})_\varepsilon \rangle_{\rho_a\rho_b\rho_c\rho_d\rho_x} -$$
$$\frac{\lambda^2}{3\cdot 4!}\langle (G_{\varepsilon,xa})_\varepsilon (G_{\varepsilon,xb})_\varepsilon (G_{\varepsilon,xc})_\varepsilon (G_{\varepsilon,xy})_\varepsilon (G_{\varepsilon,yd})_\varepsilon (G_{\varepsilon,ye})_\varepsilon (G_{\varepsilon,yf})_\varepsilon \times \quad (7.4.15)$$
$$(J_{\varepsilon,a})_\varepsilon (J_{\varepsilon,b})_\varepsilon (J_{\varepsilon,c})_\varepsilon (J_{\varepsilon,d})_\varepsilon (J_{\varepsilon,e})_\varepsilon (J_{\varepsilon,f})_\varepsilon \rangle_{\rho_a\rho_b\rho_c\rho_d\rho_e\rho_f\rho_x\rho_y}.$$

Correspondingly, the (connected) Euclidean Green's functions are given by

$$\left(G_E^{(N)}(x_1,\ldots,x_n;\varepsilon)\right)_\varepsilon = \left(\frac{\delta^N \check{Z}_\varepsilon^E[J_\varepsilon]}{\delta J_\varepsilon^{\rho_1}(x_1)\ldots\delta J_\varepsilon^{\rho_N}(x_N)}\right)_\varepsilon, \quad (7.4.16)$$

where

$$(Z_E[J_\varepsilon])_\varepsilon = \mathbf{N}_E \exp\left\{(-\check{Z}_\varepsilon^E[J_\varepsilon])_\varepsilon\right\}. \quad (7.4.17)$$

In this classical approximation we find the connected Green's functions to be

$$\left(G_E^{(2)}(\bar{x}_1,\bar{x}_2;\varepsilon)\right)_\varepsilon = (G(\bar{x}_1,\bar{x}_2;\varepsilon))_\varepsilon = \Upsilon(D,D^-)\left(\frac{r}{m}\right)^{\frac{2+D\cdot|\alpha_i|}{2}} I_{\frac{2+D|\alpha|}{2}}(mr), \quad (7.4.18)$$

see Eq.(7.2.26),

$$\left(G_E^{(4)}(\bar{x}_1,\bar{x}_2,\bar{x}_3,\bar{x}_4;\varepsilon)\right)_\varepsilon =$$
$$-\lambda\int d^D\bar{y}(G(\bar{x}_1,\bar{y};\varepsilon))_\varepsilon (G(\bar{x}_2,\bar{y};\varepsilon))_\varepsilon (G(\bar{x}_3,\bar{y};\varepsilon))_\varepsilon (G(\bar{x}_4,\bar{y};\varepsilon))_\varepsilon,$$
$$\left(G_E^{(6)}(\bar{x}_1,\bar{x}_2,\bar{x}_3,\bar{x}_4,\bar{x}_5,\bar{x}_6;\varepsilon)\right)_\varepsilon = \quad (7.4.19)$$
$$\lambda^2\int d^D\bar{x}d^D\bar{y}(G(\bar{x}_1,\bar{x}_2;\varepsilon))_\varepsilon (P(\bar{x},\bar{y},\bar{x}_1,\bar{x}_2,\bar{x}_3,\bar{x}_4,\bar{x}_5,\bar{x}_6;\varepsilon))_\varepsilon$$

where

$$(P(\bar{x},\bar{y},\{\bar{x}_i\};\varepsilon))_\varepsilon =$$
$$\sum_{(ijk)}(G(\bar{x},\bar{x}_i;\varepsilon))_\varepsilon (G(\bar{x},\bar{x}_j;\varepsilon))_\varepsilon (G(\bar{x},\bar{x}_k;\varepsilon))_\varepsilon (G(\bar{y},\bar{x}_l;\varepsilon))_\varepsilon (G(\bar{y},\bar{x}_m;\varepsilon))_\varepsilon \times \quad (7.4.20)$$
$$(G(\bar{y},\bar{x}_n;\varepsilon))_\varepsilon,$$

where the sum runs over all the following values of the triples, $(ijk) = (123),(124),(125),$ $(126),(134),(135),(136),(145),(146),(156),$ with $(lmn)$ assuming the complementary value, e.g., $(lmn) = (456)$ when $(ijk) = (123)$. Note that $ijk$ runs only over half of the possible values. This is because the expression for $P(\cdot)$ is symmetric under the interchange $\bar{x} \to \bar{y}$. In this classical approximation and to order $\lambda^2$ these are the only nonzero Green's functions.

The momentum space Green's functions, defined by

$$
\left(\widehat{G}_E^{(N)}(p_1,\ldots,p_N;\varepsilon)\right)_\varepsilon = (2\pi)^D\delta(p_1,\ldots,p_N)\times
$$
$$
\left(\int d^Dx_1\ldots\int d^Dx_N\exp[i(p_1x_1+\ldots+p_Nx_N)]G_E^{(N)}(x_1,\ldots,x_N;\varepsilon)\right)_\varepsilon = \tag{7.4.21}
$$
$$
\int d^Dx_1\ldots\int d^Dx_N\exp[i(p_1x_1+\ldots+p_Nx_N)]\left(G_E^{(N)}(x_1,\ldots,x_N;\varepsilon)\right)_\varepsilon.
$$

## 7.5.Power-counting renormalizability of $P(\varphi)_{D^-}$ scalar field theory in negative dimensions $D^-$.

Consider a free scalar field with action in negative dimensions $D^- = D(\alpha-1), \alpha < 0$ :

$$
(\mathbf{S}_{0,\varepsilon})_\varepsilon = -\frac{1}{2}\left(\int_{\mathbb{R}_x^D} d\varrho_\varepsilon(x)\,\phi_\varepsilon(x)\,P(\Box)\,\phi_\varepsilon(x)\right)_\varepsilon, \tag{7.5.1}
$$

where $\varrho = [(\varrho_\varepsilon(x))_\varepsilon] \in \mathcal{G}(\mathbb{R}^D), [(\phi_\varepsilon(x))_\varepsilon] \in \mathcal{G}(\mathbb{R}^D)$ and $P(\bullet)$ is a polinomial.

**Definition 7.5.1**. Assume that (i) $(\phi_\varepsilon(x))_\varepsilon \in \mathcal{G}(\mathbb{R}^D), \varrho \in \mathcal{G}(\mathbb{R}^D)$ and (ii) there exist Colombeau generalized function $\left(\tilde{\phi}_\varepsilon(k)\right)_\varepsilon \in \mathcal{G}(\mathbb{R}_k^D)$ such that

$$
(\phi_\varepsilon(x))_\varepsilon = \left(\int_{\mathbb{R}_k^D} d\varrho_\varepsilon(k)\,\tilde{\phi}_\varepsilon(k)\,e^{ik\cdot x}\right)_\varepsilon = (2\pi)^{-D}\int_{\mathbb{R}_k^D}(d\varrho_\varepsilon(k))_\varepsilon\left(\tilde{\phi}_\varepsilon(k)\right)_\varepsilon e^{ik\cdot x}, \tag{7.5.2}
$$

and

$$
\left(\tilde{\phi}_\varepsilon(k)\right)_\varepsilon = \left(\int_{\mathbb{R}_x^n} d\varrho_\varepsilon(x)\,\phi_\varepsilon(x)e^{-ik\cdot x}\right)_\varepsilon = \int_{\mathbb{R}_x^D}(d\varrho_\varepsilon(x))_\varepsilon(\phi_\varepsilon(x))_\varepsilon e^{-ik\cdot x}. \tag{7.5.3}
$$

Then we will say that: (1) $\left(\tilde{\phi}_\varepsilon(k)\right)_\varepsilon$ is Colombeau Fourier–Stieltjes transform of the field $(\phi_\varepsilon(x))_\varepsilon$ with weight $\varrho$ and abraviate

$$
\left(\tilde{\phi}_\varepsilon(k)\right)_\varepsilon = (\mathcal{FS})_\varrho[(\phi_\varepsilon(x))_\varepsilon](k), \tag{7.5.4}
$$

(2) $(\phi_\varepsilon(x))_\varepsilon$ is inverse Colombeau Fourier–Stieltjes transform of the field $\left(\tilde{\phi}_\varepsilon(k)\right)_\varepsilon$ with weight $(\varrho_\varepsilon)_\varepsilon$ and abraviate

$$
(\phi_\varepsilon(x))_\varepsilon = (\mathcal{FS})_\varrho^{-1}\left[\left(\tilde{\phi}_\varepsilon(k)\right)_\varepsilon\right](x). \tag{7.5.5}
$$

**Definition 7.5.2**. We will denote:

(i) the set of the Colombeau generalized functions $\left(\tilde{\phi}_\varepsilon(k)\right)_\varepsilon \in \mathcal{G}(\mathbb{R}_k^D)$ which is Colombeau
Fourier–Stieltjes transform with weight $\varrho$ by $(\mathcal{FS})_\varrho[\mathcal{G}(\mathbb{R}_k^D)]$ or by $\mathcal{G}^{(\mathcal{FS})_\varrho}(\mathbb{R}_k^D)$

(ii) the set of the Colombeau generalized functions $(\phi_\varepsilon(x))_\varepsilon \in \mathcal{G}(\mathbb{R}_x^n)$ which is inverse Colombeau-Fourier–Stieltjes transform with weight $\varrho$ by $(\mathcal{FS})^{-1}[\mathcal{G}_x(\mathbb{R}^D)]$ or $\mathcal{G}^{(\mathcal{FS})_\varrho^{-1}}(\mathbb{R}_x^D)$.

(iii) Note that we assume that in both cases (i) and (ii) the Eqs.(7.5.2)-(7.5.3) are satisfies.

**Remark 7.5.1.**Note that $\mathcal{G}^{(\mathcal{FS})_\varrho}(\mathbb{R}_k^D)\times\mathcal{G}^{(\mathcal{FS})_\varrho^{-1}}(\mathbb{R}_x^D)$ is the linear topological subspace

**Definition 7.5.3.**The free partition function $(\mathbf{Z}_{0,\varepsilon}[J_\varepsilon])_\varepsilon$ in the presence of a local source $(J_\varepsilon(x))_\varepsilon \in \mathcal{G}(\mathbb{R}_x^D)$ is

$$(\mathbf{Z}_{0,\varepsilon}[J_\varepsilon])_\varepsilon = \left(\int_{(\phi_\varepsilon)_\varepsilon \in \mathcal{G}^{(\mathcal{F}S)}{}_{\varrho}^{-1}(\mathbb{R}_x^D)} [\mathcal{D}(\phi_\varepsilon)] \exp\left\{ i[\mathbf{S}_{0,\varepsilon} + \int_{\mathbb{R}_x^D} \varrho_\varepsilon(x) J_\varepsilon(x)\phi_\varepsilon(x)] \right\}\right)_\varepsilon \triangleq$$

$$\int_{(\phi_\varepsilon)_\varepsilon \in \mathcal{G}^{(\mathcal{F}S)}{}_{\varrho}^{-1}(\mathbb{R}_x^D)} [\mathcal{D}(\phi_\varepsilon)_\varepsilon] \exp\left\{ [(\mathbf{S}_{0,\varepsilon})_\varepsilon + \int_{\mathbb{R}_x^D} (\varrho_\varepsilon(x))_\varepsilon (J_\varepsilon(x))_\varepsilon (\phi_\varepsilon(x))_\varepsilon] \right\} \triangleq \tag{7.5.6}$$

$$\int_{(\phi_\varepsilon)_\varepsilon \in \mathcal{G}^{(\mathcal{F}S)}{}_{\varrho}^{-1}(\mathbb{R}_x^D)} [\mathcal{D}(\phi_\varepsilon)_\varepsilon] e^{i(\mathbf{S}_{J_\varepsilon})_\varepsilon}.$$

Using now the definition of the $D$-dimensional Colombeau–Dirac distribution with nontrivial Colombeau–Stieltjes measure $(d\varrho_\varepsilon(k))_\varepsilon$

$$(\delta_{\varrho_\varepsilon}(k))_\varepsilon = \left(\int_{\mathbb{R}_k^D} d\varrho_\varepsilon(k) e^{ik\bullet x}\right)_\varepsilon \triangleq \int_{\mathbb{R}_k^D} (d\varrho_\varepsilon(k))_\varepsilon e^{ik\bullet x} \tag{7.5.7}$$

and Eqs.(7.5.2)-(7.5.3) we obtain

$$\begin{gathered}
(\mathbf{S}_{J_\varepsilon})_\varepsilon = \\
\frac{1}{2}\int_{\mathbb{R}_x^D} (d\varrho_\varepsilon(x))_\varepsilon \times \\
\left\{ \int_{\mathbb{R}_{k_1}^D} \frac{(d\varrho_\varepsilon(k_1))_\varepsilon}{(2\pi)^D} \int_{\mathbb{R}_{k_2}^D} \frac{(d\varrho_\varepsilon(k_2))_\varepsilon}{(2\pi)^D} e^{i(k_1+k_2)\bullet x} \times \right. \\
\left[ -\left(\left(\widetilde{\phi}_\varepsilon(k_1)\right)_\varepsilon\right) \left((f_\varepsilon(-k_2^2))_\varepsilon\right)\left(\left(\widetilde{\phi}_\varepsilon(k_2)\right)_\varepsilon\right) \right. \\
\left.\left. +\left(\widetilde{J}_\varepsilon(k_1)\right)_\varepsilon \left(\widetilde{\phi}_\varepsilon(k_2)\right)_\varepsilon + \left(\widetilde{J}_\varepsilon(k_2)\right)_\varepsilon \left(\widetilde{\phi}_\varepsilon(k_1)\right)_\varepsilon \right]\right\} = \\
= \frac{1}{2}\int_{\mathbb{R}_k^n} \frac{(d\varrho_\varepsilon(-k))_\varepsilon}{(2\pi)^D} \left[ -\left(\left(\widetilde{\phi}_\varepsilon(-k)\right)_\varepsilon\right)((f_\varepsilon(-k^2))_\varepsilon)\left(\left(\widetilde{\phi}_\varepsilon(k)\right)_\varepsilon\right) + \right. \\
\left. \left(\left(\widetilde{J}_\varepsilon(-k)\right)_\varepsilon\right)\left(\widetilde{\phi}_\varepsilon(k)\right)_\varepsilon + \left(\widetilde{J}_\varepsilon(k)\right)_\varepsilon \left(\left(\widetilde{\phi}_\varepsilon(-k)\right)_\varepsilon\right)_\varepsilon \right] \\
= \frac{1}{2}\int_{\mathbb{R}_k^D} \frac{(d\varrho_\varepsilon(-k))_\varepsilon}{(2\pi)^D} \times \\
\left[ -\left((\widetilde{\varphi}_\varepsilon(-k))_\varepsilon\right)\left(\left(f(-k^2)\right)_\varepsilon\right)\left((\widetilde{\varphi}_\varepsilon(k))_\varepsilon\right) + \frac{\left(\left(\widetilde{J}_\varepsilon(-k)\right)_\varepsilon\right)\left(\left(\widetilde{J}_\varepsilon(k)\right)_\varepsilon\right)}{(f_\varepsilon(-k^2))_\varepsilon} \right], \\
(\widetilde{\varphi}_\varepsilon)_\varepsilon = \left(\widetilde{\phi}_\varepsilon(k)\right)_\varepsilon - \frac{\left(\widetilde{J}_\varepsilon(-k)\right)_\varepsilon}{(f_\varepsilon(-k^2))_\varepsilon}, (f_\varepsilon(-k_2^2))_\varepsilon = (P(-k_2^2) + i\varepsilon)_\varepsilon.
\end{gathered} \tag{7.5.8}$$

Thus Eq.(7.5.6) becomes

$$
\begin{gathered}
(Z_{0,\varepsilon}[J_\varepsilon])_\varepsilon = \\
\left\{\int_{(\phi_\varepsilon)_\varepsilon \in \mathcal{G}^{(\mathcal{FS})_\varrho^{-1}}(\mathbb{R}_x^D)} [\mathcal{D}(\phi_\varepsilon)_\varepsilon] \times \right. \\
\left. \exp\left[-\frac{i}{2}\int_{\mathbb{R}_k^D} \frac{(d\varrho_\varepsilon(k))_\varepsilon}{(2\pi)^D}\left(\left(\widetilde{\phi}_\varepsilon(-k)\right)_\varepsilon\right)((f_\varepsilon(-k^2))_\varepsilon)\left(\left(\widetilde{\phi}_\varepsilon(k)\right)_\varepsilon\right)\right]\right\} \times \\
\exp\left[\frac{i}{2}\int_{\mathbb{R}_k^D} \frac{(d\varrho_\varepsilon(-k))_\varepsilon}{(2\pi)^D} \frac{\left(\left(\widetilde{J}_\varepsilon(-k)\right)_\varepsilon\right)\left(\left(\widetilde{J}_\varepsilon(k)\right)_\varepsilon\right)}{(f_\varepsilon(-k^2))_\varepsilon}\right] \\
= Z_0[0]\exp\left[\frac{i}{2}\int_{\mathbb{R}_k^D} \frac{(d\varrho_\varepsilon(-k))_\varepsilon}{(2\pi)^D} \frac{\left(\left(\widetilde{J}_\varepsilon(-k)\right)_\varepsilon\right)\left(\left(\widetilde{J}_\varepsilon(k)\right)_\varepsilon\right)}{(f_\varepsilon(-k^2))_\varepsilon}\right].
\end{gathered}
\tag{7.5.9}
$$

Therefore the exponent in Eq.(7.5.9) can be written as

$$
\begin{gathered}
\int_{\mathbb{R}_k^D} \frac{(d\varrho_\varepsilon(-k))_\varepsilon}{(2\pi)^D} \frac{\left(\left(\widetilde{J}_\varepsilon(-k)\right)_\varepsilon\right)\left(\left(\widetilde{J}_\varepsilon(k)\right)_\varepsilon\right)}{(f_\varepsilon(-k^2))_\varepsilon} = \\
\int_{\mathbb{R}_k^D} \frac{(d\varrho_\varepsilon(-k))_\varepsilon}{(2\pi)^D} \int_{\mathbb{R}_x^D}((d\varrho_\varepsilon(x))_\varepsilon)\int_{\mathbb{R}_y^D}((d\varrho(y))_\varepsilon)e^{ik\cdot(x-y)} \frac{\left(\left(\widetilde{J}_\varepsilon(-k)\right)_\varepsilon\right)\left(\left(\widetilde{J}_\varepsilon(k)\right)_\varepsilon\right)}{(f_\varepsilon(-k^2))_\varepsilon},
\end{gathered}
\tag{7.5.10}
$$

so that, if $(\varrho_\varepsilon(-k))_\varepsilon = (d\varrho_\varepsilon(k))_\varepsilon$, the free partition function reads

$$
\begin{gathered}
(Z_{0,\varepsilon}[J])_\varepsilon = \\
((Z_{0,\varepsilon}[0])_\varepsilon)\exp\left[\frac{i}{2}\int_{\mathbb{R}_x^D}((d\varrho_\varepsilon(x))_\varepsilon)\int_{\mathbb{R}_y^D}((d\varrho_\varepsilon(y))_\varepsilon)(J_\varepsilon(x)((G_M(x-y;\varepsilon))_\varepsilon)J_\varepsilon(y))_\varepsilon\right]
\end{gathered}
\tag{7.5.11}
$$

where

$$
(G_M(x-y;\varepsilon))_\varepsilon = \frac{1}{(2\pi)^D}\int_{\mathbb{R}_k^D}(d(\varrho_\varepsilon(k))_\varepsilon)\frac{e^{ik\cdot(x-y)}}{(f_\varepsilon(-k^2))_\varepsilon}. \tag{7.5.12}
$$

Thus, we have recovered the usual definition of the propagator as the solution of the Green equation in negative dimensions $D^- = D(\alpha-1), \alpha < 0$ :

$$
(P(\Box)\,G_M(x-y;\varepsilon))_\varepsilon = (\delta_{\varrho_\varepsilon}(x-y))_\varepsilon. \tag{7.5.13}
$$

In the following, we proceed with the conventional (perturbative) evaluation of the Euclidean Green's functions. We start from

$$
\begin{gathered}
(W_\varepsilon^E[J_\varepsilon])_\varepsilon = \exp\{(-Z_\varepsilon^E[J_\varepsilon])_\varepsilon\} = \mathbf{N}\int_{(\phi_\varepsilon)_\varepsilon \in \mathcal{G}^{(\mathcal{FS})_\varrho^{-1}}(\mathbb{R}_x^D)}[\mathcal{D}(\phi_\varepsilon)_\varepsilon] \times \\
\exp\left\{-\int_{\mathbb{R}_x^D}((d\varrho_\varepsilon(x))_\varepsilon) \times \right. \\
\left.\left[\frac{1}{2}((\partial_\mu\phi_\varepsilon)_\varepsilon)((\partial^\mu\phi_\varepsilon)_\varepsilon) + \frac{1}{2}m^2(\phi_\varepsilon^2)_\varepsilon + (V(\phi_\varepsilon))_\varepsilon - J_\varepsilon\phi_\varepsilon\right]\right\},
\end{gathered}
\tag{7.5.14}
$$

where $\mathbf{N}$ is an infinite normalization constant. The connected Euclidean Green's functions are given by

$$
\left(G_E^{(N)}(x_1,\ldots,x_n;\varepsilon)\right)_\varepsilon = \left(\frac{\delta^N Z_\varepsilon^E[J_\varepsilon]}{\delta J_\varepsilon(x_1)\ldots\delta J_\varepsilon(x_n)}\right)_\varepsilon. \tag{7.5.15}
$$

They will be calculated by perturbing in the potential $V$. Using now the standard trick, we obtain

$$(W_\varepsilon[J_\varepsilon])_\varepsilon = \mathbf{N}\exp\left(\left(-\left\langle V\left(\frac{\delta}{\delta J_\varepsilon}\right)\right\rangle\right)_\varepsilon\right)\exp\left((Z^E_{0,\varepsilon}[J_\varepsilon])_\varepsilon\right) \tag{7.5.16}$$

where

$$(Z^E_{0,\varepsilon}[J_\varepsilon])_\varepsilon = -\frac{1}{2}\langle J_\varepsilon(x)\Delta_F(x-y;\varepsilon)J_\varepsilon(y)\rangle_{xy} = -\frac{1}{2}\int_{\mathbb{R}^D_x} d^Dx\int_{\mathbb{R}^D_y} d^Dy(J_\varepsilon(x)\Delta_F(x-y;\varepsilon)J_\varepsilon(y))_\varepsilon \tag{7.5.17}$$

and

$$(\Delta_{xy;\varepsilon})_\varepsilon \triangleq (\Delta_F(x-y;\varepsilon))_\varepsilon = \frac{1}{(2\pi)^D}\int_{\mathbb{R}^D_k}(d(\varrho_\varepsilon(p))_\varepsilon)\frac{e^{ik\cdot(x-y)}}{p^2+m^2}, \tag{7.5.18}$$

where $p^2 = \sum_{\mu=1}^D p_\mu^2$. An simple algebraic rearrangement gives

$$(Z^E_\varepsilon[J_\varepsilon])_\varepsilon = -\ln\mathbf{N} + (Z^E_{0,\varepsilon}[J_\varepsilon])_\varepsilon - \ln\left(1+\exp\left[(Z^E_{0,\varepsilon}[J_\varepsilon])_\varepsilon\right]\left\{\exp\left[-\left\langle V\left(\frac{\delta}{\delta J_\varepsilon}\right)_\varepsilon\right\rangle\right]-1\right\}\exp\left[-(Z^E_{0,\varepsilon}[J_\varepsilon])_\varepsilon\right]\right) \tag{7.5.19}$$

which is ready for a perturbative expansion in the potential $V$. If we let

$$(\delta_\varepsilon)_\varepsilon = \left\{\exp\left[-\left\langle V\left(\frac{\delta}{\delta J_\varepsilon}\right)_\varepsilon\right\rangle\right]-1\right\}\exp\left[-(Z^E_{0,\varepsilon}[J_\varepsilon])_\varepsilon\right] \tag{4.5.20}$$

we obtain the expansion

$$(Z^E_\varepsilon[J_\varepsilon])_\varepsilon = -\ln\mathbf{N} + (Z^E_{0,\varepsilon}[J_\varepsilon])_\varepsilon - (\delta_\varepsilon[J_\varepsilon])_\varepsilon + \frac{1}{2}(\delta^2_\varepsilon[J_\varepsilon])_\varepsilon - \frac{1}{3}(\delta^3_\varepsilon[J_\varepsilon])_\varepsilon + \ldots \tag{7.5.21}$$

In particular for $V = \frac{\lambda}{D!}\phi^4$, we can expand $(Z^E_\varepsilon[J_\varepsilon])_\varepsilon$ in powers of the dimensionless coupling constant $\lambda$. Setting

$$(\delta_\varepsilon)_\varepsilon = \lambda(\delta_{1,\varepsilon})_\varepsilon + \lambda^2(\delta_{2,\varepsilon})_\varepsilon + \lambda^3(\delta_{3,\varepsilon})_\varepsilon\ldots \tag{7.5.22}$$

we obtain

$$\begin{aligned}(Z^E_\varepsilon[J_\varepsilon])_\varepsilon &= -\ln\mathbf{N} + (Z^E_{0,\varepsilon}[J_\varepsilon])_\varepsilon - \\ &\lambda(\delta_{1,\varepsilon}[J_\varepsilon])_\varepsilon - \lambda^2\left[(\delta_{2,\varepsilon}[J_\varepsilon])_\varepsilon - \frac{1}{2}(\delta^2_{1,\varepsilon}[J_\varepsilon])_\varepsilon\right] - \\ &\lambda^3\left[(\delta_{3,\varepsilon}[J_\varepsilon])_\varepsilon - ((\delta_{1,\varepsilon}[J_\varepsilon])_\varepsilon)((\delta_{2,\varepsilon}[J_\varepsilon])_\varepsilon) + \frac{1}{3}(\delta^3_{1,\varepsilon}[J_\varepsilon])_\varepsilon\right] + \ldots_\varepsilon\end{aligned} \tag{7.5.23}$$

By expanding the exponential in (7.5.20) we obtain

$$(\delta_{1,\varepsilon}[J_\varepsilon])_\varepsilon = -\frac{1}{4!}\exp\left((Z^E_{0,\varepsilon}[J_\varepsilon])_\varepsilon\right)\frac{\delta^4}{\delta(J^4_\varepsilon)_\varepsilon}\exp\left(-(Z^E_{0,\varepsilon}[J_\varepsilon])_\varepsilon\right) \tag{7.5.24}$$

and

$$(\delta_{2,\varepsilon}[J_\varepsilon])_\varepsilon = -\frac{1}{2(4!)^2}\exp\left((Z^E_{0,\varepsilon}[J_\varepsilon])_\varepsilon\right)\left\langle\frac{\delta^4}{\delta(J^4_{1,\varepsilon})_\varepsilon}\right\rangle\left\langle\frac{\delta^4}{\delta(J^4_{2,\varepsilon})_\varepsilon}\right\rangle\exp\left(-(Z^E_{0,\varepsilon}[J_\varepsilon])_\varepsilon\right), \tag{7.5.25}$$

etc. ...

Using the explicit form (7.5.17) for $(Z^E_{0,\varepsilon}[J_\varepsilon])_\varepsilon$ , we obtain

$$(\delta_{1,\varepsilon}[J_\varepsilon])_\varepsilon = -\frac{1}{4!}[(\langle\Delta_{xa;\varepsilon}\Delta_{xb;\varepsilon}\Delta_{xc;\varepsilon}\Delta_{xd;\varepsilon}J_{a,\varepsilon}J_{b,\varepsilon}J_{c,\varepsilon}J_{d,\varepsilon}\rangle)_\varepsilon + 6(\langle\Delta_{xx;\varepsilon}\Delta_{xa;\varepsilon}\Delta_{xb;\varepsilon}J_{a,\varepsilon}J_{b,\varepsilon}\rangle)_\varepsilon + 3(\langle\Delta^2_{xx;\varepsilon}\rangle)_\varepsilon] \tag{7.5.26}$$

where all variables $x, a, b, c, d,$ are integrated over in the relevant $\langle\ldots\rangle$. Similarly, we evaluate $(\delta_{2,\varepsilon})_\varepsilon$ in a slightly trickier fashion: we note that:

$$(\delta_{2,\varepsilon}[J_\varepsilon])_\varepsilon = -\frac{1}{2(4!)^2}\exp\left((Z^E_{0,\varepsilon}[J_\varepsilon])_\varepsilon\right)\left\langle\frac{\delta^4}{\delta(J^4_{1,\varepsilon})_\varepsilon}\right\rangle_1 \exp\left(-(Z^E_{0,\varepsilon}[J_\varepsilon])_\varepsilon\right). \tag{7.5.27}$$

By inserting $\exp\left((Z^E_{0,\varepsilon}[J_\varepsilon])_\varepsilon\right)\exp\left(-(Z^E_{0,\varepsilon}[J_\varepsilon])_\varepsilon\right)$ in the middle of (7.5.25). Next the expansion

$$\begin{gathered}\frac{\delta^4}{\delta(J^4_\varepsilon)_\varepsilon}\exp\left(-(Z^E_{0,\varepsilon}[J_\varepsilon])_\varepsilon\right) = \\ \frac{\delta^4\exp\left(-(Z^E_{0,\varepsilon}[J_\varepsilon])_\varepsilon\right)}{\delta(J^4_\varepsilon)_\varepsilon} + 4\frac{\delta^3\exp\left(-(Z^E_{0,\varepsilon}[J_\varepsilon])_\varepsilon\right)}{\delta(J^3_\varepsilon)_\varepsilon}\frac{\delta}{\delta(J_\varepsilon)_\varepsilon} + 6\frac{\delta^2\exp\left(-(Z^E_{0,\varepsilon}[J_\varepsilon])_\varepsilon\right)}{\delta(J^2_\varepsilon)_\varepsilon} \\ +4\frac{\delta\exp\left(-(Z^E_{0,\varepsilon}[J_\varepsilon])_\varepsilon\right)}{\delta(J_\varepsilon)_\varepsilon}\frac{\delta^3}{\delta(J^3_\varepsilon)_\varepsilon} + \exp\left(-(Z^E_{0,\varepsilon}[J_\varepsilon])_\varepsilon\right)\frac{\delta^4}{\delta(J^4_\varepsilon)_\varepsilon}\end{gathered} \tag{7.5.28}$$

allows us to write

$$\begin{gathered}(\delta_{2,\varepsilon}[J_\varepsilon])_\varepsilon = \frac{1}{2}(\delta^2_{1,\varepsilon})_\varepsilon - \\ \frac{1}{2(4!)}\exp\left((Z^E_{0,\varepsilon}[J_\varepsilon])_\varepsilon\right)\left\langle\left(4\frac{\delta^3\exp\left(-(Z^E_{0,\varepsilon}[J_\varepsilon])_\varepsilon\right)}{\delta(J^3_{1,\varepsilon})_\varepsilon}\frac{\delta}{\delta(J_{1,\varepsilon})_\varepsilon} + \right.\right. \\ 6\frac{\delta^2\exp\left(-(Z^E_{0,\varepsilon}[J_\varepsilon])_\varepsilon\right)}{\delta(J^2_{1,\varepsilon})_\varepsilon}\frac{\delta^2}{\delta(J^2_{1,\varepsilon})_\varepsilon} + 4\frac{\delta\exp\left(-(Z^E_{0,\varepsilon}[J_\varepsilon])_\varepsilon\right)}{\delta(J_{1,\varepsilon})_\varepsilon}\frac{\delta^3}{\delta(J^3_{1,\varepsilon})_\varepsilon} + \\ \left.\left.\exp\left(-(Z^E_{0,\varepsilon}[J_\varepsilon])_\varepsilon\right)\frac{\delta^4}{\delta(J^4_{1,\varepsilon})_\varepsilon}\right)\right\rangle_1 (\delta_{1,\varepsilon}[J_\varepsilon])_\varepsilon\end{gathered} \tag{7.5.29}$$

Comparison with the expansion (7.5.23) for $(Z^E_\varepsilon[J_\varepsilon])_\varepsilon$ shows that the "disconnected" part $\frac{1}{2}(\delta^2_{1,\varepsilon})_\varepsilon$ drops out. By disconnected we mean a contribution which can be written as the product of two or more functionals of $(J_\varepsilon)_\varepsilon$. This concept will become obvious in the diagrammatic representation. The fact that $(Z^E_\varepsilon[J_\varepsilon])_\varepsilon$ generates only connected pieces is true to all orders. For example, the order $\lambda^3$ contribution in (7.5.23) is connected: write

$$\begin{gathered}(\delta_{3,\varepsilon})_\varepsilon = -\frac{1}{3!}\left\langle\exp\left((Z^E_{0,\varepsilon}[J_\varepsilon])_\varepsilon\right)V_xV_yV_z\exp\left(-(Z^E_{0,\varepsilon}[J_\varepsilon])_\varepsilon\right)\right\rangle_{xyz} = \\ -\frac{1}{3!}\left\langle\left[\exp\left((Z^E_{0,\varepsilon}[J_\varepsilon])_\varepsilon\right)V_x\exp\left(-(Z^E_{0,\varepsilon}[J_\varepsilon])_\varepsilon\right)\right]\times\right. \\ \left[\exp\left((Z^E_{0,\varepsilon}[J_\varepsilon])_\varepsilon\right)V_y\exp\left(-(Z^E_{0,\varepsilon}[J_\varepsilon])_\varepsilon\right)\right]\times \\ \left.\left[\exp\left((Z^E_{0,\varepsilon}[J_\varepsilon])_\varepsilon\right)V_z\exp\left(-(Z^E_{0,\varepsilon}[J_\varepsilon])_\varepsilon\right)\right]\right\rangle_{xyz} - \\ \frac{1}{2}\left\langle\left[\exp\left((Z^E_{0,\varepsilon}[J_\varepsilon])_\varepsilon\right)V_x\exp\left(-(Z^E_{0,\varepsilon}[J_\varepsilon])_\varepsilon\right)\right]\times\right. \\ \left.\left[\exp\left((Z^E_{0,\varepsilon}[J_\varepsilon])_\varepsilon\right)V_yV_z\exp\left(-(Z^E_{0,\varepsilon}[J_\varepsilon])_\varepsilon\right)\right]\right\rangle_{xyz} + (\delta^c_{3,\varepsilon})_\varepsilon.\end{gathered} \tag{7.5.30}$$

In the above $(\delta^c_{2,\varepsilon})_\varepsilon$and $(\delta^c_{3,\varepsilon})_\varepsilon$ stans for the connected pieces. To arrive at this form, we have used the fact that there are only two types of "disconnectedness": all three $x,y,z$ disconnected, and only one disconnected from the other two; and there are three ways to obtain the latter possibility. The parentheses in (4.5.29) serve to shield other terms from the action of the derivative operators within them.It follows that the term appearing in the expansion of $(Z^E_\varepsilon[J_\varepsilon])_\varepsilon$ can be rewritten using (4.5.30) in the following form

$$\begin{gathered}(\delta_{3,\varepsilon})_\varepsilon - (\delta_{1,\varepsilon})_\varepsilon(\delta_{2,\varepsilon})_\varepsilon + \frac{1}{3}(\delta^3_{1,\varepsilon})_\varepsilon = (\delta^c_{3,\varepsilon})_\varepsilon + \frac{1}{3!}(\delta^3_{1,\varepsilon})_\varepsilon + (\delta_{1,\varepsilon})_\varepsilon(\delta^c_{2,\varepsilon})_\varepsilon - \\ (\delta_{1,\varepsilon})_\varepsilon\Big((\delta^c_{2,\varepsilon})_\varepsilon + \frac{1}{2}(\delta^2_{1,\varepsilon})_\varepsilon\Big) + \frac{1}{3}(\delta^3_{1,\varepsilon})_\varepsilon = (\delta^c_{3,\varepsilon})_\varepsilon.\end{gathered} \tag{7.5.31}$$

Now, the explicit evaluation of the connected part of $(\delta_{2,\varepsilon})_\varepsilon$ yields, save for the $(J_\varepsilon)_\varepsilon$-independent part,

$$\begin{gathered}(\delta^c_{2,\varepsilon}[J_\varepsilon])_\varepsilon = \frac{1}{2}\Big(\Big\langle J_{a;\varepsilon}\Delta_{ax;\varepsilon}\Big(\frac{1}{6}\Delta^3_{xy;\varepsilon} + \frac{1}{4}\Delta_{xx;\varepsilon}\Delta_{yy;\varepsilon}\Delta_{xy;\varepsilon}\Big)\Delta_{yb;\varepsilon}J_{b;\varepsilon}\Big\rangle_{xyab}\Big)_\varepsilon + \\ \frac{1}{8}\Big(\langle J_{a;\varepsilon}\Delta_{ax;\varepsilon}\Delta_{yy;\varepsilon}\Delta^2_{xy;\varepsilon}\Delta_{xb;\varepsilon}J_{b;\varepsilon}\rangle_{xyab}\Big)_\varepsilon + \\ \frac{2}{4!}\Big(\langle J_{a;\varepsilon}\Delta_{ax;\varepsilon}\Delta_{xx;\varepsilon}\Delta_{yy;\varepsilon}\Delta_{xy;\varepsilon}\Delta_{yb;\varepsilon}\Delta_{yc;\varepsilon}\Delta_{yd;\varepsilon}J_{b;\varepsilon}J_{c;\varepsilon}J_{d;\varepsilon}\rangle_{xyabcd}\Big)_\varepsilon + \\ \frac{3}{2(4!)}\Big(\langle J_{a;\varepsilon}J_{b;\varepsilon}\Delta_{ax;\varepsilon}\Delta_{bx;\varepsilon}\Delta^2_{xy;\varepsilon}\Delta_{yc;\varepsilon}\Delta_{yd;\varepsilon}J_{c;\varepsilon}J_{d;\varepsilon}\rangle_{xyabcd}\Big)_\varepsilon + \\ \frac{1}{2(3!)^2}\Big(\langle J_{a;\varepsilon}J_{b;\varepsilon}J_{c;\varepsilon}\Delta_{ax;\varepsilon}\Delta_{bx;\varepsilon}\Delta_{cx;\varepsilon}\Delta_{xy;\varepsilon}\Delta_{yd;\varepsilon}\Delta_{ye;\varepsilon}\Delta_{yf;\varepsilon}J_{d;\varepsilon}J_{e;\varepsilon}J_{f;\varepsilon}\rangle_{xyabcdef}\Big)_\varepsilon.\end{gathered} \tag{7.5.32}$$

The resulting connected Green's functions follow from Eq.(7.5.15)

$$\begin{gathered}\Big(G^{(2)}_E(x_1,x_2;\varepsilon)\Big)_\varepsilon = \Delta(x_1-x_2;\varepsilon) - \\ -\frac{\lambda}{2}\int d^Dy(\Delta(x_1-y;\varepsilon)\Delta(y-y;\varepsilon)\Delta(y-x_2;\varepsilon))_\varepsilon + \\ \frac{\lambda^2}{6}\int d^Dxd^Dy(\Delta(x_1-x;\varepsilon)\Delta^3(x-y;\varepsilon)\Delta(y-x_2;\varepsilon))_\varepsilon + \\ \frac{\lambda^2}{4}\int d^Dxd^Dy(\Delta(x_1-x;\varepsilon)\Delta^2(x-y;\varepsilon)\Delta(y-y;\varepsilon)\Delta(x-x_2;\varepsilon))_\varepsilon + \\ \frac{\lambda^2}{4}\int d^Dxd^Dy \times \\ (\Delta(x_1-x;\varepsilon)\Delta(x-x;\varepsilon)\Delta(x-y;\varepsilon)\Delta(y-y;\varepsilon)\Delta(y-x_2;\varepsilon))_\varepsilon + (O_\varepsilon(\lambda^3))_\varepsilon,\end{gathered} \tag{7.5.33}$$

$$\begin{gathered}\Big(G^{(4)}_E(x_1,x_2,x_3,x_4;\varepsilon)\Big)_\varepsilon = \\ -\lambda\int d^Dxd^Dy(\Delta(x_1-x;\varepsilon)\Delta^2(x_2-x;\varepsilon)\Delta(x_3-x;\varepsilon)\Delta(x_4-x;\varepsilon))_\varepsilon + \\ \frac{\lambda^2}{2}\int d^Dxd^Dy(\Delta^2(x-y;\varepsilon))_\varepsilon \times \\ [(\Delta(x_1-x;\varepsilon)\Delta^2(x_2-x;\varepsilon)\Delta(x_3-y;\varepsilon)\Delta(x_4-y;\varepsilon))_\varepsilon + \\ (\Delta(x_1-x;\varepsilon)\Delta(x_3-x;\varepsilon)\Delta(x_2-y;\varepsilon)\Delta(x_4-y;\varepsilon))_\varepsilon + \\ (\Delta(x_1-x;\varepsilon)\Delta(x_4-x;\varepsilon)\Delta(x_2-y;\varepsilon)\Delta(x_4-y;\varepsilon))_\varepsilon] + \\ \frac{\lambda^2}{2}\int d^Dxd^Dy(\Delta(y-y;\varepsilon)\Delta(x-y;\varepsilon))_\varepsilon \times \\ \Big[(\Delta(x_1-x;\varepsilon)\Delta(x_2-x;\varepsilon)\Delta(x_3-x;\varepsilon)\Delta(x_4-y;\varepsilon))_\varepsilon + \text{cyclic permutations}\Big] + \\ +(O_\varepsilon(\lambda^3))_\varepsilon,\end{gathered} \tag{7.5.34}$$

and finally

$$\left(G_E^{(6)}(x_1,x_2,x_3,x_4,x_5,x_6;\varepsilon)\right)_\varepsilon = \lambda^2 \int d^Dx d^Dy (\Delta(x-y;\varepsilon))_\varepsilon \times$$
$$\sum_{(ijk)} (\Delta(x_i - x;\varepsilon)\Delta^2(x_j - x;\varepsilon)\Delta(x_k - x;\varepsilon)\Delta(x_l - y;\varepsilon)\Delta(x_m - y;\varepsilon)\Delta(x_n - y;\varepsilon))_\varepsilon + \tag{7.5.35}$$
$$+(O_\varepsilon(\lambda^3))_\varepsilon,$$

where the sum in the last expression runs over the triples $(ijk) = (123),(124),(125),(126),\ (134),(135),(136),(145),(146),(156),$ with (ton) assuming the complementary value, i.e., $(lmn) = (456)$ when $(ijk) = (123)$, etc..The remaining Green's functions get no contribution to this order in $\lambda$.It is straightforward to derive the $p$-space Green's functions, using the canonical expression

$$\left(\widehat{G}_E^{(N)}(p_1,\ldots,p_N;\varepsilon)\right)_\varepsilon = (2\pi)^D \delta(p_1,\ldots,p_N) \times$$
$$\left(\int d^Dx_1 \ldots \int d^Dx_N \exp[i(p_1x_1 + \ldots + p_Nx_N)] G_E^{(N)}(x_1,\ldots,x_N;\varepsilon)\right)_\varepsilon = \tag{7.5.36}$$
$$\int d^Dx_1 \ldots \int d^Dx_N \exp[i(p_1x_1 + \ldots + p_Nx_N)] \left(G_E^{(N)}(x_1,\ldots,x_N;\varepsilon)\right)_\varepsilon$$

Finally we obtain

$$\left(\widehat{G}_E^{(2)}(p,-p;\varepsilon)\right)_\varepsilon =$$
$$\frac{1}{\left(|p|^{D|\alpha|} + \varepsilon\right)_\varepsilon (p^2 + m^2)} -$$
$$-\frac{\lambda}{2} \frac{1}{\left(|p|^{2D|\alpha|} + \varepsilon\right)_\varepsilon (p^2+m^2)^2} \int \frac{d^Dq}{(2\pi)^D} \frac{1}{\left(|q|^{D|\alpha|} + \varepsilon\right)_\varepsilon (q^2+m^2)} +$$
$$+\frac{\lambda^2}{6} \frac{1}{|p|^{2D|\alpha|}(p^2+m^2)^2} \times$$
$$\int \frac{d^Dq_1}{(2\pi)^D} \frac{d^Dq_2}{(2\pi)^D} \frac{d^Dq_3}{(2\pi)^D} \times$$
$$\frac{\delta(p - q_1 - q_2 - q_3)(2\pi)^D}{\left(|q_1|^{D|\alpha|} + \varepsilon\right)_\varepsilon \left(|q_2|^{D|\alpha|} + \varepsilon\right)_\varepsilon \left(|q_3|^{D|\alpha|} + \varepsilon\right)_\varepsilon (q_1^2+m^2)(q_2^2+m^2)(q_3^2+m^2)} + \tag{7.5.37}$$
$$+\frac{\lambda^2}{4} \frac{1}{\left(|p|^{2D|\alpha|} + \varepsilon\right)_\varepsilon (p^2+m^2)^2} \int \frac{d^Dq}{(2\pi)^D} \frac{1}{\left(|q|^{D|\alpha|} + \varepsilon\right)_\varepsilon (q^2+m^2)} \times$$
$$\int \frac{d^Dl_1}{(2\pi)^D} \int \frac{d^Dl_2}{(2\pi)^D} \times$$
$$\frac{(2\pi)^D \delta(l_1 - l_2)}{(|l_1|^{D|\alpha|} + \varepsilon)_\varepsilon (|l_2|^{D|\alpha|} + \varepsilon)_\varepsilon (l_1^2+m^2)(l_2^2+m^2)} +$$
$$+\frac{\lambda^2}{4} + \frac{1}{\left(|p|^{2D|\alpha|} + \varepsilon\right)(p^2+m^2)^2} \int \frac{d^Dq}{(2\pi)^D} \frac{1}{|q|^{D|\alpha|}(q^2+m^2)} \frac{1}{|p|^{2D|\alpha|}(p^2+m^2)^2} \times$$
$$\int \frac{d^Dl}{(2\pi)^D} \frac{1}{(|l|^{D|\alpha|} + \varepsilon)_\varepsilon (l^2+m^2)} + (O_\varepsilon(\lambda^3))_\varepsilon,$$

and

$$\left(\widehat{G}_E^{(4)}(p_1,p_2,p_3p_4;\varepsilon)\right)_\varepsilon =$$

$$\prod_{i=1}^{4}\frac{1}{\left(|p_i|^{D|\alpha|}+\varepsilon\right)_\varepsilon(p_i^2+m^2)}\left\{-\lambda+\frac{\lambda^2}{2}\int\frac{d^Dq}{(2\pi)^D}\frac{1}{\left(|q|^{D|\alpha|}+\varepsilon\right)_\varepsilon(q^2+m^2)}\times\right.$$
$$\sum\frac{1}{\left(|p_i|^{D|\alpha|}+\varepsilon\right)_\varepsilon(p_i^2+m^2)}+\frac{\lambda^2}{2}\int\frac{d^Dq_1}{(2\pi)^D}\frac{d^Dq_2}{(2\pi)^D}\times \quad (7.5.38)$$
$$\left.\frac{(2\pi)^D\sum_{\{ij\}}\delta(q_1+q_2-p_i-p_j)}{\left(|q_1|^{D|\alpha|}+\varepsilon\right)_\varepsilon(q_1^2+m^2)\left(|q_2|^{D|\alpha|}+\varepsilon\right)_\varepsilon(q_2^2+m^2)}\right\}+(O_\varepsilon(\lambda^3))_\varepsilon.$$

In the last expression, the sum $ij$ runs over $(ij) = (12),(13),(14)$ only.

**Remark 7.5.3.**Note that Fourier transform in Eq.(7.5.36) meant integration in sense of generalized function, i.e. an regularization is needed.

**Remark 7.5.4.** We aply Gel'fand regularization [40]. Note that the wollowing equality holds

$$\int_0^\infty r^\lambda\varphi(r)dr = \int_0^{r_1} r^\lambda[\varphi(r)-\varphi(0)]dr+\int_{r_1}^\infty r^\lambda\varphi(r)dr+\frac{\varphi(0)}{\lambda+1} \quad (7.5.39)$$

where $-2 < \operatorname{Re}\lambda < -1.$ Similarly as Eq.(4.7.39) one obtains

$$\int_0^\infty \frac{\varphi(r)}{(r^\lambda+\varepsilon)_\varepsilon}dr = \int_0^{r_1}\frac{[\varphi(r)-\varphi(0)]}{(r^\lambda+\varepsilon)_\varepsilon}dr+\int_{r_1}^\infty\frac{\varphi(r)dr}{(r^\lambda+\varepsilon)_\varepsilon}+\frac{\varphi(0)}{\lambda+1}. \quad (7.5.40)$$

## 7.6.Power-counting renormalizability of Einstein gravity in negative dimensions.

In the context of quantum field theory, the main obstacle against perturbative renormalizability of Einstein's theory of gravity in $D^+ = 3+1$ dimensions is well known [37].

The main problem is that the gravitational coupling constant $G_N$ is dimensionful, with a negative dimension $[G_N] = -2$ in mass units. The Feynman rules also involve the graviton propagator, which scales with the four-momentum $k_\mu = (\omega,\mathbf{k}), \mu = 0,1,2,3$ schematically as $1/k^2$, where $k = \sqrt{\omega^2-\mathbf{k}^2}$ . At increasing loop orders, the Feynman diagrams of this theory require counterterms of ever-increasing degree in curvature. The resulting theory can still be treated as an effective field theory, but it requires a UV completion. An improved UV behavior can be obtained if relativistic higher-derivative corrections are added to the Lagrangian (see [38] for a review of higher-derivative gravity). Terms quadratic in curvature not only yield new interactions (with a dimensionless coupling), they also modify the propagator. Schematically, one gets

$$\frac{1}{k^2}+\frac{1}{k^2}G_Nk^4\frac{1}{k^2}+\frac{1}{k^2}G_Nk^4\frac{1}{k^2}G_Nk^4\frac{1}{k^2}+\ldots=\frac{1}{k^2-G_Nk^4}. \quad (7.6.1)$$

At high energies, the propagator is dominated by the $1/k^4$ term. This cures the problem of UV divergences, and in fact the calculations in Euclidean signature suggest that the theory exhibits asymptotic freedom. However, this cure simultaneously produces a new pathology, which prevents this modified theory from being a solution to the problem of quantum gravity: The resummed propagator (4.4.1) has a two poles,

$$\frac{1}{k^2 - G_N k^4} = \frac{1}{k^2} - \frac{1}{k^2 - 1/G_N}. \tag{7.6.2}$$

# 7.7.Power-counting renormalizability of Hǒrava gravity in negative dimensions.

Let $P(\phi)^z_{d+1}$ be scalar quantum field theory in $(d^- + 1)$ dimensions, where $d^- \leq 0$, containing up to $2z$ spatial derivatives of the $d^- \leq 0$ dimensional spatial metric. We remind that for the $P(\phi)^z_{d^-+1}$ scalar quantum field theories each loop integral has dimension $[\kappa]^{d^{--}+z}$, while each propagator has dimension $[\kappa]^{-2z}$. To analyze the superficial degree of divergence one need only consider the one-particle-irreducible (1PI) sub-diagrams of the Feynman diagram. For each such 1PI sub-diagram the total contribution to dimensionality coming from loop integrals and internal propagators is $[\kappa]^{(d^{--}+z)L-2Iz}$, which is summarized by saying that the "superficial degree of divergence" is

$$\delta = (d^- + z)L - 2Iz = (d^- - z)L - 2(I - L)z. \tag{7.7.1}$$

Note that the quantity $I$ only counts the propagators internal to the 1PI sub-diagram. But to get $L$ loops one needs, at the very least, $I$ internal propagators. So for any 1PI Feynman diagram we certainly have

$$\delta \leq (d^- - z)L. \tag{7.7.2}$$

Consequently, if one picks $d^- \leq 0$ then

$$\delta < 0, \tag{7.7.3}$$

and the worst divergence one can possibly encounter is logarithmic. This observation is enough to guarantee that the non-normal-ordered $P(\phi)^z_{d^{--}+1}$ is power-counting renormalizable, and to render the normal-ordered : $P(\phi)^z_{d^{--}+1}$ : power-counting finite. Furthermore if one takes $d^- \leq 0$ this discussion is sufficient to render $P(\phi)^z_{d^{--}+1}$ (with or without normal ordering) power-counting finite.

Turning our attention now to a $d^- \leq 0$ variant of Hořava gravity in $(d^- + 1)$ dimensions, (containing up to $2z$ spatial derivatives of the $d^-$ dimensional spatial metric), one obtains the same power-counting for the loop integrals and the propagators — the difference now lies in the graviton self-interaction vertices. While the vertices for the scalar field theory carried no factors of momentum, for Hořava gravity and its variants the graviton self-interaction vertices arise from a perturbative action of the form

$$S \sim \int \{\dot{h}^2 + P(\nabla^{2z}, h)\}\ dt\, d^{d^-}x, \tag{7.7.4}$$

where $P(\nabla^{2z}, h)$ is now an infinite-order polynomial in the graviton field $h$, which contains up to $2z$ spatial derivatives.

In contrast to the scalar self-interaction vertices, the graviton self-interaction vertices thus contain up to $2z$ factors of momentum. If these are external momenta they do not contribute to the superficial degree of divergence. However internal momenta, and for any 1PI Feynman diagram with $V$ vertices there can be up to $2zV$ factors of internal momenta, do contribute to the superficial degree of divergence. Consequently we now have the inequality

$$\delta \leq (d^- + z)L + 2z(V - I) = (d^- - z)L + 2z(V + L - I). \tag{7.7.5}$$

But as always, Euler's theorem for graphs implies

$$V + L - I = 1 \quad (7.7.6)$$

so that

$$\delta \leq (d^- - z)L + 2z. \quad (7.7.7)$$

For $|d^-| \geq z$ one simply has

$$\delta \leq 0. \quad (7.7.8)$$

# 8.The solution cosmological constant problem

## 8.1.Einstein-Gliner-Zel'dovich vacuum with tiny Lorentz invariance violation.

We assume now that:

(i) Poincaré group of momentum space is deformed at some fundamental high-energy cutoff $\Lambda_*$ [9],[10].

(ii) The canonical quadratic invariant $\|p\|^2 = \eta^{ab}p_a p_b$ collapses at high-energy cutoff $\Lambda_*$ and being replaced by the non-quadratic invariant:

$$\|p\|^2 = \frac{\eta^{ab}p_a p_b}{(1 + l_{\Lambda_*}p_0)}. \quad (8.1.1)$$

(iii) The canonical concept of Minkowski space-time collapses at a small distances $l_{\Lambda_*} = \Lambda_*^{-1}$ to fractal space-time with Hausdorff-Colombeau negative dimension and therefore the canonical Lebesgue measure $d^4x$ being replaced by the Colombeau-Stieltjes measure with negative Hausdorff-Colombeau dimension $D^-$ :

$$(d\eta(x,\varepsilon))_\varepsilon = (v_\varepsilon(s(x))d^4x)_\varepsilon, \quad (8.1.2)$$

where

$$(v_\varepsilon(s(x)))_\varepsilon = \left(\left(|s(x)|^{|D^-|} + \varepsilon\right)^{-1}\right)_\varepsilon,$$
$$s(x) = \sqrt{x_\mu x^\mu}, \quad (8.1.3)$$

see subsection VI.3.

(iv) The canonical concept of momentum space collapses at fundamental high-energy cutoff $\Lambda_*$ to fractal momentum space with Hausdorff-Colombeau negative dimension and therefore the canonical Lebesgue measure $d^3\mathbf{k}$, where $\mathbf{k} = (k_x, k_y, k_z)$ being replaced by the Hausdorff-Colombeau measure

$$d^{D^+,D^-}\mathbf{k} \triangleq \frac{\Delta(D^-)d^{D^+}\mathbf{k}}{\left(|\mathbf{k}|^{|D^-|} + \varepsilon\right)_\varepsilon} = \frac{\Delta(D^+)\Delta(D^-)p^{D^+-1}dp}{(p^{|D^-|} + \varepsilon)_\varepsilon}, \quad (8.1.4)$$

where $\Delta(D^\pm) = \dfrac{2\pi^{D^\pm/2}}{\Gamma(D^\pm/2)}$ and $p = |\mathbf{k}| = \sqrt{k_x + k_y + k_z}$.

**Remark 8.1.1.** Note that the integral over measure $d^{D^+,D^-}\mathbf{k}$ is given by formula(6.3.32). Thus vacuum energy density $\varepsilon(D^+,D^-,\mu_{\mathbf{eff}},p_*)$ for free quantum fields is

$$\varepsilon(D^+,D^-,\mu_{\mathbf{eff}},p_*) = \varepsilon(\mu_{\mathbf{eff}}) + \varepsilon(\mu_{\mathbf{eff}},p_*) + \breve{\varepsilon}(D^+,D^-,\mu_{\mathbf{eff}},p_*). \quad (8.1.5)$$

Here the quantity $\varepsilon(\mu_{\mathbf{eff}})$ is given by formula

$$\varepsilon(\mu_{\mathbf{eff}}) = \frac{1}{2(2\pi\hbar)^3}\int_0^{\mu_{\mathbf{eff}}} d\mu f(\mu) \int\limits_{\|\mathbf{k}\|\leq\mu} \sqrt{\mathbf{k}^2+\mu^2}\, d^3\mathbf{k} =$$

$$K\int_0^{\mu_{\mathbf{eff}}} d\mu f(\mu) \int\limits_{p\leq\mu} \sqrt{p^2+\mu^2}\, p^2 dp = K\int_0^{\mu_{\mathbf{eff}}} d\mu f(\mu) \int\limits_0^{\mu} \sqrt{p^2+\mu^2}\, p^2 dp \quad (8.1.6)$$

where $K = \frac{2\pi}{(2\pi\hbar)^3}, c = 1$. The quantity $\varepsilon(\mu_{\mathbf{eff}}, p_*)$ is given by formula

$$\varepsilon(\mu_{\mathbf{eff}}, p_*) = \frac{1}{2(2\pi\hbar)^3}\int_0^{\mu_{\mathbf{eff}}} d\mu f(\mu) \int\limits_{\mu<\|\mathbf{k}\|<p_*} \sqrt{\mathbf{k}^2+\mu^2}\, d^3\mathbf{k} =$$

$$K\int_0^{\mu_{\mathbf{eff}}} d\mu f(\mu) \int\limits_{\mu<\|\mathbf{k}\|<p_*} \sqrt{p^2+\mu^2}\, p^2 dp. \quad (8.1.7)$$

The quantity $\breve{\varepsilon}(D^+, D^-, \mu_{\mathbf{eff}}, p_*)$ (since Eq.(1.1.18) holds) is given by formula

$$\breve{\varepsilon}(D^+, D^-, \mu_{\mathbf{eff}}, p_*) =$$

$$K'\int_0^{\mu_{\mathbf{eff}}} d\mu f(\mu) \times \quad (8.1.8)$$

$$\int\limits_{\|\mathbf{k}\|\geq p_*}\left[\frac{\mu^2 l_{\Lambda_*}}{1-\mu^2 l_{\Lambda_*}^2} + \frac{1}{\sqrt{1-\mu^2 l_{\Lambda_*}^2}}\sqrt{\frac{\mu^4 l_{\Lambda_*}^2}{1-\mu^2 l_{\Lambda_*}^2} + (|\mathbf{k}|^2+\mu^2)}\right] d^{D^+,D^-}\mathbf{k},$$

where $K' = \frac{1}{2(2\pi\hbar)^3}, c = 1$.

**Remark 8.1.2**. We assume now that $\mu^2 l_{\Lambda_*}^2 \ll 1, \mu^4 l_{\Lambda_*}^2 \ll 1$ and therefore from Eq.(8.1.8)
we obtain

$$\varepsilon(D^+, D^-, \mu_{\mathbf{eff}}, p_*) =$$

$$K' l_\Lambda \int_0^{\mu_{\mathbf{eff}}} f(\mu)\mu^2 d\mu \int\limits_{\|\mathbf{k}\|\geq p_*} d^{3,D_-}\mathbf{k} + K'\int_0^{\mu_{\mathbf{eff}}} d\mu f(\mu) \int\limits_{\|\mathbf{k}\|\geq p_*} \sqrt{\mathbf{k}^2+\mu^2}\, d^{D^+,D^-}\mathbf{k}. \quad (8.1.9)$$

From Eq.(8.1.9) and Eq.(8.1.4) we obtain

$$\varepsilon(D^+, D^-, \mu_{\mathbf{eff}}, p_*) =$$

$$K' l_\Lambda \int_0^{\mu_{\mathbf{eff}}} f(\mu)\mu^2 d\mu \int\limits_{\|\mathbf{k}\|\geq p_*} d^{D^+,D^-}\mathbf{k} + K'\int_0^{\mu_{\mathbf{eff}}} d\mu f(\mu) \int\limits_{\|\mathbf{k}\|\geq p_*} \sqrt{\mathbf{k}^2+\mu^2}\, d^{D^+,D^-}\mathbf{k} =$$

$$\left(K' l_\Lambda \Delta(D^+)\Delta(D^-)\int_0^{\mu_{\mathbf{eff}}} f d\mu(\mu)\mu^2\right)\int_{p_*}^{\infty} \frac{p^{D^+-1}dp}{(p^{|D^-|}+\varepsilon)_\varepsilon} +$$

$$+K'\Delta(D^+)\Delta(D^-)\int_0^{\mu_{\mathbf{eff}}} d\mu f(\mu)\int_{p_*}^{\infty} \frac{\sqrt{p^2+\mu^2}\, p^{D^+-1}dp}{(p^{|D^-|}+\varepsilon)_\varepsilon} = \quad (8.1.10)$$

$$\left(K' l_\Lambda \Delta(D^+)\Delta(D^-)\int_0^{\mu_{\mathbf{eff}}} f(\mu)\mu^2 d\mu\right)\int_{p_*}^{\infty} p^{D^-+D^+-1}dp +$$

$$+K'\Delta(D^+)\Delta(D^-)\int_0^{\mu_{\mathbf{eff}}} d\mu f(\mu)\int_{p_*}^{\infty} \sqrt{p^2+\mu^2}\, p^{D^-+D^+-1}dp.$$

**Remark 8.1.2.**We assume now that:

$$D^{-}+D^{+}+2 \leq -6. \tag{8.1.11}$$

Note that

$$\begin{gathered}\int_0^{\mu_{\mathbf{eff}}} d\mu f(\mu)\int_{p_*}^{\infty}\sqrt{p^2+\mu^2}\,p^{D^-+D^+-1}dp = \int_0^{\mu_{\mathbf{eff}}} d\mu f(\mu)\int_{p_*}^{\infty}\sqrt{1+\frac{\mu^2}{p^2}}\,p^{D^-+D^+}dp = \\ \int_0^{\mu_{\mathbf{eff}}} f(\mu)d\mu\int_{p_*}^{\infty} p^{D^-+D^+}dp+\frac{1}{2}\int_0^{\mu_{\mathbf{eff}}} f(\mu)\mu^2 d\mu\int_{p_*}^{\infty} p^{D^-+D^+-1}dp - \\ -\frac{1}{8}\int_0^{\mu_{\mathbf{eff}}} f(\mu)\mu^4 d\mu\int_{p_*}^{\infty} p^{D^-+D^+-3}dp+O(p_*^{D^-+D^+-4}) = \\ \frac{p_*^{D^-+D^++1}}{D^-+D^++1}\int_0^{\mu_{\mathbf{eff}}} f(\mu)d\mu+\frac{p_*^{D^-+D^+}}{2(D^-+D^+)}\int_0^{\mu_{\mathbf{eff}}} f(\mu)\mu^2 d\mu - \\ -\frac{p_*^{D^-+D^+-1}}{8(D^-+D^+-1)}\int_0^{\mu_{\mathbf{eff}}} f(\mu)\mu^4 d\mu+O(p_*^{D^-+D^+-4}).\end{gathered} \tag{8.1.12}$$

Thus finally we obtain

$$\begin{gathered}\varepsilon(D^+,D^-,\mu_{\mathbf{eff}},p_*) = \\ \frac{K'p_*^{D^-+D^++1}}{D^-+D^++1}\int_0^{\mu_{\mathbf{eff}}} f(\mu)d\mu+\left([K'l_\Lambda\Delta(D^+)\Delta(D^-)+0.5]\int_0^{\mu_{\mathbf{eff}}} f(\mu)\mu^2 d\mu\right)\frac{p_*^{D^-+D^+}}{D^-+D^+} - \\ -\frac{K'p_*^{D^-+D^+-2}}{8(D^-+D^+-1)}\int_0^{\mu_{\mathbf{eff}}} f(\mu)\mu^4 d\mu+O(p_*^{D^-+D^+-4}).\end{gathered} \tag{5.1.13}$$

**Remark 8.1.3.** Note that (see Eqs.(1.2.12)):

$$\begin{gathered}\widetilde{\varepsilon}(\mu_{\mathbf{eff}},p_*) = \varepsilon(\mu_{\mathbf{eff}})+\varepsilon(\mu_{\mathbf{eff}},p_*) = \\ \frac{1}{4}p_*^4\int_0^{\mu_{\mathbf{eff}}} f(\mu)d\mu+\frac{1}{4}p_*^2\int_0^{\mu_{\mathbf{eff}}} f(\mu)\mu^2 d\mu+\left(C_1-\frac{1}{8}\ln p_*\right)\int_0^{\mu_{\mathbf{eff}}} f(\mu)\mu^4 d\mu + \\ +\frac{1}{8}\int_0^{\mu_{\mathbf{eff}}} f(\mu)\mu^4(\ln\mu)d\mu-\left(\frac{1}{p_*^2}\right)\frac{1}{32}\int_0^{\mu_{\mathbf{eff}}} f(\mu)\mu^6 d\mu+O\left(\int_0^{\mu_{\mathbf{eff}}} f(\mu)\mu^8\right)p_*^{-5}.\end{gathered} \tag{8.1.14}$$

From Eq.(8.1.5),Eq.(8.1.13) and Eq.(8.1.14) finally we obtain

$$\begin{gathered}\varepsilon(D^+,D^-,\mu_{\mathbf{eff}},p_*) = \varepsilon(\mu_{\mathbf{eff}})+\varepsilon(\mu_{\mathbf{eff}},p_*)+\breve{\varepsilon}(D^+,D^-,\mu_{\mathbf{eff}},p_*) = \\ \frac{1}{4}p_*^4\int_0^{\mu_{\mathbf{eff}}} f(\mu)d\mu+\frac{1}{4}p_*^2\int_0^{\mu_{\mathbf{eff}}} f(\mu)\mu^2 d\mu+\left(C_1-\frac{1}{8}\ln p_*\right)\int_0^{\mu_{\mathbf{eff}}} f(\mu)\mu^4 d\mu + \\ +\frac{1}{8}\int_0^{\mu_{\mathbf{eff}}} f(\mu)\mu^4(\ln\mu)d\mu-\left(\frac{1}{p_*^2}\right)\frac{1}{32}\int_0^{\mu_{\mathbf{eff}}} f(\mu)\mu^6 d\mu+O\left(\int_0^{\mu_{\mathbf{eff}}} f(\mu)\mu^8\right)p_*^{-5} + \\ +O(p_*^{D^-+D^++2}).\end{gathered} \tag{8.1.15}$$

The pressure $p(D^+,D^-,\mu_{\mathbf{eff}},p_*)$ for free scalar quantum field is

$$p(D^+,D^-,\mu_{\mathbf{eff}},p_*) = p(\mu_{\mathbf{eff}})+p(\mu_{\mathbf{eff}},p_*)+\breve{p}(D^+,D^-,\mu_{\mathbf{eff}},p_*). \tag{8.1.16}$$

Here the quantity $p(\mu_{\mathbf{eff}})$ is given by formula

$$p(\mu_{\mathbf{eff}}) = \frac{K}{3}\int_0^{\mu_{\mathbf{eff}}} d\mu f(\mu) \int\limits_{\|p\|<\mu} \frac{p^4}{\sqrt{p^2+\mu^2}} dp. \tag{8.1.17}$$

The quantity $p(\mu_{\mathbf{eff}}, p_*)$ is given by formula

$$p(\mu_{\mathbf{eff}}, p_*) = \frac{K}{3}\int_0^{\mu_{\mathbf{eff}}} d\mu f(\mu) \int\limits_{\mu \le \|p\| \le p_*} \frac{p^4}{\sqrt{p^2+\mu^2}} dp. \tag{8.1.18}$$

The quantity $\breve{p}(D^+, D^-, \mu_{\mathbf{eff}}, p_*)$ is given by formula

$$\breve{p}(D^+, D^-, \mu_{\mathbf{eff}}, p_*) \simeq \frac{K'}{3}\int_0^{\mu_{\mathbf{eff}}} d\mu \int\limits_{\|p\|>p_*} f(\mu) \frac{p^4}{\sqrt{p^2+\mu^2}} dp, \tag{8.1.19}$$

where $K' = \dfrac{1}{2(2\pi\hbar)^3}, c = 1.$

**Remark 8.1.4**.Note that (see Eqs.(1.2.12)):

$$\begin{gathered}
\widetilde{p}(\mu_{\mathbf{eff}}, p_*) = p(\mu_{\mathbf{eff}}) + p(\mu_{\mathbf{eff}}, p_*) = \\
\frac{1}{12} p_*^4 \int\limits_0^{\mu_{\mathbf{eff}}} f(\mu) d\mu - \frac{1}{12} p_*^2 \int\limits_0^{\mu_{\mathbf{eff}}} f(\mu)\mu^2 d\mu + \left(C_2 + \frac{1}{8}\ln p_*\right) \int\limits_0^{\mu_{\mathbf{eff}}} f(\mu)\mu^4 d\mu - \\
-\frac{1}{8} \int\limits_0^{\mu_{\mathbf{eff}}} f(\mu)\mu^4 (\ln\mu) d\mu + \left(\frac{5}{p_*^2}\right)\frac{1}{32} \int\limits_0^{\mu_{\mathbf{eff}}} f(\mu)\mu^6 d\mu + O\left(\int\limits_0^{\mu_{\mathbf{eff}}} f(\mu)\mu^8\right) p_*^{-5}.
\end{gathered} \tag{8.1.20}$$

From Eq.(8.1.15),Eq.(8.1.19) and Eq.(8.1.20) similarly as above finally we get

$$\begin{gathered}
p(D^+, D^-, \mu_{\mathbf{eff}}, p_*) = \\
\frac{1}{12} p_*^4 \int\limits_0^{\mu_{\mathbf{eff}}} f(\mu) d\mu - \frac{1}{12} p_*^2 \int\limits_0^{\mu_{\mathbf{eff}}} f(\mu)\mu^2 d\mu + \left(C_2 + \frac{1}{8}\ln p_*\right) \int\limits_0^{\mu_{\mathbf{eff}}} f(\mu)\mu^4 d\mu - \\
-\frac{1}{8} \int\limits_0^{\mu_{\mathbf{eff}}} f(\mu)\mu^4 (\ln\mu) d\mu + \left(\frac{5}{p_*^2}\right)\frac{1}{32} \int\limits_0^{\mu_{\mathbf{eff}}} f(\mu)\mu^6 d\mu + O\left(\int\limits_0^{\mu_{\mathbf{eff}}} f(\mu)\mu^8\right) p_*^{-5} + \\
+O(p_*^{D^- + D^+ + 2}).
\end{gathered} \tag{8.1.21}$$

**Remark 8.1.5**.We assume now that:

$$\int\limits_0^{\mu_{\mathbf{eff}}} f(\mu) d\mu = \int\limits_0^{\mu_{\mathbf{eff}}} f(\mu)\mu^2 d\mu = \int\limits_0^{\mu_{\mathbf{eff}}} f(\mu)\mu^4 d\mu = 0. \tag{8.1.22}$$

From Eq.(8.1.15),Eq.(8.1.21) and Eq.(8.1.22) finally we get

$$\begin{gathered}
\varepsilon \triangleq \varepsilon(D^+, D^-, \mu_{\mathbf{eff}}, p_*) = \frac{1}{8} \int\limits_0^{\mu_{\mathbf{eff}}} f(\mu)\mu^4 (\ln\mu) d\mu + O(p_*^{-2}), \\
p \triangleq (D^+, D^-, \mu_{\mathbf{eff}}, p_*) = -\frac{1}{8} \int\limits_0^{\mu_{\mathbf{eff}}} f(\mu)\mu^4 (\ln\mu) d\mu + O(p_*^{-2}).
\end{gathered} \tag{8.1.23}$$

**Remark 8.1.5**. The fine tuning assumed by (5.1.22) is a problematic in order to obtain the mass distribution $f(\mu)$ wich gives an observed value of $\varepsilon$.

**Remark 8.1.6**. Note that the Eq.(5.1.23) can be obtained without fine-tuning

(8.1.22) which was ussumed in Zel'dovich paper [1].

In order to obtain Eq.(8.1.23) ander strictly weaker conditions we assume now that:

(i)

$$|f(\mu,p_*)| = |f_{s.m.}(\mu,p_*) + f_{g.m.}(\mu,p_*)| = \left(\mu^*_{\mathbf{eff}}\right)^{-n}, \tag{8.1.24}$$

where $\mu^*_{\mathbf{eff}} = \mu_{\mathbf{eff}}(p_*)$, $n = n(p_*) > 0$ is an parametr, $f_{s.m.}(\mu,p_*)$ corresponds to standard matter and where $f_{g.m.}(\mu,p_*)$ corresponds to physical ghost matter,see Eq.(1.2.2).

(ii)

$$\begin{gathered} I_1 = p_*^4 \int\limits_0^{\mu^*_{\mathbf{eff}}} f(\mu,p_*)d\mu \approx 0, I_2 = p_*^2 \int\limits_0^{\mu^*_{\mathbf{eff}}} f(\mu,p_*)\mu^2 d\mu \approx 0, \\ I_3 = \ln p_* \int\limits_0^{\mu^*_{\mathbf{eff}}} f(\mu,p_*)\mu^4 d\mu \approx 0. \end{gathered} \tag{8.1.25}$$

(iii)

$$\left| \, I_1 + I_2 + I_3 \, \right| \ll \left| \int\limits_0^{\mu^*_{\mathbf{eff}}} f(\mu,p_*)\mu^4(\ln\mu)d\mu \right|. \tag{8.1.26}$$

Finally we get

$$\begin{gathered} \varepsilon \triangleq \varepsilon(D^+,D^-,\mu_{\mathbf{eff}},p_*) = \frac{1}{8} \int\limits_0^{\mu^*_{\mathbf{eff}}} f(\mu,p_*)\mu^4(\ln\mu)d\mu + O(p_*^{-2}), \\ p \triangleq (D^+,D^-,\mu_{\mathbf{eff}},p_*) = -\frac{1}{8} \int\limits_0^{\mu^*_{\mathbf{eff}}} f(\mu,p_*)\mu^4(\ln\mu)d\mu + O(p_*^{-2}). \end{gathered} \tag{8.1.26}$$

## 8.2. Zeropoint energy density corresponding to a non-singular Gliner cosmology.

We assume now that

$$\begin{gathered} \int\limits_0^{\mu^*_{\mathbf{eff}}} f(\mu,p_*)d\mu \approx 0, \int\limits_0^{\mu^*_{\mathbf{eff}}} f(\mu,p_*)\mu^4 d\mu < 0, \int\limits_0^{\mu^*_{\mathbf{eff}}} f(\mu,p_*)\mu^2 d\mu > 0, \\ p_* \gg \mu^*_{\mathbf{eff}}. \end{gathered} \tag{8.2.1}$$

From Eq.(8.1.15),Eq.(8.1.21) and (8.2.1) we obtain

$$\varepsilon \triangleq \varepsilon\left(D^{+}, D^{-}, \mu_{\mathbf{eff}}^{*}, p_{*}\right) =$$

$$\frac{1}{4} p_{*}^{2} \int_{0}^{\mu_{\mathbf{eff}}} f(\mu, p_{*}) \mu^{2} d\mu - \left(C_{1} - \frac{1}{8} \ln p_{*}\right) \left| \int_{0}^{\mu_{\mathbf{eff}}^{*}} f(\mu, p_{*}) \mu^{4} d\mu \right| +$$

$$+\frac{1}{8} \int_{0}^{\mu_{\mathbf{eff}}^{*}} f(\mu, p_{*}) \mu^{4} (\ln \mu) d\mu - \left(\frac{1}{p_{*}^{2}}\right) \frac{1}{32} \int_{0}^{\mu_{\mathbf{eff}}^{*}} f(\mu, p_{*}) \mu^{6} d\mu \tag{8.2.2}$$

$$+O\left(\int_{0}^{\mu_{\mathbf{eff}}^{*}} f(\mu, p_{*}) \mu^{8}\right) p_{*}^{-5} +$$

$$+O(p_{*}^{D^{-}+D^{+}+2}),$$

and

$$p \triangleq p\left(D^{+}, D^{-}, \mu_{\mathbf{eff}}^{*}, p_{*}\right) =$$

$$-\frac{1}{12} p_{*}^{2} \int_{0}^{\mu_{\mathbf{eff}}^{*}} f(\mu, p_{*}) \mu^{2} d\mu - \left(C_{2} + \frac{1}{8} \ln p_{*}\right) \left| \int_{0}^{\mu_{\mathbf{eff}}^{*}} f(\mu, p_{*}) \mu^{4} d\mu \right| -$$

$$-\frac{1}{8} \int_{0}^{\mu_{\mathbf{eff}}^{*}} f(\mu, p_{*}) \mu^{4} (\ln \mu) d\mu + \left(\frac{5}{p_{*}^{2}}\right) \frac{1}{32} \int_{0}^{\mu_{\mathbf{eff}}^{*}} f(\mu, p_{*}) \mu^{6} d\mu + \tag{8.2.3}$$

$$O\left(\int_{0}^{\mu_{\mathbf{eff}}^{*}} f(\mu, p_{*}) \mu^{8}\right) p_{*}^{-5} +$$

$$+O(p_{*}^{D^{-}+D^{+}+2})$$

correspondingly.From Eq.(8.2.2) and Eq.(8.2.3) we obtain

$$3p+\varepsilon =$$

$$-\frac{1}{4}p_*^2\int\limits_0^{\mu_{\text{eff}}^*} f(\mu,p_*)\mu^2 d\mu - \left(3C_2+\frac{3}{8}\ln p_*\right)\left|\int\limits_0^{\mu_{\text{eff}}^*} f(\mu,p_*)\mu^4 d\mu\right| -$$

$$-\frac{3}{8}\int\limits_0^{\mu_{\text{eff}}^*} f(\mu,p_*)\mu^4(\ln\mu)d\mu + \left(\frac{5}{p_*^2}\right)\frac{3}{32}\int\limits_0^{\mu_{\text{eff}}^*} f(\mu,p_*)\mu^6 d\mu +$$

$$\frac{1}{4}p_*^2\int\limits_0^{\mu_{\text{eff}}^*} f(\mu,p_*)\mu^2 d\mu - \left(C_1-\frac{1}{8}\ln p_*\right)\left|\int\limits_0^{\mu_{\text{eff}}^*} f(\mu,p_*)\mu^4 d\mu\right| +$$

$$+\frac{1}{8}\int\limits_0^{\mu_{\text{eff}}^*} f(\mu,p_*)\mu^4(\ln\mu)d\mu - \left(\frac{1}{p_*^2}\right)\frac{1}{32}\int\limits_0^{\mu_{\text{eff}}^*} f(\mu,p_*)\mu^6 d\mu = \qquad (8.2.4)$$

$$-\frac{1}{4}\ln p_*\left|\int\limits_0^{\mu_{\text{eff}}^*} f(\mu,p_*)\mu^4 d\mu\right| - (3C_2+C_1)\left|\int\limits_0^{\mu_{\text{eff}}^*} f(\mu,p_*)\mu^4 d\mu\right| -$$

$$\frac{1}{4}\int\limits_0^{\mu_{\text{eff}}^*} f(\mu,p_*)\mu^4(\ln\mu)d\mu +$$

$$+\left(\frac{5}{p_*^2}\right)\frac{1}{16}\int\limits_0^{\mu_{\text{eff}}^*} f(\mu,p_*)\mu^6 d\mu < 0.$$

Therefore under conditions (5.2.1) the inequality

$$-2\varepsilon < 3p+\varepsilon < 0 \qquad (8.2.5)$$

corresponding to Gliner non-singular cosmology [2],[4] is satisfied.

# 8.3. Zeropoint energy density in models with supermassive physical ghost fields.

We assume now that:

(i) ghost fields corresponding to massive spin-2 particle with mass $m_2$ and to massive scalar particle with mass $m_0$ appears (see section 5.1) as real physical fields in action (5.1.1)

**Remark 8.3.1.** Note that their unphysical behavior may be restricted to arbitrarily high-energy cutoff $\Lambda_*$ by an appropriate limitation on the renormalized masses $m_2$ and $m_0$.

Actually, it is only the massive spin-two excitations of the field which give the problem with unitarity and thus require a very large mass (see subsection II.2).

(ii) Poincaré group is deformed at some fundamental high-energy cutoff $\Lambda_*$

$$\Lambda_* = \Lambda_*(m_0,m_2) \ll m_0c^2 < m_2c^2. \qquad (8.3.1)$$

The canonical quadratic invariant $\|p\|^2 = \eta^{ab}p_ap_b$ collapses at high-energy cutoff $\Lambda_*$ and being replaced by the non-quadratic invariant:

$$\|p\|^2 = \frac{\eta^{ab}p_a p_b}{(1+l_{\Lambda_*}p_0)}. \tag{8.3.2}$$

(iii) The canonical concept of Minkowski space-time collapses at a small distances to fractal space-time with Hausdorff-Colombeau negative dimension and therefore the canonical Lebesgue measure $d^4x$ being replaced by the Colombeau-Stieltjes measure

$$(d\eta(x,\varepsilon))_\varepsilon = (v_\varepsilon(s(x))d^4x)_\varepsilon, \tag{8.3.3}$$

where

$$(v_\varepsilon(s(x)))_\varepsilon = \left(\left(|s(x)|^{|D^-|}+\varepsilon\right)^{-1}\right)_\varepsilon, s(x) = \sqrt{x_\mu x^\mu}, \tag{8.3.4}$$

(iv) we assume now that

$$f(\mu) = f_{s.m.}(\mu) + f_{g.m.}(\mu), \tag{8.3.5}$$

where $f_{s.m.}(\mu)$ corresponds to standard matter and where $f_{g.m.}(\mu)$ corresponds to physical ghost matter.

**Remark 8.3.2.** We assume now that

$$|f(\mu)| = \begin{cases} O(\mu^{-n}), n>1 & m_0 c \ll \mu_{\mathbf{eff}}^{(1)} \leq \mu \leq \mu_{\mathbf{eff}}^{(2)} \ll m_2 c \\ 0 & \mu > \mu_{\mathbf{eff}}^{(2)} \end{cases} \tag{8.3.6}$$

Thus vacuum energy density $\varepsilon\left(D^+,D^-,\mu_{\mathbf{eff}}^{(1)},\mu_{\mathbf{eff}}^{(2)}\right)$ for free quantum fields is

$$\varepsilon\left(D^+,D^-,\mu_{\mathbf{eff}}^{(1)},\mu_{\mathbf{eff}}^{(2)}\right) = \varepsilon\left(\mu_{\mathbf{eff}}^{(1)},\mu_{\mathbf{eff}}^{(2)},\right) + \breve{\varepsilon}\left(D^+,D^-,\mu_{\mathbf{eff}}^{(1)},\mu_{\mathbf{eff}}^{(2)}\right). \tag{8.3.7}$$

Here the quantity $\varepsilon\left(\mu_{\mathbf{eff}}^{1},\mu_{\mathbf{eff}}^{2},\right)$ is given by formula

$$\begin{aligned}\varepsilon\left(\mu_{\mathbf{eff}}^{(1)},\mu_{\mathbf{eff}}^{(2)}\right) &= \frac{1}{2(2\pi\hbar)^3}\int_{\mu_{\mathbf{eff}}^{(1)}}^{\mu_{\mathbf{eff}}^{(2)}} d\mu f(\mu) \int_{\|\mathbf{k}\|\leq\sqrt{\mu}} \sqrt{\mathbf{k}^2+\mu^2}\, d^3\mathbf{k} = \\ &= K\int_{\mu_{\mathbf{eff}}^{(1)}}^{\mu_{\mathbf{eff}}^{(2)}} d\mu f(\mu) \int_{p\leq\sqrt{\mu}} \sqrt{p^2+\mu^2}\, p^2 dp,\end{aligned} \tag{8.3.8}$$

where $K = \frac{2\pi}{(2\pi\hbar)^3}, c = 1$. The quantity $\breve{\varepsilon}\left(D^+,D^-,\mu_{\mathbf{eff}}^{(1)},\mu_{\mathbf{eff}}^{(2)}\right)$ is given by formula

$$\begin{aligned}&\breve{\varepsilon}\left(D^+,D^-,\mu_{\mathbf{eff}}^{(1)},\mu_{\mathbf{eff}}^{(2)}\right) = \\ &K'\int_{\mu_{\mathbf{eff}}^{(1)}}^{\mu_{\mathbf{eff}}^{(2)}} d\mu f(\mu) \times \\ &\int_{\|\mathbf{k}\|>\sqrt{\mu}} \left[\frac{\mu^2 l_\Lambda}{1-\mu^2 l_\Lambda^2} + \frac{1}{\sqrt{1-\mu^2 l_{\Lambda_*}^2}}\sqrt{\frac{\mu^4 l_{\Lambda_*}^2}{1-\mu^2 l_{\Lambda_*}^2} + (|\mathbf{k}|^2+\mu^2)}\right] d^{D^+,D^-}\mathbf{k},\end{aligned} \tag{8.3.9}$$

where $K' = \frac{1}{2(2\pi\hbar)^3}, c = 1.$

**Remark 8.3.2**. We assume now that $\mu^2 l_{\Lambda_*}^2 < 1$, and therefore from Eq.(8.3.9) we obtain

$$\varepsilon\left(D^+,D^-,\mu^{(1)}_{\mathbf{eff}},\mu^{(2)}_{\mathbf{eff}}\right) \simeq$$

$$K' l_\Lambda \int_{\mu^{(1)}_{\mathbf{eff}}}^{\mu^{(2)}_{\mathbf{eff}}} d\mu f(\mu)\mu^2 \int\limits_{\|\mathbf{k}\|>\sqrt{\mu}} d^{3,D^-}\mathbf{k} + K' \int_{\mu^{(1)}_{\mathbf{eff}}}^{\mu^{(2)}_{\mathbf{eff}}} d\mu f(\mu) \int\limits_{\|\mathbf{k}\|>\sqrt{\mu}} \sqrt{\mathbf{k}^2+\mu^2}\, d^{D^+,D^-}\mathbf{k}. \tag{8.3.10}$$

From Eq.(8.3.10) and Eq.(8.1.4) we obtain

$$\begin{aligned}
&\varepsilon\left(D^+,D^-,\mu^{(1)}_{\mathbf{eff}},\mu^{(2)}_{\mathbf{eff}}\right) \simeq \\
&K' l_\Lambda \int_{\mu^{(1)}_{\mathbf{eff}}}^{\mu^{(2)}_{\mathbf{eff}}} d\mu f(\mu)\mu^2 \int\limits_{\|\mathbf{k}\|>\sqrt{\mu}} d^{D^+,D^-}\mathbf{k} + K' \int_{\mu^{(1)}_{\mathbf{eff}}}^{\mu^{(2)}_{\mathbf{eff}}} d\mu f(\mu) \int\limits_{\|\mathbf{k}\|>\sqrt{\mu}} \sqrt{\mathbf{k}^2+\mu^2}\, d^{D^+,D^-}\mathbf{k} = \\
&K'\Delta(D^+)\Delta(D^-) l_\Lambda \int_{\mu^{(1)}_{\mathbf{eff}}}^{\mu^{(2)}_{\mathbf{eff}}} d\mu f(\mu)\mu^2 \left[\int_{\sqrt{\mu}}^{\infty} \frac{p^{D^+-1}dp}{(p^{|D^-|}+\varepsilon)_\varepsilon}\right] + \\
&+K'\Delta(D^+)\Delta(D^-) \int_{\mu^{(1)}_{\mathbf{eff}}}^{\mu^{(2)}_{\mathbf{eff}}} d\mu f(\mu) \left[\int_{\sqrt{\mu}}^{\infty} \frac{\sqrt{p^2+\mu^2}\, p^{D^+-1}dp}{(p^{|D^-|}+\varepsilon)_\varepsilon}\right] = \\
&K'\Delta(D^+)\Delta(D^-) l_\Lambda \int_{\mu^{1}_{\mathbf{eff}}}^{\mu^{2}_{\mathbf{eff}}} d\mu f(\mu)\mu^2 \left[\int_{\sqrt{\mu}}^{\infty} p^{D^-+D^+-1}dp\right] + \\
&+K'\Delta(D^+)\Delta(D^-) \int_{\mu^{1}_{\mathbf{eff}}}^{\mu^{2}_{\mathbf{eff}}} d\mu f(\mu) \left[\int_{\sqrt{\mu}}^{\infty} \sqrt{p^2+\mu^2}\, p^{D^-+D^+-1}dp\right].
\end{aligned} \tag{8.3.11}$$

Note that

$$\begin{aligned}
\sqrt{p^2+\mu^2} &= \mu\sqrt{1+\frac{p^2}{\mu^2}} = \mu\left(1+\frac{1}{2}\frac{p^2}{\mu^2}-\frac{1}{8}\frac{p^4}{\mu^4}+\frac{1}{16}\frac{p^6}{\mu^6}+\ldots\right) = \\
&= \mu+\frac{1}{2}\frac{p^2}{\mu}-\frac{1}{8}\frac{p^4}{\mu^3}+\frac{1}{16}\frac{p^6}{\mu^5}+\ldots
\end{aligned} \tag{8.3.12}$$

By inserting the Eq.(8.3.12) into the Eq.(8.3.8) we get

$$\begin{aligned}
&\varepsilon\left(\mu^{(1)}_{\mathbf{eff}},\mu^{(2)}_{\mathbf{eff}}\right) = \\
&K\int_{\mu^{(1)}_{\mathbf{eff}}}^{\mu^{(2)}_{\mathbf{eff}}} d\mu f(\mu) \int\limits_{p\leq\sqrt{\mu}} \left(\mu+\frac{1}{2}\frac{p^2}{\mu}-\frac{1}{8}\frac{p^4}{\mu^3}+\frac{1}{16}\frac{p^6}{\mu^5}+\ldots\right)p^2dp = \\
&K\int_{\mu^{(1)}_{\mathbf{eff}}}^{\mu^{(2)}_{\mathbf{eff}}} d\mu f(\mu) \left[\int\limits_0^{\sqrt{\mu}} \left(\mu p^2+\frac{1}{2}\frac{p^4}{\mu}-\frac{1}{8}\frac{p^6}{\mu^3}+\frac{1}{16}\frac{p^8}{\mu^5}+\ldots\right)dp\right] = \\
&K\int_{\mu^{(1)}_{\mathbf{eff}}}^{\mu^{(2)}_{\mathbf{eff}}} d\mu f(\mu) \left[\mu\frac{p^3}{3}+\frac{1}{2}\frac{p^5}{5\mu}-\frac{1}{8}\frac{p^7}{7\mu^3}+\frac{1}{16}\frac{p^9}{9\mu^5}+\ldots\right]_0^{\sqrt{\mu}} = \\
&K\int_{\mu^{(1)}_{\mathbf{eff}}}^{\mu^{(2)}_{\mathbf{eff}}} d\mu f(\mu) \left[\mu\frac{\mu^{\frac{3}{2}}}{3}+\frac{1}{2}\frac{\mu^{\frac{5}{2}}}{5\mu}-\frac{1}{8}\frac{\mu^{\frac{1}{2}}}{7\mu^3}+\frac{1}{16}\frac{\mu^{\frac{9}{2}}}{9\mu^5}+\ldots\right] = \\
&K\int_{\mu^{(1)}_{\mathbf{eff}}}^{\mu^{(2)}_{\mathbf{eff}}} f(\mu)d\mu\left[\frac{1}{3}\mu^{\frac{5}{2}}+\frac{1}{10}\mu^{\frac{3}{2}}-\frac{1}{56}\mu^{\frac{1}{2}}+\frac{1}{144}\mu^{-\frac{1}{2}}+\ldots\right] = \\
&K\int_{\mu^{(1)}_{\mathbf{eff}}}^{\mu^{(2)}_{\mathbf{eff}}} f(\mu)d\mu\left[\frac{1}{3}\mu^{\frac{5}{2}}+\frac{1}{10}\mu^{\frac{3}{2}}-\frac{1}{56}\mu^{\frac{1}{2}}+\frac{1}{144}\mu^{-\frac{1}{2}}\right]+o\left(\left(\mu^{(1)}_{\mathbf{eff}}\right)^{-n+1/2}\right).
\end{aligned} \tag{8.3.13}$$

The pressure $p\left(D^{+},D^{-},\mu_{\mathbf{eff}}^{(1)},\mu_{\mathbf{eff}}^{(2)}\right)$ for free quantum fields is

$$p\left(D^{+},D^{-},\mu_{\mathbf{eff}}^{(1)},\mu_{\mathbf{eff}}^{(2)}\right)=p\left(\mu_{\mathbf{eff}}^{(1)},\mu_{\mathbf{eff}}^{(2)},\right)+\breve{p}\left(D^{+},D^{-},\mu_{\mathbf{eff}}^{(1)},\mu_{\mathbf{eff}}^{(2)}\right). \quad (8.3.14)$$

Here the quantity $p\left(\mu_{\mathbf{eff}}^{(1)},\mu_{\mathbf{eff}}^{(2)},\right)$ is given by formula

$$\begin{aligned} p\left(\mu_{\mathbf{eff}}^{(1)},\mu_{\mathbf{eff}}^{(2)}\right) &= \frac{1}{2(2\pi\hbar)^{3}}\int_{\mu_{\mathbf{eff}}^{(1)}}^{\mu_{\mathbf{eff}}^{(2)}} d\mu f(\mu)\int\limits_{\|\mathbf{k}\|\leq\sqrt{\mu}}\frac{\|\mathbf{k}\|^{2}}{\sqrt{\mathbf{k}^{2}+\mu^{2}}}d^{3}\mathbf{k} = \\ &= \frac{K}{3}\int_{\mu_{\mathbf{eff}}^{(1)}}^{\mu_{\mathbf{eff}}^{(2)}} d\mu f(\mu)\int\limits_{p\leq\sqrt{\mu}}\frac{p^{4}}{\sqrt{p^{2}+\mu^{2}}}dp. \end{aligned} \quad (8.3.15)$$

The quantity $\breve{p}\left(D^{+},D^{-},\mu_{\mathbf{eff}}^{(1)},\mu_{\mathbf{eff}}^{(2)}\right)$ is given by formula

$$\breve{p}\left(D^{+},D^{-},\mu_{\mathbf{eff}}^{(1)},\mu_{\mathbf{eff}}^{(2)}\right) \simeq \frac{K'}{3}\int_{\mu_{\mathbf{eff}}^{(1)}}^{\mu_{\mathbf{eff}}^{(2)}} d\mu f(\mu)\int\limits_{\|p\|>\sqrt{\mu}}\frac{\|\mathbf{k}\|^{2}}{\sqrt{\mathbf{k}^{2}+\mu^{2}}}d^{D^{+},D^{-}}\mathbf{k}, \quad (8.3.16)$$

where $K'=\dfrac{1}{2(2\pi\hbar)^{3}}, c=1$.Note that

$$\begin{aligned} &\frac{1}{\sqrt{p^{2}+\mu^{2}}}=\mu^{-1}\left(\sqrt{1+\frac{p^{2}}{\mu^{2}}}\right)^{-1}= \\ &\mu^{-1}\left(1-\frac{1}{2}\frac{p^{2}}{\mu^{2}}+\frac{3}{8}\frac{p^{4}}{\mu^{4}}-\frac{5}{16}\frac{p^{6}}{\mu^{6}}+\ldots\right)= \\ &=\frac{1}{\mu}-\frac{1}{2}\frac{p^{2}}{\mu^{3}}+\frac{3}{8}\frac{p^{4}}{\mu^{5}}-\frac{5}{16}\frac{p^{6}}{\mu^{7}}+\ldots \end{aligned} \quad (8.3.17)$$

By inserting Eq.(8.3.17) into Eq.(8.3.15) we get

$$\begin{aligned} &p\left(\mu_{\mathbf{eff}}^{(1)},\mu_{\mathbf{eff}}^{(2)}\right)= \\ &\frac{K}{3}\int_{\mu_{\mathbf{eff}}^{(1)}}^{\mu_{\mathbf{eff}}^{(2)}} d\mu f(\mu)\int\limits_{p\leq\sqrt{\mu}}\left[\frac{1}{\mu}-\frac{1}{2}\frac{p^{2}}{\mu^{3}}+\frac{3}{8}\frac{p^{4}}{\mu^{5}}-\frac{5}{16}\frac{p^{6}}{\mu^{7}}+\ldots\right]p^{4}dp= \\ &\frac{K}{3}\int_{\mu_{\mathbf{eff}}^{(1)}}^{\mu_{\mathbf{eff}}^{(2)}} d\mu f(\mu)\int\limits_{p\leq\sqrt{\mu}}\left[\frac{p^{4}}{\mu}-\frac{1}{2}\frac{p^{6}}{\mu^{3}}+\frac{3}{8}\frac{p^{8}}{\mu^{5}}-\frac{5}{16}\frac{p^{10}}{\mu^{7}}+\ldots\right]dp= \\ &\frac{K}{3}\int_{\mu_{\mathbf{eff}}^{(1)}}^{\mu_{\mathbf{eff}}^{(2)}} d\mu f(\mu)\left[\frac{p^{5}}{5\mu}-\frac{1}{2}\frac{p^{7}}{7\mu^{3}}+\frac{3}{8}\frac{p^{9}}{9\mu^{5}}-\frac{5}{16}\frac{p^{11}}{10\mu^{7}}+\ldots\right]_{0}^{\sqrt{\mu}}= \\ &\frac{K}{3}\int_{\mu_{\mathbf{eff}}^{(1)}}^{\mu_{\mathbf{eff}}^{(2)}} d\mu f(\mu)\left[\frac{\mu^{\frac{5}{2}}}{5\mu}-\frac{1}{2}\frac{\mu^{\frac{7}{2}}}{7\mu^{3}}+\frac{3}{8}\frac{\mu^{\frac{9}{2}}}{9\mu^{5}}-\frac{5}{16}\frac{\mu^{\frac{11}{2}}}{10\mu^{7}}+\ldots\right]= \\ &\frac{K}{3}\int_{\mu_{\mathbf{eff}}^{(1)}}^{\mu_{\mathbf{eff}}^{(2)}} d\mu f(\mu)\left[\frac{1}{5}\mu^{\frac{3}{2}}-\frac{1}{14}\mu^{\frac{1}{2}}+\frac{1}{24}\mu^{-\frac{1}{2}}-\frac{1}{32}\mu^{-\frac{3}{2}}+\ldots\right]= \\ &\frac{K}{3}\int_{\mu_{\mathbf{eff}}^{(1)}}^{\mu_{\mathbf{eff}}^{(2)}} d\mu f(\mu)\left[\frac{1}{5}\mu^{\frac{3}{2}}-\frac{1}{14}\mu^{\frac{1}{2}}+\frac{1}{24}\mu^{-\frac{1}{2}}-\frac{1}{32}\mu^{-\frac{3}{2}}\right]+o\left(\left(\mu_{\mathbf{eff}}^{(1)}\right)^{-n-1/2}\right). \end{aligned} \quad (8.3.18)$$

# 9.Discussion and conclusions

We will now briefly review the canonical assumptions that are made in the usual formulation of the cosmological constant problem.

# 9.1.The canonical assumptions:

**1**.**The physical dark matter**.

Dark matter is a hypothetical form of matter that is thought to account for approximately 85% of the matter in the universe, and about a quarter of its total energy density. The majority of dark matter is thought to be non-baryonic in nature, possibly being composed of some as-yet undiscovered subatomic particles.Its presence is implied in a variety of astrophysical observations, including gravitational effects that cannot be explained unless more matter is present than can be seen. For this reason, most experts think dark matter to be ubiquitous in the universe and to have had a strong influence on its structure and evolution. The name dark matter refers to the fact that it does not appear to interact with observable electromagnetic radiation, such as light, and is thus invisible (or 'dark') to the entire electromagnetic spectrum, making it extremely difficult to detect using usual astronomical equipment.Because dark matter has not yet been observed directly, it must barely interact with ordinary baryonic matter and radiation. The primary candidate for dark matter is some new kind of elementary particle that has not yet been discovered, in particular, weakly-interacting massive particles (WIMPs), or gravitationally-interacting massive particles (GIMPs).Many experiments to directly detect and study dark matter particles are being actively undertaken, but none has yet succeeded.

**2**.The total effective cosmological constant $\lambda_{\mathbf{eff}}$ is on at least the order of magnitude of the vacuum energy density generated by zero-point fluctuations of the standard particle fields.

**3**.Canonical QFT is an effective field theory description of a more fundamental theory, which becomes significant at some high-energy scale $\Lambda_*$.

**4**.The vacuum energy-momentum tensor is Lorentz invariant.

**5**.The Moller-Rosenfeld approach [35],[36] to semiclassical gravity by using an expectation value for the energy-momentum tensor is sound.

**6**.The Einstein equations for the homogeneous Friedmann-Robertson-Walker metric accurately describes the large-scale evolution of the Universe.

**Remark 9**.**1**.**1**.Note that obviously there is a strong inconsistency between Assumptions **2** and **3**: the vacuum state cannot be Lorentz invariant if modes are ignored above some high-energy cutoff $\Lambda_*$, because a mode that is high energy in one reference frame will be low energy in another appropriately boosted frame. In this paper Assumption 3 is not used and this contradiction is avoided.

**Remark 9.1.2.**Note that also, Assumptions 1,3,4 and 5 is modifed, which we denote as

Assumptions 4 and 5 respectively.

## 9.2.Modified assumptions

**1**$'$.The physical dark matter.

**2**$'$.The total effective cosmological constant $\lambda_{\mathbf{eff}}$ is on at least the order $|\mu_{\mathbf{eff}}|^{-n+5}\ln|\mu_{\mathbf{eff}}|$ of

magnitude of the*renormalized* vacuum energy density generated by zero-point fluctuations of standard particle fields and ghost particle fields,see subsection I.2.

**4**$'$.The vacuum energy-momentum tensor is not Lorentz invariant.

## 9.3.The physical ghost matter and dark matter nature

In the contemporary quantum field theory, a ghost field, or gauge ghost is an unphysical state in a gauge theory. Ghosts are necessary to keep gauge invariance in theories where the local fields exceed a number of physical degrees of freedom.For example in quantum electrodynamics, in order to maintain manifest Lorentz invariance, one uses a four component vector potential $A_{\mu}(x)$, whereas the photon has only two polarizations. Thus, one needs a suitable mechanism in order to get rid of the unphysical degrees of freedom. Introducing fictitious fields, the ghosts, is one way of achieving this goal. Faddeev-Popov ghosts are extraneous fields which are introduced to maintain the consistency of the path integral formulation. Faddeev-Popov ghosts are sometimes referred to as "good ghosts".

"Bad ghosts" represent another, more general meaning of the word "ghost" in theoretical physics: states of negative norm,or fields with the wrong sign of the kinetic term, such as Pauli-Villars ghosts, whose existence allows the probabilities to be negative thus violating unitarity.

(**IX.1**) In contrary with standard Assumption1 in the case of the new approach introduced

in this paper we assume that:

(**IX.1.1**.a) The ghosts fields and ghosts particles with masses at a scale less then an fixed scale $m_{\mathbf{eff}}$ really exist in the universe and formed dark matter sector of the universe,in particular:

(**IX.1.1**.b) these ghosts fields gives additive contribution to a full zero-point fluctuation (i.e.

also to effective cosmological constant $\lambda_{\mathbf{eff}}$ [5],see subsection I.2).

(**IX.1.1**.c) Pauli-Villars renormalization of zero-point fluctuations (see subsection I.2) is no longer considered as an intermediate mathematical construct but obviously has rigorous physical meaning supported by assumption (**I**.a-b).

(**IX.1.2**) The physical dark matter formed by ghosts particles;
(**IX.1.3**) The standard model fields do not to couple directly to the ghost sector in the ultraviolet region of energy at a scale less then an fixed large energy scale $\Lambda_*$,in particular:
(**IX.1.3**.a) The "bad" ghosts fields with masses at a scale less then an fixed scale $m_{\mathbf{eff}}$, where $m_{\mathbf{eff}}c^2 \ll \Lambda_*$,cannot appear in any effective physycal lagrangian which contain also
the standard particles fields.
In additional though not necessary we assume that:
(**IX.1.4**) The "bad" ghosts fields with masses at a scale $m_*$,where $m_*c^2 \gg \Lambda_*$can appear
in any effective physycal lagrangian which contain also the standard particles fields,in particular:
(**IX.1.4**.a) Pauli-Villars finite renormalization with masses of ghosts fields at a scale $m_*$ of
the S-matrix in QFT (see subsection II.I-2) is no longer considered as an intermediate mathematical construct but obviously has rigorous physical meaning supported by assumption (**IX**).
(**IX.1.4**.b) If the "bad" ghosts fields coupled to matter directly, it gives rise to small and controlable violetion of the unitarity condition.
**Remark 9.3.1.**We emphazize that in universe standard matter coupled with a *physical* ghost matter has the equation of state [3]:

$$\varepsilon_{\mathbf{vac}}(\mu_{\mathbf{eff}}) = -p(\mu_{\mathbf{eff}}) = \frac{1}{8}\int_0^{\mu_{\mathbf{eff}}} f(\mu)\mu^4(\ln\mu)d\mu = \frac{c^4\lambda_{\mathbf{vac}}}{8\pi G}, \qquad (9.3.1)$$

where

$$|f(\mu)| = \begin{cases} O(\mu^{-n}), n > 1 & \mu \leq \mu_{\mathbf{eff}} \\ 0 & \mu > \mu_{\mathbf{eff}} \end{cases} \qquad (9.3.2)$$

and where $\mu_{\mathbf{eff}} = m_{\mathbf{eff}}c$ (see subsection I.2,Eq.(1.2.16)) and therefore gives rise to a de Sitter phase of the universe even if bare cosmological constant $\lambda = 0$.

## 9.4.Different contributions to $\lambda_{\mathbf{eff}}$

The total effective cosmological constant $\lambda_{\mathbf{eff}}$ is on at least the order of magnitude of the vacuum energy density generated by zero-point fluctuations of *standard* particle fields.

Assumption 2 is well justified in the case of the traditional approach, because the contribution from zero-point fluctuations is on the order of $1$ in Planck units and no other known contributions are as large thus, assuming no significant cancellation of terms (e.g. fine tuning of the bare cosmological constant $\lambda$), the total $\lambda_{\mathbf{eff}}$ should be at least on the order of the largest contribution [14].

(**9.4**) In contrary with standard Assumption1 in the case of the new approach introduced
in this paper we assume that:
(**9.4.1**) For simplisity though not necessary bare cosmological constant $\lambda = 0$.

(**9.4.2**) The total effective cosmological constant $\lambda_{\mathbf{eff}}$ depend only on mass distribution $f(\mu)$ and constant $\mu_{\mathbf{eff}}$ where $\mu_{\mathbf{eff}}c^2 < \Lambda_*$ but cannot depend on large energy scale $\Lambda_*$

**Remark 9.4.1.**Note that in subsection 9.1 we pointed out that under Assumption 1 if bare cosmological constant $\lambda = 0$ the total cosmological constant $\lambda_{\mathbf{vac}}$ is on at least the order $\sim |\mu_{\mathbf{eff}}|^{-n+5}$ of magnitude of the *renormalized* vacuum energy density generated by
zero-point fluctuations of standard particle fields and ghost particle fields

$$\begin{aligned} \varepsilon_{\mathbf{vac}}\left(\mu^*_{\mathbf{eff}}\right) &= \frac{1}{8}\int\limits_0^{\mu^*_{\mathbf{eff}}} f(\mu,\Lambda_*)\mu^4(\ln\mu)d\mu + O(\Lambda_*^{-2}), \\ p_{\mathbf{vac}}\left(\mu^*_{\mathbf{eff}}\right) &= -\frac{1}{8}\int\limits_0^{\mu^*_{\mathbf{eff}}} f(\mu,\Lambda_*)\mu^4(\ln\mu)d\mu + O(\Lambda_*^{-2}), \end{aligned} \tag{9.4.1}$$

where $\mu^*_{\mathbf{eff}} = \mu_{\mathbf{eff}}\,(\Lambda_*)$

# 9.5. Effective field theory and Lorentz invariance violetion

To prevent the vacuum energy density from diverging,the traditional approach also assumes that performing a high-energy cutoff is acceptable. This type of regularization is a common step in renormalization procedures, which aim to eventually arrive at a physical, cutoff-independent result. However,in the case of the vacuum energy density, the result is inherently cutoff dependent, scaling quartically with the cutoff $\Lambda_*$.

**Remark 9.5.1**. By restricting to modes with particle energy a certain cutoff energy $\omega_{\mathbf{k}} \leq \Lambda_*$a finite, regularized result for the energy density can be obtained. The result is proportional to $\Lambda_*^4$.Any other fields will contribute similarly, so that if there are $n_b$ bosonic
fields and $n_f$ fermionic fields, the density scales with $(n_b - 4n_f)\ \Lambda_*^4$. Typically, the cutoff is
taken to be near $= 1$ in Planck units (i.e.the Planck energy), so the vacuum energy gives a
contribution to the cosmological constant on the order of at least unity according to Eq. (9.6.5). Thus we see the extreme fine-tuning problem: the original cosmological constant $\lambda$ must cancel this large vacuum energy density $\varepsilon_{\mathbf{vac}} \simeq 1$ to a precision of 1 in $10^{120}$ -but not completely- to result in the observed value $\lambda_{\mathbf{eff}} = 10^{-120}$[5].

**Remark 9.5.2**. As it pointed out in this paper that a high energy theory, i.e. QFT in fractal space-time with Hausdorff-Colombeau negative dimension would not display the zero-point fluctuations that are characteristic of QFT, and hence that the divergence caused by oscillations above the corresponding cutoff frequency is unphysical. In this case, the cutoff $\Lambda_*$ is no longer an intermediate mathematical construct, but instead a physical scale at which the smooth,continuous behavior of QFT breaks down.

Poincaré group of the momentum space is deformed at some fundamental high-energy cutoff $\Lambda_*$The canonical quadratic invariant $\|p\|^2 = \eta^{ab}p_a p_b$ collapses at high-energy cutoff $\Lambda_*$ and being replaced by the non-quadratic invariant:

$$\|p\|^2 = \frac{\eta^{ab}p_a p_b}{(1 + l_{\Lambda_*}p_0)}. \tag{9.5.1}$$

**Remark 9.5.3**. In contrary with canonical approuch the total effective cosmological constant $\lambda_{\mathbf{eff}}$ depend only on mass distribution $f(\mu)$ and constant $\mu_{\mathbf{eff}} = m_{\mathbf{eff}}c$ but cannot depend on large energy scale $\sim \Lambda_*$.

## 9.6. Semiclassical Moller-Rosenfeld gravity

Assumption 5 means that it is valid to replace the right-hand side of the Einstein equation $T_{\mu\nu}$ with its expectation $\langle T_{\mu\nu}\rangle$. It requires that either gravity is not in fact quantum, and the Moller-Rosenfeld approach is a complete description of reality, or at least a valid approximation in the weak field limit.The usual argument states that the vacuum state $|0\rangle$ should be locally Lorentz invariant so that observers agree on the vacuum state. This means that the expectation value of the energy-momentum tensor on the vacuum, $\langle 0|\widehat{T}_{\mu\nu}|0\rangle$,must be a scalar multiple of the metric tensor $g_{\mu\nu}$ which is the only Lorentz invariant rank $(0,2)$ tensor. By using Moller-Rosenfeld approach the Einstein field equations of general relativity, a term representing the curvature of spacetime $R_{\mu\nu}$ is related to a term describing the energy-momentum of matter $\langle 0|\widehat{T}_{\mu\nu}|0\rangle$, as well as the cosmological constant $\lambda$ and metric tensor $g_{\mu\nu}$ reads:

$$R_{\mu\nu} - \frac{1}{2}R^{\upsilon}_{\upsilon}g_{\mu\nu} + \lambda g_{\mu\nu} = 8\pi\langle 0|\widehat{T}_{\mu\nu}|0\rangle. \tag{9.6.1}$$

The $\widehat{T}_{00}$ component is an energy density, we label $\langle 0|\widehat{T}_{\mu\nu}|0\rangle = \varepsilon_{\mathbf{vac}}$,so that the vacuum contribution to the right-hand side of Eq.(9.4.1) can be written as

$$8\pi\langle 0|\widehat{T}_{\mu\nu}|0\rangle = 8\pi\varepsilon_{\mathbf{vac}}g_{\mu\nu}. \tag{9.6.2}$$

Subtracting this from the right-hand side of Eq.(9.4.1) and grouping it with the cosmological constant term replaces with an "effective" cosmological constant [5]:

$$\lambda_{\mathbf{eff}} = \lambda + 8\pi\varepsilon_{\mathbf{vac}}. \tag{9.6.3}$$

Note that in flat spacetime, where $g_{\mu\nu} = diag(-1,+1,+1,+1)$, Eq.(9.4.2) implies $\varepsilon_{\mathbf{vac}} = -p_{\mathbf{vac}}$, where $p_{\mathbf{vac}} = \langle 0|\widehat{T}_{ii}|0\rangle$ for any $i = 1,2,3$ is the pressure. Obviously this implies that if the energy density is positive as is usually assumed, then the pressure must be negative, a conclusion which extends to any metric $g_{\mu\nu}$ with a $(-1,+1,+1,+1)$ signature.

**Remark 9.6.1.**In this paper we assume that the vacuum state $|0\rangle$ should be locally invariant under modified Lorentz boost (1.1.18) but not locally Lorentz invariant. Obviously this assumption violate the Eq.(9.6.2). However modified Lorentz boosts (1.1.18) becomes Lorentz boosts for a sufficiently small energies and therefore in IR region one obtain in a good aproximation

$$8\pi\langle 0|\widehat{T}_{\mu\nu}|0\rangle \approx 8\pi\varepsilon_{\mathbf{vac}}g_{\mu\nu} \tag{9.6.4}$$

and

$$\lambda_{\mathbf{eff}} \approx \lambda + 8\pi\varepsilon_{\mathbf{vac}}. \tag{9.6.5}$$

Thus Moller-Rosenfeld approach holds in a good approximation.

## 9.7. Quantum gravity at energy scale $\Lambda \leq \Lambda_*$Controlable violetion of the unitarity condition.

Gravitational actions which include terms quadratic in the curvature tensor are renormalizable. The necessary Slavnov identities are derived from Becchi-Rouet-Stora (BRS) transformations of the gravitational and Faddeev-Popov ghost fields. In general, non-gauge-invariant divergences do arise, but they may be absorbed by nonlinear renormalizations of the gravitational and ghost fields and of the BRS transformations [13].The geneic expression of the action reads

$$I_{sym} = -\int d^4x\sqrt{-g}\,(\alpha R_{\mu\nu}R^{\mu\nu} - \beta R^2 + 2\kappa^{-2}R), \qquad (9.7.1)$$

where the curvature tensor and the Ricci is defined by $R^{\lambda}_{\mu\alpha\nu} = \partial_\nu \Gamma^{\lambda}_{\mu\alpha}$ and $R_{\mu\nu} = R^{\lambda}_{\mu\lambda\nu}$ correspondingly, $\kappa^2 = 32\pi G.$The convenient definition of the gravitational field variable in terms of the contravariant metric density reads

$$\kappa h^{\mu\nu} = g^{\mu\nu}\sqrt{-g} - \eta^{\mu\nu}. \qquad (9.7.2)$$

Analysis of the linearized radiation shows that there are eight dynamical degrees of freedom in the field. Two of these excitations correspond to the familiar massless spin-2 graviton. Five more correspond to a massive spin-2 particle with mass $m_2$. The eighth corresponds to a massive scalar particle with mass $m_0$. Although the linearized field energy of the massless spin-2 and massive scalar excitations is positive definite, the linearized energy of the massive spin-2 excitations is negative definite. This feature is characteristic of higher-derivative models, and poses the major obstacle to their physical interpretation.

In the quantum theory, there is an alternative problem which may be substituted for the negative energy. It is possible to recast the theory so that the massive spin-2 eigenstates of the free-fieid Hamiltonian have positive-definite energy, but also negative norm in the state vector space.These negative-norm states cannot be excluded from the physical sector of the vector space without destroying the unitarity of the **S** matrix. The requirement that the graviton propagator behave like $p^{-4}$ for large momenta makes it necessary to choose the indefinite-metric vector space over the negative-energy states.The presence of massive quantum states of negative norm which cancel some of the divergences due to the massless states is analogous to the Pauli-Villars regularization of other field theories. For quantum gravity, however, the resulting improvement in the ultraviolet behavior of the theory is sufficient only to make it renormalizable,but not finite.

**Remark 9.7.1.(I)**The renormalizable models which we have considered in this paper many years mistakenly regarded only as constructs for a study of the ultraviolet problem of quantum gravity. The difficulties with unitarity appear to preclude their direct acceptability as canonical physical theories in locally Minkowski space-time. In canonical case they do have only some promise as phenomenological models.

(**II**) However, for their unphysical behavior may be restricted to *arbitrarily large energy scales* $\Lambda_*$ mentioned above by an appropriate limitation on the renormalized masses $m_2$ and $m_0$. Actually, it is only the massive spin-two excitations of the field which give the trouble with unitarity and thus require a very large mass. The limit on the mass $m_0$ is determined only by the observational constraints on the static field.

# 10.Conclusion

We argue that a solution to the cosmological constant problem is to assume that there exists hidden physical mechanism which cancel divergences in canonical $QED_4, QCD_4$, Higher-Derivative-Quantum-Gravity, etc. In fact we argue that corresponding supermassive Pauli-Villars ghost fields,etc.really exists. New theory of elementary particles which contains hidden ghost sector is proposed. Zel'dovich hypotesis [1] we suggest that physics of elementary particles is separated into low/high energy ones the standard notion of smooth spacetime is assumed to be altered at a high energy cutoff scale $\Lambda_*$ and a new treatment based on QFT in a fractal spacetime with negative dimension is used above that scale.This would fit in the observed value of the dark energy needed to explain the accelerated expansion of the universe if we choose highly symmetric masses distribution below that scale $\Lambda_*$, i.e.,
$f_{s.m}(\mu) \approx f_{g.m}(\mu), \mu \leq \mu_{\mathbf{eff}}, \mu_{\mathbf{eff}}c^2 < \Lambda_*$

# ACKNOWLEDGMENTS

To rewiever

# Refferences